\documentclass[reqno,article]{amsart}
\RequirePackage{amsthm,amsmath,amsfonts,amssymb}
\usepackage[dvipdfmx]{graphicx}
\usepackage{bm}
\usepackage{mathrsfs}

\newtheorem{theorem}{Theorem}[section]
\newtheorem{proposition}{Proposition}[section]% 

\newtheorem{lemma}{Lemma}[section]

\theoremstyle{remark}
	\newtheorem{remark}{Remark}[section]
\newtheorem{definition}{Definition}[section]

\newcommand\DN{\newcommand}

\DN\B{\mathbf{B}}
\DN\bs{\bigskip}
\DN\Ixs{I (x,x,s)}
\DN\Ixxs{I (x,x',s)}
\DN\Ixxxs{I (x,x,s)}
\DN\SUMklm{\sum_{k=0}^{L-1}\sum_{m=0}^k}
\DN\SECT{\subsection}%\DN\SECT{\section}
\DN\iotaN{\iota _{N}} 
\DN\iotaNx{\iota _{N}(x)} \DN\iotaNy{\iota _{N}(y)} 
 \DN\iotaNz{\iota _{N}(z)} 
\DN\yhat{\check{y}}
\DN\xhat{\check{x}}
\DN\yhatyhat{y}

\DN\xyhat{x - \yhat }

\DN\yhaT{y}
 \DN\zhat{\check{z}}
\DN\xy{\lvert x - y \rvert }
\DN\QR{(Q,R)}
\DN\ellell{\ell}
\DN\Kakk{\mathfrak{K}[\mathit{a}_{\qq }^+]}
\DN\Stwoqm{\overline{\Sr }^m } %{\Rm _{\rr } }
\DN\Ssi{\mathfrak{S} _{\mathrm{s.i.}}}

\DN\nnNNN{\nn \in \NNN }

\DN\xs{ x ,\sss }
\DN\xxs{\mathbf{x},\sss }
\DN\tbm{ \widetilde{ b }^m}

\DN\Ha{\mathfrak{H}}
\DN\Ham{\Ha ^{[m]}}
\DN\Hpqrm{\Ha _{\pqr }^{[m]}}
\DN\HpqrmC{\Ha _{\pqr }^{[m], \circ }}
 \DN\HIiC{ \HanC \backslash \LLL \Hmnc \rRRR } 
\DN\HmnP{\Pitwo (\Hmn )} \DN\HmnPc{\Pitwo (\Hmnc )}
\DN\Hmnc{\HanC \cap \IIm }
\DN\HanC{\Han ^{\circ }}
\DN\Han{\Ha _{\nn }}
\DN\Hanm{\Han ^{[m]} }

\DN\NNN{\mathsf{N}} \DN\nn{\mathsf{n}}
\DN\NNNthree{\NNN _3 }

\DN\IIm{\mathfrak{I}_{\mathsf{m}}}
\DN\SmSS{\Rm \times \SSS }
\DN\mm{\mathsf{m}}
\DN\Pitwo{\prod}
\DN\pqr{p,q,r}
\DN\fg{\widetilde{b}^{m}}

\DN\Ka{\mathfrak{K}[\mathbf{a}]}
\DN\ak{a_{q}}
\DN\akk{a_{q+1}}
\DN\qq{q}

\DN\pP{P}
\DN\pPF{ P _{\{ \Fs \} } }
\DN\Ft{\{ \mathcal{F}_t \} }
\DN\WRNz{W_{\mathbf{0}} (\RN )}
\DN\WRN{W (\RN )}
\DN\PBr{\mathcal{W}_{\mathbf{0}}}

	\DN\capamu{\capa ^{\mu }}
	\DN\capamuA{\capa ^{\muA }}
\DN\capa{\mathrm{Cap}_1}
\DN\capaone{\mathrm{Cap}_1^{\muairybeta ^{[1]} }}

\DN\XB{(\mathbf{X} , \BBB )}
\DN\ASSUMP{{\As{SIN}, \As{AC} for $ \mu_2 $, \As{IFC}}, and \As{MF}}
\DN\satASSUMP{{\As{SIN$ _{\mathbf{A}}$}, \As{AC} for $ \mu_2 $, \As{IFC}}, and \As{MF}}
\DN\ASSUMPt{{\As{SIN}, \As{AC} for $ \muAt $, \As{IFC}}, and \As{MF}}
\DN\ms{\medskip}

\DN\done{\textsf{done}}
\DN\Xmstar{\mathbf{X} ^{m,*}}
\DN\Xmmstar{\mathbf{X} ^{m+1,*}}

\DN\mmm{\mathsf{m}}

\DN\yyone{[\yY , \yY + 1]}
\DN\Yone{-1} 	\DN\YYone{ \min\{ \yY + 1 , -1 \}}

\DN\gdelta{( g , \delta )}
\DN\V{V }
\DN\Vi{\V _{ N ,i}}
\DN\Vone{\V _{ N , 1}}
\DN\Vtwo{\V _{ N , 2}}

\DN\None{[-2 N^{2/3}, -1]}
\DN\Ninfty{[-2 N^{2/3}, \infty) }
\DN\Noneopen{[-2 N^{2/3}, -1)}

\DN\Ntheta{( N , \theta )}
\DN\FNvV{\FN (\vv ^{-1} (\vV ) )}
\DN\FNuU{\FN (\uu ^{-1} (\uU ) )}
\DN\FNvVi{\FN ^{(\iii )}(\vv ^{-1} (\vV ) )}
\DN\FNvVz{\FN ^{(0)}(\vv ^{-1} (\vV ) )}
\DN\FNvVone{\FN ^{(1)}(\vv ^{-1} (\vV ) )}
\DN\FNvVtwo{\FN ^{(2)}(\vv ^{-1} (\vV ) )}
\DN\QM{4^{-1}}
\DN\zzz{^{(0)}}
\DN\Fone{^{(1)}}
\DN\Ftwo{^{(2)}}

\DN\xX{x}
\DN\yY{y}
\DN\xXx{x}
\DN\yYy{y}
\DN\xXX{x}
\DN\yYY{y}
\DN\uU{u} 
\DN\vV{v}
\DN\vVV{\lvert v \rvert ^}

\DN\rR{R}
\DN\yy{y}
\DN\zz{z}

\DN\iii{i}
\DN\partialtwo{\partial_2}

\DN\SV{\mathit{SV}}
\DN\dN{D_N}
\DN\dNz{\dN^{[0]}}
\DN\dNone{\dN^{[1]}}

\DN\uu{\mathsf{u}} \DN\vv{\mathsf{v}}

\DN\qqq{\mathsf{q}}
\DN\rrr{\mathsf{r}}

\DN\pppW{\mathsf{p} (\theta )}\DN\qqqW{\qqq (\theta )}
\DN\pppWW{\mathsf{p} '(\theta )}\DN\qqqWW{\qqq '(\theta )}

\DN\FN{F_N}
\DN\FNu{F_{N,\uu }}
\DN\FNv{F_{N,\vv }}

\DN\half{2^{-1}}
\DN\oneNell{1+N^{-1}\ell } \DN\oneN{1+N^{-1}} 
\DN{\aNp}{\alpha_{N , 1 }}
\DN{\aNm}{\alpha_{N , -1 }}
\DN{\aNl}{\alpha_{N , \ell }}
\DN{\aNll}{\alpha_{N , -\ell }}
\DN{\aNz}{\alpha_{N , 0 }}
\DN{\aN }{\alpha _N }

\DN\ga{\gamma }
\DN\gaN{\ga _{N }}
\DN\gaNz{\ga _{N , 0}}
\DN\gaNp{\ga _{N , 1 }} \DN\gaNm{\ga _{N , -1 }}
\DN\gaNl{\ga _{N ,\ell }}
\DN\gaNll{\ga _{N ,-\ell }}

\DN\gaga{\prod_{\lpm }\gaNl ^{1/4} }

\DN\deltal{\delta _{\ell }}
\DN\thetal{\theta _{\ell }}

\DN\lpm{\ell = \pm 1}
\DN\fz{ \mathsf{f}} \DN\gz{\mathsf{g}}
\DN\fell{ \mathsf{f}_{\ell }} \DN\gell{\mathsf{g}_{\ell }}
\DN\Azz{{\mathsf{f}^{-1/4}}}
\DN\Al{{\mathsf{f}_{\ell }^{-1/4}}}
\DN\All{{\mathsf{f}_{-1}^{-1/4}}}
\DN\Alll{{\mathsf{f}_{1}^{-1/4}}{\mathsf{f}_{-1}^{-1/4}}}

\DN\BB{\mathsf{B}}

\DN\Bell{\BB _{\ell }}
\DN\Bellp{\BB _{1}} 
\DN\Bellm{\BB _{-1} }

\DN\Bellkm{\Bell ^{\KM }}

\DN\Bellkmh{\wB _{\ell }^{\KM }}
\DN\wB{\widehat{\BB }}

\DN\wBell{\wB _{\ell }}
\DN\wBellm{\wB _{-\ell }}

\DN\sSS{\mathsf{S}}
\DN\Sell{\sSS _{\ell }}
\DN\wSell{\wS _{\ell }}
\DN\wS{\widehat{\sSS }}

\DN\CC{\mathsf{C}}
\DN\CCC{\mathbb{C}}
\DN\KM{ k m }

\DN\ct{\mathsf{c}} \DN\ckm{\ct ^{\KM }} \DN\cln{\ct ^{\LN }}
\DN\Ckm{\CC ^{\KM }}
\DN\Cln{\CC ^{\LN }}
\DN\Clkm{\CC _{\ell }^{\KM }}
\DN\Cz{\Ckm \Cln }

\DN\Cone{\CCC _{1}}
\DN\Cm{\CCC _{-1}}
\DN\Cl{\CCC _{\ell }}
\DN\CCl{\mathsf{C}_{\ell }}

\DN\SSSSp{\mathbb{S}_1-\mathbb{S}_{-1}}
\DN\SSSSpt{\widehat{\mathbb{S}}_1-\widehat{\mathbb{S}}_{-1}}

\DN\h{\uu }

\DN\gp{\mathfrak{q}_+}%\gp{\mathfrak{g}_+}
\DN\gm{\mathfrak{q}_-}%{\mathfrak{g}_-}

\DN\deltaz{\delta }

\DN\sumL{\sum_{k=0}^{L-1}\sum_{m=0}^k}
\DN\sumLL{{\sum_{l=0}^{L-1}\sum_{n=0}^l}}
\DN\IL{\mathbf{I}(L)}
\DN\DDN{\mathsf{D}(N )}

 \DN\PsiN{\Psi _N (x)}
\DN\PsiH{\widehat{\Psi }_N (x)}
\DN\psiNellH{\widehat{\psi }_{N,\ell }}
\DN\psiNpH{\widehat{\psi }_{N,1}}
\DN\psiNmH{\widehat{\psi }_{N,-1}} 
\DN\psiNzH{\widehat{\psi }_{N,0}} 

\DN\llll{l} \DN\n{n}\DN\LN{ln}

\DN\vertx{\lvert x \rvert }
\DN\verty{\lvert y \rvert }
\DN\xF{\vertx ^{3/2}}
\DN\yF{\verty ^{3/2}}

\DN\yrone{\verty \ge r+1 }

\DN\xvee{1 \vee \vertx } 
\DN\yvee{1 \vee \verty } 
	\DN\zvee{ 1 \vee \vert z \vert }
\DN\yhatvee{ 1 \vee \vert \yhat \rvert } 

\DN\yoneF{(\yvee )^{-1/4}}
\DN\ythreeF{(\yvee )^{-3/4}}
\DN\zoneF{(\zvee )^{-1/4}}

\DN\Pms{\mathbf{P}_{\sbm }^{m , * } }
\DN\Pones{\mathbf{P}_{\sbm }^{1 , * } }
\DN\Pmsh{\mathbf{P}_{\sbm , h }^{m , * }}

\DN\JJ{\mathsf{J}}
\DN\II{\mathsf{I}}
\DN\Ja{\JJ _1} \DN\JaN{\Ja ^N }
\DN\Jb{\JJ _2} \DN\JbN{\Jb ^N }
\DN\Jc{\JJ _3} \DN\JcN{\Jc ^N }
\DN\Jp{\JJ _p} \DN\JpN{\Jp ^N}

\DN\HHa{\mathsf{H}_1} 
\DN\HHb{\mathsf{H}_2}
\DN\HHc{\mathsf{H}_3}
\DN\HHp{\mathsf{H}_p}

\DN\RRpr{\mathbb{R}_{\pp ,\rr }}
\DN\RRprCs{\RRpr ^{m , \circ } (\sss )}
\DN\RRprsm{\RRpr ^m (\sss )}

\DN\pp{p}\DN\rr{r}
\DN\simpr{\sim _{\pp ,\rr }}
\DN\lpathX{\lpath (\mathfrak{X}) }
\DN\OFPFs{(\Omega ,\mathcal{F}, \PPs , \{ \mathcal{F}_t \} )}
\DN\PPs{\PP _{\mathfrak{s}}}
\DN\PP{\mathsf{P}}
\DN\HImnC{\LLL \Hmnc \rRRR } 
\DN\LLL{\langle}\DN\rRRR{\rangle}

\DN\UNkm{\UN ^{km}}
\DN\UNz{\UN ^{00}}
\DN\UN{\Upsilon_N }

\DN\OL{\underline{\mathcal{O}} }
\DN\overBRN {\overline{\mathcal{B}( \mathbb{R} ^{\mathbb{N}} )}_0 }

\DN\1{\eqref{:Y1g}}%{\eqref{:Y1g}} 
\DN\2{ \OL ( x ^{-1})}
\DN\3{\OL (N^{-1/3}) } 
\DN\4{\OL (N^{-1}) } 
\DN\5{ \OL (N^{-2/3} x ^{-1})}
\DN\6{ \OL (N^{-2/3}\vertx ^{-1/2}) }
\DN\8{ \OL (N^{-1} x )}
\DN\9{\mathfrak{H}_{(p+1,q,r) }^{[1]}}
\DN\0{\prod_{\lpm }\Ai (\aNl ) }
\DN\XXX{\mathbb{X}}
\DN\7{ (\Omega , \mathcal{M}, \{ \XXX _t \}_{t\in[0,\infty)},\{ \PPPs \}_{\sss \in \SSS } ) }
\DN\BS{\mathcal{B}(\SSS ) }\DN\SBS{ (\SSS , \BS ) }
\DN\XP{(\mathbb{X},\{ \PPPs \}_{\sss \in \SSS } )}
\DN\PPP{\mathbb{P}}
\DN\PPPs{\mathbb{P}_{\sss }}

\DN\RA{\Rightarrow }

\DN\rrY{\int_{\verty <r}\frac{\rhohat (y)}{-y}dy }
 \DN\sbm{\mathbf{s}}
\DN\bbm{\mathbf{b}}
\DN\Hbm{\mathbf{H}}

\DN\Fs{F_{\sbm }}
\DN\Fsm{F_{\sbm }^m }
\DN\Wnex{W_{\mathrm{NE}}}

\DN\Supr{\sup_{\vertx \le r }}

\DN\Kxy{\frac{1}{2^{\frac{2}{3}}}\Kairyfour \Big(\frac{x}{2^{\frac{2}{3}}}, \frac{y}{2^{\frac{2}{3}}}\Big)}
\DN\At{\mathcal{A}(t)}
\DN\Az{\mathcal{A}(0)}

\DN\aAA{\mathfrak{A}}
\DN\bBB{\mathfrak{B}}
\DN\fFF{\mathfrak{F}}
\DN\kKK{\mathfrak{K}}
\DN\oOO{\mathfrak{O}}
\DN\NN{\mathcal{R}}
\DN\Lmu{L^2(\SSS ,\mu )}
\DN\Lmuairy{L^2(\SSS ,\muairy )}

\DN\Lloc{L_{\mathrm{loc}}} %¶'«

\DN\Ws{\mathcal{W}_{\sbm } }
\DN\Wsh{\mathcal{W}_{\sbm , h }}
\DN\bXm{\bbm _{\mathbf{X} }^{m }}
\DN\bXone{\bbm _{\mathbf{X} }^{1}}
\DN\bXmi{\bbm _{\mathbf{X} }^{m,i }}
\DN\bXonei{\bbm _{\mathbf{X} }^{1,1 }}

\DN\bYm{\mathbf{Y}^m} %{\bbm _{\mathbf{Y} }^{m }}

\DN\bYone{\bbm _{\mathbf{Y} }^{1}}
\DN\bYonei{\bbm _{\mathbf{Y} }^{1,1 }}

\DN\kk{{q}}%\DN\kk{\mathsf{k}}
\DN\nnnn{\nn } %\DN\nn{\mathsf{n}}

\DN\vq{\varphi _{\q }}
\DN\vr{\varphi _{r }}
\DN\chin{\chi _{\nn }}
\DN\chinn{\chi _{\nn +1}}

	 \DN\K{\mathfrak{K}^{\kappa }}
	\DN\Kq{\K _{q}} 

\DN\BBB{\mathbf{B}}
	\DN\E{\mathbb{E}}
\DN\zN{\{ 0 \} \cup \mathbb{N}}

\DN\Emu{\mathcal{E}^{\mu }}
\DN\domu{\dom ^{\mu }}

\DN\Ai{\mathrm{Ai}}
\DN\Det{\mathrm{Det}}
\DN\supp{{\mathrm{supp}}\ }

\DN\lref[1]{Lemma~\ref{#1}}
\DN\tref[1]{Theorem~\ref{#1}}
\DN\pref[1]{Proposition~\ref{#1}}
\DN\sref[1]{Section~\ref{#1}}
\DN\dref[1]{Definition~\ref{#1}}
\DN\rref[1]{Remark~\ref{#1}} 
\DN\coref[1]{Corollary~\ref{#1}}
\DN\eref[1]{Example~\ref{#1}}
\numberwithin{equation}{section}

\newcounter{Const} 
\DN\Ct{\refstepcounter{Const}c_{\theConst}}%\label{#1}
\DN\cref[1]{c_{\ref{#1}}}	

\DN\map[3]{#1\!:\!#2\!\to\!#3}
\DN\ot{\otimes} 
\DN\ts{\times }
\DN\qdet{\mathrm{qdet}}

\DN\limi[1]{\lim_{#1\to\infty}} 	
\DN\limz[1]{\lim_{#1\to 0}}
\DN\limsupi[1]{\limsup_{#1\to\infty}} 	
\DN\liminfi[1]{\liminf_{#1\to\infty}} 	
\DN\limsupz[1]{\limsup_{#1\to 0}} 	
\DN\liminfz[1]{\liminf_{#1\to 0}} 	

\DN\PD[2]{\frac{\partial#1}{\partial#2}}
\DN\ODod[2]{\frac{ d #1}{ d #2}}
\DN\Rd{\mathbb{R} ^d}
\DN\Rm{\mathbb{R} ^m}

\DN\supN{\sup_{N \in\mathbb{N}}}
\DN\SupN{\sup_{2\le N \in \mathbb{N}}}

\DN\rA{\rho _{\beta }}
\DN\rAx{\rho _{\beta , x }}

\DN\rairybeta{\rho _{ \beta }}
\DN\rairybetax{\rho _{\beta ,\, x }}
\DN\rairybetaxx{\rho _{ \beta ,\, X_t^i }}
\DN\rairybetazero{\rho _{ \beta ,\, 0 }}
\DN\rairyone{\rho _{ 1}}
\DN\rairyonex{\rho _{\Ai ,\, 1,\, x }}
\DN\rairy{\rho _{ 2}}
\DN\rairyx{\rho _{ 2,\, x }}
\DN\rairyfour{\rho _{ 4}}
\DN\rairyfourx{\rho _{ 4,\, x }}
\DN\Kairybeta{K _{ \beta }}
\DN\Kairybetax{K _{ \beta ,\, x }}

\DN\Kairyone{K _{ 1}}
\DN\Kairyonex{K _{ 1,\, x }}
\DN\Kairy{K _{ 2}}
\DN\Kairyx{K _{ 2,\, x }}
\DN\Kairyfour{K _{ 4}}
\DN\Kairyfourx{K _{ 4,\, x }}

\DN\rNone{\rho ^{N ,\, 1}}
\DN\rNtwo{\rho ^{N ,\, 2}}
\DN\rNnk{\rho ^{N ,\, n+k}}
\DN\rhohat{\hat{\varrho }}
\DN\rN{\rho ^{N ,\, n}}

\DN\sS{S}
\DN\Sk{\mathbb{R} ^{k}}
\DN\Sr{\sS _{r}}
\DN\Srr{\sS _{r + 1}}
\DN\SR{\sS _{ \rR }}
\DN\SRR{\sS _{ \rR +1 }}

\DN\SSS{\mathfrak{S}} 
\DN\SSSz{\mathfrak{S}_0} 
\DN\SSSzz{\mathfrak{S}} 
\DN\SSSk{\SSS ^{[k]}}
\DN\SSSsi{\SSS _{\mathrm{s.i.}}}
\DN\SSSsip{\SSS _{\mathrm{s.i.}}^{+f}}
\DN\RNN{\RN _{>} }
\DN\SoneSSS{\mathbb{R} \ts \SSS }
\DN\RN{\mathbb{R} ^{\mathbb{N}}}

\DN\OFPF{(\Omega ,\mathcal{F}, P ,\{ \mathcal{F}_t \} )}

\DN\dlog{\mathrm{d}} 
\DN\dmu{\dlog ^{\mu }}
\DN\dmuA{\dlog ^{\muA }}

\DN\sss{\mathfrak{s}} 
\DN\dom{\mathcal{D}} 
\DN\di{\dom _{\circ }}

\DN\ulab{\mathfrak{u} }
\DN\lab{\mathfrak{l} }

\DN\upath{\ulab _{\mathrm{path}}}
\DN\lpath{\lab _{\mathrm{path}}}

\DN\muA{\muairybeta }
\DN\muAt{\mu _{\beta , \mathfrak{t}}}
\DN\muAx{\mu _{\beta , x }}
\DN\muAxone{\mu _{\beta , x }^{[1]}}
\DN\muAxx{\mu _{\beta , \mathbf{x}}}
\DN\muAm{\muairybeta ^{[m]}} 
\DN\muAone{\muairybeta ^{[1]}} 

\DN\muairy{\mu_{ 2}}
\DN\muairyone{\mu_{ 1}}
\DN\muairyfour{\mu_{ 4}}
\DN\muairybeta{\mu_{ \beta }}

\DN\muairybetax{\mu_{ \beta ,x}} %¶'«
\DN\mairybeta{m_{\beta }}%¶'«
\DN\musin{\mu_{\mathrm{sin}, 2}} %'±'ê'͐¶'«

\DN\tauh{\tau _{h}} %'±'ê'͐¶'«
\DN\Srm{\SSS _r^m} %'±'ê'͐¶'«
\DN\muk{\mu ^{[k]}} %'±'ê'͐¶'«
\DN\mux{\mu _{x}} %'±'ê'͐¶'«
\DN\murkm{\mu _{r,k}^{m}} %'±'ê'͐¶'«
\DN\murkk{\mu _{r,k+1}} %'±'ê'͐¶'«
\DN\muone{\mu ^{[1]}} %¶'« 
\DN\mum{\mu ^{[m]}} %¶'« 
\DN\muN{\mu ^{N }} %¶'« 

\DN\MbN{M_{\beta}^N}
\DN\TailS{\mathcal{T}(\mathfrak{S})}
\DN\WTsol{W(\SSSsde )} %¶'«
\DN\SSSsde{\mathbf{S}_{\mathrm{sde} }}
\DN\Bt{\mathscr{B}_t}
\DN\zti{0 \le t < \infty }
\DN\Psbf{\mathbf{P}_{\sbm }}
\DN\Pssf{\mathbb{P}_{\mathfrak{s}}}
\DN\mut{\mu _{\mathfrak{t}}}
\DN\mutl{\mu _{\mathfrak{t}}^{\lab }}
\DN\mmutl{(m\mu _{\mathfrak{t}})^{\lab }}
\DN\Pmut{\mathbb{P}_{\mut }}
\DN\Pmu{\mathfrak{P}_{\mu }}

\DN\Ps{P_{\sbm }} %¶'« 
\DN\tzT{t \in [0,T]} %¶'« 

\DN\Lambdar{\Lambda _r }
\DN\Lambdarm{\Lambdar ^m }

\DN\As[1]{\thetag{\textbf{#1}}}

\DN\Om{\Omega}

\DN\Omh{\Om ^{\rm h}}
\DN\OmA{\Om _2} 
\DN\OmAd{\Om ^{\delta }}
\DN\OmAdz{\Om ^{\deltaz }}

\DN\OmAe{\OmA ^{\varepsilon }}
\DN\OmAde{\OmA ^{\delta , \varepsilon }}
\DN\Ome{\Om _1^{\varepsilon}}
\DN\Omei{\Om _i^{\varepsilon}}
\DN\OmAeone{\OmA ^{\varepsilon }} 
\DN\OmAetwo{\OmA ^{\varepsilon'}} 

\DN\bmu{\mathbf{u}}

\begin{document}

\title{Infinite-dimensional stochastic differential equations arising from Airy random point fields}

\maketitle
\markboth{Hirofumi Osada , Hideki Tanemura }{ISDEs arising from Airy random point fields}

\begin{center}
Revised version: 2024/6/18
\end{center}

\bs

\begin{center}
\textsf{Hirofumi Osada$ ^*$, \quad Hideki Tanemura{$^ \dagger$}}
\end{center}

\bs\noindent 
\begin{flushright}
$ {}^*$ Chubu University, \texttt{osada@isc.chubu.ac.jp}\\
$ ^\dagger$ Keio University, \texttt{tanemura@math.keio.ac.jp}
\end{flushright}

\bs 

% \author[1]{\fnm{Hirofumi} \sur{Osada}}\email{osada@isc.chubu.ac.jp}
% 
% \author*[2]{\fnm{Hideki} \sur{Tanemura}}\email{tanemura@math.keio.ac.jp}
% \affil[1]{ \orgname{Chubu University}, \orgaddress{\street{1200 Matsumoto}, \city{Kasugai}, \postcode{487-8501}, \state{Aichi}, \country{Japan}}}
% 
% 
% \affil*[2]{\orgdiv{Department of Mathematics}, \orgname{Keio University}, \orgaddress{\street{3141 Hiyoshi, Kohoku-ku}, \city{Yokohama}, \postcode{223-8522}, \state{Kanagawa}, \country{Japan}}}

\noindent {\bf Abstract: }
The Airy$_{\beta }$ random point fields ($ \beta = 1,2,4$) are random point fields emerging as the soft-edge scaling limits of eigenvalues of Gaussian random matrices. 
We construct the unlabeled diffusion reversible with respect to the Airy$_{\beta }$ random point field for each $ \beta = 1,2,4$. 
We identify the infinite-dimensional stochastic differential equations (ISDEs) describing the labeled stochastic dynamics for the unlabeled diffusion mentioned above. 
We prove the existence and pathwise uniqueness of strong solutions of these ISDEs. 
Furthermore, the solution of the ISDE is the limit of the solutions of the stochastic differential equations describing the dynamics of the $ N $-particle system in the soft-edge limit.
We thus establish the construction of the stochastic dynamics whose unlabeled dynamics are reversible with respect to the Airy random point fields. 

When $ \beta=2 $, the solution is equal to the stochastic dynamics defined by the space-time correlation functions obtained by Pr\"ahofer--Spohn, Johansson, Katori--Tanemura, and Corwin--Hammond, among others. 
We develop a new method whereby these ISDEs have unique, strong solutions. We expect that our approach is valid for other soft-edge scaling limits of stochastic dynamics arising from the random matrix theory.

\bs 


\maketitle

	\tableofcontents

%%%%%%%%%%%%%%%%%%%%%%%%%%%%%%%%%%%%%%%%%%%
\section{Introduction}\label{s:1}%%%%%%%%%%%%%%
%%%%%%%%%%%%%%%%%%%%%%%%%%%%%%%%%%%%%%%%%%%
Gaussian ensembles are introduced as random matrices 
with independent elements of Gaussian random variables, 
under the constraint that the joint distribution 
is invariant under conjugation with appropriate unitary matrices. 
The ensembles are divided into classes according to whether
their elements are real, complex, or real quaternion and 
their invariance under orthogonal 
(Gaussian orthogonal ensemble, GOE), unitary 
(Gaussian unitary ensemble, GUE), and unitary symplectic 
(Gaussian symplectic ensemble, GSE) conjugation.

The distribution of eigenvalues of the ensembles 
with size $N \times N $ is given by 
\begin{align}
m_{\beta }^{N }(d\mathbf{x} ^{N })= 
\frac{1}{Z}
\Big( \prod_{i<j}^{N } \lvert x^i-x^j \rvert ^\beta \Big) 
\exp \Big(-\frac{\beta}{4} 
\sum_{k=1}^{N}
 \lvert x^k \rvert ^2
\Big) 
d\mathbf{x} ^{N} 
\label{:11z}
,\end{align}
where $ \mathbf{x} ^{N }=(x^1,\ldots,x^{N })$ and $ d\mathbf{x} ^{N }=dx^1\cdots dx^{N } $. 
The GOE, GUE, and GSE correspond to $\beta=1, 2$, and $4$, respectively.
The probability density coincides with the Boltzmann factor normalized by the partition function $Z$ for a log-gas system at three specific values of the inverse temperature, namely $\beta=1, 2$, and $4$. 
The measures $ m_{\beta }^{N } $ still make sense for any $ 0 < \beta < \infty $ 
and are examples of log-gasses \cite{For10}. 

The eigenvalues of the GOE/GUE/GSE are spread on an interval of width of order $\sqrt{N}$. 
Let $ \nu _{\beta }^{N } $ be the distribution of 
$ N ^{-1} \sum \delta_{x^i/\sqrt{N }} $ under $ m_{\beta }^{N }(d\mathbf{x} ^{N }) $. 
Wigner's celebrated semi-circle law states that the sequence $ \{ \nu _{\beta }^{N } \} $ weakly 
converges to the nonrandom $ \sigma_{\rm semi} (x) dx $ in the space of probability measures on 
$ (\mathbb{R} , \mathcal{B}(\mathbb{R} ))$. Here, $ \sigma_{\rm semi} $ is defined as 
\begin{align} & \label{:11a}
\sigma_{\rm semi} (x)=\frac{1}{2\pi}\sqrt{4 - x^2 }{\bf 1}_{(-2,2)} (x)
.\end{align}

There exist two typical thermodynamic scalings in \eqref{:11z}, namely bulk and soft-edge scalings. 

At a point $ x $ in the interval $(-2\sqrt{N}, 2\sqrt{N})$, the density of eigenvalues is of order $\sqrt{N}$, 
and the appropriate scaling to zoom at the point is $x \mapsto x/\sqrt{N}$, which is called the bulk scaling. 
For example, at the origin the scaled random point field (RPF) $\mu_{{\rm bulk}, \beta}^N$ has the labeled density $m_{{\rm bulk}, \beta}^N$ such that
\begin{align}\label{:11b}& 
m_{\mathrm{bulk},\beta }^{N }(d\mathbf{x} ^{N })= 
\frac{1}{Z}
\Big(\prod_{i<j}^{N }\lvert x^i-x^j \rvert ^\beta \Big) 
\exp\Big(-\frac{\beta}{4N }\sum_{k=1}^{N } \lvert x^k \rvert ^2 \Big) 
d\mathbf{x} ^{N } 
.\end{align}

At the endpoint $\{2\sqrt{N}\}$ of the interval, the density of eigenvalues is of order $o(\sqrt{N})$, $N\to\infty$, and the appropriate scaling is known to be $x \mapsto xN^{-1/6} + 2\sqrt{N}$, which is called the soft edge scaling. 
This scaling yields the RPF $\mu_\beta^{N } $ with the labeled density $ \mairybeta ^{N } $ such that 
\begin{align}\label{:11c}&
\mairybeta ^{N }(d\mathbf{x} ^{N })=\frac{1}{Z}
\Big( \prod_{i<j}^{N }\lvert x^i-x^j \rvert ^\beta \Big) 
\exp\Big(-\frac{\beta}{4}\sum_{k=1}^{N } 
\left\lvert 2\sqrt{N }+\frac{x^k}{N ^{1/6} }\right\rvert ^2 \Big)
d\mathbf{x} ^{N } 
.\end{align}
The scaling $N^{-1/6}$ is related the fluctuation with the variance $N^{1/3}$. This kind of fluctuation can be seen in various stochastic models such as TASEP, last passage percolation, and the longest increasing subsequence, which are related to KPZ universality \cite{BDJ, joh.03, KPZ}.

Suppose $ \beta=2 $. The limit RPF $ \musin $ of the finite particle system \eqref{:11b} is then 
the determinantal RPF with the $ n $-point correlation functions 
$ \rho _{\mathrm{sin},2}^{n} $ defined as 
\begin{align}\label{:11d}
\rho _{\mathrm{sin},2}^{n} (\mathbf{x} ^{n})=\det [K_{\sin, 2}(x^i,x^j) ]_{i,j=1}^{n}
.\end{align}
Here, $ K_{\sin, 2} $ is a continuous kernel such that, for $ x \not= y $, 
\begin{align} & \notag %\label{:11f} 
K_{\sin, 2}(x,y) 
= \frac{\sin ( \pi(x-y) ) }{\pi (x-y)}
.\end{align}
%
%	\end{document}
% 
The limit RPF $ \muairy =\mu_{{\rm Ai},2}$ of the finite particle system \eqref{:11c} is also 
the determinantal RPF with the $ n $-point correlation functions 
 $ \rho _{2}^{n}=\rho _{\Ai ,2}^{n} $ defined as 
\begin{align} \label{:11e}&
\rho _{2}^{n} (\mathbf{x} ^{n})=\det [K_{2}(x^i,x^j) ]_{i,j=1}^{n} 
.\end{align}
Here, $K_{2} =K_{{\rm Ai}, 2}$ is the continuous kernel given by, for $ x \not= y $, 
%
%	\end{document}
\begin{align} \label{:11g}&
K_{2}(x,y)=
\frac{{\Ai }(x) {\Ai }'(y)-{\Ai }'(x) {\Ai }(y)}{x-y}
,\end{align}
where we set ${\Ai }'(x)=d {\Ai }(x)/dx$ and denote by ${\Ai }(\cdot)$ the Airy function 
%
% \end{document}
\begin{align}\label{:11h}&
{\mathrm{Ai} }(z)= \frac{1}{2 \pi} 
 \int _{\mathbb{R} } 
dk 
 e^{\sqrt{-1}(z k+k^3/3)}, \quad z\in\mathbb{R} 
.\end{align}
See \eqref{:20o} for similar expressions in terms of the quaternion determinant for $ \beta=1,4 $. 
We obtain from \eqref{:11b} the associated stochastic dynamics $ \mathbf{X}^{N }=( X^{N ,i})_{i=1}^{N }$
%
%\end{document}
according to the stochastic differential equation (SDE): 
\begin{align}\label{:11i}&
dX_t^{N ,i}=dB_t^i + 
\frac{\beta }{2} \sum_{j=1,j\not= i}^{N } 
\frac{1}{X_t^{N ,i} - X_t^{N ,j} } dt - 
\frac{\beta }{4N } X_t^{N ,i} dt 
\quad (i=1,\ldots,N )
.\end{align}
For $ \beta=1,2$, and $4$, 
\eqref{:11i} was introduced by Dyson and is referred to as the equation for Dyson's Brownian motions. Taking $ N \to \infty $, we obtain the infinite-dimensional SDE (ISDE) 
\begin{align} & \notag %\label{:11j}
dX_t^i=dB_t^i + 
\frac{\beta }{2} \sum_{j=1, \, j\not= i}^{\infty} 
\frac{1}{X_t^i - X_t^j } dt \quad (i\in\mathbb{Z})
.\end{align}
This ISDE (with $ \beta =2$) is often called Dyson's model (in infinite dimensions) and was introduced 
by Spohn at the heuristic level. 
Spohn \cite{Spo87} constructed the associated unlabeled dynamics 
as an $ L^2 $ Markovian semi-group by introducing a Dirichlet form related to $\mu_{\sin, 2}$. 
%	\end{document}

In \cite{o.tp,o.isde,o.rm}, not only for $ \beta=2$ but also for $ \beta=1,4$, 
the first author constructed the $ \mu _{\mathrm{sin},\beta } $-reversible unlabeled diffusion 
$ \mathbb{X}=\sum_{i\in\mathbb{Z}} \delta _{X^i}$. 
He proved that tagged particles $ X_t^i$ never collide with one another and that the associated 
labeled system $ \mathbf{X} =(X^i)_{i\in\mathbb{Z}}$ solves the ISDE 
%	\end{document}
\begin{align}\label{:11k}&
dX_t^i=dB_t^i + 
\frac{\beta }{2} \lim_{r\to\infty }\sum_{\lvert X_t^i - X_t^j \rvert < r, \ j\not= i}^{\infty} 
\frac{1}{X_t^i - X_t^j } dt \quad (i\in\mathbb{Z})
.\end{align}
We therefore have infinitely many, non-colliding paths describing the motions of a limit particle system 
as an $ \mathbb{R}^{\mathbb{Z}}$-valued diffusion process. 
We remark that 
because $ \mu _{\mathrm{sin},\beta } $ is translation invariant, only conditional convergence is possible
for the sum in \eqref{:11k}. The shape of the limit SDEs is quite sensitive owing to the conditional convergence. 
In fact, the ISDE describing Ginibre interacting Brownian motions, which is the two-dimensional counterpart of \eqref{:11k}, has multiple expressions \cite{o.isde}. 

In soft-edge scaling, we see from \eqref{:11c} that the $ N $-particle dynamics 
$ \mathbf{X} ^{N }=( X^{N ,i})_{i=1}^{N }$ is given by 
\begin{align}\label{:11l}&
dX_t^{N ,i}=dB_t^i + 
\frac{\beta }{2} \sum_{j=1,\, j\not= i}^{N } 
\frac{1}{X_t^{N ,i} - X_t^{N ,j} } dt 
- \frac{\beta }{2 } 
\Big( N ^{1/3} + \frac{1}{2N ^{1/3}}X_t^{N ,i} \Big) dt 
.\end{align}
No simple guess of the limit SDE is possible because of the divergence of the second and third terms on the right-hand side as $ N \to \infty $. 

The purpose of this paper is to detect and solve the limit ISDE at the soft-edge scaling. 
In fact, from \eqref{:11l}, we derive the ISDE of $ \mathbf{X} =(X^i)_{i\in\mathbb{N} } $ such that 
\begin{align}\label{:12t}&
dX_t^i=dB_t^i 
+ \frac{\beta}{2} \lim_{r\to\infty} \Big( \Big( 
 \sum_{ \lvert X_t^j \rvert <r ,\, j \not= i}\frac{1}{X_t^i -X_t^j } \Big) - \rrY \Big) dt 
\tag{{\textsf{Airy}}}
.\end{align}
Here, we set 
\begin{align}\label{:12u}&
 \rhohat (x)=\frac{1_{(-\infty , 0)}(x)}{\pi} \sqrt{-x} 
.\end{align}
We remark that the index $ i $ can be taken in $\mathbb{N}$ rather than in $\mathbb{Z}$, because we are looking at the edge of the spectrum, see \eqref{:20e} below. 
We prove that \eqref{:12t} has a family of unique strong solutions $ \mathbf{X} $ 
starting at $ \mathbf{s} = \lab (\sss ) \in \RNN $ for $ \muairybeta $-a.s.\,$ \sss =\sum_{i=1}^{\infty} \delta_{X_0^i}$ in \tref{l:23}. 
Here, $ \map{\lab }{ \SSSsip }{ \RNN}$ is the injection defined by \eqref{:22x}, where $ \SSSsip $ and $\RNN$ are the sets defined by \eqref{:20d} and \eqref{:12f}, respectively. 
 We thus have a non-equilibrium solution in the sense that it starts at a point 
 $ \mathbf{s}= \lab (\sss )$ in $ \mathbb{R} ^{\mathbb{N}} $ for $ \muA $-a.s.\,$ \sss $. 
% (see \cite{fot} for the definition of quasi-everywhere). 

We give the outline of the derivation of \eqref{:12t}. 
Considering the inverse of the soft-edge scaling 
\begin{align*}&
x \mapsto x N ^{-1/6} + 2\sqrt{N }=\sqrt{N }( x N ^{-2/3} + 2 ) 
,\end{align*}
we set 
\begin{align}\label{:12w}& 
\rhohat ^{N } (x)=N ^{1/3} \sigma_{\rm semi} ( x N ^{-2/3} + 2 ) 
.\end{align}
From \eqref{:12w}, we take $ \rhohat ^{N } $ as the first approximation of 
the one-point correlation function $ \rairybetax ^{N ,1} $ 
of the reduced Palm measure $ \muairybetax ^{N } $ of the $ N $-particle RPF $ \muairybeta ^{N } $. 
We see from \eqref{:12w} that 
\begin{align}\label{:12x} & 
 \int _{\mathbb{R}} \rhohat ^{N } (y) dy=N 
.\end{align}
The constant $ N ^{1/3}$ in \eqref{:12w} is chosen for \eqref{:12x}. 
A simple calculation shows that 
\begin{align}
\label{:12z}&
 \rhohat ^{N } (x)=%% N ^{1/3} \sigma (\frac{x}{N ^{2/3}} + 2 )=
 \frac{1_{(-4N ^{2/3} , 0)}(x)}{\pi} \sqrt{-x\left(1+\frac{x}{4N ^{2/3}}\right)} 
,\\\label{:12y}&
\lim_{N \to \infty} \rhohat ^{N } (x)=\rhohat (x)\quad \text{ uniformly on each compact set}
.\end{align}
The key point of the derivation is that 
\begin{align}\label{:12b}&
 N ^{1/3}=\int_{\mathbb{R}} \frac{\rhohat ^{N } (y) }{-y} dy 
.\end{align}
Equations \eqref{:12z} and \eqref{:12y} justify the appearance of $ \rhohat (x) $ in \eqref{:12t}. Indeed, as $ N \to \infty $, 
\begin{align*} 
&dX_t^i \approx dB_t^i 
+ \frac{\beta}{2} \Big( \Big(
 \sum_{ j=1,\, j \not= i}^{N } 
\frac{1}{X_t^i -X_t^j }\Big ) - N ^{1/3} \Big) dt 
\\ \notag & \approx 
dB_t^i 
+ \frac{\beta}{2} \lim_{r\to\infty} \Big( \Big( 
 \sum_{ \lvert X_t^j \rvert <r ,\, j \not= i}\frac{1}{X_t^i -X_t^j } \Big) 
-\int_{\verty <r}\frac{ \rhohat ^{N } (y) }{-y}dy \Big) dt
 &&\text{by }\eqref{:12b} 
\\ \notag & \approx 
dB_t^i 
+ \frac{\beta}{2} \lim_{r\to\infty} \Big( \Big( 
 \sum_{ \lvert X_t^j \rvert <r ,\, j \not= i}\frac{1}{X_t^i -X_t^j } \Big) 
-\int_{\verty <r}\frac{\rhohat (y)}{-y}dy 
\Big) dt 
 &&\text{by }\eqref{:12y} 
.\end{align*}
We thus obtain \eqref{:12t}.

Systems described by ISDEs of the form 
\begin{align} & \notag %\label{:123j}
dX_t^i=dB_t^i -\frac{1}{2}\nabla\Phi(X_t^i) 
 -\frac{1}{2}\sum_{j=1,\, j\not= i}^{\infty} 
\nabla \Psi(X_t^i, X_t^j) dt \quad (i\in\mathbb{Z})
\end{align}
are called interacting Brownian motions in infinite dimensions. 
They describe infinitely many Brownian particles moving in $ \Rd $ with free potential 
$ \Phi=\Phi (x)$ and interaction potential $ \Psi=\Psi (x,y)$. 
The study of interacting Brownian motions in infinite dimensions was initiated by 
Lang \cite{lang.1,lang.2} and continued by Fritz \cite{Fr}, the second author \cite{T2}, and others. 
In their works, interaction potentials $ \Psi (x,y)= \Psi (x-y)$ are 
smooth functions with compact supports or functions which exponentially decay at infinity. 
Such a restriction on $ \Psi $ excludes logarithmic potentials. 

%We use a general theory presented in \cite{o.isde} to solve \eqref{:12t} in \tref{l:22}. 

In \tref{l:21}, we construct an unlabeled diffusion reversible with respect to the Airy$_{\beta }$ random point field for each $ \beta = 1,2,4$. 
We also prove that none of the tagged particles $ \mathbf{X} =\{(X_t^i)_{i\in\mathbb{N}}\}_{t\in [0,\infty)}$ of the unlabeled diffusion $ \mathbb{X}_t=\sum_i \delta_{X_t^i}$ collide with one another for all $ t $. That is, 
\begin{align} & \notag %\label{:12c}
P(X_t^i\not= X_t^j \text{ for all } 0 \le t < \infty , i\not=j )=1
.\end{align}
Hence, we always label the particles in such a way that $ X_t^i > X_t^j $ for all $ i < j $. 
In \tref{l:21}, $ \mathbb{X}_0 \in \SSSsip $ a.s.. 
Then, a right-most particle exists and does not explode from \tref{l:21} \thetag{3}. 
Therefore, there exists the right-most particle denoted by $ X_t^1 $. 
% We note that a right-most particle $ X_t^i $ exists and does not explode from \tref{l:21} \thetag{3}. 
%Therefore, we denote the right-most particle denoted by $ X_t^1 $. 

The labeled process $ (X_t^i)_{i\in\mathbb{N}}$ is thus an $ \RNN $-valued process, 
where 
\begin{align}\label{:12f} \RNN =
\{ (s^i)_{i\in\mathbb{N} }\in \RN ; s^i > s^j \ (i<j) \}
.\end{align}
Here $ \RNN \subset \RN $, and we define on $ \RN $ the metric 
\begin{align}&\label{:12g}
\mathsf{dist} ( \mathbf{x},\mathbf{y}) = \sum_{ i =1}^{\infty} 2^{-i} \min \{ 1, \lvert x^i-y^i\rvert \} 
.\end{align}
 
% In \tref{l:22}, we indentify and solve the ISDE describing the unlabeled diffusion. 
% The solution at this stage has the usual meaning; that is, a pair of infinite-dimensional processes $ (\mathbf{X} , \mathbf{B})$ satisfying \eqref{:12t}. 
% Such a solution is called a weak solution. 

In \tref{l:22}, we identify and solve the ISDE describing the labeled dynamics. 
The solution at this stage is a weak solution; that is, a pair of infinite-dimensional processes $ (\mathbf{X} , \mathbf{B})$ satisfying \eqref{:12t}.

In \tref{l:23} and \tref{l:24}, we refine the results in \tref{l:22} much further stating that \eqref{:12t} has a family of unique, strong solutions. 
%
% Here, as usual, a strong solution means that $ \mathbf{X} $ is a function of the (given) Brownian motion $ \BBB $. 
%
Here, loosely speaking, a strong solution of \eqref{:12t} with \eqref{:22s} starting at ${\sbm }$ means that, for a given Brownian motion $ \mathbf{B} $, there exists a function $ \mathbf{X} =\Fs (\mathbf{B} ) $ of $ \mathbf{B} $ such that the pair $ (\mathbf{X} , \mathbf{B} )$ satisfies 
 \eqref{:12t}, \eqref{:22s}, and $\mathbf{X} _0=\sbm $. 
We also call $ \Fs $ itself a strong solution. 
See Definitions \ref{d:55}--\ref{d:53} or \cite[Section 3.1]{o-t.tail} for the precise formulation and related concepts.

The uniqueness of the solutions of the ISDE yields several important consequences, such as 
the uniqueness of Dirichlet forms \cite{o-t.core,o-t.sm,k-o-t.udf} and martingale problems, 
the convergence of $ N$-particle dynamics in \eqref{:11l} to the solution of \eqref{:12t} \cite{k-o.fpa} (see \tref{l:27}), and the dynamical universality of Airy RPFs in the sense of \cite{k-o.du}.

In \tref{l:25}, we obtain a Girsanov formula for the solution of \eqref{:12t}. 
As an application, we deduce that each tagged particle has a semi-martingale manner, 
and the local property of the particle trajectories is similar to that of Brownian motions.

For $ \beta=2$, natural infinite-dimensional stochastic dynamics 
have been constructed using the extended Airy kernel \cite{FNH99, NKT03, KNT04, KT07b, CH14}. 
The identity between our construction and the construction based on 
 space-time correlation functions for $ \beta=2 $ was proved in \cite[Theorem 2.2]{o-t.sm}, which 
is an important application of the uniqueness of the weak solutions of \eqref{:12t}. 
See \tref{l:26}.

The Airy process $ \At $ corresponds to the top particle of the dynamics defined by the extended Airy kernel. 
This process was introduced by Pr\"{a}hofer--Spohn \cite{p-spohn} and has attracted much attention. 
 From the uniqueness of the weak solutions $ \mathbf{X} = (X^i)_{i\in\mathbb{N}}$ of \eqref{:12t}, we deduce that the right-most particle $ X_t^1$ is treated with the Airy process $ \At $. 

When $ \beta=1,4$, no construction of the stochastic dynamics 
by means of space-time correlation functions is known, even for finite systems. 
Our construction based on stochastic analysis is valid for $ \beta=1,4 $ as well as $ \beta=2 $.

In \cite{o-t.tail}, we developed a general theory on the existence of a family of unique, strong solutions and the pathwise uniqueness of solutions for ISDEs of interacting Brownian motions. 
Because of the appearance of the divergent term $ N^{1/3}$ in \eqref{:11l}, it is a challenging problem to verify the assumptions necessary for the general theory. 

We solve \eqref{:12t} through an equivalent 
yet more refined representation of the ISDE 
\begin{align}\notag 
dX_t^i=dB_t^i 
& + \frac{\beta}{2} \lim_{r\to\infty} \Big( 
 \sum_{ \lvert X_t^j \rvert <r ,\, j \not= i}\frac{1}{X_t^i -X_t^j } 
- \int_{\verty < r}\frac{\rairybetaxx ^{1}(y)}{X_t^i - y}dy 
\Big) dt 
\\ \label{:12v} &
+ 
 \frac{\beta}{2} \lim_{r\to\infty} \Big( 
\int_{\verty < r}\frac{\rairybetaxx ^{1}(y)}{X_t^i - y}dy - \rrY \Big) dt 
.\end{align}
Here, $ \rairybetax ^{1} $ is the one-point correlation function of 
the reduced Palm measure $ \muairybetax $ conditioned at $ x $,
where the reduced Palm measure is defined in (\ref{:37a}). 
Note that the second term on the right-hand side is neutral in the sense of \eqref{:13a}, 
and the coefficient of 
the third term can be regarded as the $-\frac{1}{2}$ multiple of the derivative of {\it the free potential} $\Phi_\beta$ in \tref{l:47}.

The most crucial requirement for analyzing \eqref{:12v} is a global and uniform estimate of the variance of the tail part of the logarithmic derivative such that, for each $ r \in \mathbb{N}$, 
\begin{align}\label{:13a}&
\limi{s} \SupN \Supr 
E^{\muairybetax ^{N }} \Big[\Big\lvert \sum_{\lvert x- y^i \rvert \geq s} \frac{1}{x- y^i}
 - 
\int_{\lvert x-y \rvert \geq s} \frac{\rairybetax ^{N ,1} (y)}{x-y} dy \Big\rvert ^{2}\Big]=0 
.\end{align}
Here, $ \mathfrak{y}=\sum_i \delta_{y^i}$ is a configuration. We prove \eqref{:13a} in \pref{l:3S}. 

Another main requirement is to find a constant 
$ \Ct \label{;13b} $ independent of $ N \in \mathbb{N}$ such that 
\begin{align}&%\notag & 
 \lvert \rairybeta ^{N ,1} (x) - \rhohat ^{N } (x) \rvert \le \cref{;13b} 
 \Big( \frac{1}{ \xvee } +\frac{{\bf 1}(\beta\not=2)}{ (\xvee )^{1/4}} \Big) 
,\ x \in [ -2N ^{2/3}, \infty )
\label{:13b} 
.\end{align}
Here $ \rairybeta ^{N ,1}$ is the one-point correlation function of $ \muairybeta ^{N }$. We prove \eqref{:13b} in \pref{l:35}.

For \eqref{:13a} and \eqref{:13b}, we use the estimates of Plancherel--Rotach \cite{PR29} and Deift et al. \cite{DKMVZ}. The estimate \eqref{:13b} is spatially global, and we thus finely tune these estimates and stitch them together. We shall carry out many calculations to prove \eqref{:13b} and related estimates. 

We utilize the general theory for solving ISDEs and related topics developed in 
\cite{k-o.fpa,k-o.du,k-o-t.udf,k-o-t.ifc,o.dfa,Osa04,o.tp,o.isde,o.rm,o.rm2,o-t.core,o-t.sm,o-t.tail}. 
See \cite{o.ams1,o.ams2} for explanatory papers on the theory.

The remainder of the paper is organized as follows. 
In \sref{s:2}, we define the problems and state our main theorems (Theorems \ref{l:21}--\ref{l:25}). 
 In \sref{s:3}, we investigate various estimates of finite particle systems approximating Airy RPFs. 
 In \sref{s:40}, we prove the main theorems, namely Theorems \ref{l:21}--\ref{l:25}. 
 In \sref{s:8} (Appendix 1), we collect estimates of the Airy function $ \mathrm{Ai}(x)$ and related quantities. 
We devote \sref{s:Ap2} (Appendix 2) to some estimates for the normalized oscillator functions $ \psi_{N } $. 
In Sections \ref{s:Ap3} (Appendix 3), we prove \pref{l:35}, \lref{l:38}, and \lref{l:3R}.

\section{Main results} \label{s:2}
%%%%%%%%%%%%%%%%%%%%%%%%%%%%%%%%%%%%%%%%%%%
This section defines the problem and states the main theorems. 
Let $ \delta _{a} $ denote the delta measure at $ a $, and
\begin{align}&\notag 
\SSS=\{ \sss=\sum _i \delta _{s^i}\, ;\, 
\sss ( K ) < \infty \text{ for all compact sets } K \subset \mathbb{R} \} 
.\end{align}
We endow $ \SSS $ with the vague topology, under which $ \SSS $ is homeomorphic to a complete separable metric space. It is known that a complete separable metric space is homeomorphic to a Borel subset of a compact metric space. Thus, $ \SSS $ is a Lusin space in the sense of \cite[p.14]{c-f}. 
%
%By definition, a Polish space is a separable completely metrizable topological space; that is, a space homeomorphic to a complete metric space that has a countable dense subset. 

A probability measure $ \mu $ on $ (\SSS ,\mathcal{B}(\SSS ) )$ is called a random point field ( RPF) on $ \mathbb{R} $. 
A symmetric locally integrable function 
$ \map{\rho ^n }{\mathbb{R} ^n}{[0,\infty ) } $ is called 
the $ n $-point correlation function of $ \mu $ 
on $ \SSS $ with respect to the Lebesgue measure if 
\begin{align}& \label{:20a}
\int_{A_1^{k_1}\ts \cdots \ts A_m^{k_m}} 
\rho ^n (\mathbf{x} ^n) d\mathbf{x} ^n 
=\int _{\SSS } \prod _{i=1}^{m} 
\frac{\sss (A_i) ! }
{(\sss (A_i) - k_i )!} d\mu
 \end{align}
for any $m$ disjoint bounded measurable subsets 
$ A_1,\ldots,A_m \subset \mathbb{R}$ and any $m$ natural numbers 
$ k_1,\ldots,k_m $ satisfying $ k_1+\cdots + k_m=n $. 
It is known that 
$ \{ \rho ^n \}_{n \in \mathbb{N}}$ determines the measure $ \mu $ under a weak condition. In particular, determinantal RPFs generated by given kernels and reference measures are uniquely given \cite{Sos00,ST03}. 

Let $ \muairy $ be the RPF whose correlation functions are given by \eqref{:11e}. 
The RPFs $\mu_{ \beta } $ ($ \beta=1,4$) 
are defined similarly through the $ n $-point correlation functions $ \rairybeta ^{n} $ given by 
\begin{align}\label{:20o}&
\rairybeta ^{n} (\mathbf{x} ^n)=\qdet [K_{\beta }(x^i,x^j) ]_{i,j=1}^{n} 
.\end{align}
Here, $ \qdet $ denotes the quaternion determinant defined by \eqref{:32x}, 
and $ \Kairybeta $ denotes quaternion-valued kernels defined by \eqref{:31v} for $\beta=1$, and \eqref{:31w} for $\beta=4$. 
We call $\mu_{ \beta }$ the $ \mathrm{Airy}_\beta$ RPFs. 

Let 
$ \SSSsi=\{ \sss \in \SSS \, ;\, \, \sss (\{ x \} ) \le 1 \text{ for all }x \in \mathbb{R} 
 ,\, \sss (\mathbb{R} )= \infty \} $ 
 be the set of simple infinite configurations, consisting that is of an infinite number of single-point measures. 
Let $ \SSSsip $ be the subset of $ \SSSsi $ with only a finite number of point measures in 
$\mathbb{R} ^+=[0,\infty)$. That is, 
\begin{align} \label{:20d} &
\SSSsip =\{ \sss \in \SSSsi \, ;\, \, \sss (\mathbb{R} ^+ )< \infty \} 
.\end{align}
It is clear that $ \SSSsip $ is a Borel subset of $ \SSS $. Thus, $ \SSSsip $ is a Lusin space under the relative topology. It is well known that 
\begin{align}\label{:20e}&
\muairybeta (\SSSsip )=1 \quad \text{ for } \beta=1,2,4 
.\end{align}

Using \eqref{:20d} and \eqref{:20e}, for $ \muairybeta $-a.s.\ $ \sss=\sum_{i \in \mathbb{N}} \delta_{s^i}$, we label $ \{ s^i \} $ in such a way that $ s^i>s^j $ for all $ i<j $. We define the map 
$ \map{\lab }{\SSSsip }{\RNN }$ by %$ \lab (\sss )=(s^1,s^2,\ldots ) $. 
 \begin{align}\label{:22x}&
 \lab (\sss )=(s^1,s^2,s^3, \ldots )
 .\end{align}
Here, $ \RNN $ is as in \eqref{:12f}. We write $ \lab (\mathfrak{s}) = (\lab ^i (\mathfrak{s}))_{i \in \mathbb{N}} = (s^i)_{i \in \mathbb{N}} $. 
It is easy to see that $ \lab = (\lab ^i )_{i\in\mathbb{N} }$ is $ \mathcal{B}(\SSSsip )/\mathcal{B}(\RNN ) $-measurable. 
Indeed, we have 
\begin{align*}&
\lab ^1 (\sss ) = \min \{ t \in \mathbb{R} ; \sss((t,\infty )) = 0 \} 
,\\&
\lab ^i (\sss ) = \min \{ t \in \mathbb{R} ; (\sss - \sum_{j=1}^{i-1}\delta_{\lab ^j (\sss )}) ((t,\infty )) = 0 \} 
 \quad \text{ for $ i \ge 2$}
.\end{align*}
From this, we see that $ \lab ^i $ is $ \mathcal{B}(\SSSsip )/\mathcal{B}(\mathbb{R} ) $-measurable for each $ i \in \mathbb{N}$. 
This yields that $ \lab = (\lab ^i )_{i\in\mathbb{N} }$ is $ \mathcal{B}(\SSSsip )/\mathcal{B}(\RNN ) $-measurable. 
We remark that $ \lab $ is not continuous from the vague topology on $ \SSSsip $ to the topology on $ \RNN $ 
defined by the distance in \eqref{:12g}. Indeed, if 
$ \sss (k) = \delta_{k} + \sum_{i=2}^{\infty} \delta_{-i} $, then 
\begin{align*}&
\limi{k} \sss (k) = \sum_{i=2}^{\infty} \delta_{-i}, \quad 
\\&
\limi{k} \lab (\sss (k)) = \limi{k} (k, -2,-3,\ldots ) \ne (-2,-3,\ldots ) = \lab ( \sum_{i=2}^{\infty} \delta_{-i} )
.\end{align*}
We call $ \lab $ the canonical label of the Airy$_\beta$ RPF.

For a complete separable metric space $A$ with a metric $ \mathsf{dist} $, we denote by $W(A)$ the set of $A$-valued continuous paths defined on $[0,\infty)$. We set 
\begin{align*}&
\mathsf{dist}_W ( w_1 , w_2 ) = \sum_{r=1}^{\infty} 
2^{-r} \min \{ 1, \sup_{t \in [0,r]} \mathsf{dist}(w_1(t) , w_2(t) ) \} 
.\end{align*}
Then $W(A)$ is a complete separable metric space with the metric $ \mathsf{dist}_W $. 
 Clearly, $ w_n $ converges to $ w $ with respect to the metric $ \mathsf{dist}_W $ if and only if $ w_n(t)$ converges to $ w (t)$ uniformly in $ t $ with respect to the metric $ \mathsf{dist}$ on each compact interval $ [0,r]$ for all $ r \in \mathbb{N}$. 
%See \cite[ttt]{IW} for the proof in the case of $ A = \mathbb{R}^d $ with the euclidean distance. The general case can be proved similarly. 
See \cite[Chap.I, Sect.4]{IW} for some properties of the Borel $ \sigma $-field $ \mathcal{B}(W(A)) $ in case of $ A = \mathbb{R}^d$. 

Recall that $ \SSS $ is homeomorphic to a complete separable metric space with a metric $ \mathfrak{dist}$. Hence, $ W (\SSS )$ is homeomorphic to a complete separable metric space with the metric $ \mathfrak{dist}_W $. 

For the canonical label $ \map{\lab }{\SSSsip }{\RNN }$ defined by \eqref{:22x}, 
we write $ \lab = (\lab ^i)_{i \in \mathbb{N}}$, that is, 
$ \lab (\mathfrak{s}) = (\lab ^i (\mathfrak{s}))_{i \in \mathbb{N}} = (s^i)_{i \in \mathbb{N}} $ 
for $ \mathfrak{s} $ as in \eqref{:22x}. 

For $ \mathfrak{w} \in W ( \SSSsip ) $ and $ i \in \mathbb{N}$, let 
\begin{align}\notag &
\tau^1 (\mathfrak{w}) = \sup \{ t \ge 0 ; 
\sup_{u \in [0,t]} \lvert \lab^1 (\mathfrak{w}(u)) \rvert < \infty 
 \} , 
\\& \notag %\label{:20X} %\quad 
 \tau ^i (\mathfrak{w}) = \sup \{ t \in [0,\tau ^{i-1} (\mathfrak{w}) ) ; 
\sup_{u \in [0,t]} \lvert \lab^i (\mathfrak{w}(u)) \rvert < \infty 
\}
\quad \text{ for } i \ge 2 
.\end{align}
Thus, $ \tau ^1$ is the explosion time of the right-most particle. 
We have $ \tau ^{i+1} \le \tau ^i \le \tau ^{1}$ for all $ i \in \mathbb{N}$ by construction. 
We set 
\begin{align}\label{:20z}&
\tau (\mathfrak{w}) = \inf \{ \tau ^i (\mathfrak{w}) ; i \in \mathbb{N} \} 
.\end{align}
Then $ \lab^i (\mathfrak{w}(t))$ is continuous on $ [0, \tau )$ for all $ i \in \mathbb{N}$ 
because $ \mathfrak{w} \in W ( \SSSsip ) $. 
For the sake of completeness, we give a proof of this fact in \lref{l:28}. 
We note that no particle can enter from $ -\infty $ before time $ \tau (\mathfrak{w}) $. 

% Let $ \tau $ be as in \eqref{:20z}. 
Let $ W_{\mathrm{NE}}(\SSSsip ) $ be the subset of $ W(\SSSsip ) $ defined by
\begin{align}&\notag % \label{:20i}%
W_{\mathrm{NE}}(\SSSsip ) =\{ \mathfrak{w}\in W(\SSSsip ) \, ;\, 
 \tau (\mathfrak{w}) = \infty \} 
.\end{align}
For each $ \mathfrak{w} \in W_{\mathrm{NE}}(\SSSsip ) $, particles neither enter nor explode by definition. 
From \lref{l:28}, $ W_{\mathrm{NE}}(\SSSsip ) $ is a Borel measurable set. 

If $ \mathfrak{w} \in W_{\mathrm{NE}}(\SSSsip )$, then $ \tau (\mathfrak{w}) = \infty$ and 
$ t \mapsto \lab ^i (\mathfrak{w}(t))$ is continuous in $ [0,\infty)$ for all $ i \in \mathbb{N}$. 
Hence, we have 
% We write $ w^i = \lab ^i (\mathfrak{w}) $, where 
% $ w^i (t) = \lab ^i (\mathfrak{w}(t) ) $. Then we have $ w ^i \in W (\mathbb{R})$ and 
\begin{align}& \notag %
W_{\mathrm{NE}}(\SSSsip )=\{ \mathfrak{w}\in W(\SSSsip ) ; 
\mathfrak{w} (t) =\sum_{i=1}^{\infty} \delta_{ \lab ^i (\mathfrak{w}(t))} \forall t \in [0,\infty), 
 \lab ^i (\mathfrak{w}) \in W (\mathbb{R}) \} 
.\end{align}
% Let $ \mathbf{w} = (w ^i)_{i\in \mathbb{N}}$. 
% Note that $ w ^i \in W (\mathbb{R})$ for each $ i \in \mathbb{N}$. 
% Hence, $ \mathbf{w} \in W (\RN )$. 
Thus, we relate $ \mathfrak{w} \in W_{\mathrm{NE}} ( \SSSsi ) $ 
 to the labeled path $ \{ (\lab ^i (\mathfrak{w}(t) ) )_{i\in\mathbb{N}}\}_{t\in[0,\infty)} \in W(\mathbb{R} ^{\mathbb{N}} )$. 
% by $ \mathbf{w} = (\lab ^i (\mathfrak{w}) )_{i\in\mathbb{N}}$. 
We denote this correspondence by 
\begin{align}\label{:20j}
\map{\lpath }{W_{\mathrm{NE}}(\SSSsip ) }{ W (\RNN )} , \quad 
\mathfrak{w} \mapsto 
\lpath (\mathfrak{w}) 
.\end{align}
We note that $ t \mapsto \lpath (\mathfrak{w})(t) \in \RNN $ is continuous in $ [0,\infty ) $. 
% As we see above, $ t \mapsto \lab ^i (\mathfrak{w}(t)) $ is continuous in $ [0,\infty ) $ for all $ i \in \mathbb{N}$. Hence, $ t \mapsto \lpath (\mathfrak{w})(t) \in \RNN $ is continuous in $ [0,\infty ) $. 
% because $ \lab (\mathfrak{w}(t)) = (\lab ^i (\mathfrak{w}(t)))_{i\in \mathbb{N}} $, 
% $ \lab ^i (\mathfrak{w}(t))$ is continuous for each $ i \in \mathbb{N}$ in $ t \in [0,\tau (\mathfrak{w}))$, 
% and $ \tau (\mathfrak{w}) = \infty $. 
We prove that $ \lpath $ is $ \mathcal{B}(W_{\mathrm{NE}}(\SSSsip ) ) / \mathcal{B}(W (\RNN )) $-measurable in \lref{l:29}. We call $ \lpath $ a path-label map. 
% We can construct a similar map $ \lpath $ for general RPFs supported on $ \SSSsi $. 
% For the Airy RPFs $ \muairybeta $, the construction is much simpler than the general case because $ \muairybeta $ is supported on $ \SSSsip $. 

For $ A \subset \mathbb{R}$, we set $ \pi_{A}(\sss )=\sss (\cdot \cap A)$. We say a function $ f : \SSS \to \mathbb{R} $ is local if $ f $ is $\sigma[\pi_K]$-measurable for some compact set $ K $. 
For a local function $ f $, there is a function 
$\check{f}=\check{f}_O$ on $ {\cup _{k=0}^\infty O^k}$ 
with bounded open set $O\supset K$, such that, for each $ k \in \mathbb{N}$, the restriction $\check{f} $ on $O^k$ is symmetric in $ (s^1,\dots,s^k) $ and $ f(\sss)=\check{f}(s^1,\dots,s^k) $ 
if $\pi_O(\sss)=\sum_{i=1}^k \delta_{s^i} $.

We say a local function $\map{f}{\SSS }{\mathbb{R}}$ is smooth if $\check{f}$ is smooth. 
Let $\di $ be the set of all local, smooth functions on $\SSS $. 
For $f,g\in \di $, we set $\mathbb{D} [f,g] : \SSS \to \mathbb{R} $ according to 
\begin{align}& \label{:20b}
\mathbb{D} [f,g] (\sss )= \frac{1}{2}\sum_{i} 
\frac{\partial _i\check{f}(\sbm )}{\partial s^i} 
\frac{\partial _i\check{g}(\sbm )}{\partial s^i}
,\end{align}
where $ \sss=\sum_{i}\delta _{s^i} $ and $ \sbm =(s^i)$. 
Let $ (\Emu , \di ^{\mu }) $ be the bilinear form defined as 
\begin{align} & \notag %\label{:20c}
\Emu (f,g)=\int_{\SSS } \mathbb{D}[f,g] d\mu ,\quad 
 \di ^{\mu }= \{ f \in \di \, ;\, \mathcal{E}^{\mu}(f,f)< \infty ,\ f \in L^2(\SSS ,\mu ) \}
.\end{align}

We say that a non-negative symmetric bilinear form $ (\Emu , \di ^{\mu }) $ is closable on $ L^2(\SSS ,\mu ) $ if $ \limi{n}\Emu (f_n , f_n ) = 0 $ for any $ \Emu $-Cauchy sequence 
$ f_n\in \di ^{\mu }$ such that $ \limi{n} \| f_n \| _{L^2(\SSS ,\mu ) } = 0 $. 
If $ (\Emu , \di ^{\mu }) $ is closable on $ L^2(\SSS ,\mu ) $, then there exists 
a closed extension of $ (\Emu , \di ^{\mu }) $. 
The smallest closed extension $ (\Emu , \domu ) $ of $ (\Emu , \di ^{\mu }) $ is called the closure of $ (\Emu , \di ^{\mu }) $. 
We refer to \cite{fot} for detail. 

We say a closed non-negative symmetric bilinear form $ (\Emu , \domu ) $ on $ \Lmu $ is 
a Dirichlet form if every $ u \in \domu $ satisfies 
\begin{align*}& 
v := \min \{ 1 , \max\{ 0 , u \} \} \in \domu \text{ and }
\Emu (v,v) \le \Emu (u,u)
.\end{align*}
In general, a Dirichlet form is not necessarily symmetric. In the present paper, a Dirichlet form means a symmetric Dirichlet form. 

We recall the concept of the $ 1 $-capacity \cite{c-f}. 
Let $ \mathcal{O} $ be the family of all open subsets of $ \SSS $. 
For $ \aAA \in \mathcal{O} $, let 
$ \mathcal{L}_{\aAA ,1}=\{ f \in \domu \, ;\, f \ge 1 \ \mu \text{-a.e.\,on }\aAA \} $. 
We set $ \mathcal{O}_0=\{ \aAA \in \mathcal{O} \, ;\, \mathcal{L}_{\aAA ,1} \ne \emptyset \} $. 
Let $ \Emu _1 = \Emu + (\cdot,*)_{\Lmu }$. 
For an open set $ \aAA \in \mathcal{O} $, we set 
%For an open set $ \aAA $, we define
\begin{align*}&
\capamu (\aAA ) = 
\begin{cases}
\inf \{ \Emu _1(f,f); f \in \mathcal{L}_{\aAA ,1}\}, & \aAA \in \mathcal{O}_0 
\\
\infty , & \aAA \notin \mathcal{O}_0 
.\end{cases}
\end{align*}
For an arbitrary set $ \bBB \subset \SSS $, we set 
\begin{align}& \label{:20p}
\capamu (\bBB ) = \inf \{ \capamu (\aAA ) ; \aAA \in \mathcal{O},\, \aAA \supset \bBB \} 
.\end{align}
An increasing sequence of closed sets $ \{ \fFF_{k} \} $ is called an $ \Emu $-nest if, 
for any compact set $ \kKK $, $ \limi{k} \capamu (\kKK \backslash \fFF_{k}) = 0 $. 
%A nest $ \{ \fFF_{k} \} $ is called a compact nest if each $ \fFF_{k} $ is a compact set. 
%
A function $ f $ is called $ \Emu $-quasi-continuous if for any $ \epsilon >0 $, there exists an open set $ \oOO $ with $ \capamu (\oOO ) < \epsilon $ such that $ f \vert _{\SSS \backslash \oOO }$ is finite and continuous. 
We call a subset $ \mathfrak{N} \subset \SSS $ an $ \Emu $-polar set if 
there exists an $ \Emu $-nest $ \{ \fFF_{k} \} $ such that 
$ \mathfrak{N} \subset \cap_k ( \SSS \backslash \fFF_{k} )$. 
% \cite{c-f}. 
% See \cite{c-f} for the definitions of \lq\lq nest'', \lq\lq quasi-continuous'', and \lq\lq a polar set''. 
Although these concepts depend on the Dirichlet form $ (\Emu , \domu )$ on $ \Lmu $, we often suppress it from the notation.

A Dirichlet form $ (\Emu , \domu ) $ on $ L^2(\SSS ,\mu ) $ is called quasi-regular if: \\ 
\thetag{1} there exists an $ \Emu $-nest $ \{ \fFF_{k} , k \ge 1 \} $ consisting of compact sets; 
\\\thetag{2} 
there exists $ \Emu + ( \cdot , *)_{ L^2(\SSS ,\mu ) }$-dense subset of $ \domu $ 
whose elements have $ \Emu $-quasi-continuous $ \mu $-version; 
\\\thetag{3} 
there exists $ \{ f_k , k \ge 1 \} \subset \domu $ having $ \Emu $-quasi-continuous $ \mu $-versions $ \{ \tilde{f}_k , k \ge 1 \} \subset \domu $ and 
an $ \Emu $-polar set $ \mathfrak{N} \subset \SSS $ such that 
$ \{ \tilde{f}_k , k \ge 1 \} \subset \domu $ 
separates the points on $ \SSS \backslash \mathfrak{N} $. 

A symmetric form $ (\Emu , \domu ) $ on $ L^2(\SSS ,\mu ) $ is called strongly local if $ \Emu ( f , g ) = 0 $ for any $ f , g \in \domu $ \ such that $ f $ is constant on a neighborhood of the support of $ g $. 
This definition is natural but not invariant under quasi-homeomorphisms. 
Hence we quote the second definition of strong locality from \cite[Remark 2.4.4]{c-f}. A symmetric form $ (\Emu , \domu ) $ on $ L^2(\SSS ,\mu ) $ is strongly local if $ \Emu (f,g) = 0 $ whenever $ f,g \in \domu $ with $ (f-c) g = 0 $ $ \mu $-a.e.\,on $ \SSS $ for some constant $ c $. 
The two definitions are equivalent when the symmetric form is dominated by a regular Dirichlet form \cite[Theorem 2.4.3]{c-f}. 
The point is that the second definition is invariant under quasi-homeomorphisms. 
We adopt the second definition in the present paper. 
It is noteworthy that there is a concept of the generalized strong locality of quasi-regular Dirichlet forms introduced in \cite[Theorem 5.1]{kuwae.98}, which is useful to apply the Beurling-Deny decomposition for quasi-regular Dirichlet forms. 

A diffusion process is a family of Markov processes with continuous sample paths and has the strong Markov property. We say a diffusion process is conservative if it has an infinite lifetime. 
% See \sref{s:D} for the definition of a conservative diffusion (also \cite{cf,fot}). 
%
In general, a diffusion process has a cemetery point. We always consider a conservative diffusion process; thus, we will formulate diffusion processes without a cemetery point below. 

We recall that $ \SSS $ is a Lusin space, that is, a topological space homeomorphic to a Borel subset of a compact metric space. 
Let $ (\Omega , \mathcal{M} )$ be a measurable space. 
Let $ \XXX = \{ \XXX _t \} $ be an $ \SSS $-valued stochastic process on $ (\Omega , \mathcal{M} )$. 

A family $ \{ \mathcal{M}_t; t \ge 0 \} $ of sub-$ \sigma $-fields of $ \mathcal{M} $ is called a filtration if $ \mathcal{M}_s \subset \mathcal{M}_t $ for $ 0 \le s \le t $. 
A filtration $ \{ \mathcal{M}_{t} \} $ is called an admissible filtration for a stochastic process 
$ \XXX =\{ \XXX _t\} $ on $ ( \Omega , \mathcal{M} ) $ if $ \XXX _t $ is $ \mathcal{M}_t/\BS $-measurable for all $ t $. We call $ \{ \mathcal{M}_t\} $ a right continuous filtration if 
$ \mathcal{M}_t = \bigcap_{\epsilon > 0 }\mathcal{M}_{t+\epsilon} $. 

We now introduce the concept of a conservative diffusion process on $\SSS$. 
\begin{definition}\label{d:A91}
A quadruplet $ \7 $ is called a conservative diffusion process on $\SBS $ if it satisfies \As{D1}--\As{D5}. 
\\
\As{D1} For each $ \mathfrak{s} \in \SSS $, $ (\Omega , \mathcal{M}, \{ \XXX _t \}_{t\in[0,\infty)}, \PPPs ) $ is a stochastic process with state space $ \SBS $ and time parameter set $ [0,\infty)$. 
\\
\As{D2} For each $ t \ge 0 $ and $ B \in \BS $, $ \PPPs (\XXX _t \in B)$ is $ \BS $-measurable as a function of $ \sss \in \SSS $. 
\\
\As{D3} For each $ \omega \in \Omega $, the sample path $ t \mapsto \XXX _t(\omega ) \in \SSS $ is continuous on $ [0,\infty )$. 
\\
\As{D4} $ \PPPs (\XXX _0 = \sss ) = 1$ for every $ \sss \in \SSS $. 
\\
\As{D5} A right continuous admissible filtration $ \{ \mathcal{M}_t \} $ exists such that 
% the following hold: 
\begin{align}&\notag %\label{:A921}&
\PPP _{\mu } (\{ \XXX _{\sigma + s }\in B \} \vert \mathcal{M}_{\sigma } ) = 
\PPP _{\XXX _{\sigma }} (\XXX _s \in B ) \quad \PPP _{\mu}\text{-a.s}
\end{align}
for any $ \{ \mathcal{M}_t \} $-stopping time $ \sigma $, probability measure $ \mu $ on $ \SBS $, $ s \ge 0 $, and $ B \in \BS $. 
\end{definition}
We use the term diffusion to represent a diffusion process. 
We omit $ (\Omega , \mathcal{M} )$ from $ \7 $ and write $ \XP $, where $ \mathbb{X}= \{ \XXX _t \}_{t\in[0,\infty)}$. 

For a Dirichlet form $ (\mathcal{E}, \domu )$ on $ \Lmu $, we say a conservative diffusion $ \XP $ is associated with $ (\mathcal{E}, \domu )$ on $ \Lmu $ if $ \mathbb{E}_{ \sss }[f ( \mathbb{X}_t)]=T_t f ( \sss )$ $ \mu $-a.e $ \sss \in \SSS $ for any $ f \in \Lmu $ and for all $ t \in [0,\infty ) $. 
Here, $ \mathbb{E}_{\sss }$ is the expectation with respect to $ \PPPs $ and $ T_t $ is the Markovian $ L^2 $-semi-group associated with $ (\mathcal{E}, \domu )$ on $ \Lmu $. 
Furthermore, $ \XP $ is said to be properly associated with $ (\mathcal{E}, \domu )$ on $ \Lmu $ if, in addition, $ \mathbb{E}_{ \sss }[f( \mathbb{X}_t)]$ is a quasi-continuous version of $ T_t f ( \sss )$ for any $ f \in \Lmu $ and $ t \in [0,\infty)$ \cite{c-f}. 

We say $ \XP $ is (resp.\,properly) associated with a closable form $ (\mathcal{E}, \dom _0 )$ on $ \Lmu $ if $ \XP $ is (resp.\,properly) associated with the closure of $ (\mathcal{E}, \dom _0 )$ on $ \Lmu $. 
%
% We can extend the state space $ \SSSzz $ to $ \SSS $ by taking $ \mathbb{P}_{\sss } (\XXX _t= \sss \text{ for all }t)$ for $ \sss \in \SSSzz ^c$. 

For a quasi-regular Dirichlet form $ (\Emu , \domu ) $ on $ L^2(\SSS ,\mu ) $, an associated conservative diffusion is always properly associated, and is unique q.e.\,in the sense that, if two conservative diffusions $ \XP $ and $ (\mathbb{X}', \{ \PPPs ' \}_{\sss \in \SSS }) $ are associated with $ (\Emu , \domu ) $ on $ L^2(\SSS ,\mu ) $, then $ \PPPs = \PPPs '$ for all $ \sss \in \SSS \backslash \mathfrak{N} $ for some $\mathfrak{N}$ such that $ \capa ^{\mu} (\mathfrak{N} ) = 0 $.

A diffusion $\XP $ is called $ \mu $-reversible if $\XP $ is $ \mu $-symmetric and $ \mu $ is an invariant probability measure of $ \XP $. We note that a $ \mu $-reversible diffusion is conservative. 

Let $ \Lambda $ denote the Poisson RPF whose intensity is the Lebesgue measure. 
If $ \mu=\Lambda $, then the bilinear form 
$(\mathcal{E}^{\Lambda },\di ^{\Lambda })$ 
is closable on $ L^2 (\SSS ,\Lambda )$, and the closure of $(\mathcal{E}^{\Lambda },\di ^{\Lambda })$ is a quasi-regular Dirichlet form. The associated diffusion 
$ \mathbb{B}_t^{\sss }=\sum_{i\in\mathbb{N}}\delta_{B_t^i+s^i}$ is the 
$ \SSS $-valued Brownian motion starting at $ \sss=\sum_i \delta_{s^i}$. 
Here, $\{ B^i \}_{i\in\mathbb{N}}$ denotes independent copies of 
the one-dimensional standard Brownian motion \cite{o.dfa}. 
It is thus natural to ask, if we replace $ \Lambda $ by $ \muairybeta $, whether the forms 
$ (\mathcal{E}^{\muairybeta },\di ^{\muairybeta })$ are 
 still closable on $ L^2(\SSS ,\muairybeta ) $ and the associated diffusions exist. 

\begin{theorem}	\label{l:21}
 Assume $ \beta=1,2,4 $. The following then hold. \\
\thetag{1} 
The bilinear form $ (\mathcal{E}^{\muairybeta },\di ^{\muairybeta })$ is closable on $ L^2(\SSS ,\muairybeta ) $. 
\\
\thetag{2} 
The closure $ (\mathcal{E}^{\muairybeta },\dom ^{\muairybeta })$ of 
$ (\mathcal{E}^{\muairybeta },\di ^{\muairybeta }) $ on $ L^2(\SSS ,\muairybeta ) $ 
is a strongly local, quasi-regular Dirichlet form. 
There exists a $\muairybeta $-reversible diffusion $(\mathbb{X}, \{\mathbb{P}_{\sss }\}_{\sss \in \SSS } )$ 
properly associated with 
$ (\mathcal{E}^{\muairybeta },\dom ^{\muairybeta })$ on $ L^2(\SSS ,\muairybeta ) $. 
\\
\thetag{3} 
Let $\mathbb{P}_{\muairybeta} =\int_{\SSS}\mathbb{P}_{\sss} \muairybeta (d\sss ) $. Then 
$ \mathbb{P}_{\muairybeta } (\mathbb{X}\in \Wnex (\SSSsip ) )=1 $. 
 \end{theorem}

\begin{remark}\label{r:21} \thetag{1}
For $ \XP = \7 $ in \tref{l:21}, we can take $ \Omega = \Wnex (\SSSsip ) $, $ \mathcal{B}( \Wnex (\SSSsip ) ) \subset \mathcal{M} $, and $ \mathbb{X}_t(\mathfrak{w}) = \mathfrak{w}(t) $. 
Furthermore, $ \{\mathcal{M}_t\} $ is an admissible filtration satisfying 
$ \mathcal{B}(W_{\mathrm{NE}}(\SSSsip ) )_t \subset \mathcal{M}_t $. 
Here
\begin{align}\label{:21z}&
 \mathcal{B}(W_{\mathrm{NE}}(\SSSsip ) )_t = \sigma [\mathfrak{w}(u); u \le t ]
%, \mathfrak{w} \in \mathcal{B}(W_{\mathrm{NE}}(\SSSsip ) ) ] 
\end{align}
and $ \map{\mathfrak{w}(u)}{W_{\mathrm{NE}}(\SSSsip )}{\SSSsip }$ is such that 
$ \mathfrak{w} \mapsto \mathfrak{w}(u)$. 

\noindent 
\thetag{2} 
Let $ \lab $ be the canonical label defined by \eqref{:22x}. Let $ \lpath $ be as in \eqref{:20j}. 
% Recall that $ \lpath (\mathfrak{w}) (0) = \lab (\mathfrak{w}(0))$. 
Let 
\begin{align}\label{:22y}&
\text{$ \mathbf{P}_{\sbm }=\mathbb{P}_{\sss} \circ \lpath ^{-1}$, \quad $ \sbm=\lab (\sss ) $}
.\end{align}
% Let $ \mathbf{P}_{\mu_{\beta}\circ\mathfrak{l}^{-1}}=\mathbb{P}_{\muairybeta } \circ \lpath ^{-1} $ and $ \mathbf{X} =\lpath (\mathbb{X})$. 
Then, from \tref{l:21} \thetag{3}, Fubini's theorem, and \eqref{:22y} 
\begin{align}\label{:22a}
% \mathbf{P}_{\mu_{\beta}\circ\mathfrak{l}^{-1}} (W(\RNN ))=
 \mathbf{P}_{\sbm } (W(\RNN ))
=1 
\quad \text{ for $ \mu_{\beta}\circ\mathfrak{l}^{-1}$-a.s.\,$ \sbm $}
.\end{align}
We easily see that $ (\Omega , \mathcal{M}, \{ \lab (\XXX _t) \}_{t\in[0,\infty)},\{ \PPPs \}_{\sbm \in \RNN } ) $ is a conservative diffusion process with state space $ \RNN $ 
because $ \lab = (\lab ^i)_{i\in\mathbb{N}}$ is injective on $ \SSSsip $ and $t \mapsto \lpath (\mathbb{X}) (t)$ is continuous by \eqref{:22y} and \eqref{:22a}, and 
\begin{align*}&
 (\lab ^i (\mathbb{X}_t))_{i\in\mathbb{N}} = 
\lab (\mathbb{X}_t) = \lpath (\mathbb{X}) (t)
.\end{align*}
From \lref{l:29}, we see that $ \{\mathcal{M}_t\} $ is an admissible filtration for $ \lpath \circ \mathbb{X} (t) = \lab (\XXX _t) $. 
Furthermore, $ (\mathbf{X}, \{ \PPPs \}_{\sbm \in \RN }) $ is $ \muA \circ \lab ^{-1}$-reversible because $ \XP $ is $ \muA $-reversible and 
\begin{align*}&
\int _{\SSSsip }\mathbb{E}_{\sss }[f (\lab (\mathbb{X}_t)) ] g (\lab (\sss )) \muA (d\sss ) = 
\int_{\RNN } \mathbb{E}_{\sbm } [f (\mathbf{X}_t)] g (\sbm ) \muA \circ \lab ^{-1} (d\sbm )
.\end{align*}
\end{remark}

We next specify the ISDE for $( \mathbf{X} , \mathbf{P}_{\sbm }) $. % for $ \muA \circ \lab ^{-1} $-a.s.\,$ \mathbf{s}$. 
We consider ISDE \eqref{:12t}, $ i \in \mathbb{N}$, 
\begin{align}\tag{\ref{:12t}}
& dX_t^i=dB_t^i 
+ \frac{\beta}{2} \lim_{r\to\infty} \Big( 
\sum_{ \lvert X_t^j \rvert <r ,\, j \not= i}\frac{1}{X_t^i -X_t^j } 
-\int_{\vertx <r}\frac{\rhohat (x)}{-x}dx \Big) dt 
%\quad ( i \in \mathbb{N})
\end{align}
under the condition that there exist subsets $\Hbm $ and 
$\mathbf{S}_{\rm sde}$ such that 
\begin{align} & \notag %\label{:22p}
\Hbm \subset\mathbf{S}_{\rm sde} \subset \RNN 
, \quad %\\& \notag % \label{:22q}
\mu_{\beta}\circ\mathfrak{l}^{-1}(\Hbm ) =1 
.\end{align}
Here, $ \Hbm $ is the set of the initial starting points of solutions and $ \mathbf{S}_{\rm sde} $ is the set where the drift of the ISDE is well defined. 
We impose on $ \mathbf{X}$ to satisfy 
%We consider the condition on a solution $ \mathbf{X} $ of \thetag{Airy} such that 
\begin{align} \label{:22s}
& \mathbf{X} _0 \in \Hbm , \quad \mathbf{X} \in \WTsol 
.\end{align}
The solution $ \XB $ in \tref{l:22} below is a weak solution of \eqref{:12t} for some 
 $\mathbf{S}_{\rm sde}$ in the sense of Definition 3.1 in \cite[p. 1151]{o-t.tail}: 
\begin{definition}\label{d:21}
A weak solution of \eqref{:12t} is 
 an $\mathbb{R} ^{\mathbb{N}} \ts \mathbb{R} ^{\mathbb{N}}$-valued stochastic process $ \XB $ 
defined on a probability space $ (\Omega , \mathcal{F}, P )$ 
with a filtration $ \{ \mathcal{F}_t \}_{t \ge 0 } $ such that

\noindent 
\thetag{i} $ \mathbf{X}=(X^i)_{i \in \mathbb{N}} $ is an $ \SSSsde $-valued continuous process 
%
%Furthermore, $ \mathbf{X}$ is 
adapted to $ \{ \mathcal{F}_t \}_{t \ge 0 }$, that is, $ \mathbf{X}_t $ is 
$ \mathcal{F}_t / \Bt $-measurable for each $ \zti $, where 
\begin{align}& \label{:50k}
\Bt = \sigma [ \mathbf{w} _s ; 0\le s \le t ,\, \mathbf{w} \in W (\RN ) ] 
.\end{align}
\thetag{ii} $ \mathbf{B} = (B^i)_{i \in \mathbb{N}}$ is an $ \RN $-valued $ \{ \mathcal{F}_t \} $-Brownian motion with $ \mathbf{B}_0 = \mathbf{0}$, 
\\
\thetag{iii} 
the family of measurable $ \{ \mathcal{F}_t \}_{t \ge 0 } $-adapted process defined by 
\begin{align*}&
\Psi ^i (t , \omega ) :=
 \frac{\beta}{2} \lim_{r\to\infty} \Big( 
\sum_{ \lvert X _t^j (\omega ) \rvert <r ,\, j \not= i}\frac{1}{X_t^i (\omega ) -X_t^j (\omega ) } 
-\int_{\vertx <r}\frac{\rhohat (x)}{-x}dx \Big) 
%\beta \limi{s} 
% \Big( \sum_{\lvert x-y^i\rvert <s}\frac{1}{x-y^i} - \int_{\verty <s} \frac{\rhohat (y)}{-y}dy\Big) 
\end{align*}
belongs to $ \mathcal{L}^1 $. 
Here $ \mathcal{L}^{1} $ is the set of all measurable $ \{ \mathcal{F}_t \}_{t \ge 0 } $-adapted 
processes $ \alpha $ such that $ E[ \int_0^T \lvert \alpha (t,\omega) \rvert dt ] < \infty $ for all $ T $. 
Here we can and do take a predictable version of $\Psi ^i (t , \omega ) $ (see pp. 45-46 in \cite{IW}). 
\\
\thetag{iv} with probability one, the process $ \XB $ satisfies for all $ t $ and $ i \in \mathbb{N}$
\begin{align}\notag &% \label{:50m}&
X_t^i - X_0^i = B_t^i 
 + 
\int_0^t 
 \frac{\beta}{2} \lim_{r\to\infty} \Big( 
\sum_{ \lvert X _u^j \rvert <r ,\, j \not= i}\frac{1}{X_u^i -X_u^j } 
-\int_{\vertx <r}\frac{\rhohat (x)}{-x}dx \Big) 
du 
.\end{align}
\end{definition}

We say $ \mathbf{X}$ is a weak solution of \eqref{:12t} if $ (\mathbf{X},\mathbf{B})$ is a weak solution of \eqref{:12t} for some Brownian motion $ \mathbf{B}$. 

In the next theorem, we specify the ISDE describing $ \mathbf{X} =\lpath (\mathbb{X})$ and 
prove that the distribution of the right-most particle coincides with the $ \beta $-Tracy--Widom distribution for 
 $ \beta=1,2,4$. 
We refer to \cite[pp. 92-94]{AGZ10} for the definition of the Tracy--Widom distributions. 
\begin{theorem}	\label{l:22} 
Assume $ \beta=1,2,4 $. Let $ \mathbf{X} =\lpath (\mathbb{X})$. 
Let $ \mathbb{P}_{\sss } $ be as in \tref{l:21}. 
\\ \thetag{1} 
For $ \muairybeta $-a.s.\,$\sss $, $ \mathbf{X} =( X ^i)_{i\in\mathbb{N}} $ under $ \mathbb{P}_{\sss } $ satisfies $ \mathbf{X} _0=\lab ( \sss ) $ and 
\begin{align}\tag{\ref{:12t}}
& dX_t^i=dB_t^i 
+ \frac{\beta}{2} \lim_{r\to\infty} \Big( 
\sum_{ \lvert X_t^j \rvert <r ,\, j \not= i}\frac{1}{X_t^i -X_t^j } 
-\int_{\vertx <r}\frac{\rhohat (x)}{-x}dx \Big) dt 
\quad ( i \in \mathbb{N})
\end{align}
with an $ \mathbb{R} ^{\mathbb{N}} $-valued standard Brownian motion $ \BBB =\{B^i\}_{i\in\mathbb{N} }$. 
\\
\thetag{2} The distribution of $ X_t^1 $ under 
$ \mathbf{P}_{\muairybeta \circ \lab ^{-1}}$ is the $ \beta $-Tracy--Widom distribution. 
\end{theorem}

We next study whether the process $ \mathbf{X} =\lpath (\mathbb{X})$ is a family of unique strong solutions starting at each $ \sbm \in \Hbm $. 
For this purpose, we introduce a sequence of conditions and concepts. 
We consider the following assumptions. 

%Let $\overline{\mathcal{B}(\mathbb{R} ^{\mathbb{N}} ) } $ be the completion of $\mathcal{B}(\mathbb{R} ^{\mathbb{N}} )$ with respect to $P\circ \mathbf{X} _0^{-1}$. 

\smallskip
\smallskip

\noindent 
\As{AC} \quad \ $ \mathbb{P}_{\muairybeta} ( \mathbb{X}_t \in \cdot) $ is absolutely continuous with respect to $ \muairybeta $ for all $ t >0$. 

\noindent 
\As{SIN$ _{\mathbf{A}}$} \ 
$ \mathbb{P}_{\muairybeta} (\lpath (\mathbb{X}) \in W(\RNN ))=1 $.

\smallskip
\noindent 
We further consider a milder condition than \As{SIN$ _{\mathbf{A}}$}: \\
\As{SIN} \quad 
$ \mathbb{P}_{\muairybeta } ( \mathbb{X}\in \Wnex (\SSSsi))=1$.

\begin{remark}\label{r:23a} 
\thetag{1}
If the solution $ \mathbf{X} =\lpath (\mathbb{X})$ comes from \tref{l:22}, then 
\As{SIN$ _{\mathbf{A}}$} and \As{AC} are clearly satisfied. 
Indeed, \As{SIN$ _{\mathbf{A}}$} follows from \eqref{:22a} and 
\As{AC} is obvious because $ \mathbb{P}_{\muairybeta} ( \mathbb{X}_t \in \cdot) $ is associated with the symmetric Dirichlet form with invariant probability measure $ \muairybeta $ and thus 
$ \mathbb{P}_{\muairybeta} ( \mathbb{X}_t \in \cdot)= \muairybeta $. 
\\\thetag{2} The name {\bf SIN} comes from {\bf n}on explosion and {\bf n}on entering paths with values 
{\bf i}nfinite configurations with {\bf s}imple points. 
Furthermore, {\bf SIN}$ _{\mathbf{A}}$ implies the state space is related to the {\bf A}iry RPFs as we see in \eqref{:20e}. 
Clearly, \As{SIN$ _{\mathbf{A}}$} implies \As{SIN}. 
The name {\bf AC} comes from the {\bf a}bsolute {\bf c}ontinuity. 
\\
\thetag{3} In Theorems \ref{l:23} and \ref{l:24}, we shall prove the existence of strong solutions satisfying 
\As{SIN$ _{\mathbf{A}}$} and the uniqueness of weak solutions under the constraints of \As{SIN}. 
\end{remark}

Let $ \XB $ be a weak solution of \eqref{:12t} with \eqref{:22s} under $ \Ps $. 
Let $ \Xmstar =(X^{i})_{i> m}^{\infty}$. 
Let $\sbm =(s^i)_{i=1}^{\infty} \in \Hbm $ and $ m \in \mathbb{N}$. 
Let $ \sbm ^m $ be the first $ m $-components of $ \sbm $. 
We consider the finite-dimensional SDE \eqref{:22u} of $ \mathbf{Y} ^{m}=(Y^{m,i})_{i=1}^m $ 
such that 
\begin{align}\notag 
 dY_t^{1,1} &= dB_t^1+ 
 \frac{\beta}{2} \lim_{r\to\infty} \Big( 
\sum_{ \lvert X_t^j \rvert <r ,\, j > 1 }\frac{1}{Y_t^{1,1} -X_t^j } 
-\int_{\vertx <r}\frac{\rhohat (x)}{-x}dx \Big) dt 
,\\ \notag 
 dY_t^{m,i} &= dB_t^i+ \frac{\beta}{2} 
 \Big( 
\sum_{ j \not= i}^m 
\frac{1}{Y_t^{m,i} - Y_t^{m,j} }\Big) dt 
\\ & \notag \ 
 + \frac{\beta}{2} \lim_{r\to\infty} \Big( 
\sum_{ \lvert X_t^j \rvert <r ,\, j > m }\frac{1}{Y_t^{m,i} -X_t^j } 
-\int_{\vertx <r}\frac{\rhohat (x)}{-x}dx \Big) dt 
, \quad m \ge 2 ,\\ \label{:22u} 
\mathbf{Y} _0^{m}=&\sbm ^m , \quad 
(\mathbf{Y} _t^{m}, \Xmstar ) \in \WTsol 
.\end{align}

We set 
$ \widetilde{P}^m=\Ps \circ (\BBB ^m , \Xmstar )^{-1} $. 
Let $ \mathcal{B}_t (W_0(\Rm ) \ts W (\mathbb{R} ^{\mathbb{N}} )) $ be 
the sub $ \sigma $-field of $ \mathcal{B}(W_0(\Rm )\ts W (\mathbb{R} ^{\mathbb{N}} )) $ 
generated by $ \{(\mathbf{v}_s,\mathbf{w}_s); 0 \le s \le t \} $. 
We set $ \mathcal{B}_t (W (\Rm ) ) $ similarly. 
Let $ \mathcal{C}^m $ and $ \mathcal{C}_t^m $ be the completion of 
$ \mathcal{B}(W_0(\Rm )\ts W (\mathbb{R} ^{\mathbb{N}} )) $ and 
$ \mathcal{B}_t (W_0(\Rm )\ts W (\mathbb{R} ^{\mathbb{N}} )) $ 
with respect to $ \widetilde{P}^m $, respectively. 

We call $\mathbf{Y} ^m$ a {\em strong solution} of \eqref{:22u} for $ \XB $ under $ \Ps $ if 
$ (\mathbf{Y} ^m , \BBB ^m , \Xmstar ) $ satisfies \eqref{:22u} and 
there exists a $ \mathcal{C}^m $-measurable function 
$ \Fsm : W_{\bf{0}}(\mathbb{R} ^m)\times W(\mathbb{R} ^{\mathbb{N}} )\to W(\mathbb{R} ^m)$ 
such that $ \Fsm $ is $ \mathcal{C}_t^m/ \mathcal{B}_t (W (\Rm ) ) $-measurable for each $ t $, 
and $ \Fsm $ satisfies 
\begin{align*}&
\mathbf{Y} ^m=\Fsm (\mathbf{B} ^m, \Xmstar ) \quad \text{ $ \Ps $-a.s}
.\end{align*}
Additionally, we call $\Fsm $ a strong solution of \eqref{:22u}. 
We say that SDE \eqref{:22u} has a {\em unique strong solution} for $(\mathbf{X} , \mathbf{B} ) $ 
under $ \Ps $ if there exists a strong solution $\Fsm $ such that 
$ \widehat{\mathbf{Y} }^m=\Fsm (\mathbf{B} ^m, \Xmstar ) $ for $ \Ps $-a.s.\,for any solution $(\widehat{\mathbf{Y} }^m, \mathbf{B} ^m, \Xmstar )$ of \eqref{:22u} under $ \Ps $ (see pp. 1154--1156 in \cite{o-t.tail}). 

We introduce the IFC condition of a weak solution $(\mathbf{X} , \mathbf{B} )$ as follows. 
 \smallskip 

\noindent \As{IFC} the SDE \eqref{:22u} has a unique strong solution $ \Fsm $ 
for $(\mathbf{X} , \mathbf{B} )$ under $\Ps $ for $P\circ \mathbf{X} _0^{-1}$-a.s.\,$\sbm $ for all $ m \in \mathbb{N} $. 

\begin{remark}\label{r:23c}
{\bf IFC} is an abbreviation of an {\bf i}nfinite system of {\bf f}inite-dimensional stochastic differential equations with {\bf c}onsistency. 
We write $ \Fsm = (\Fs ^{m,1},\ldots , \Fs ^{m,m})$. 
Using the uniqueness in \As{IFC}, we see that 
% \begin{align}\label{:23y}&
% \Fsm (\mathbf{B} ^m, \Xmstar ) = (\Fs ^{m+1, i} (\mathbf{B} ^{m+1}, \Xmmstar ) )_{i=1}^m 
% .\end{align}
\begin{align}&\notag %\label{:23y}&
\Fsm (\mathbf{B} ^m, \Xmstar ) = 
(\Fs ^{m+1, 1} (\mathbf{B} ^{m+1}, \Xmmstar ),\ldots,\Fs ^{m+1, m} (\mathbf{B} ^{m+1}, \Xmmstar ) ) 
.\end{align}
This consistency is critical in the proof of \tref{l:23}. From the equation above, we can reduce the infinite-dimensional problem to the one in finite-dimensional issue with random environment-type \cite{o-t.tail}. 
See \cite{k-o-t.ifc} for a sufficient condition for \As{IFC}. 
\end{remark}

We now state the definition of strong solution of \eqref{:12t}, which is analogous to 
Definition 1.6 in \cite[p. 163]{IW}. 
Let $ \PBr $ be the distribution of an $ \RN $-valued Brownian motion 
$ \mathbf{B} $ with $ \mathbf{B}_0 = \mathbf{0}$. 
%
%Let $ \WRNz = \{ \mathbf{w} \in \WRN \, ; \mathbf{w} _0 =\mathbf{0} \}$. 
Clearly, $ \PBr (\WRNz )= 1 $. 

Let $ \Bt $ be in \eqref{:50k}. 
Let $ \Bt (\PBr )$ be the completion of 
$ \sigma [ \mathbf{w} _s ; 0\le s \le t ,\, \mathbf{w} \in \WRNz ] $ with respect to $ \PBr $. 
Let $ \mathcal{B}(\PBr ) $ be the completion of $ \mathcal{B}(\WRNz ) $ with respect to $ \PBr $. 

\begin{definition}[a strong solution starting at $ \mathbf{s}$] \label{d:55} 
A weak solution $ \mathbf{X}$ of \eqref{:12t} starting at $ \sss $ 
with an $ \RN $-valued $ \Ft $-Brownian motion 
$\mathbf{B}=\mathbf{B}(\omega ) $ is called a strong solution starting at $ \mathbf{s}$ defined on $ \OFPF $ 
if $ \mathbf{X}_0=\mathbf{s}$ $ P$-a.s.\, and if there exists a function 
$ \map{\Fs }{\WRNz }{\WRN }$ such that 
$ \Fs $ is $ \mathcal{B}(\PBr ) /\mathcal{B}(\WRN )$-measurable, and that 
$ \Fs $ is $ \Bt (\PBr ) /\Bt $-measurable for each $ t $, and 
that $ \Fs $ satisfies 
\begin{align}\notag &% \label{:50q}&
\mathbf{X} = \Fs (\mathbf{B}) \quad \text{ $ P$-a.s}
.\end{align}
Also, we call $ \Fs $ itself a strong solution starting at $ \mathbf{s}$. 
\end{definition}

We do not impose any measurability of $\mathbf{s}\mapsto \Fs (\omega ) $ in \dref{d:55}. We have only imposed the measurability of $ \omega \mapsto \Fs (\omega ) $, given that $ \mathbf{s} $ is fixed. Considering these, we introduce the assumption on a family of strong solutions $ \{ \Fs \} $ starting at $ \mathbf{s} $. 

\smallskip
\noindent 
\As{MF} \quad 
$ \sbm \mapsto P( \Fs (\BBB ) \in A) $ is $\overBRN $-measurable
for any $A\in \mathcal{B}(W( \mathbb{R} ^{\mathbb{N}} ))$, where $\overBRN $ is the completion of 
$\mathcal{B}(\mathbb{R} ^{\mathbb{N}} )$ with respect to $P\circ \mathbf{X} _0^{-1}$. 
\smallskip

 \begin{remark}\label{r:23cc} % \thetag{1} 
The name \As{MF} comes from the assumption that $ P( \Fs (\BBB ) \in A)$ is a {\bf m}easurable {\bf f}unction in $ \mathbf{s}$. Because $ \Fs = F _{\lab (\sss )}$ is constructed for fixed a.s.\,$ \sss $, it is not apparent that $ P( \Fs (\BBB ) \in A)$ is 
$\overBRN $-measurable as a function in $ \mathbf{s}$. 
In \tref{l:23}, the dynamics are associated with the Dirichlet form, the relation 
$ \ulab ( \mathbf{s}) = \sss $ holds, where $ \ulab (\mathbf{s}) = \sum_i \delta_{s^i}$ for $ \mathbf{s} = (s^i)_i $, and 
\begin{align}\label{:23z}&
P( \Fs (\BBB ) \in A) = \mathbb{P}_{\sss } (\lpath (\mathbb{X}) \in A )
.\end{align}
The measurability of $ \mathbb{P}_{\sss } (\lpath (\mathbb{X}) \in A )$ in $ \sss $ is obvious from the Dirichlet form theory. 
Indeed, the transition probability $ \mathbb{P}_{\sss }$ 
is a quasi-continuous function in the starting points $ \sss $ (cf.\,\cite{fot,c-f}). 
Thus, the measurability of $ P( \Fs (\BBB ) \in A)$ in $ \mathbf{s}$ holds from the identity \eqref{:23z} and the measurability of the map $ \ulab $. 
 \end{remark}

For a family of strong solutions $ \{ \Fs \} $ satisfying \As{$\mathbf{MF}$}, we set 
\begin{align}\label{:5Az}&
\pPF = \int \pP ( \Fs (\mathbf{B} ) \in \cdot ) \pP \circ \mathbf{X}_0^{-1} (d\mathbf{s}) 
.\end{align}
We need \As{MF} to give meaning to the right-hand side of \eqref{:5Az}. 
Following \cite[p.1157]{o-t.tail}, we recall the concept of \lq\lq a family of unique strong solutions $ \{ \Fs \} $.'' %
\begin{definition}\label{dfn:43}
For a condition \As{$ \bullet $}, 
we say \eqref{:12t} has a family of unique strong solutions $ \{ \Fs \} $ starting at $ \mathbf{s}$ for 
$ P \circ \mathbf{X}_0^{-1} $-a.s.\,$ \mathbf{s}$ 
under the constraints of \As{$\mathbf{MF}$} and \As{$ \bullet $} 
if %$ \{ \Fs \} $ satisfies \As{$\mathbf{MF}$} and 
$ \{ \Fs \} $ satisfies \As{$\mathbf{MF}$} and $\pPF $ satisfies \As{$ \bullet $}. Furthermore, 
\thetag{i} and \thetag{ii} are satisfied.

\noindent\thetag{i} 
For any weak solution $ (\hat{\mathbf{X}},\hat{\mathbf{B}})$ under $ \hat{P}$ of \eqref{:12t} with 
$$ 
\hat{P}\circ \hat{\mathbf{X}}_0^{-1} \prec P \circ \mathbf{X}_0^{-1} 
$$
 satisfying \As{$ \bullet $}, it holds that, for $\hat{P} \circ \hat{\mathbf{X}}_0^{-1}$-a.s. $ \mathbf{s}$, 
\begin{align}\notag &%\label{:5Ab}&
\hat{\mathbf{X}}=\Fs (\hat{\mathbf{B}}) \quad 
\text{ $ \hat{P}_{\mathbf{s}} $-a.s.}
,\end{align}
where $ \hat{P}_{\mathbf{s}} = \hat{P}(\cdot \vert \hat{\mathbf{X}}_0={\mathbf{s}})$. 

\noindent 
\thetag{ii} 
For an arbitrary $ \RN $-valued $ \{\mathcal{F}_t\} $-Brownian motion $ \mathbf{B}$
defined on $ \OFPF $ with $ \mathbf{B}_0=\mathbf{0}$, 
$ \Fs (\mathbf{B})$ is a strong solution of \eqref{:12t} satisfying \As{$ \bullet $} 
starting at $ \mathbf{s}$ 
for $ P \circ \mathbf{X}_0^{-1}$-a.s.\,$ \mathbf{s}$. 
\end{definition}

\begin{definition}[uniqueness in law]\label{d:53} \thetag{i}
We say that the uniqueness in law of weak solutions on $ \Hbm $ for \eqref{:12t} 
holds if whenever $ \mathbf{X}$ and $ \mathbf{X}'$ are two weak solutions 
whose initial distributions coincide and are supported on $ \Hbm $, then the laws of the processes 
$ \mathbf{X}$ and $ \mathbf{X}'$ on the space $ \WRN $ coincide. 
If this uniqueness in law holds for an initial distribution $ \delta _{\mathbf{s}}$, then we say 
the uniqueness in law of weak solutions for \eqref{:12t} starting at $ \mathbf{s}$ holds. 
\\\thetag{ii}
We say that the uniqueness in law of weak solutions on $ \Hbm $ for \eqref{:12t} 
holds under the constraints \As{$ \bullet $} if whenever $ \mathbf{X}$ and $ \mathbf{X}'$ are two weak solutions satisfying \As{$ \bullet $} whose initial distributions coincide and are supported on $ \Hbm $, then the laws of the processes $ \mathbf{X}$ and $ \mathbf{X}'$ on the space $ \WRN $ coincide. 
If this uniqueness in law under the constraints \As{$ \bullet $} holds for an initial distribution $ \delta _{\mathbf{s}}$, then we say the uniqueness in law of weak solutions for \eqref{:12t} starting at $ \mathbf{s}$ holds under the constraints of \As{$ \bullet $} . 
\end{definition}

\begin{theorem}\label{l:23} Assume $ \beta=2$. The following then hold. 

\noindent 
\thetag{1} 
$ \mathbf{X} = \lpath (\mathbb{X}) $ under $ \mathbb{P}_{\sss } $ 
is a family of unique strong solutions of \eqref{:12t} with \eqref{:22s} 
starting at $\sbm = \lab (\sss )$ for $\mu_2 $-a.s.\,$ \sss $ under the constraints of 
\ASSUMP. 
\\\thetag{2} 
The uniqueness in law of weak solutions of \eqref{:12t} with \eqref{:22s} holds 
under the constraints of {\As{SIN}, \As{AC} for $ \mu_2 $, and \As{IFC}}. 
\\\thetag{3} 
$ \mathbf{X} = \lpath (\mathbb{X}) $ under $ \mathbb{P}_{\muairy } $ satisfies \As{SIN$ _{\mathbf{A}}$}.
\end{theorem} 

Let $ \TailS= \bigcap_{\rR =1}^{\infty} \sigma [\pi _{\SR ^c}]$ be the tail $ \sigma $-field of $ \SSS $. 
Let $ \muAt (\cdot )=\muairybeta (\cdot \vert \TailS ) (\mathfrak{t} ) $ 
be the regular conditional probability of $ \muairybeta $ with respect to $ \TailS $. 
We note that $ \muAt $ exists because $ \SSS $ is homeomorphic to a complete separable metric space \cite[Section 3.3]{o-t.tail}. 

By construction, we have a tail decomposition such that 
\begin{align}\label{:23A}&
\muairybeta (\cdot ) = \int_{\SSS } \muAt (\cdot ) \muairybeta (d\mathfrak{t})
.\end{align}
From \cite[Lemma 14.2]{o-t.tail}, we deduce that, 
for each $ \mathfrak{H}$ such that $ \muairybeta (\mathfrak{H}) = 1 $, 
there exists a subset of $\mathfrak{H}$ denoted by the same symbol $\mathfrak{H}$ 
satisfying $\muairybeta (\mathfrak{H})=1$ 
and that the following hold for all $\mathfrak{t},\mathfrak{u}\in\mathfrak{H}$. 
\begin{align}&\notag %\label{:23a}&
\muAt ( \mathfrak{A} ) \in \{0,1\}
\quad \text{ for all } \mathfrak{A} \in \TailS 
,\\ \notag &%\label{:23b}&
\muAt (\{\mathfrak{v} \in \SSS ; 
\muAt = \mu_{\beta ,\mathfrak{v}} \})=1
, \\ \label{:23c}&
\muAt \text{ and } \mu_{\beta ,\mathfrak{u}} 
\text{ are mutually singular on } \TailS \text{ if }
 \muAt \not= \mu_{\beta ,\mathfrak{u}} 
.\end{align}
Using Lemmas 9.2 and 12.1 in \cite{o-t.tail}, we find that, for $ \muairybeta $-a.s.\,$ \mathfrak{t}$, 
Theorems \ref{l:21} and \ref{l:22} hold if we substitute 
$\mu_{\beta, \mathfrak{t}} $ with $\mu_\beta $. 
We denote the associated unlabeled process by 
$(\mathbb{X}, \{\mathbb{P}_{\sss}^{\mathfrak{t}}\}_{\sss \in \SSS})$. 

\begin{theorem}	\label{l:24}
Let $\beta=1,4$. For $ \muairybeta $-a.s.\,$ \mathfrak{t} \in \SSS $, the following hold. 
\\\thetag{1} 
$ \mathbf{X} =\lpath (\mathbb{X})$ under $ \mathbb{P}_{\sss}^{\mathfrak{t}} $ is 
a family of unique strong solutions of \eqref{:12t} with \eqref{:22s} starting at $\sbm=\mathfrak{l}(\sss)$ 
 for $ \muAt $-a.e.\,$ \sss $ under the constraints of 
\ASSUMPt. 
\\
\thetag{2} 
The uniqueness in law of weak solutions of \eqref{:12t} with \eqref{:22s} holds 
under the constraints of 
\As{SIN}, \As{AC} for $ \muAt $, and \As{IFC}. 
\\\thetag{3} 
$ \mathbf{X} = \lpath (\mathbb{X}) $ under 
$ \mathbb{P}_{ \mu_{\beta, \mathfrak{t}}}^{\mathfrak{t}} $ % for $ \muA $-a.s.\,$ \mathfrak{t}$ 
satisfy \As{SIN$ _{\mathbf{A}}$}. 
\end{theorem}

% %%%%%%%%%%%%TH 2.5%%%%%%%%%
We turn to a Girsanov formula. 
For $T\in(0,\infty)$, let $ \Ws ^{m }$ be the distribution of 
$ \mathbb{R}^{m }$-valued Brownian motion 
$ \{ (B_t^i+s^i)_{i=1}^{m } \}_{t\in[0,T]} $ starting at 
$ \sbm ^m =(s^i)_{i=1}^m\in\mathbb{R} ^{m }$, $m \in \mathbb{N} \cup \{\infty\}$. 
Let $ \mathbf{P}_{\sbm }=\mathbb{P}_{\sss} \circ \lpath ^{-1}$, where $ \sbm=\lab (\sss ) $, as before. 
It is then clear that the probability measure $ \mathbf{P}_{\sbm }$ 
 is {\em not} absolutely continuous with respect to $\Ws ^{\infty} $. 
Hence, we formulate a Girsanov-type formula for a finite number of particles 
$ \mathbf{X} ^{m }=\{(X_t^i)_{i=1}^m \}_{t\in[0,T]} $ for each $ m \in \mathbb{N}$.

We set $ \Xmstar =\{ (X_t^{i})_{i= m +1}^{\infty} \}_{\tzT } $. Here, we omit $ T $ from the notation. 
In \tref{l:25}, $ \mathbf{X} ^{m }$ is the restriction of the original $ \mathbf{X} ^{m }$ to $ [0,T]$. 
Let $ \Pms $ be the regular conditional probability 
of $ \mathbf{P}_{\sbm } (\mathbf{X} ^{m } \in \, \cdot \, ) $ conditioned at $ \Xmstar $ such that 
\begin{align}&\label{:25x}
 \Pms (\ \cdot \ ) =\mathbf{P}_{\sbm } (\mathbf{X} ^{m } \in \cdot \ \vert \Xmstar )
.\end{align}
Let $\mathbf{x} ^m= (x^i)_{i=1}^m \in \mathbb{R} ^m $. Let $\map{ \bXm }{\mathbb{R} ^{m }\ts [0,T]}{\mathbb{R} ^{m }}$ be the vector $ \bXm=(\bXmi )_{i=1}^m$ such that for $ i=m = 1$ 
\begin{align*}&
\bXonei (\mathbf{x} ^1,t)= \frac{\beta}{2}
\lim_{r\to\infty} 
\Big( \sum_{ j= 2 ,\, \lvert X_t^j \rvert <r}^{\infty} \frac{1}{x^1 -X_t^j } -
 \int_{\verty <r}\frac{\rhohat (y)}{-y}dy \Big)\Big) 
\end{align*}
and for $ 2 \le i \le m $ 
\begin{align}\notag %\label{:25y}
\bXmi &(\mathbf{x} ^m,t)=
\frac{\beta}{2} \Big( 
\sum_{j=1 ,\, j\not= i}^{m } \frac{1}{x^i -x^j} + 
\lim_{r\to\infty} 
\Big( \sum_{ j= m +1 ,\, \lvert X_t^j \rvert <r}^{\infty} \frac{1}{x^i -X_t^j } -
 \int_{\verty <r}\frac{\rhohat (y)}{-y}dy \Big)\Big) 
.\end{align}
\begin{theorem} \label{l:25}
Assume $ \beta=1,2,4 $ and $T\in(0, \infty)$. 
For $ h \in \mathbb{N}$, let $ \map{\tauh }{W(\mathbb{R} ^{m})}{\mathbb{R} \cup\{ \infty \} } $
 be the stopping time with respect to the canonical filtration such that 
\begin{align}&\label{:25a}
\tauh ( \mathbf{w} ^m ) =
\inf \Big\{ t \wedge T \, ;\, h \le \int _0^t \lvert 
\bXm ( \mathbf{w} _u^m,u) \rvert ^2 du \Big\} 
.\end{align}
For $ \muairybeta \circ \lab ^{-1} $-a.s.\,$ \sbm $, let 
$ \Pmsh =
\Pms \circ ( \mathbf{w} _{\cdot\wedge \tauh }^m)^{-1} $ and 
$ \Wsh ^{m }=\Ws ^{m } \circ ( \mathbf{w} _{\cdot \wedge \tauh }^m)^{-1}$. 
Then, for $ \mathbf{P}_{\sbm }$-a.s.\,$ \Xmstar $, $\Pmsh $ 
is absolutely continuous with respect to $ \Wsh ^{m } $ with Radon--Nikodym density 
\begin{align} \label{:25c}& 
\frac
{d \Pmsh }
{d \Wsh ^{m }}=\exp \Big( 
 \int_0^{\cdot \wedge \tauh } \bXm ( \mathbf{w} _u^m,u)d \mathbf{w} _u^m - 
 \frac{1}{2}\int_0^{\cdot \wedge \tauh } \lvert \bXm ( \mathbf{w} _u^m,u) \rvert ^2 du 
 \Big) 
.\end{align}
Furthermore, 
we obtain for $ \mathbf{P}_{\sbm }$-a.s.\,$ \Xmstar $ 
\begin{align}\label{:25d}&
\text{$ \lim_{ h \to \infty } \tauh ( \mathbf{w} ^m )=T $ for $ \Pms $-a.s.\,$ \mathbf{w} ^m $}
.\end{align}
\end{theorem}

We present applications of the uniqueness results in \tref{l:23} and \tref{l:24}. 

For $ \beta=2 $, the infinite-dimensional stochastic dynamics 
have been constructed using the space-time correlation functions \cite{FNH99, p-spohn, joh.03,NKT03, KNT04,CH14}. 
Let $ \mathbb{K}_{\mathrm{Ai}}(s,x,t,y) $ be the extended Airy kernel defined as 
\begin{align} \notag 
\mathbb{K}_{\mathrm{Ai}}(s,x,t,y) =\left\{
\begin{array}{cc}
\displaystyle{\int_{-\infty}^0 du \ 
e^{(t-s)u/2}\Ai(x-u)\Ai(y-u)},
&\mbox{if $s\le t$}\cr
&\cr
\displaystyle{-\int_0^\infty du \ 
e^{(t-s)u/2}\Ai(x-u)\Ai(y-u)},
&\mbox{if $s>t$}
.\end{array}
\right 
.\end{align}
The determinantal process $\mathbb{Y}=\{ \mathbb{Y}_t \} $ 
with the extended kernel $ \mathbb{K}_{\mathrm{Ai}}$ is 
an $\SSS $-valued stochastic process such that, for any $M\in \mathbb{N} $, 
$\mathbf{f} =(f_1,f_2,\dots,f_M) \in C_0(\mathbb{R} )^M$ and a sequence of times 
${\bf t}=(t_1,t_2,\dots,t_M)$ with $0 < t_1 < \cdots < t_M < \infty$,
if we set $\chi_{t_m}(x)=e^{f_{m}(x)}-1, 1 \leq m \leq M$,
the moment generating function of a multi-time distribution, 
$$
{\Psi}^{\bf t}[\mathbf{f} ]
\equiv \E \Big[\exp \Big( \sum_{m=1}^{M} 
\int_{\mathbb{R} } f_m(x) \mathbb{Y}_{t_m}(dx) \Big) \Big],
$$
is given by a Fredholm determinant
\begin{align} \notag 
{\Psi}^{\bf t}[\mathbf{f}] 
= \mathop{\Det}_
{ %\substack
{(s,t)\in \{t_1, t_2, \dots, t_M \}^2, (x,y)\in \mathbb{R} ^2}}
\Big[\delta_{s t} \delta(x-y)
+ \mathbb{K}_{\mathrm{Ai}} (s, x, t, y) \chi_{t}(y) \Big] 
.\end{align}
Here $ \delta_{s t}$ denotes the Kronecker delta and $ \delta (\cdot )$ is the delta measure at $ \cdot $.
The reader may refer to \cite{KT07b} for the details.

Let $ \mathbb{Q}_{\sss }$ be the distribution of 
a determinantal process with the extended Airy kernel $ \mathbb{K}_{\mathrm{Ai}} $ 
 starting at $ \sss $ \cite{KT09}. 
Then there exists $ \SSS _{0}$ such that $ \muairy (\SSS _{0})=1 $ and that 
$ \{\mathbb{Q}_{\sss }\}_{\sss \in\SSS _{0}}$ is the family of transition probabilities of a diffusion process with state space $ \SSS _{0}$ \cite{o-t.sm}. 
There are thus two completely different approaches for constructing infinite-dimensional stochastic dynamics related to the Airy RPF with $ \beta=2 $. 
In \cite{o-t.sm}, we proved that these two infinite-dimensional stochastic dynamics coincide. 

%%%%%%%%%%%%%%%%%%%%%%%%%TH 2.5%%%%%%%%%%%%%%%%
\begin{theorem}[{\cite[Theorem 2.2]{o-t.sm}}]	\label{l:26} 
Let $ \mathbb{Q}_{\muairy }=\int_{\SSS _{0}} \mathbb{Q}_{\sss }d \muairy $. Then, 
\begin{align} &\notag %\label{:24a}
\mathbb{Q}_{\muairy } ( \mathbb{Y} \in \Wnex (\SSSsip ) )=1 
.\end{align}
The labeled process $ \mathbf{Y} =\lpath (\mathbb{Y})$ under $ \mathbb{Q}_{\sss } $ 
is a strong solution of \eqref{:12t} with \eqref{:22s} starting at 
$\sbm=\lab (\sss)$ for $ \mu $-a.s.\,$ \sss$ satisfying 
\satASSUMP. 
Moreover, 
\begin{align}&\notag 
\quad \quad \mathbb{P}_{\sss }=\mathbb{Q}_{\sss } 
\quad \text{ for $ \muairy $-a.s.\,$ \sss$}
.\end{align}
 \end{theorem}

For $ N \in\mathbb{N} $ and $\sss $ with $\lab(\sss)=(s^i)_{i\in \mathbb{N}} $, 
we set $\lab^{N }(\sss) =(s^i)_{i=1}^N \in \mathbb{R} ^{N }$. 
Let ${\mathbf{X} }^{N } $ be the solution of \eqref{:11l} with ${\mathbf{X} }_0^{N }=\lab^{N }(\sss) $. 
For $m\le N $, let $ \mathbf{X} ^{N ,m}=(X^{N , i})_{i=1}^m $ be the first $ m $-components of $ \mathbf{X} ^{N }$. Let $ \mathbb{X}^{N } = \sum_{i=1}^N \delta_{X^{N ,i }}$. 
Let $ \mathbb{P}_{\muairybeta ^{N }}^{N } $ be the distribution of the unlabeled process 
$ \mathbb{X}^N $ such that $ \mathbb{X}_0^N = \muairybeta ^{N }$ in law. 
\begin{theorem}[{\cite[{Corollary 2.3}]{o-t.sm}, \cite[Proposition 5.1]{k-o.fpa}}]	\label{l:27}
The following hold. 
\\
\thetag{1} 
$ \mathbb{P}_{\muairybeta ^{N }}^{N } $ weakly converges to 
$ \mathbb{P}_{\muairybeta } $ in $ W (\SSS )$ as $N \to\infty$.
\\\thetag{2} 
Under $ \mathbb{P}_{\muairybeta ^{N }}^{N } $ and $ \mathbb{P}_{\muairybeta } $, 
the first $ m $-labeled processes $ \mathbf{X} ^{N ,m }$ converge to 
$ \mathbf{X} ^{m } $ in distribution in $ W(\mathbb{R} ^{m })$ as $N \to\infty$ for each $ m \in \mathbb{N}$. 
\end{theorem}
For $ \beta = 2$, \tref{l:27} was proved in \cite{o-t.sm} and, for $ \beta = 1,2,4$, \tref{l:27} was proved in \cite{k-o.fpa}. The methods in \cite{o-t.sm} and \cite{k-o.fpa} are completely different.

\begin{remark}\label{r:24} 
\thetag{1} \tref{l:25} implies that the local properties 
of tagged particles of $ \mathbf{X} $ are the same as those of Brownian motion $ \BBB $. 
In particular, each tagged particle $ X_t^i$ ($i\in\mathbb{N} $) is non-differentiable and 
$ \alpha $-H\"{o}lder continuous for any $ \alpha < 1/2 $ in $ t $. 
\\\thetag{2} 
The unlabeled process $ \mathbb{X}$ is reversible, and the time parameter can thus be extended 
from $ [0,\infty)$ to $ \mathbb{R}$ by the stationarity of the time shift. 
Assume $ \beta =2 $ and consider the process given by the extended Airy kernel. 
Let $ \At $ be the top particle. Thus far, we have denoted this by $ X_t^1$. 
$ \At $ is usually called the Airy process. 
We consider the case of a time stationary Airy process. 
Therefore, the distribution of $ \At $ is independent of $ t $. 
In \cite[Conjecture 1.5]{joh.03}, Johansson conjectured that 
$ H(t)=\At -t^2$ almost surely has a unique maximum point in $ [-T,T]$. 
We have an immediate affirmative answer for this conjecture from \tref{l:25} and \tref{l:26}. 
In fact, this is the case of the Brownian path, and 
the Airy process is absolutely continuous on the time interval $ [-T,T]$ 
with respect to Brownian motions starting from the distribution of $ \Az $ at time $ -T $. 
We remark that the conjecture referred to above 
was solved by Corwin--Hammond \cite{CH14} and H\"agg \cite{hag.08} using a different method. 
\end{remark}

%{\em\sc Discussion. } Let $ \mu = \muairy $ for $ \beta = 2$ and $ \muAt $ for $ \beta = 1,4$. The condition \As{AC} for $ \mu $ is satisfied if $ \mathbf{X}$ is given by \tref{l:22} and the distribution of the associated unlabeled dynamics $ \mathbb{X}_0$ has a density with respect to $ \mu $. It is a challenging problem to prove $ \mathbf{X} $ satisfies \As{AC} for $ \mu $ when the distribution of $ \mathbb{X}_0$ is singular to $ \mu $ and, in particular, $ \mathbb{X}_0$ concentrates on a delta mass $ \delta_{\sss }$. This problem is related to a construction of non-equilibrium dynamics. 

We prove the continuity of $ t \mapsto \lab^i (\mathfrak{w}(t)) \in \mathbb{R}$ for each $ i \in \mathbb{N}$. 
Let $ \SSSsip $ and $ \tau $ be as in \eqref{:20d} and \eqref{:20z}, respectively. 
\begin{lemma} \label{l:28}
For each $ \mathfrak{w} \in W ( \SSSsip ) $, 
$ (\lab ^1(\mathfrak{w}(t) ) , \ldots , \lab ^k (\mathfrak{w}(t)) ) $ is continuous in $ [0,\tau )$ 
for each $ k \in \mathbb{N}$. 
In particular, $\lab^i (\mathfrak{w}(t))$ is continuous in $ [0,\tau )$ for each $ i \in \mathbb{N}$. 
\end{lemma}
\begin{proof}
Let $ a \in [0,\tau )$ be fixed. 
Recall that $ t \mapsto \mathfrak{w}(t) $ is continuous at $ a $. 
Then there exist $ b , c $ such that $ b < a < c $ and $ Q < R $ such that 
\begin{align*}&
\mathfrak{w} (t) (\QR ) = k ,\quad \mathfrak{w}(t)([R , \infty )) = 0 
,\quad \text{ for all } t \in (b,c)
.\end{align*}
Thus, at time $ t=a$, we have exactly $ k $ particles in $ \QR $ and no particle in $ [R,\infty)$. 
Let $ \pi_{\QR} (\sss ) = \sss (\cdot \cap \QR )$ is the projection. Let 
\begin{align*}&
\SSS (Q,R,k) = \{ \sss \in \pi_{\QR}( \SSSsip ) ; \sss (\QR ) = k \} ,
\\&
\mathbb{R}_{<} (Q,R,k) = \{ \mathbf{s}=(s^i)_{i=1}^k \in \QR ^k\, ;\, s^i > s^{i+1} \text{ for all } i = 1\ldots, k-1 \} 
.\end{align*}
We pose $ \SSS (Q,R,k) $ the relative topology of $ \SSS $. 
Because we fix the number of particles, % as $ k $, 
$ \SSS (Q,R,k) $ and $ \mathbb{R}_{<} (Q,R,k) $ are homeomorphic via the map 
$ \mathfrak{s} \mapsto (\lab^i(\mathfrak{s}))_{i=1}^k$. Hence, 
\begin{align*}&
t \mapsto (\lab ^1(\mathfrak{w}(t) ) , \ldots , \lab ^k (\mathfrak{w}(t)) ) 
\end{align*}
is continuous in $ t \in (b,c)$. From this we obtain the first claim. 
The second claim is clear from the first claim. 
\end{proof}

Let 
$ \mathcal{B}(W_{\mathrm{NE}}(\SSSsip ) )_t = \sigma [\mathfrak{w}(u); u \le t]$ be as in \eqref{:21z}. 
Let $ \mathcal{B}(W (\RNN ))_t $ be the sub-$ \sigma $-field of $ \mathcal{B}(W (\RNN )) $ such that 
\begin{align*}&
 \mathcal{B}(W (\RNN ))_t = \sigma [\mathbf{w}(u); u \le t ]
.\end{align*}
Here, for a fixed $ u $, $ \map{\mathbf{w}(u)}{W (\RNN ))}{\RNN }$ is such that $ \mathbf{w} \mapsto \mathbf{w}(u)$. 

As for the measurability of $ \lpath $, we have the following. 
\begin{lemma} \label{l:29} \thetag{1}
$ \lpath $ is $ \mathcal{B}(W_{\mathrm{NE}}(\SSSsip ) ) / \mathcal{B}(W (\RNN )) $-measurable. 
\\\thetag{2} 
$ \lpath $ is $ \mathcal{B}(W_{\mathrm{NE}}(\SSSsip ) )_t / \mathcal{B}(W (\RNN ))_t $-measurable for each $ t \ge 0 $. 
\end{lemma}
\begin{proof}
For a complete separable metric space $ A $, 
let $ \mathscr{C}(W (A)) $ be the totality of the Borel cylinder sets of $ W (A)$. 
By definition, $ B \in \mathscr{C}(W (A)) $ if 
\begin{align}& \notag 
B = \{ w \in W (A); (w (t_1),\ldots, w(t_n)) \in C \} 
\end{align}
for some $ t_1,\ldots t_n \in [0,\infty )$ and $ C \in \mathcal{B}(A^n) $, where $ n \in \mathbb{N}$. 
It is known that 
\begin{align}&\label{:29f}
\sigma [\mathscr{C}(W (A)) ] = \mathcal{B} (W (A)) 
.\end{align}
See \cite[p.16]{IW} for the proof in the case of $ A = \mathbb{R}^d$; the proof of the general case is the same, and we omit it. 

Recall that $( \lpath (\mathfrak{w}) (t_1) , \ldots , \lpath (\mathfrak{w}) (t_n)) = (\lab (\mathfrak{w}(t_1)), \ldots, \lab (\mathfrak{w}(t_n)))$ and that $\lab $ is a $ \mathcal{B}(\SSSsip )/\mathcal{B}(\RNN ) $-measurable map. Then, it is clear that 
\begin{align}&\notag % \label{:29g}
\{ \mathfrak{w} \in W (\SSSsip )\, ;\, ( \lpath (\mathfrak{w}) (t_1) , \ldots , \lpath (\mathfrak{w}) (t_n)) \in 
\mathbf{C}
\} 
\in \mathscr{C}( W (\SSSsip )) 
\end{align}
for any $ t_1,\ldots t_n \in [0,\infty )$ and $ \mathbf{C} \in \mathcal{B}((\RNN )^n) $, where $ n \in \mathbb{N}$. Hence, 
\begin{align}& \label{:29h}
\lpath ^{-1} (\mathscr{C}(W (\RNN )) ) \subset \mathscr{C} (W (\SSSsip )) 
.\end{align}
Using \eqref{:29f} and \eqref{:29h} and noting that $ W_{\mathrm{NE}}(\SSSsip ) $ is a Borel subset of $ W (\SSSsip ) $, we obtain \thetag{1}. 
The proof of \thetag{2} is the same as the proof of \thetag{1}. Hence we omit it. 
\end{proof}

% \begin{remark}\label{r:28} 
% Then, $ \{\mathbf{P }_{\sbm }\}_{\sbm \in \RNN }$ is 
% a $ \muairybeta \circ \lab ^{-1}$-reversible diffusion with state space $ \RNN $ 
% because $ \lab $ is injective and \tref{l:21} \thetag{2} holds. 
% \end{remark}
%%%%%%%%%%%%%%%%%%%%%%%%%%%%%%%%%%%%%%
\section{Finite particle approximation of Airy RPFs}\label{s:3}
%%%%%%%%%%%%%%%%%%%%%%%%%%%%%%%%%%%%%%%%%%%

Let $ \muairybeta ^{N }$ be the RPF whose labeled density 
$ \mairybeta ^{N } $ is given by \eqref{:11c} as before. 
The sequence $ \{\muairybeta ^{N }\}_{N \in \mathbb{N}}$ approximates the measure $\muairybeta $ and 
plays an important role in the proof of the main theorems. In this section, we collect estimates of 
correlation functions and quantities concerning $ \mairybeta ^{N } $. 

\subsection{Quaternion determinantal point fields and their kernels $ \beta=1,4$ } \label{s:31}

Let $ \beta=1,4$. 
The random point fields $\mu_\beta$ and $ \muairybeta ^{N }$, $N\in\mathbb{N} $ 
are quaternion determinantal random point fields by definition. 
In this subsection, we give the precise definition of the associated kernels.

We first recall the standard quaternion notation for $ 2\ts 2 $ matrices (see \cite[Ch.\ 2.4]{Meh04}, \cite{Dy.70}), 
\begin{align} &\notag 
\mathbf{1}=
\begin{bmatrix}
1& 0 \\ 0&1
\end{bmatrix},\quad 
\mathbf{e}_1=
\begin{bmatrix}
\sqrt{-1} & 0 \\ 0&-\sqrt{-1}
\end{bmatrix},\quad 
\mathbf{e}_2=
\begin{bmatrix}
0&1\\-1&0
\end{bmatrix},\quad 
\mathbf{e}_3=
\begin{bmatrix}
0&\sqrt{-1} \\\sqrt{-1} &0
\end{bmatrix}
.\end{align}
A quaternion $ q $ is represented by 
$ q=q^{(0)}\mathbf{1} + q^{(1)}\mathbf{e}_1 + 
 q^{(2)}\mathbf{e}_2 + q^{(3)}\mathbf{e}_3 $, 
where $ q^{(i)} \in \mathbb{C} $. There is a natural identification between the $ 2\ts 2 $ complex matrices and the quaternions given by 
\begin{align} &\notag 
\begin{bmatrix}
a&b \\ c&d
\end{bmatrix}=
\frac{1}{2}(a+d)\mathbf{1} -
\frac{\sqrt{-1} }{2}(a-d)\mathbf{e}_1 +
\frac{1}{2}(b-c)\mathbf{e}_2 -
\frac{\sqrt{-1} }{2}(b+c)\mathbf{e}_3 
.\end{align}
Let $ q=q^{(0)}\mathbf{1} + \sum_{i=1}^3 q^{(i)}\mathbf{e}_i $. 
A quaternion $ q $ is called complex scalar if $ q^{(i)}=0 $ for 
$ i=1,2,3 $. We often identify a complex scalar quaternion 
$ q=q^{(0)}\mathbf{1}$ by the complex number $ q^{(0)} $. 
A scalar quaternion $ q=q^{(0)}\mathbf{1}$ is regarded as a real number if $ q^{(0)} \in \mathbb{R} $.

Let $ \Kairy $ and $ \Ai $ be as in \eqref{:11g} and \eqref{:11h}, respectively. Let 
\begin{align}\label{:31t}&
J_1(x,y)=K_{2}(x,y)+\frac{1}{2}
\Ai(x) \Big(1-\int_y^\infty \Ai(u)du \Big) 
,\\ \label{:31u}& 
J_4(x,y)=\Kairy(x,y) - \frac{1}{2} \Ai(x)\int_y^\infty \Ai(u)du 
.\end{align}
Let $ \partial_i = \partial/\partial{x^i}$, $i\in \mathbb{N} $. 
We then define quaternion kernels $ \Kairyone $ and $ \Kairyfour $ as 
\begin{align} \label{:31v}&
\Kairyone (x,y)
 =\begin{bmatrix} 
J_1(x,y) & -\partialtwo J_1(x,y) \\
\int_y^x J_1(u,y)du - \frac{1}{2}{\rm sign}(x-y) & J_1(y,x)
\end{bmatrix}
,\\ \label{:31w}&
 \Kxy =\frac{1}{2}
\begin{bmatrix} 
J_4(x,y) & -\partialtwo J_4(x,y) \\ 
\int_y^x J_4(u,y) du & J_4(y,x)
\end{bmatrix}
.\end{align}
(See \cite[Subsection 3.9.3]{AGZ10}.) The scaling $2^{\frac{2}{3}}$ in \eqref{:31w} is chosen to obtain \eqref{:34a} below.

Let $\widehat{H}_{k }$, $\phi_{k } $, $\psi_{k }^{ N } $, and $ \psi_{N } $, where 
$ k \in \zN $ and $ N\in \mathbb{N}$, be the Hermite polynomials, the normalized oscillator functions, and their modifications defined as 
\begin{align}\label{:31y}&
\widehat{H}_{k }(x)=(-1)^{k } e^{x^2/2}\frac{d^{k }}{dx^{k }}e^{-x^2/2} 
,\quad 
\phi_{k }(x)=
(\sqrt{2\pi}{k }!)^{-1/2}
% \frac{1}{\sqrt{\sqrt{2\pi}{k }!}}
e^{-x^2/4}\widehat{H}_{k }(x) 
,\\\label{:31a} &
\psi_{k }^{ N }(x)= N ^{1/12}\phi_{k } ( 2\sqrt{ N } + N ^{-1/6}x ) 
,\quad 
\psi_{N } (x)=\psi_{N }^{N } (x) 
.\end{align}
The following properties of $ \psi _{N }$ $ (N \in\mathbb{N} ) $ are well known (see \cite[Lemma 3.2.5]{AGZ10}). 
\begin{align}
&
 \mbox{ $\psi_{N }(x-2N ^{2/3})$ is even (resp. odd) if $N $ is even (resp. odd)},
\label{:31b}
\\ \label{:31c} 
&\int_{\mathbb{R} }\psi_{2N -1}(x)dx =0, \quad \int_{\mathbb{R} }\psi_{2N }(x)dx \neq 0 ,\quad 
\lim_{N \to\infty}\int_{\mathbb{R} }\psi_{2N }(x)dx= 2 
.\end{align}
The correlation kernel $\Kairy^{N }$ of $\muairy^{N }$ is given by 
\begin{align}\label{:31f}&
\Kairy^{N }(x,y)=\frac{1}{N ^{1/3}} \sum_{ k=0}^{N -1} 
\psi_{k }^N(x)\psi_{k }^N(y)
.\end{align}

For a proposition $\mathsf{A}$, we write
$ {\bf 1}(\mathsf{A})= 1$ if $\mathsf{A} $ is true, and 
$ {\bf 1}(\mathsf{A})= 0 $ otherwise. 
For each integrable real-valued function $ f $ on $\mathbb{R} $, we set 
\begin{align}\label{:31g}&
(\varepsilon f)(x)=\int_{\mathbb{R} }\frac{1}{2}{\rm sign}(x-y)f(y)dy=\frac{1}{2}\int_\mathbb{R} f(y)dy -\int_x^\infty f(y)dy 
.\end{align}
Let $ J_1^{N }(x,y) $ and $ J_4^{N }(x,y) $ be complex valued functions such that 
\begin{align}\notag %\label{:31i}
J_1^{N }(x,y) &=
\Kairy^{N }(x,y)+\frac{1}{2}\psi_{N -1}^{N } (x)\varepsilon \psi_{N }(y)
+\frac{\psi_{N -1}^{N } (x)}
{\int_{\mathbb{R} }\psi_{N -1}^{N } (\zz )d\zz } 
{\bf 1}(\mbox{$N $ is odd})
,\\ \label{:31j} 
J_4^{N }(x,y) & 
= 
 K_2^{2N} (x,y) + \frac{\psi_{2N}(x) \psi_{2N}(y)}{(2N)^{1/3} } + 
\frac{\sqrt{2N +1}}{2(2N )^{1/2}}
 \psi_{2N }(x)\varepsilon \psi_{2N +1}^{2N }(y) 
.\end{align}
We easily see that
\begin{align}\label{:31p}
\partialtwo J_{\beta }^{N }(x,y) &= - \partialtwo J_{\beta }^{N }(y,x) 
,\quad 
 \int_y^x J_{\beta }^{N }(u,y) du = -\int_x^y J_{\beta }^{N }(u,x) du 
.\end{align}
The correlation kernels $\Kairybeta^{N }$ of $\mu_{\beta}^{N }$ for $\beta=1, 4$ are defined as 
(see \cite[Subsection 3.9.2]{AGZ10})
%%%%%%%%%%%%%%
\begin{align}\notag %\label{:31k}&
\Kairyone^N (x,y)&=
\begin{bmatrix} 
J_1^{N }(x,y) & -\partialtwo J_1^{N }(x,y)
\\
 \int_y^x J_1^{N }(\zz ,y) d\zz - \frac{1}{2}{\rm sign}(x-y) & J_1^{N }(y,x)
\end{bmatrix}
,\\ \label{:31l} 
\frac{1}{2^{\frac{2}{3}}}
\Kairyfour^{N }\Big( \frac{x}{2^{\frac{2}{3}}}, \frac{y}{2^{\frac{2}{3}}}\Big)
&= 
\frac{1}{2}
\begin{bmatrix} 
J_4^{N }(x,y) & -\partialtwo J_4^{N }(x,y)
\\
\int_y^x J_4^{N }(\zz ,y) d\zz & J_4^{N }(y,x)
\end{bmatrix}
.\end{align}
%We easily see $ \Kairybeta^{N } (x,x)$ are scalar quaternions regarded as positive numbers. %Hence, $ \rairybeta ^{N ,1} (x)=\Kairybeta^{N } (x,x)^{(0)}=\Kairybeta^{N } (x,x) $. 
To simplify the notation, we set 
\begin{align}\label{:31n}&
L_{\beta }^{N } (x,y)=
\Kairybeta ^{N }(y,x)\Kairybeta ^{N }(x,y)
,\\&\label{:31o}
\MbN (x,y,z)=
\Kairybeta ^{N }(y,z)
\Kairybeta ^{N }(z,x)
\Kairybeta ^{N }(x,y)
+
\Kairybeta ^{N }(y,x)
\Kairybeta ^{N }(x,z)
\Kairybeta ^{N }(z,y)
.\end{align}
%
% Let $ q=q^{(0)}\mathbf{1} +q^{(1)}\mathbf{e}_1 +q^{(2)}\mathbf{e}_2 +q^{(3)}\mathbf{e}_3 $. 
% 
\begin{lemma} \label{l:31} 
Let $ \beta=1,4$. Then $ \Kairybeta^{N } (x,x)$, 
$ L_{\beta }^{N } (x,y) $, and $M_{\beta }^{N}(x,y,z)$
 are permutation invariant and scalar quaternions regarded as real numbers. 
\end{lemma}
\begin{proof} 
\lref{l:31} easily follows from \eqref{:31p}. 
\end{proof}

\subsection{Estimates of correlation functions} \label{s:32}

In \sref{s:32}, we present the $ n$-point correlation function 
 $\rairybeta ^{N ,n }$ of $\muairybeta^{N }$ and 
investigate the property of $\rairybeta ^{N ,n }$ as $ N \to \infty $. 

For $ q=q^{(0)}\mathbf{1} + \sum_{i=1}^3 q^{(i)}\mathbf{e}_i $, let 
$ {q}^{\dagger}=q^{(0)}\mathbf{1}- \sum_{i=1}^3 q^{(i)}\mathbf{e}_i $. 
A quaternion matrix $ A=[a_{ij}]$ is called self-dual if $ a_{ij}= a_{ji}^{\dagger} $ for all $ i,j $. 
For a self-dual $ n\ts n $ quaternion matrix $ A=[a_{ij}]$, 
we set 
\begin{align} \label{:32x}&
\qdet A =\sum_{\sigma \in \mathfrak{S}_n }
\mathrm{sign} [\sigma ] 
\prod _{i=1}^{L(\sigma )}
 ( a_{\sigma_i(1)\sigma_i(2)}a_{\sigma_i(2)\sigma_i(3)}
\cdots a_{\sigma_i(\ell _i)\sigma_i(1)} )^{(0)} 
.\end{align}
Here, $ \sigma=\sigma_1\cdots\sigma_{L(\sigma )}$ 
denotes a decomposition of $ \sigma $ into 
the cyclic permutations $ \{\sigma_i\} $ with disjoint indices. 
We write $ \sigma_i=(\sigma_i(1),\ldots,\sigma_i(\ell_i )) $, 
where $ \ell_i $ is the length of the cyclic permutation $ \sigma_i $. 
It is known that the right-hand side is well defined (see \cite[Section 5.1]{Meh04}).

%	\end{document}
Let $\Kairybeta ^{N }$ be as for \eqref{:31f} and \eqref{:31l}. 
It is known (see, for example, \cite{Meh04, AGZ10, For10}) that $[\Kairybeta ^{N }(X^i,X^j)]_{i,j=1,\ldots,n} $ are self dual for $\beta = 1,4$ and 
\begin{align}\label{:32y}&
\rairy ^{N ,n } (\mathbf{x} ^n)=\det [\Kairy ^{N }(x^i,x^j)]_{i,j=1,\ldots,n}
,\\\label{:32z}&
\rairybeta ^{N ,n } (\mathbf{x} ^n)=
\qdet [\Kairybeta ^{N }(x^i,x^j)]_{i,j=1,\ldots,n}
\quad \text{ for $ \beta=1,4$}
.\end{align}
%	\end{document}
%
The next equality follows from the Christoffel--Darboux formula and 
properties of Hermite polynomials. (Refer, for instance, to \cite[Lemma 3.9.9]{AGZ10}). 
\begin{align}\label{:Y1k}
\Kairy^{N }(x,y) &
=\frac{\psi_{N }(x)\psi'_{N }(y)-\psi'_{N }(x)\psi_{N }(y)}{x-y}
-\frac{1}{2N ^{1/3}}\psi_{N }(x)\psi_{N }(y) 
.\end{align}
\begin{lemma} \label{l:32} Let $ \beta=1,2,4$. 
Then, it holds for each $ n \in \mathbb{N}$, $ i=1\ldots n $ that 
 \begin{align}&\notag 
 \limi{N } \rairybeta ^{N ,n}=\rairybeta ^{n} ,\
 \limi{N } \partial _{i} \rairybeta ^{N ,n}=\partial _{i} \rairybeta ^{n} 
 \text{ uniformly on each compact set} 
 .\end{align}
% \begin{align}&\notag 
% \limi{N }\Supr \lvert \rairybeta ^{N ,n} (x) - \rairybeta ^{n} (x) \rvert = 0 ,\ 
% \limi{N }\Supr \lvert \partial _{i} \rairybeta ^{N ,n}(x) - \partial _{i} \rairybeta ^{n} (x) \rvert =0
%% \ \text{ uniformly on each compact set} 
% .\end{align}
%
\end{lemma}
\begin{proof}
Note that $ \psi_{N } $ and $ \Ai $ are entire functions and 
$ \psi_{N } $ converge to $ \Ai $ uniformly on each compact set as $ N \to \infty $ 
(see \cite[p 134, Lemma 3.7.2]{AGZ10}). 
Then, using \eqref{:Y1k}, we see 
 that $ \Kairy ^{N }$ and $ \partial _{i} \Kairy ^{N } $ ($ i=1,2$) 
converge to $ \Kairy $ and $ \partial _{i} \Kairy $ ($ i=1,2$) 
 uniformly on each compact set. 
The same convergence also holds for $ \beta=1,4 $ from the above combined with the definitions of 
$ \Kairybeta $ and $ \Kairybeta ^{N }$. 
Hence, we obtain \lref{l:32} from \eqref{:20o} and \eqref{:32z}. 
\end{proof}

Let $ \rhohat ^{N } $ and $ \rhohat $ be as in \eqref{:12z} and \eqref{:12u}, respectively. 
It is noteworthy that $ \rhohat^N (x) = \rhohat (x) = 0 $ for all $ x \in [0,\infty )$ by definition. 
We shall take $ \rhohat ^{N } $ and $ \rhohat $ as the main terms of the approximations for the one-point correlation functions 
$ \rairybeta ^{N ,1}$ and $ \rairybeta ^{1}$, as we see in \eqref{:13b} and \eqref{:34a}. 
\begin{lemma} \label{l:33} 
Let $ 0 < \varepsilon < 1/6$. 
We then have a constant $ \Ct \label{;37} > 0$ such that 
\begin{align}\label{:33a}&
 \lvert \rhohat^N (x) - \rhohat (x) \rvert \le \frac{\cref{;37}}{\xvee } 
 \quad \text{ for each $ x \in [ - N^{\varepsilon }, 0 )$}
,\\ \label{:33b}&
\limi{N } \limi{s} \Supr 
\Big\lvert 
 \int_{\verty <s} (\rhohat ^{N }(y)-\rhohat (y))
 \Big ( \frac{1} {x-y} - \frac{1}{-y}\Big) dy \Big\rvert =0 
\end{align}
for each $ r \in \mathbb{N}$. 
The integral in \eqref{:33b} is understood in the sense of Cauchy principal value. 
\end{lemma}

\begin{proof}
From \eqref{:12u} and \eqref{:12z}, we obtain \eqref{:33a}. 
Using \eqref{:12u} and \eqref{:12z}, we see 
\begin{align}\notag 
(\rhohat ^{N } )'(y) &=
\begin{cases}
 \frac{-1}{2\pi } \Big(\frac{2y}{4N^{2/3}} + 1 \Big)
\frac{1}{\sqrt{-y ( 1+ \frac{y}{4N ^{2/3}})}}
, & y \in (-4N^{2/3}, 0) 
,\\ 
 0 
, & y \in (-\infty -4N^{2/3}) \cup ( 0 , \infty ) 
,\end{cases}
\\\notag 
\rhohat ' (y) &= 
\begin{cases}
\frac{-1}{2\pi \sqrt{-y}} %,\quad 
& y \in (-\infty , 0) 
,\\
0& y \in (0,\infty )
.\end{cases}
\end{align}
%		and $ (\rhohat ^{N } )'(y) = \rhohat ' (y) = 0 $ for $ 0 < y < \infty $. 
From these, we see for each $ r \in \mathbb{N}$
\begin{align}&\label{:33g}
\limi{N} \sup_{\verty < r + 1 , y \ne 0 } \lvert (\rhohat ^{N } )'(y) - \rhohat ' (y) \rvert = 0 
.\end{align}
Hence, we obtain from \eqref{:12y} and \eqref{:33g} 
\begin{align} \label{:33h}&
\limi{N } \sup_{\vertx < r} 
\Big\lvert 
 \int_{ \verty < r + 1 } (\rhohat ^{N }(y)-\rhohat (y))
 \Big ( \frac{1} {x-y} - \frac{1}{-y}\Big) dy \Big\rvert = 0 
.\end{align}
Note that $ \lvert \frac{1} {x-y} - \frac{1}{-y} \rvert \le r/ (r - \verty )^2 $ for 
$\vertx < r $, $ r + 1 \le \verty $. 
Hence using this, \eqref{:12u}, and \eqref{:12z}, we apply the Lebesgue convergence theorem to deduce 
\begin{align}\notag &
\limi{N } \limi{s} \sup_{\vertx < r} 
\Big\lvert 
 \int_{ r + 1 \le \verty <s} (\rhohat ^{N }(y)-\rhohat (y)) 
 \Big ( \frac{1} {x-y} - \frac{1}{-y}\Big) dy \Big\rvert 
\\ \notag \le & 
\limi{N } \sup_{\vertx < r} \int_{ r + 1 \le \verty } 
\Big\lvert 
(\rhohat ^{N }(y)-\rhohat (y))
 \Big ( \frac{1} {x-y} - \frac{1}{-y}\Big) \Big\rvert dy 
\\ \label{:33j} \le & 
\int_{ r + 1 \le \verty } 
\limi{N } 
\Big\lvert 
(\rhohat ^{N }(y)-\rhohat (y))
 \frac{r}{ (r - \verty )^2} \Big\rvert dy = 0 
.\end{align}
Combining \eqref{:33h} and \eqref{:33j}, we obtain \eqref{:33b}. 
\end{proof}

\begin{lemma} \label{l:34} 
Let $ \beta=1,2,4$. 
There then exists a positive constant $ \Ct \label{;35} > 0 $ such that 
\begin{align} 
\label{:34a} 
\lvert \rairybeta ^{1} (x) - \rhohat (x) \rvert 
& \le \cref{;35}\Big(\frac{1}{\xvee } +\frac{{\bf 1}(\beta=4)}{ ( \xvee ) ^{1/4}}\Big) 
\text{ for all } x \in \mathbb{R}
,\\ \label{:34b}
\rairybeta ^1 (x) & 
\le \cref{;35} \sqrt{ \xvee } 
\quad \text{ for all }x \in \mathbb{R} 
.\end{align}
\end{lemma}

\begin{proof}
From \eqref{:83a} and \eqref{:12u}, we have 
\begin{align*}&
\rho_{ 2}^{1}(x)= {\Ai }'(x) ^2 - x {\Ai }(x)^2 , \quad \rhohat (x)=\frac{1_{(-\infty , 0)}(x)}{\pi} \sqrt{-x} 
.\end{align*}
Applying the estimates in \lref{l:81} to the right-hand side yields \eqref{:34a} for $\beta=2$. 

The correlation functions $\rairybeta^{1}$, $ \beta=1,4$, are represented by a quaternion determinant with kernels in \eqref{:31v} and \eqref{:31w}. Then, we see from \eqref{:31t} and \eqref{:31u} that 
\begin{align} &\notag 
\rairyone ^{1} (x)=
\rairy ^{1} (x) + \frac{1}{2} 
{\Ai }(x)( 1-\int_x^\infty {\Ai }(u)du ) 
,\\ &\notag 
\rairyfour ^{1} (x)=
\frac{1}{2^{1/3}} \rairy ^{1} (2^{2/3}x) - 
\frac{1}{2^{2/3}}{\Ai }(2^{2/3}x) \int_x^\infty {\Ai }(2^{2/3}u)du 
.\end{align}
Hence, \eqref{:34a} for $\beta=1,4$ follows from the case $ \beta=2$ 
with Lemmas \ref{l:81} and \ref{l:82}. 

We readily obtain \eqref{:34b} from \eqref{:34a} and \eqref{:12u}. 
\end{proof}

The proof of the next proposition is long, and we leave it until \sref{s:X}. 
\begin{proposition} \label{l:35} 
Let $ \beta=1,2,4$. 
Then, \eqref{:13b} holds for all $ N \in \mathbb{N}$. 
\end{proposition}

Using \pref{l:35}, we obtain the following. 
\begin{lemma} \label{l:36} 
Let $ \beta=1,2,4$. 
We then have a constant $ \Ct \label{;36} >0 $ such that for all $ N \in \mathbb{N}$ 
\begin{align} \label{:36b} 
\rairybeta ^{N ,1} (x) 
& \le \cref{;36} \sqrt{ \xvee } 
 \quad \text{ for all } x \in \mathbb{R} 
.\end{align}
\end{lemma}
\begin{proof}
We obtain \eqref{:36b} for $ x \in [-2N ^{2/3},\infty)$ from \eqref{:12z} and \eqref{:13b} immediately. 
Because of the symmetry of $ \rairybeta ^{N ,1} (x) $ 
around $ x=- N ^{2/3}$, we deduce \eqref{:36b} for all $ x \in \mathbb{R} $. 
\end{proof}

\subsection{Correlation functions of reduced Palm measures}\label{s:33}

The reduced Palm measure $ \muAx $ of $ \muA $ conditioned at $ x $ is defined as 
\begin{align}&\label{:37a}
\muAx ( A ) = \muA \big( A+ \delta _{x} \big\vert \sss ( \{ x \} )\ge 1 \big) 
.\end{align}
The $ n $-point correlation function $ \rAx ^{n}$ 
of the reduced Palm measure $ \muAx $ conditioned at $ x $ is then given by 
\begin{align}\label{:37b}&
 \rAx ^{n} (\mathbf{y} ^n )= {\rA ^{1} (x)}^{-1} {\rA ^{n+1}(x, \mathbf{y} ^n)}
\quad \text{ for a.e.\,$ x $}
.\end{align} 
We note that $\rA ^1$ is positive and $ \rA ^{n}$ are continuous and thus $ \rAx ^{n} (\mathbf{y} ^n ) $ become continuous. 
Let $ \beta=1,2,4 $. From \lref{l:32}, \eqref{:37b}, and $ \rairybeta ^{1} (x)>0$, we obtain 
for each $ n , r , s\in \mathbb{N}$ and $ i=1,\ldots, n$ 
\begin{align}& \notag 	
\limi{N } \Supr \sup_{\lvert y \rvert \le s }
\lvert \rairybetax ^{N ,n} ( y ) - \rairybetax ^{n} ( y ) \rvert = 0 
, \\&\label{:37c} 
\limi{N } \Supr \sup_{\lvert y \rvert \le s }
\lvert \partial_{i} \rairybetax ^{N , n } ( y ) - \partial _{i} \rairybetax ^{n} ( y ) \rvert =0 
.\end{align}
\begin{lemma} \label{l:37}
Let $\beta=1,2,4$. Then, for each $ r \in\mathbb{N} $, there exists $ \Ct \label{;39} >0$ such that 
\begin{align} \label{:37f} & 
\Supr \lvert \rairybetax ^1 (y)-\rairybeta ^1 (y)\rvert
 \le 
\cref{;39} \Big( \frac{1}{(\yvee ) ^{3/2} } + 
\frac{{\bf 1}(\beta\not=2)}{ (\yvee ) ^{1/4}} \Big) 
 \text{ for all $ y \in \mathbb{R} $}
.\end{align}
\end{lemma}
\begin{proof}
We set 
\begin{align*}&
L_{\beta } (x,y)=\Kairybeta (y,x)\Kairybeta (x,y) 
.\end{align*}
Then, by a direct calculation, we see that $ L_{\beta } (x,y) $ is a complex scalar quaternion 
regarded as a positive number. We therefore set 
$ L_{\beta } (x,y)=[L_{\beta } (x,y) ]^{(0)} $. 
From the relations \eqref{:20o} and \eqref{:37b}, we deduce that 
\begin{align}\label{:86a}&
\rairybeta ^{1} (y) - \rairybetax ^{1} (y)=
{\rairybeta ^{1} (x)^{-1}}
 { L_{\beta } (y,x) }
.\end{align}
Because $ \rairybeta ^{1} $ are continuous and positive, it remains to 
control $ L_{\beta } (x,y) $. 

Suppose that $\beta=2$. 
Then, \eqref{:37f} immediately follows from \eqref{:84a} (see \lref{l:84} below) and \eqref{:86a}. 

Suppose that $\beta=1$. We then see that, by \eqref{:31v} and \eqref{:31t} 
\begin{align} \notag &
 L_{1} (x,y)=J_1(x,y)J_1(y,x)
 -\partialtwo J_1(x,y) 
\Big( \int_x^y J_1(u,x)du +\frac{{\rm sign}(x-y)}{2} \Big) 
,\\ \label{:86b}&
\partialtwo J_1(x,y)= 
\partialtwo K_{2}(x,y) + \half {\Ai(x)\Ai(y)}
.\end{align}
From \lref{l:81}, \lref{l:82}, \lref{l:84}, and \lref{l:85}, for each $ r \in \mathbb{N}$, there exists a constant $ \Ct \label{;86}$ such that
\begin{align} &\notag 
\Supr \lvert J_1(x,y)J_1(y,x) \rvert \le \cref{;86} \frac{1}{\yvee }, \quad 
\Supr \lvert \int_x^y J_1(u,x)du \rvert \le \cref{;86} 
,\\ &\label{:86c}
\Supr \partialtwo J_1(x,y) \le \cref{;86} \frac{1}{(\yvee )^{1/4} }
\quad \text{ for all } y \in \mathbb{R}
.\end{align}
Combining \eqref{:86a}--\eqref{:86c}, we obtain \eqref{:37f} for $ \beta=1 $. 

Suppose $\beta=4$. We then see that by \eqref{:31w} 
\begin{align} \notag 
\frac{1}{2^{\frac{4}{3}}} L_{4}
 \Big(\frac{x}{2^{\frac{2}{3}}},\frac{y}{2^{\frac{2}{3}}}\Big) 
=& \frac{1}{4}\Big(
 J_4(x,y)J_4(y,x)-\partialtwo J_4(x,y) \int_x^y J_4(u,x)du 
\Big)
.\end{align}
Hence, we obtain \eqref{:37f} similarly as in the case $\beta=1$.
\end{proof}

\begin{lemma} \label{l:38}
Let $\beta=1,2,4$. 
Then, for each $ r \in\mathbb{N} $, there exists $ \Ct \label{;3P}>0$ independent of $N \in\mathbb{N} $ 
such that %, for all $y\in [-2N ^{2/3}, \infty) $ with $ \yrone $, 
\begin{align} \label{:38a}&
 \Supr 
\lvert \rairybetax^{N ,1}(y)-\rairybeta^{N ,1}(y) \rvert 
\le \cref{;3P} 
\Big( \frac{1}{(\yvee )^{3/2}} + 
\frac{{\bf 1} (\beta\not=2) }
{(\yvee )^{1/4}} 
\Big) 
\quad \text{ for all } x \in \mathbb{R}
.\end{align}
\end{lemma}
We shall prove \lref{l:38} in Subsection \ref{s:Y1}.

%%%%%%%%%%%%%%%%%%%%%%%%%%%%%%%%%%%%%%%%%%%%%%%%
We use the following lemma to prove \tref{l:21} and \tref{l:22}. 
%%%%%%%%%% LEMMA %%%%%%%%%%%%%%%%%%%%%%%%%%%%%%%%
\begin{lemma} \label{l:39} 
Let $\beta=1,2,4$. 
For each $ r \in \mathbb{N}$, 
\begin{align} & \label{:39z}
\limi{N } \limi{s} \Supr \Big\lvert \int_{\verty <s} 
\frac{\rairybetax ^{N , 1}(y)- \rairybetax ^{1}(y) - 
\rhohat ^{N }(y) + \rhohat (y) 
 } {x-y} dy \Big\rvert 
=0 
.\end{align}
Here the integral in \eqref{:39z} is understood in the sense of Cauchy principal value. 
\end{lemma}
%%%%%%%%%% LEMMA %%%%%%%%%%%%%%%%%%%%%%%%%%%%%%%%
\begin{proof}
From \lref{l:32} and \eqref{:37c}, we see 
\begin{align}\label{:39a}&
\limi{N} \Supr \sup_{\verty < r + 1 }
\lvert (
\rairybetax ^{N , 1} - \rairybetax ^{1} - \rhohat ^{N } + \rhohat )' (y)
 \rvert = 0 
.\end{align}
We deduce from \lref{l:32}, \eqref{:37c}, and \eqref{:39a} that 
\begin{align}\label{:39d}&
\limi{N } \Supr \Big\lvert \int_{ \verty < r + 1 } 
\frac%{\rairybetax ^{N , 1}(y)- \rairybetax ^{1}(y) - (\rairybeta ^{N , 1}(y)- \rairybeta ^{1}(y)) } 
{\rairybetax ^{N , 1}(y)- \rairybetax ^{1}(y) - \rhohat ^{N }(y) + \rhohat (y) }
{x-y}
dy \Big\rvert 
=0
.\end{align}

We divide the numerator in \eqref{:39d} as follows. 
\begin{align}\label{:39f}&
\rairybetax ^{N , 1} - \rairybetax ^{1} - \rhohat ^{N } + \rhohat 
 = 
 \rairybeta ^1 - \rairybetax ^{1} +
 \rhohat - \rairybeta ^1+ 
\rairybetax ^{N , 1} -\rairybeta ^{N ,1} 
+ 
\rairybeta ^{N ,1} -\rhohat ^{N }
.\end{align}
Using \lref{l:37} and \lref{l:34}, we immediately have for $ \yrone $
\begin{align}\label{:39g}&
\Supr 
\Big\lvert 
\frac{ \big(
 \rairybeta ^1 - \rairybetax ^{1} + \rhohat - \rairybeta ^1\big)(y) }{x-y}
\Big\rvert 
\le 
\frac{ 2(\cref{;35} +\cref{;39})} {\verty ^{1/4}\lvert x-y \rvert }
.\end{align}

From \eqref{:13b} and \lref{l:38}, we have for $ y \in [-2N^{2/3}, \infty)$ with $ \yrone $ 
\begin{align}\label{:39p}&
\Supr 
\Big\lvert 
\frac{ (
\rairybetax ^{N , 1} -\rairybeta ^{N ,1} + 
\rairybeta ^{N , 1} - \rhohat ^{N }) (y) } {x-y} 
\Big\rvert 
\le 
\frac{ \cref{;13b} + 2 \cref{;35}} {\verty ^{1/4}\lvert x-y \rvert }
.\end{align}
From \lref{l:32}, \eqref{:37c}, and \eqref{:39f}--\eqref{:39p}, we use the Lebesgue convergence theorem to obtain 
\begin{align}\notag &
\limi{N } \limi{s} \Supr 
\Big\lvert 
 \int_{ [-2N^{2/3}, \infty)\cap \{ r + 1 \le \verty < s \} }
\frac%{\rairybeta ^{N , 1}(y)- \rhohat ^{N }(y) } 
{\rairybetax ^{N , 1}(y)- \rairybetax ^{1}(y) - \rhohat ^{N }(y) + \rhohat (y) }
{x-y} 
dy 
\Big\rvert 
\\ \label{:39q} \le & 
 \int_{ r + 1 \le \verty } 
\limi{N } 
\Supr 
\Big\lvert 
\frac
{\rairybetax ^{N , 1}(y)- \rairybetax ^{1}(y) - \rhohat ^{N }(y) + \rhohat (y) }
{x-y} 
\Big\rvert dy 
 = 0 
.\end{align}
Note that $ \rairybetax ^{N , 1}(y) = \rairybeta ^{N , 2}(x,y)/\rairybeta ^{N , 1}(x)$ and that 
$\rairybeta ^{N ,1} $, $\rairybeta ^{N ,2} $, and $\rhohat ^{N } $ are symmetric around $ x=- 2 N ^{2/3}$. Using these, we see from \eqref{:39f}--\eqref{:39p} that 
\begin{align}\notag &
\limi{N } \limi{s} \Supr 
\Big\lvert 
 \int_{ ( -\infty ,-2N^{2/3}] \cap \{ r + 1 \le \verty < s \} }
\frac
{\rairybetax ^{N , 1}(y)- \rairybetax ^{1}(y) - \rhohat ^{N }(y) + \rhohat (y) }
{x-y} 
dy 
\Big\rvert 
\\ \label{:39r} \le & 
\limi{N } 
 \int_{-\infty }^{-2N^{2/3}} 
\Supr 
\Big\lvert 
\frac
{\rairybetax ^{N , 1}(y)- \rairybetax ^{1}(y) - \rhohat ^{N }(y) + \rhohat (y) }
{x-y} 
\Big\rvert dy 
 = 0 
.\end{align}
Putting \eqref{:39q} and \eqref{:39r} together, we obtain 
\begin{align}\label{:39j}&
\limi{N } \limi{s} \Supr 
\Big\lvert \int_{ r + 1 \le \verty <s} 
\frac{\rairybeta ^{N , 1}(y)- \rairybeta ^{1}(y) - (\rhohat ^{N }(y)-\rhohat (y))} {x-y}
dy \Big\rvert 
=0
.\end{align}

From \eqref{:39d} and \eqref{:39j}, we deduce \lref{l:39}. 
\end{proof}

\subsection{Key estimate} \label{s:34} 
Let $ \Kairybeta^{N } $, $ L_{\beta }^{N } $, and $ \MbN $ be as in 
\eqref{:31f} and \eqref{:31l}--\eqref{:31o}. 
As before, we regard $ \Kairybeta^{N } (x,x)$, $ L_{\beta }^{N } (x,y)$, and 
$ \MbN (x,y,z)$ as real numbers by \lref{l:31}. 
We set 
\begin{align} \notag &
I^{N }_{\beta,1}(x,s)=\int_{\lvert x-\yy \rvert >s} 
\frac{ \big\lvert L_{\beta }^{N } (x,\yy ) \big\lvert }
{\lvert x-\yy \rvert ^2}d\yy 
,\\ \notag &
 I^{N }_{\beta,2}(x,s)=\int_{%\substack
 {\lvert x-\yy \rvert >s ,\, \lvert x-\zz \rvert >s}}
\frac{\big\lvert L_{\beta }^{N } (\yy ,\zz )
\big\rvert }
{\lvert x-\yy \rvert \lvert x-\zz \rvert } d\yy d\zz 
, \\ \label{:3Rz}
& I^{N }_{\beta,3}(x,s)=\int_{%\substack
{\lvert x-\yy \rvert >s ,\, \lvert x-\zz \rvert >s} }
\frac{\big\lvert \MbN (x,y,z)\big\rvert }
{\lvert x-\yy \rvert \lvert x-\zz \rvert } d\yy d\zz 
. \end{align}
\begin{lemma} \label{l:3R}
Let $ \beta=1,2,4$. For each $ 1\le k \le 3 $, we obtain 
\begin{align}&\label{:3Rb}
\lim_{s\to\infty}\SupN \Supr 
I^{N }_{\beta,k}(x,s)=0 
.\end{align}
\end{lemma}
The proof of \lref{l:3R} is long and left until \sref{s:Z}. 
\begin{proposition}[\textrm{Key estimate}]
\label{l:3S} 
Let $ \beta=1,2,4$. 
For each $ r\in\mathbb{N} $, \eqref{:13a} holds. 
\end{proposition}
%%%%%%%%%%%%%%%%%LEMMA%%%%%%%%%%%%%%%%%%

\begin{proof}
We set $ \mathfrak{y}=\sum_i \delta _{y^i}$ and 
\begin{align}\label{:3Sx}&
w_{\beta ,s}^{N } (x,\mathfrak{y})=\sum_{\lvert x-y^i\rvert \geq s} \frac{1}{x-y^i}
 - 
\int_{\lvert x-y \rvert \geq s} \frac{\rairybetax ^{N ,1} (y)}{x-y} dy 
.\end{align}
From the standard calculation of correlation functions, we deduce that 
\begin{align}\notag &
\mathrm{Var}^{\muairybetax ^{N }} [w_{\beta ,s}^{N } (x,\cdot )] = 
\int_{\lvert x-\yy \rvert >s} \frac{ \rairybetax ^{N ,1} (\yy ) }{(x-\yy )^2}d\yy 
\\ \label{:3Sb}
& +
 \int_{%\substack
 {\lvert x-\yy \rvert >s ,\, \lvert x-\zz \rvert >s} } \frac{ \rairybetax ^{N ,2} (\yy ,\zz )} 
{(x-\yy )(x-\zz )} d\yy d\zz 
 - 
\Big(\int_{\lvert x-\yy \rvert >s} \frac{ \rairybetax ^{N ,1} (\yy ) }{(x-\yy ) }d\yy \Big)^2
.\end{align}
From the relation \eqref{:31n}, \eqref{:31o}, \eqref{:32y}, \eqref{:32z}, and \eqref{:37b}, we see that 
\begin{align}\label{:3Sc}
 \rairybetax ^{N ,1} (\yy )
=&
 \frac{ \rairybeta ^{N ,2} (x, \yy ) } 
{ \rairybeta ^{N ,1} (x) }
= \rairybeta ^{N ,1} (\yy ) - \frac{L_{\beta }^{N } (x,\yy ) } { \rairybeta ^{N ,1} (x) }
,\\ \notag 
 \rairybetax ^{N ,2} (\yy ,\zz ) =&
 \frac{ \rairybeta ^{N ,3} (x, \yy ,\zz ) } 
{ \rairybeta ^{N ,1} (x) } = 
 \rairybeta ^{N ,1} (\yy ) \rairybeta ^{N ,1} (\zz ) 
- L_{\beta }^{N } (\zz ,\yy ) 
\\ \label{:3Sd} 
- &\frac{ L_{\beta }^{N } (x,\zz ) \rairybeta ^{N ,1} (\yy )}{\rairybeta ^{N ,1} (x)}
- \frac{ L_{\beta }^{N } (x,\yy ) \rairybeta ^{N ,1} (\zz )}{\rairybeta ^{N ,1} (x)}
+\frac{\MbN (x,y,z)}
{\rairybeta ^{N ,1} (x)}
.\end{align}

Using \eqref{:3Sd}, we have 
\begin{align}\notag &
 \int_{%\substack
 {\lvert x-\yy \rvert >s ,\, \lvert x-\zz \rvert >s} }
 \frac{ \rairybetax ^{N ,2} (\yy ,\zz )} 
{(x-\yy )(x-\zz )} d\yy d\zz 
\\\notag =&
 \int_{%\substack
 {\lvert x-\yy \rvert >s ,\, \lvert x-\zz \rvert >s} }
\frac{1}{(x-\yy )(x-\zz )} 
\Big( 
 \rairybeta ^{N ,1} (\yy ) \rairybeta ^{N ,1} (\zz ) 
- L_{\beta }^{N } (\zz ,\yy ) 
\\\notag &\quad \quad 
- \frac{ L_{\beta }^{N } (x,\zz ) \rairybeta ^{N ,1} (\yy )}{\rairybeta ^{N ,1} (x)}
- \frac{ L_{\beta }^{N } (x,\yy ) \rairybeta ^{N ,1} (\zz )}{\rairybeta ^{N ,1} (x)}
+\frac{\MbN (x,y,z)}
{\rairybeta ^{N ,1} (x)}
\Big) 
 d\yy d\zz 
\\\notag =&
\Big(\int_{\lvert x-\yy \rvert >s} \frac{ \rairybeta ^{N ,1} (\yy ) }{x-\yy }d\yy \Big)^2
		% \rairybeta ^{N ,1} (\yy ) \rairybeta ^{N ,1} (\zz ) 
%		- L_{\beta }^{N } (\zz ,\yy ) 
 - \int_{%\substack
 {\lvert x-\yy \rvert >s ,\, \lvert x-\zz \rvert >s} } \frac{ L_{\beta }^{N } (\yy ,\zz ) } {(x-\yy )(x-\zz )} d\yy d\zz 
\\ \notag & \quad \quad 
- 2 \Big(
\int_{\lvert x-\yy \rvert >s} \frac{ \rairybeta ^{N ,1} (\yy ) }{x-\yy }d\yy 
\Big)\Big(
\int_{%\substack
 {\lvert x-\yy \rvert >s} } \frac{ L_{\beta }^{N } (x,\yy ) } {\rairybeta ^{N ,1} (x)(x-\yy )} d\yy 
\Big)
\\\label{:3Se}&\quad \quad 
+ \frac{1}{\rairybeta ^{N ,1} (x) }
 \int_{%\substack
 {\lvert x-\yy \rvert >s ,\, \lvert x-\zz \rvert >s} }
 \frac{\MbN (x,y,z)}
 {(x-\yy )(x-\zz )} d\yy d\zz 
.\end{align}
Using \eqref{:3Sc}, we have 
\begin{align} &\notag 
\Big(\int_{\lvert x-\yy \rvert >s} \frac{ \rairybetax ^{N ,1} (\yy ) }{(x-\yy ) }d\yy \Big)^2
=\Big(
\int_{\lvert x-\yy \rvert >s} \frac{ \rairybeta ^{N ,1} (\yy ) }{x-\yy }d\yy 
-
\int_{%\substack
 {\lvert x-\yy \rvert >s} } \frac{ L_{\beta }^{N } (x,\yy ) } {\rairybeta ^{N ,1} (x)(x-\yy )} d\yy 
\Big)^2
\\&\notag =
\Big(
\int_{\lvert x-\yy \rvert >s} \frac{ \rairybeta ^{N ,1} (\yy ) }{x-\yy }d\yy 
\Big)^2 + \Big(
\int_{%\substack
 {\lvert x-\yy \rvert >s} } \frac{ L_{\beta }^{N } (x,\yy ) } {\rairybeta ^{N ,1} (x)(x-\yy )} d\yy 
\Big)^2
\\\label{:3Sf}&
-2 
\Big(
\int_{\lvert x-\yy \rvert >s} \frac{ \rairybeta ^{N ,1} (\yy ) }{x-\yy }d\yy 
\Big)\Big(
\int_{%\substack
 {\lvert x-\yy \rvert >s} } \frac{ L_{\beta }^{N } (x,\yy ) } {\rairybeta ^{N ,1} (x)(x-\yy )} d\yy 
\Big)
.\end{align} 
Putting \eqref{:3Se} and \eqref{:3Sf} into \eqref{:3Sb}, we obtain 
\begin{align} \notag &
\mathrm{Var}^{\muairybetax ^{N }}[w_{\beta ,s}^{N } (x,\cdot )] 
= 
\int_{\lvert x-\yy \rvert >s} \frac{ \rairybeta ^{N ,1} (\yy ) }{(x-\yy )^2}d\yy 
 - \frac{1}{\rairybeta ^{N ,1} (x) }
 \int_{\lvert x-\yy \rvert >s} \frac{ L_{\beta }^{N } ( x , \yy ) }{(x-\yy )^2}d\yy 
\\ \notag 
& % - \frac{1}{\rairybeta ^{N ,1} (x) }
% \int_{\lvert x-\yy \rvert >s} \frac{ L_{\beta }^{N } ( x , \yy ) }{(x-\yy )^2}d\yy 
%\\ \notag &
 - \int_{%\substack
 {\lvert x-\yy \rvert >s ,\, \lvert x-\zz \rvert >s} } \frac{ L_{\beta }^{N } (\yy ,\zz ) } {(x-\yy )(x-\zz )} d\yy d\zz 
 - \Big(%	 \frac{1}{\rairybeta ^{N ,1} (x) }
 \int_{\lvert x-\yy \rvert >s } \frac{ L_{\beta }^{N } ( x , \yy ) } {\rairybeta ^{N ,1} (x) ( x-\yy ) } d\yy 
\Big)^2
\\ & 
+ \frac{1}{\rairybeta ^{N ,1} (x) }
 \int_{%\substack
 {\lvert x-\yy \rvert >s ,\, \lvert x-\zz \rvert >s} }
 \frac{\MbN (x,y,z)}
 {(x-\yy )(x-\zz )} d\yy d\zz 
 \label{:3Sg} 
.\end{align}
Neglecting the forth term in the right-hand side of \eqref{:3Sg} and putting $ I^{N }_{\beta,k}$ in \lref{l:3R} into \eqref{:3Sg}, we obtain 
\begin{align}& \label{:3Sh} 
 \mathrm{Var}^{\muairybetax ^{N }} [w_{\beta ,s}^{N } (x,\cdot ) ] \le 
\int_{\lvert x-\yy \rvert >s} \frac{\rairybeta ^{N ,1} (\yy )}{( x-\yy )^2}d\yy 
+ \frac{I^{N }_{\beta,1}(x,s)}{\rairybeta ^{N ,1} (x) } + I^{N }_{\beta,2}(x,s) 
+\frac{ I^{N }_{\beta,3}(x,s)}{\rairybeta ^{N ,1} (x) } 
.\end{align}
From \eqref{:36b}, we obtain 
\begin{align}\label{:3Sk}&
\limi{s} \SupN \Supr 
\int_{ \lvert x-\yy \rvert >s} \frac{\rairybeta ^{N ,1} (\yy )}{ ( x-\yy )^2} d\yy=0 
.\end{align}
Recall that $ \rairybeta ^{N ,1} $ are positive and continuous, and converge to $\rairybeta ^{1}$ uniformly on each compact set as $ N \to \infty $ by \lref{l:32}. Note that $\rairybeta ^{1}$ is locally uniformly positive. Hence, 
\begin{align}\label{:3Si}&
\inf _{ N \in \mathbb{N}} \inf_{\vertx \le r} \rairybeta ^{N ,1} (x) > 0 
\quad \text{ for each }r\in\mathbb{N} 
.\end{align}
We therefore deduce from \eqref{:3Sh}--\eqref{:3Si} and \lref{l:3R} that 
\begin{align} \label{:3Sl}&
\limi{s} \SupN \Supr 
\mathrm{Var}^{\muairybetax ^{N }} [ \lvert w_{\beta ,s}^{N } (x,\cdot ) \rvert ^{2}]=0 
.\end{align}
Note that $ E^{\muairybetax ^{N }} [w_{\beta ,s}^{N } (x,\cdot )]=0 $. Then 
\eqref{:13a} follows from \eqref{:3Sx} and \eqref{:3Sl}. 
\end{proof}

\section{Proof of the main theorems}\label{s:40}
In this section, we prove the main theorems (Theorems \ref{l:21}--\ref{l:25})
\subsection{Proof of Theorem \ref{l:21}: Unlabeled diffusions} \label{s:4}
For two measures $ \nu _1$ and $ \nu _2 $ on a measurable space 
$ (\Omega , \mathcal{B})$, we write $ \nu _1 \le \nu _2 $ 
if $ \nu _1(A)\le \nu _2(A)$ for all $ A\in\mathcal{B}$. 
For $r>0$, we set $ \Sr =\{ s \in \mathbb{R} \,;\, \lvert s \rvert < r \}$. 
For Borel measurable functions $ \map{\Phi }{\mathbb{R} }{\mathbb{R} \cup \{\infty \}} $ and 
$ \map{\Psi }{\mathbb{R} \ts \mathbb{R} }{\mathbb{R} \cup \{\infty \}} $ with $ \Psi (x,y)=\Psi (y,x) $, let
\begin{align} & \notag %
 \mathcal{H}_{r} (\mathfrak{x})=
\sum_{x^i\in \Sr } \Phi ( x^i ) + \sum_{x^i, x^j\in \Sr , i < j } \Psi ( x^i, x^j) 
, \quad \mathfrak{x}=\sum _i \delta _{x^i} 
.\end{align} 
Let $ \Lambdar $ be the Poisson RPF whose intensity is $ 1_{\Sr }dx $. 
Let $\Lambdarm=\Lambdar (\cdot \cap \Srm )$, where $ \Srm=\{\sss\in\SSS ;\sss(\Sr )=m \}$. 

A RPF $ \mu $ is said to be a $ (\Phi , \Psi ) $-quasi Gibbs measure if 
there exists a sequence of measures $ \{\mu_{r,k}^{m}\} $ on $ \SSS $ such that, for each 
$ r , k , m \in \mathbb{N}$, 
\begin{align}&\notag %\label{:41p}&
\mu_{r,k}^{m} \le \mu_{r,k+1}^{m} \text{ for all }k , \quad 
\limi{k} \mu_{r,k}^{m} = \mu (\cdot \cap \Srm ) \text{ weakly}
\end{align}
and the regular conditional probabilities 
\begin{align} &\notag %\label{:qg4} 
 \mu_{r,k,\sss}^m (d\mathfrak{x})= 
\mu_{r,k}^{m} (\pi _{ \Sr } \in d\mathfrak{x} \vert \ 
 \pi _{ \Sr^c }(\mathfrak{x})=\pi _{ \Sr^c }(\sss) ) 
\end{align}
satisfy for $\mu_{r,k}^{m}$-a.e.\,$\sss $ 
\begin{align}&\notag %\label{:qg2}
\cref{;2y}^{-1} 
e^{-\mathcal{H}_{r}(\mathfrak{x})} \Lambdarm (d\mathfrak{x}) 
 \le 
\mu_{r,k,\sss}^m (d\mathfrak{x}) \le 
\cref{;2y} 
 e^{-\mathcal{H}_{r}(\mathfrak{x})} \Lambdarm (d\mathfrak{x}) 
% \quad \text{ for $ \mu_r^m $-a.e.\! $ \sss \in \SSS $}
.\end{align}
Here, $\Ct \label{;2y}=\cref{;2y}(r,k,m,\pi _{\Sr ^c}(\sss))$ is a positive constant \cite{o.rm, o.rm2}. 
% \marginpar{q-Gibbs の定義を整理する}

%We say a probability measure $ \mu $ satisfies the assumptions \As{A1} and \As{A2} if: 
We consider the assumptions \As{A1} and \As{A2} for a RPF $ \mu $: 

\noindent 
\As{A1} $ \mu $ is a $ (\Phi , \Psi )$-quasi Gibbs measure such that 
$\Phi$ and $\Psi$ are locally bounded from below and that there exist 
upper semi-continuous $ (\hat{\Phi} , \hat{\Psi} )$ 
and positive constants $\Ct \label{;4y}$ and $\Ct \label{;4z}$ satisfying
$ \cref{;4y}^{-1} \hat{\Phi}(x) \le \Phi(x) \le \cref{;4y} \hat{\Phi}(x) $ and $
\cref{;4z}^{-1} \hat{\Psi}(x,y) \le \Psi(x,y) \le \cref{;4z}\hat{\Psi}(x,y) $. 
\\
\As{A2} \ $ \mu $ has a locally bounded $ n $-point correlation function $ \rho ^{n}$ 
with respect to the Lebesgue measure for each $ n \in \mathbb{N}$. 
\begin{lemma} %[{\cite[Lemma 2.1, Corollary 2.1]{o.rm}}]
 \label{l:41}
Let $ \mu $ be a RPF satisfying \As{A1} and \As{A2}. We then obtain 
\\ \thetag{1}
$ (\Emu ,\di ^{\mu } ) $ is closable on $ L^2(\SSS ,\mu )$. 
\\ \thetag{2} 
The closure $ (\Emu ,\domu ) $ of $ (\Emu ,\di ^{\mu } ) $ on $ \Lmu $ is a strongly local, quasi-regular Dirichlet form. 
There exists a $ \mu $-reversible diffusion $ \{\mathbb{P}_{\sss }\}_{\sss \in\SSS _{\mu }} $ properly associated with $ (\Emu ,\domu ) $ on $ L^2(\SSS ,\mu )$. 
\end{lemma}
\begin{proof} 
The closability of $ (\Emu ,\di ^{\mu } ) $ on $ L^2(\SSS ,\mu )$ 
follows from \cite[Lemma 3.6]{o.rm}. We thus obtain \thetag{1}. 
The quasi-regularity of $ (\Emu ,\domu ) $ on $ L^2(\SSS ,\mu )$ follows from 
\cite[Theorem 1]{o.dfa}. 
%\cite[Lemma 2.1, Corollary 2.1]{o.rm}. 

We next prove the strong locality of $ (\Emu ,\domu ) $ on $ L^2(\SSS ,\mu )$. 
For $ f,g \in \di ^{\mu } $, we set 
\begin{align}& \notag % \label{:20B}
\mathbb{D} _{\rR } [f,g] (\sss )= \frac{1}{2}\sum_{s^i \in \SR } 
\frac{\partial _i\check{f}(\sbm )}{\partial s^i} 
\frac{\partial _i\check{g}(\sbm )}{\partial s^i}
,\end{align}
where $ \sss=\sum_{i}\delta _{s^i} $, $ \sbm =(s^i)$, and 
$ \check{f}$ and $ \check{g}$ are defined before \eqref{:20b}. 
We set %
\begin{align*}&
\Emu _{\rR } (f,g) = \int_{\SSS } \mathbb{D}_{\rR} [f,g] d\mu 
.\end{align*}
Then $ (\Emu _{\rR } ,\di ^{\mu })$ is closable on $ \Lmu $ (see \cite[Lemmas 3.4,\,3.5]{o.rm}). 
Let $ (\underline{\mathcal{E}} _{\rR }^{\mu },\underline{\mathcal{D}}_{\rR }^{\mu }) $ 
be the closure of $ (\Emu _{\rR } , \di ^{\mu })$ on $ \Lmu $. 
Then $ (\underline{\mathcal{E}} _{\rR }^{\mu },\underline{\mathcal{D}}_{\rR }^{\mu }) $ is increasing. 
Let $ (\underline{\mathcal{E}}^{\mu },\underline{\mathcal{D}}^{\mu }) $ be the increasing limit, that is, 
\begin{align*}&
 \underline{\mathcal{E}} _{\rR }^{\mu } (f,f) \uparrow \underline{\mathcal{E}}^{\mu }(f,f), 
\rR \to \infty 
,\end{align*}
Then by definition 
\begin{align*}&
\underline{\mathcal{D}}^{\mu } = \{ f \in \cap_{\rR =1 }^{\infty} \underline{\mathcal{D}}_{\rR }^{\mu } 
\, ;\, 
\limi{\rR } \underline{\mathcal{E}} _{\rR }^{\mu } ( f , f ) < \infty 
\} 
\end{align*}
and $ (\underline{\mathcal{E}}^{\mu },\underline{\mathcal{D}}^{\mu }) $ is a closed form which is an extension of $ (\Emu ,\domu ) $ such that 
\begin{align*}&
\Emu (f,g) = \underline{\mathcal{E}}^{\mu }(f,g) \text{ for all } f,g \in \domu , \quad 
\domu \subset \underline{\mathcal{D}}^{\mu }
%(\underline{\mathcal{E}}^{\mu },\underline{\mathcal{D}}^{\mu }) \le (\Emu ,\domu )
.\end{align*}
Hence, the strong locality of $ (\Emu ,\domu ) $ follows from that of 
$ (\underline{\mathcal{E}}^{\mu },\underline{\mathcal{D}}^{\mu }) $. 
Furthermore, the strong locality of $ (\underline{\mathcal{E}}^{\mu },\underline{\mathcal{D}}^{\mu })$
follows from that of $ (\underline{\mathcal{E}} _{\rR }^{\mu },\underline{\mathcal{D}}_{\rR }^{\mu }) $ for each $ \rR \in \mathbb{N}$, which is clear by construction. 

We have thus seen that $ (\Emu ,\domu ) $ is a strongly local, quasi-regular Dirichlet form on $ \Lmu $. 
Recall that $ \SSS $ is a Lusin space. 
Note that $ 1 \in \domu $. Hence, $ (\Emu ,\domu ) $ is conservative. 
Clearly, $ (\Emu ,\domu ) $ is $ \mu $-symmetric. 
Hence from the general theory of Dirichlet forms \cite{c-f}, there exists a $ \mu $-reversible diffusion $ \{\mathbb{P}_{\sss }\}_{\sss \in\SSS _{\mu }} $ properly associated with $ (\Emu ,\domu ) $ on $ L^2(\SSS ,\mu )$. 
We thus complete the proof. 
\end{proof}

We shall check \As{A1} and \As{A2} for $ \muairybeta $. 
We introduce \As{A1$\bullet $} to prove \As{A1}. 
Condition \As{A1$\bullet $} consists of three conditions \As{A1a}, \As{A1b}, and \As{A1c} below. 
These conditions guarantee that $ \mu $ has a good finite particle approximation 
$\{\muN \}_{N \in\mathbb{N}}$, which enables us to prove the quasi-Gibbs property of $ \mu $. 

\noindent 
\As{A1a} 
$\mu$ is an RPF on $ \mathbb{R} $ with correlation functions $\{\rho^n \}_{n\in\mathbb{N} }$. 
There exists a sequence of RPFs $\{ \muN \}_{N \in\mathbb{N}}$ on 
$ \SSS $ satisfying the following. 

\noindent \thetag{1} 
The $ n $-point correlation functions $ \rN $ of $ \muN $ satisfy 
\begin{align} \label{:41a} &
\lim _{N \to \infty } \rN (\mathbf{x} )=\rho ^n (\mathbf{x} ) 
\quad \text{ for all }\mathbf{x} \in \mathbb{R} ^n \text{ for all $ n \in \mathbb{N}$}
,\\ \label{:41b} 
&\sup_{ N \in \mathbb{N}} \sup_{\mathbf{x} \in \Sr ^n} 
\rN (\mathbf{x} ) \le 
\big( \cref{;41} n ^{\cref{;41b} } \big) ^n 
\quad \text{ for all $ n, r\in \mathbb{N}$} 
,\end{align}
where $ 0 < \Ct (r) < \infty \label{;41}$ and $ 0 < \Ct ( r ) < 1 \label{;41b}$. 

\noindent \thetag{2} 
$ \muN (\sss (\mathbb{R}) \le n_{N } )=1 $ for some $ n_{N } \in \mathbb{N}$. 

\noindent \thetag{3} 
$ \muN $ is a $ (\Phi ^{N }, -\beta \log \lvert x-y \rvert ) $-canonical Gibbs measure 
for some $ \beta \in (0 , \infty ) $. 

\noindent \As{A1b} %\thetag{4} 
There exists a sequence 
$ \{ \mmm ^{N }_{\infty} \}_{N \in \mathbb{N}} $ in $ \mathbb{R} $ such that 
\begin{align}\label{:41d}&
\limi{N } \big(\Phi ^{N } (x) - \mmm ^{N }_{\infty} x \big) = \Phi (x) \quad \text{ for a.e.\! $ x $,} 
\quad 
\\ \label{:41e} &
\inf_{N \in \mathbb{N}} \inf _{ \vertx < r } \big(\Phi ^{N } (x) - \mmm ^{N }_{\infty} x \big) 
 > -\infty 
\quad \text{ for each } r \in \mathbb{N} 
.\end{align}

\noindent 
\As{A1c} There exists a sequence $ \{ \mmm ^{N }_r \}_{N , r \in \mathbb{N} }$ in $ \mathbb{R} $ such that 
\begin{align} \label{:41f}&
\limi{r} \mmm ^{N }_r=\mmm ^{N }_{\infty} 
\quad \text{ for all } N \in \mathbb{N},
\\\label{:41F}&
\supN \lvert \mmm ^{N }_r \rvert < \infty 
\quad \, \text{ for all } r \in \mathbb{N} 
,\\\label{:41g}&
\limi{r} \supN 
\Big\| 
\Big(\beta \sum_{r \le \lvert s^i \rvert < \infty } \frac{1}{s^i} \Big) + 
(\mmm ^{N }_{\infty} - \mmm ^{N }_r ) \Big\|_{L^2(\SSS , \muN ) } 
=0 
,\end{align}
where $ \sss=\sum_i \delta_{s^i}$. Moreover, $ \rNone $ satisfy 
\begin{align} \label{:41h}&
\supN \Big( \int_{1\le \vertx } \frac{1 }{ \vertx ^{2}} \, \rNone (x ) dx \Big) < \infty 
.\end{align}
\begin{lemma}[{\cite[Theorem 2.2]{o.rm2}}]\label{l:42} 
Let $\beta\in (0,\infty)$. 
Assume $ \mu $ is a RPF satisfying \As{A1$\bullet $} and \As{A2}. 
Then, $ \mu $ is a $ (\Phi , -\beta \log \lvert x-y \rvert ) $-quasi-Gibbs measure. 
\end{lemma}

We now apply \lref{l:42} to Airy RPFs $ \muairybeta $. 
We achieve this through a sequence of lemmas, in which we assume $\beta=1,2$, or $4$. 

\begin{lemma} \label{l:43} 
Let $ \muairybeta ^{N }$ be the RPF on $ \mathbb{R}$ 
whose labeled distribution is given by 
$ \mairybeta ^{N } $ in \eqref{:11c}. Set $ n_{N }=N $ and 
$ \Phi ^{N } (x)=\frac{\beta }{4} N ^{-1/3}x^2 + \beta N ^{1/3} x $.
Then, $\mu_\beta$ satisfies \As{A1a}. 
\end{lemma}

\begin{proof}
\eqref{:41a} follows from \lref{l:32}. 
\eqref{:41b} follows from \eqref{:32y} and \eqref{:32z}. Indeed, the determinantal structure 
of correlation functions gives a bound in \eqref{:41b}. We thus obtain \thetag{1}. 
Moreover, \thetag{2} and \thetag{3} of \As{A1a} are obvious from \eqref{:11c}. 
\end{proof}
\begin{lemma} \label{l:44} \thetag{1} 
The following limit exists uniformly on each compact set. 
\begin{align}\label{:44a}&
u_{\beta}(x)=\limi{s} 
\int_{\verty < s} \Big( \frac{\rairybetax ^{1}(y)}{x-y} - \frac{\rhohat (y)}{-y} \Big) dy 
.\end{align}
\thetag{2} 
Let $ u_{\beta}^{N } $ $( N \in \mathbb{N})$ be the continuous functions defined as 
\begin{align} \label{:44z}&
u_{\beta}^{N }(x)=
\int_{\mathbb{R}}\frac{\rairybetax ^{N , 1}(y)}{x-y}dy-N ^{1/3} -\frac{N ^{-1/3}}{2}x 
.\end{align}
Then, $ u_{\beta}^{N }$ converges to $ u_{\beta} $ uniformly on each compact set.
\end{lemma}
\begin{proof}
Equation \eqref{:44a} follows from \eqref{:34a} and \eqref{:37f}. Indeed, for $ s > \vertx + 1 $ 
\begin{align}& \notag 
\sup_{s \le t < \infty } 
\Big\lvert \int_{s \le \verty < t } \frac{\rairybetax ^{1}(y)}{x-y} - 
\frac{\rhohat (y)}{-y}dy \Big\rvert 
\\ \notag \le & 
\sup_{s \le t < \infty} 
\int_{s\le \verty < t } 
\Big\lvert \frac{\rairybetax ^{1}(y)}{x-y} - \frac{\rairybeta ^{1}(y)}{x-y} \Big\rvert 
+
\Big\lvert \frac{\rairybeta ^{1}(y)}{x-y} - \frac{\rhohat (y)}{x-y} \Big\rvert 
+ 
\Big\lvert \frac{\rhohat (y)}{x-y} - \frac{\rhohat (y)}{-y} \Big\rvert 
dy 
\\ \label{:44b} 
\le & 
\int_{s\le \verty < \infty }
 \frac{ \cref{;39} + \cref{;35} }{\lvert x-y \rvert \verty ^{1/4}} + \frac{ \rhohat (y)\vertx }{\lvert x-y \rvert \verty }
 dy 
= \mathcal{O} (s^{-1/4}) \quad (s\to\infty)
.\end{align}
Here, we used $\displaystyle{ \rhohat (y)= (1/\pi ) 1_{(-\infty,0]}\sqrt{-y} }$ in the last line. 
From \eqref{:44b}, we deduce \eqref{:44a}. We thus obtain \thetag{1}. 
%%%
From \eqref{:12b}, \eqref{:44a}, and \eqref{:44z}, we see that 
\begin{align}\notag &
\lvert u_{\beta}^{N }(x) - u_{\beta}(x) \rvert \\ \notag=&
\Big\lvert \Big( \int_{\mathbb{R}} 
\frac{\rairybetax ^{N , 1}(y)}{x-y}dy - 
\int_{\mathbb{R}} \frac{\rhohat ^{N }(y)}{-y} dy - \frac{N ^{-1/3}}{2}x \Big) 
 - \limi{s} \Big( \int_{\verty < s} 
\frac{\rairybetax ^{1}(y)}{x-y} - \frac{\rhohat (y)}{-y} dy \Big) \Big\rvert 
\\ \notag=& 
\Big\lvert \limi{s} 
\Big( \int_{\verty <s} 
\frac{\rairybetax ^{N , 1}(y)- \rairybetax ^{1}(y)} {x-y} - 
 \frac{\rhohat ^{N }(y)-\rhohat (y)}{-y}dy \Big) - 
\frac{N ^{-1/3}}{2}x \Big\rvert 
\\ \notag \le & 
 \limi{s} \Big\lvert 
 \int_{\verty <s} 
\frac{\rairybetax ^{N , 1}(y)- \rairybetax ^{1}(y) - (\rhohat ^{N }(y)-\rhohat (y))} {x-y} dy 
\Big\rvert 
\\ \notag & % \quad \quad \quad \quad \quad 
+ \limi{s} \Big\lvert 
 \int_{\verty <s} (\rhohat ^{N }(y)-\rhohat (y))\Big( \frac{1} {x-y} - \frac{1}{-y}\Big) dy \Big\rvert 
+ 
\frac{N ^{-1/3}}{2} \vertx 
.\end{align}
Applying \eqref{:33b} and \lref{l:39} to the last two lines, we obtain \thetag{2}. 
\end{proof}
For each $ r \in \mathbb{N}\cup\{ \infty \} $ and $ N \in \mathbb{N}$, 
let $ \mmm _{r}^{N }$ be the constant such that 
\begin{align}\label{:45z} &
\mmm _{r}^{N }=\beta \int_{\verty < r } \frac{\rairybetazero ^{N ,1} (y)}{-y} dy 
%\quad \text{ for $ r \in \mathbb{N}\cup\{ \infty \} $}
.\end{align}
\begin{lemma} \label{l:45}
Let $ \Phi ^{N } $ be as in \lref{l:43}. 
Let $ u_{\beta} $ be as in \lref{l:44} and set $ \Phi (x)=- \beta u_{\beta}(0) x $. 
Then, \As{A1b} holds. 
\end{lemma}
\begin{proof}
By definition, $ \Phi ^{N } (x)=\frac{\beta }{4} N ^{-1/3}x^2 + \beta N ^{1/3} x $. Then, 
\begin{align} \notag 
\limi{N } \big( \Phi ^{N } (x) - \mmm ^{N }_{\infty} x \big) 
&= \limi{N } \Big( \beta N ^{1/3} x - \mmm ^{N }_{\infty} x \Big) 
\\ \label{:45b}
= \lim_{N \to \infty}& - \beta u_{\beta}^{N }(0) x=- \beta u_{\beta}(0) x 
= \Phi (x) 
.\end{align}
Here, we used \lref{l:44} \thetag{2}, \eqref{:45z}, and 
$ \Phi (x)=- \beta u_{\beta}(0) x $ in the last line. 
Using \lref{l:44} \thetag{2}, we find that the convergence in \eqref{:45b} 
takes place uniformly on each compact set. Hence, \As{A1b} follows from \eqref{:45b}. 
\end{proof}

\begin{lemma} \label{l:46}
 \As{A1c} holds with $\rNone=\rairybeta ^{N ,1} $ and 
 $ \mmm ^{N }_{r} $ and $ \mmm ^{N }_{\infty} $ as in \eqref{:45z}. 
\end{lemma}
\begin{proof}
Note that $ \int_{\mathbb{R}} \rairybetazero ^{N ,1} (y) dy=N -1 $ because of \eqref{:20a}. 
Hence, \eqref{:41f} follows from 
\begin{align}& \notag %\label{:46a}
\lvert \mmm ^{N }_{\infty} - \mmm ^{N }_{r} \rvert =
\Big\lvert \beta \int_{r \le \verty } \frac{\rairybetazero ^{N ,1} (y)}{-y} dy \Big\rvert 
 \le 
\frac{\beta (N -1) }{r} \to 0 \quad \quad \text{as } r \to \infty 
.\end{align}
Recall that $\rairybeta ^{N ,1} (0) > 0 $ and 
$ \rairybetazero ^{N ,1} (x)=\rairybeta ^{N ,2} (x,0)/ \rairybeta ^N (0) $ and that 
the same hold for $ \rairybeta ^{1}$ and $ \rairybetazero ^{1} $. 
Hence, by \lref{l:32}, we have the uniform convergence 
of $ \rairybetazero ^{N ,1} $ and 
$ \partial _1 \rairybetazero ^{N ,1} $ to $ \rairybetazero ^{1} $ and 
$ \partial _1 \rairybetazero ^{1}$ on each compact set, respectively. 
This implies \eqref{:41F}. 
We obtain \eqref{:41g} from \pref{l:3S} and \eqref{:45z} easily. 
Using \eqref{:36b}, we obtain \eqref{:41h}. We thus deduce \As{A1c}. 
\end{proof}

\begin{theorem}	\label{l:47} 
Set 
$ \Phi_\beta(x)= \beta\int_x^0 u_{\beta}(y) dy $ and $ \Psi_\beta(x,y)= -\beta \log \lvert x-y \rvert $, where 
 $ \beta =1,2,4 $ and $u_\beta$ are as in \eqref{:44a}. 
Then, $ \muairybeta $ satisfies \As{A1} and \As{A2}. 
\end{theorem}
\begin{proof}
Let $\Phi (x) =-\beta u_{\beta}(0) x $. Then, $ \muairybeta $ satisfies 
\As{A1$\bullet $} for $ ( \Phi , \Psi_\beta ) $ from \lref{l:43}, \lref{l:45}, and \lref{l:46}. 
Hence, $ \muairybeta $ is a $ ( \Phi , \Psi_\beta ) $-quasi Gibbs measure. 
Note that 
$ \Phi_\beta=\Phi + \Phi_\beta -\Phi $ and that $ \Phi_\beta - \Phi $ is a continuous function. 
Then, $ \muairybeta $ is a $ ( \Phi_\beta , \Psi_\beta ) $-quasi Gibbs measure. 
Furthermore, $ ( \Phi_\beta , \Psi_\beta ) $ is upper semi-continuous. 
Hence, $ \muairybeta $ satisfies \As{A1}. 
\As{A2} is obvious from local boundedness of the kernels $ K_{\beta } $ combined with 
\eqref{:11e} and \eqref{:20o}. 
\end{proof}

From \lref{l:41} and \tref{l:47}, $ (\mathcal{E}^{\muairybeta },\di ^{\muairybeta })$ is closable on $ L^2(\SSS ,\muairybeta ) $. 
Let $ (\mathcal{E}^{\muairybeta },\dom ^{\muairybeta }) $ be the closure of 
$ (\mathcal{E}^{\muairybeta },\di ^{\muairybeta }) $ on $ L^2(\SSS ,\muairybeta ) $. 
Then, from \lref{l:41}, $ (\mathcal{E}^{\muairybeta },\dom ^{\muairybeta }) $ is a quasi-regular Dirichlet form on $ L^2(\SSS ,\muairybeta ) $. 

Let $ \capamuA $ be the $ 1 $-capacity given by the Dirichlet form 
$ (\mathcal{E}^{\muairybeta },\dom ^{\muairybeta } ) $ on $ L^2(\SSS ,\muairybeta ) $. 
See \eqref{:20p} for the definition of $ \capamuA $. 

Let $ \NN (t)=({1}/{\sqrt{2\pi }}) \int_t^{\infty} e^{- x ^2/2} dx $ be a scaled complementary error function. 
We consider the following assumption. 
% We make an assumption. 

\smallskip 

\noindent\As{A3} \; 
(1) $ \capamuA $ satisfies the following. 
\begin{align}\label{:48a}
&\capamuA (( \SSSsi )^c)=0 
.\end{align}
\noindent (2) \; 
For each $r, T \in \mathbb{N} $, 
\begin{align}&\label{:48b} 
\int_{\mathbb{R} } \NN \Big( \frac{\vertx -r}{\sqrt{T}} \Big) \rA ^{1} (x) dx < \infty 
\end{align}
and there exists a $ T>0 $ such that for each $ R>0 $, 
\begin{align}\label{:48c}&
\liminf_{r\to \infty} \NN \Big(\frac{r}{\sqrt{(r+R )T}}\Big) \, 
\Big( \int_{\vertx \le {r+R}} \rA ^1 (x)dx \Big) = 0
.\end{align}

%%%%%%%%%%%%%%%%%LEMMA%%%%%%%%%%%%%%%%%%
\begin{lemma} \label{l:48} Let $ \beta=1,2,4$. Then $ \muairybeta $ satisfies \As{A3}. 
\end{lemma}
%%%%%%%%%%%%%%%%%%LEMMA%%%%%%%%%%%%%%%%%%
\begin{proof}
By \tref{l:47}, $ \muairybeta $ is a quasi-Gibbs measure with upper semi-continuous potentials. 
Recall that the correlation functions $\rho_{ \beta}^k(\mathbf{x} ^k)$ are given
by the determinant with elements $K_{ 2}(x^i, x^j)$ if $\beta=2$,
and the quaternion determinant with elements $K_{\beta}(x^i, x^j)$ if $\beta=1,4$. 
Note that the kernels $K_{\beta}$ are locally Lipschitz continuous. 
We then deduce \eqref{:48a} from \cite[Theorem 2.1]{Osa04}. 
The properties \eqref{:48b} and \eqref{:48c} are derived from the asymptote of $\rho_{\beta}^1(x)=\mathcal{O}(\vertx ^{1/2})$ given in \eqref{:34b}. Collecting these, we see that $ \muairybeta $ satisfies \As{A3}. 
\end{proof}

Let $(\mathbb{X} , \{\mathbb{P}_{\sss } \} ) $ be a $ \muA $-reversible diffusion associated with 
 $ ( \mathcal{E}^{\muairybeta },\dom ^{\muairybeta }) $ on $ L^2(\SSS ,\muairybeta ) $. 
Then, \eqref{:48a} means that the particles never collide with one another as follows. 
\begin{align} \label{:48r}&
\mathbb{P}_{\sss } (\mathbb{X}_t \in \SSSsi 
 \text{ for all }0 \le t < \infty )=1 \quad \text{ for q.e.\! }
\sss \in \SSS 
.\end{align}
Here, q.e.\! means quasi-everywhere. We refer to \cite{fot} for the concept of quasi-everywhere. 
We see from \eqref{:48c} that no labeled particle ever explodes \cite[Lemma 10.2]{o-t.tail}; that is, 
\begin{align}\label{:48s} 
\mathbb{P}_{\sss } ( \sup_{\tzT }\lvert X_t^i\rvert < \infty \text{ for all }
T , i \in \mathbb{N} )=1 \quad \text{ for q.e.\! }\sss \in \SSS 
.\end{align}
In \cite[Lemma 10.2]{o-t.tail}, we derived \eqref{:48s} from \As{A1}--\As{A4} in \cite{o-t.tail}. 
The assumptions \As{A2}--\As{A4} in \cite{o-t.tail} correspond to \As{A1}--\As{A3} in the present paper. 
\As{A1} in \cite{o-t.tail} consists of uniform ellipticity and the existence of logarithmic derivative. 
To prove \cite[Lemma 10.2]{o-t.tail}, we do not need the existence of a logarithmic derivative, only uniform ellipticity. In the present paper, uniform ellipticity is obvious, and thus we obtain \eqref{:48s} 
because \As{A1}--\As{A3} in the present paper are satisfied. 

From \eqref{:48r} and \eqref{:48s}, we can and do take a state space $ \SSS _{\muA } $ of the unlabeled diffusion $ ( \mathbb{X}, \{ \mathbb{P}_{\sss }\} )$ in such a way that $\SSS _{\muA } \subset \SSSsi$ 
and that \eqref{:48r} and \eqref{:48s} hold for all $ \sss \in \SSS _{\muA }$. 

Note that \eqref{:48a} is equivalent to \eqref{:48r} and that 
\eqref{:48c} is a sufficient condition for \eqref{:48s}. 
Hence, all subsequent arguments follow from \eqref{:48r} and \eqref{:48s} 
instead of \eqref{:48a} and \eqref{:48c}. 

Let $ \mathbb{P}_{\muairybeta } = \int_{\SSS } \mathbb{P}_{\sss } d\muA $. 
We consider the following assumptions. 

\smallskip 
\noindent
\As{SIN} \qquad\ 
$ \mathbb{P}_{\muairybeta } ( \mathbb{X}\in \Wnex (\SSSsi))=1$. 

\noindent
\As{NBJ} \qquad 
$ \mathbb{P}_{\muairybeta } 
(m_{r,T}( \lpath (\mathbb{X}))<\infty)=1$ \quad \mbox{for each $r, T \in \mathbb{N} $}. 
\smallskip 

\noindent Here, for $ \mathbf{X}=(X^n)_{n\in\mathbb{N}} \in W (\RN )$, we set 
\begin{align*}&
m_{r,T}(\mathbf{X} )=\inf \{m\in \mathbb{N} ; \{ \min_{ t\in [0,T]} \lvert X_t^n\rvert \} 
> r \mbox{ for all $n\in\mathbb{N} $ such that $n>m$}\}
.\end{align*}

We note that, using \As{SIN}, we can construct the labeled dynamics $ \lpath (\mathbb{X}) $ in \As{NBJ}. 
Using Fubini's theorem, we see that \As{SIN} and \As{NBJ} hold for $ \mathbb{P}_{\sss }$ 
for $ \muA $-a.s.\,$ \sss $. 
Because $ \mathbb{P}_{\sss } (\cdot)$ is quasi-continuous in $ \sss $, we can refine this for q.e.\,$ \sss $. 

%%%%%%%%%%%%%%%%%LEMMA%%%%%%%%%%%%%%%%%%
\begin{lemma} \label{l:49} Let $ \beta=1,2,4$. 
Then, $ \mathbb{P}_{\muairybeta }$ satisfies \As{SIN} and \As{NBJ}. 
\end{lemma}
%%%%%%%%%%%%%%%%%%LEMMA%%%%%%%%%%%%%%%%%%
\begin{proof}
From \tref{l:47} and \lref{l:48}, $ \muairybeta $ satisfies \As{A1}--\As{A3}. 
Applying \cite[Lemmas 10.2 and 10.3]{o-t.tail}, we deduce 
\As{SIN} and \As{NBJ} from \As{A1}--\As{A3}. 
\end{proof}

\noindent 
{\em Proof of \tref{l:21}}: 
\As{A1} and \As{A2} follow from \tref{l:47}. 
Hence, the assumptions of \lref{l:41} are fulfilled. 
In particular, \thetag{1} and \thetag{2} of \tref{l:21} follow from \lref{l:41}. 

Recall that $ \SSSsip \subset \SSSsi $ and $ \muairybeta (\SSSsip ) =1 $. 
Hence, $ \mathbb{P}_{\muairybeta } ( \mathbb{X}_t \in \SSSsip ) = 1 $ for all $ t \ge 0 $. 
Because $ \mathbb{X}$ is a continuous process, this implies 
\begin{align}\label{:4Xa}&
\mathbb{P}_{\muairybeta } ( \mathbb{X}\in W (\SSSsip ))=
 \mathbb{P}_{\muairybeta } ( \mathbb{X}_t \in \SSSsip \text{ for all } t \ge 0 ) = 1 
.\end{align}
From \lref{l:49}, we have 
\begin{align}&\label{:4Xb}
 \mathbb{P}_{\muairybeta } ( \mathbb{X}\in \Wnex (\SSSsi ))=1
.\end{align}
Thus, particles neither enter nor explode. 
Hence, from \eqref{:4Xa} and \eqref{:4Xb}, 
\begin{align*}&
\mathbb{P}_{\muairybeta } ( \mathbb{X}\in \Wnex (\SSSsip ))=1 
,\end{align*}
which implies \thetag{3} of \tref{l:21}. 
\qed

%%%%%%%%%%%%%%%%%% SECTION 5 %%%%%%
\subsection{Proof of \tref{l:22}: ISDEs of Airy interacting Brownian motions}\label{s:5}

In \sref{s:5}, we specify the ISDEs describing Airy interacting Brownian motions. 
The critical concept for this is the logarithmic derivative of $ \muA $ \cite{o.isde}. 
Let $ \rA ^{1}$ be the one-point correlation function of $ \muA $ with respect to the Lebesgue measure. Let $ \muAx $ be the reduced Palm measure of $ \muA $ defined by \eqref{:37a}. 
Let $ \muAone $ be the reduced one-Campbell measure of $ \muA $ on $ \mathbb{R} \times \SSS $ defined as 
\begin{align}& \label{:50b} 
\muAone (A\ts B )=\int_{A} \rA ^{1} ( x ) \muAx (B) d x 
.\end{align}
We call $ \dmuA \in \Lloc ^1 (\muAone ) $ the logarithmic derivative of $ \muA $ if $ \dmuA $ satisfies 
\begin{align}\label{:50c}&
\int_{\SoneSSS } \dmuA f d \muAone =- \int_{\SoneSSS }\nabla _{x}f d \muAone 
\quad \text{ for all } f \in C^{\infty}_{0} (\mathbb{R} )\ot \di 
.\end{align}
%The reader can refer to \cite{bdo} for the reason and background for calling it logarithmic derivative. 
The reader can refer to \cite{bdo} for the background for the logarithmic derivative. 

We begin by calculating 
the logarithmic derivative of the finite particle approximation $ \{ \muairybeta ^{N } \} $. 
Let $ \mairybeta ^{N }$ be the density of $ \{ \muairybeta ^{N } \} $ given by \eqref{:11c}. 
Then 
\begin{align}\notag &
\nabla _{x^1} \log \mairybeta ^{N } (x^1, \ldots, x^{N}) 
\\\notag 
=& \nabla _{x^1} \log 
\frac{1}{Z}
\Big( \prod_{i<j}^{N }\lvert x^i-x^j \rvert ^\beta \Big) 
\exp\Big(-\frac{\beta}{4}\sum_{k=1}^{N } 
\left\lvert 2\sqrt{N }+\frac{x^k}{N ^{1/6} }\right\rvert ^2 \Big)
\\\label{:51o} = &
\beta \Big( 
 \sum_{k=2}^{N}\frac{1}{x^1-x^k} - N ^{1/3}-\frac{N ^{-1/3}}{2} x^1
\Big) 
.\end{align}
Taking $ x = x^1$ and $ \mathfrak{y}=\sum_{i=2}^N \delta_{x^i}$, we regard 
$ \nabla _{x^1} \log \mairybeta ^{N } (x^1, \ldots, x^{N}) $ as a function on $ \mathbb{R}\ts \SSS $. 
From \eqref{:51o} and integration by parts, we easily deduce that the logarithmic derivative 
$ \dlog^{\muairybeta ^{N }} $ of $ \muairybeta ^{N } $ is 
\begin{align}\label{:51p}&
 \dlog^{\muairybeta ^{N }} (x, \mathfrak{y} )=\beta \Big( 
 \sum_{i=1}^{N -1}\frac{1}{x-y^i} - N ^{1/3}-\frac{N ^{-1/3}}{2}x
\Big) 
.\end{align}
Here $ \mathfrak{y} = \sum_{i=1}^{N-1} \delta_{y^i}$. 

\begin{remark}\label{r:51}
The name\lq\lq logarithmic derivative" comes from \eqref{:51o} when $ N < \infty $. 
We note that \eqref{:51o} does not make sense for $ N = \infty $ because there exists no infinite product of the Lebesgue measure on $ \mathbb{R}^{\infty}$. We can still justify the logarithmic derivative by introducing the one-Campbell measure instead of the infinite product of the Lebesgue measure. 
\end{remark}

For $ \mathfrak{y} \in \SSS $, we write $ \mathfrak{y}=\sum_{i} \delta_{y^i}$. We set 
\begin{align}\label{:51q}&
\mathfrak{g} _{\beta ,s}^{N }(x, \mathfrak{y} )= 
\sum_{\lvert x-y^i\rvert <s}\frac{1}{x-y^i} - \int_{\lvert x-y \rvert <s} \frac{\rairybetax ^{N , 1}(y)}{x-y}dy 
,\\
\label{:51r}&
\mathfrak{g} _{\beta ,s}(x, \mathfrak{y} )=
 \sum_{\lvert x-y^i\rvert <s}\frac{1}{x-y^i} - \int_{\lvert x-y \rvert <s}\frac{\rairybetax ^{1}(y)}{x-y}dy 
.\end{align}
Let $ u_{\beta} $, $ u_{\beta}^{N } $, and $ w_{\beta ,s}^{N } $ 
be as in \eqref{:44a}, \eqref{:44z} and \eqref{:3Sx}, respectively. 
We see from \eqref{:51p} that 
\begin{align} \label{:51s}&
\dlog^{\muairybeta ^{N }} (x, \mathfrak{y} )=
\beta ( u_{\beta}^{N }(x) + 
\mathfrak{g} _{\beta ,s}^{N }(x, \mathfrak{y} ) + w_{\beta ,s}^{N }(x, \mathfrak{y} ) 
 ) 
.\end{align}

\begin{lemma} \label{l:51}
Assume $\beta=1,2,4$. Then, 
\begin{align} \label{:51a} & 
\lim_{N \to\infty}u_{\beta}^{N }(x)=
u_{\beta}(x) 
\mbox{ in $ \Lloc ^{2} ({\mathbb{R} },dx)$}
,\\\label{:51b}& 
\lim_{N \to\infty}\mathfrak{g} _{\beta ,s}^{N }(x, \mathfrak{y} ) 
=
\mathfrak{g} _{\beta ,s}(x, \mathfrak{y} ) 
\mbox{ in $ \Lloc ^{2} ( \muairybeta ^{[1]} )$ for any $s>0$}, 
\\ & 
\lim_{s \to\infty} \SupN 
\int_{[-r,r]\times \SSS } \lvert w^{N }_{\beta ,s}(x, \mathfrak{y} )\rvert ^{2} 
d\muairybeta ^{N ,[1]}=0 
\text{ for all } r \in \mathbb{N} 
\label{:51c}
.\end{align}
\end{lemma}
%%%%%%%%%%%%%%%%%%%%%%%%%%%%%%%%%%%%%%%%%%%%%%%%%%%%%%
\begin{proof}
Equation \eqref{:51a} is obvious because 
$ u_{\beta}^{N }$ converge to $ u_{\beta} $ uniformly on each compact set by \lref{l:44}. 
From \eqref{:51q} and \eqref{:51r}, we have 
\begin{align} \notag 
\int_{\SSS } \mathfrak{g} _{\beta ,s}^{N }(x, \mathfrak{y} ) - 
\mathfrak{g} _{\beta ,s}(x, \mathfrak{y} ) 
\, \muAxone (d \mathfrak{y})
&= -
\int_{\lvert x-y \rvert <s} \frac{\rairybetax ^{N , 1}(y)}{x-y}dy + 
\int_{\lvert x-y \rvert <s} \frac{\rairybetax ^{ 1}(y)}{x-y}dy 
\\& \label{:51d}
=
\int_{\lvert x-y \rvert <s} \frac{\rairybetax ^{ 1}(y) - \rairybetax ^{N , 1}(y)}{x-y} dy 
.\end{align}
Note that $ \rairybetax ^{ 1}(x)=\rairybetax ^{N , 1}(x)= 0 $. Then, from this and \eqref{:37c}, 
we see that for each $ r $
\begin{align}& \notag %\label{:51e} 
\limi{N} \Big( \sup_{x \in \Sr \atop 0 < \lvert x-y \rvert < s }
 \Big\lvert \frac{\rairybetax ^{ 1}(y) - \rairybetax ^{N , 1}(y)}{x-y} \Big\rvert \Big) 
\\ \notag 
 =& 
\limi{N} \Big( \sup_{x \in \Sr \atop 0 < \lvert x-y \rvert < s }
 \Big\lvert 
\frac{ \int_x^y \big( \partial \rairybetax ^{ 1}(z) - \partial \rairybetax ^{N , 1}(z) \big) dz }{x-y}
 \Big\rvert \Big) 
=0 
.\end{align}
Combining this with \eqref{:51d}, we obtain \eqref{:51b}. Eq.\,\eqref{:51c} follows from \eqref{:13a} and \eqref{:3Sx}. 
\end{proof}

\begin{theorem}\label{l:52}
For each $\beta=1,2,4 $, $ \muairybeta $ has a logarithmic derivative $ \dlog^{\muairybeta } $ 
% in $ L_{\mathrm{loc}}^2 (\muairybeta ^{[1]} ) $
 given by 
\begin{align} \label{:52a}&
 \dlog^{\muairybeta } (x,\mathfrak{y}) 
= \beta ( u_{\beta}(x) + \limi{s}\mathfrak{g} _{\beta ,s}(x, \mathfrak{y} ) ) 
.\end{align}
Furthermore, $ \dlog^{\muairybeta } (x,\mathfrak{y}) $ also satisfies 
 in $ L_{\mathrm{loc}}^2 (\muairybeta ^{[1]} ) $ 
\begin{align}\label{:52b}&
\dlog^{\muairybeta }(x,\mathfrak{y})= \beta \limi{s} 
\Big( \sum_{\lvert x-y^i\rvert <s}\frac{1}{x-y^i} - \int_{\verty <s} \frac{\rhohat (y)}{-y}dy\Big) 
.\end{align}
\end{theorem}
\begin{proof}
We use \cite[Theorem 45]{o.isde} to prove the existence of the logarithmic derivative 
$ \dlog^{\muairybeta } $ in $ L_{\mathrm{loc}}^p (\muairybeta ^{[1]} ) $ for any $ 1< p < 2$. 
\cite[Theorem 45]{o.isde} requires six assumptions, namely, 
\thetag{4.1}, \thetag{4.2}, \thetag{4.15}, \thetag{4.29}--\thetag{4.31} in \cite{o.isde}. 
To distinguish the labels in the present paper, 
we denote these assumptions by \thetag{I}--\thetag{VI}, respectively. 

\thetag{I} is the requirement such that 
$ \limi{N } \rairybeta ^{N ,n}=\rairybeta ^{n} $ uniformly on each compact set. 
Then, \thetag{I} follows from \lref{l:32}. 
\thetag{II} is the same as \eqref{:41b}, which we have already proved in \lref{l:43}. 
\thetag{III} follows from \eqref{:51a} in \lref{l:51}. 
\thetag{IV}, \thetag{V}, and \thetag{VI} follow from 
\eqref{:51s}, \eqref{:51b}, and \eqref{:51c} in \lref{l:51}, respectively. 

Thus, all the assumptions of \cite[Theorem 45]{o.isde} are satisfied. 
Hence, by \cite[Theorem 45]{o.isde}, 
the logarithmic derivative $ \dlog^{\muairybeta } \in \Lloc ^{p} ( \muairybeta ^{[1]} ) $ 
for any $ 1 < p < 2 $ exists and satisfies \eqref{:52a}. 

Note that all the convergences in \lref{l:51} are in $ \Lloc ^{2} ( \muairybeta ^{[1]} ) $ or $ \Lloc ^{2} ({\mathbb{R} },dx)$. 
Then, $ \dlog^{\muairybeta } \in \Lloc ^{2} ( \muairybeta ^{[1]} )$ from \eqref{:51s}. 
Moreover, $ \{ \mathfrak{g} _{\beta ,s}(x, \mathfrak{y} ) \}_s $ is bounded in $ \Lloc ^{2} ( \muairybeta ^{[1]} ) $. 
The convergence of $ \limi{s}\mathfrak{g} _{\beta ,s}(x, \mathfrak{y} ) $ 
follows from the martingale convergence theorem (see the proof of Theorem 45 in \cite{o.isde}). 
Hence, it takes place in $ \Lloc ^{2} ( \muairybeta ^{[1]} ) $. 
We have thus strengthened the meaning of \eqref{:52a} from $ \Lloc ^{p} ( \muairybeta ^{[1]} ) $ 
with $ 1 < p < 2$ to $ \Lloc ^{2} ( \muairybeta ^{[1]} ) $.

From \eqref{:34b} and \eqref{:37f}, we obtain 
\begin{align}& \label{:52e}
\limi{s}\Big( \int_{\verty <s}\frac{\rairybetax ^{1}(y)}{x-y}dy 
- \int_{\lvert x-y \rvert <s}\frac{\rairybetax ^{1}(y)}{x-y}dy \Big) = 0 
.\end{align}
Then, using \eqref{:44a}, \eqref{:51r}, and \eqref{:52e}, we obtain 
\begin{align} \notag &
 u_{\beta}(x) + \limi{s}\mathfrak{g} _{\beta ,s}(x, \mathfrak{y} ) 
\\ \notag 
= &
 \limi{s} \Big( 
\int_{\verty < s} \Big( \frac{\rairybetax ^{1}(y)}{x-y} - \frac{\rhohat (y)}{-y} \Big) dy + 
 \sum_{\lvert x-y^i\rvert <s}\frac{1}{x-y^i} - 
\int_{\lvert x-y \rvert <s}\frac{\rairybetax ^{1}(y)}{x-y}dy \Big) 
\\\label{:52f} 
=& 
 \limi{s}\Big( \sum_{\lvert x-y^i\rvert <s}\frac{1}{x-y^i} - \int_{\verty <s} \frac{\rhohat (y)}{-y}dy \Big) 
.\end{align}
Combining \eqref{:52a} and \eqref{:52f} implies \eqref{:52b}. 
\end{proof}
%
% We quote a general result concerning the ISDE representation of the labeled dynamics $ \mathbf{X} $. 

% We make an assumption. 
We consider the following assumption.

\smallskip 
\noindent 
\As{A4} $ \muA $ has a logarithmic derivative $ \dmuA $ in the sense of \eqref{:50c}. 

\smallskip 
\noindent 
The following lemma is a special case of \cite[Theorem 26]{o.isde} with a slight modification. 
\begin{lemma}\label{l:53} %% 
Assume \As{A1}--\As{A4}. 
Let $ \{\mathbb{P}_{\sss }\}_{\sss \in\SSS _{\muA }} $ be as in \lref{l:41}. 
Let $ \widetilde{\lab }$ be a label. Then, for $\muA \circ \widetilde{\lab }^{-1}$-a.s.\,$ \sbm $, 
the labeled path 
$ \mathbf{X}=(X^i)_{i\in\mathbb{N}} $ under $ \mathbb{P}_{\sss } \circ \widetilde{\lab }_{\mathrm{path}}^{-1} $ satisfies 
\begin{align} \label{:53a}&
dX_t^i=dB_t^i + \frac{1}{2} \dmuA (X_t^i,\mathbb{X}^{i\diamondsuit }_t)dt \quad (i\in \mathbb{N}) 
\quad \text{ and }\quad \mathbf{X} _0=\sbm 
.\end{align}
Here, $ \B= (B^i)_{i\in\mathbb{N}} $ is an $\mathbb{R} ^{\mathbb{N}}$-valued Brownian motion, and 
$ \mathbb{X}^{i\diamondsuit }_t=\sum_{j\not=i} \delta_{X_t^j}$. 
\end{lemma}
%%%%%%%%%%%%%%

\noindent 
{\em Proof of \tref{l:22}}: 
From \tref{l:47} and \lref{l:48}, we see that \As{A1}--\As{A3} hold.
By \tref{l:52}, \As{A4} holds with the logarithmic derivative given by \eqref{:52b}. 
Hence, \tref{l:22} \thetag{1} follows directly from \lref{l:53}. 
\tref{l:22} \thetag{2} is obvious because the Tracy--Widom distribution is equal to the distribution of the right-most particle under $ \muairybeta$. 
\qed

\subsection{Proof of \tref{l:23} and \tref{l:24}: Unique strong solutions}\label{s:6} 
In \sref{s:6}, we prove \tref{l:23} and \tref{l:24}. 
We use the general theory for the existence and uniqueness of strong solutions to ISDEs developed in \cite{o-t.tail}. 
To state the results of \cite{o-t.tail}, we prepare a set of assumptions and concepts. 

Let $ \TailS= \bigcap_{\rR =1}^{\infty} \sigma [\pi _{\SR ^c}]$ be the tail $ \sigma $-field of $ \SSS $ as before. We consider the following assumption. 

\smallskip
\noindent 
\As{TT} $ \muA $ is tail trivial, that is, 
$ \muA (\mathfrak{A}) \in \{ 0,1 \} $ for all $ \mathfrak{A}\in \TailS $. 
\smallskip

Let $ \mathbf{a}=\{ \ak \}_{\qq \in\mathbb{N}} $ be a sequence of 
increasing sequences $ \ak = \{ \ak (r) \}_{r\in\mathbb{N}} $ 
 of natural numbers such that $ \ak (r) < \akk (r)$ for all $ r , \qq \in \mathbb{N}$. 
 We set for $ \mathbf{a} = \{\ak \}_{\qq \in\mathbb{N}} $ 
\begin{align}\label{:95w}&
\Ka = \bigcup_{\qq =1}^{\infty} \mathfrak{K}[\ak ], \quad 
\mathfrak{K}[\ak ] =\{ \sss \, ;\, 
\sss (\Sr ) \le \ak (r) \text{ for all } r \in \mathbb{N} \} 
.\end{align}
By construction, $ \mathfrak{K}[\ak ] \subset \mathfrak{K}[\akk ]$ 
for all $ \qq \in\mathbb{N}$. It is well known that $ \mathfrak{K}[\ak ]$ is a compact set 
 in $ \SSS $ for each $ \qq \in \mathbb{N}$. We assume 

\smallskip 
\noindent \As{B1} \quad 
 $ \muA (\Ka ) = 1 $. 

\smallskip 
\noindent 
% Such a sequence $ \mathbf{a}$ with $ \muA (\Ka ) = 1 $ always exists for any RPF $ \muA $ \cite[Lemma 2.6]{o.dfa}. 
We set $ \ak ^+ (r) = \ak (r+1)$ and $ \mathbf{a}^+ =\{ \ak ^+ \}_{\qq \in \mathbb{N}} $. 
Clearly, $ \mathfrak{K}[\mathbf{a} ] \subset \mathfrak{K}[ \mathbf{a}^+ ]$.

Let $ \dlog^{\muairybeta } $ be the logarithmic derivative of $ \muairybeta $ given by \eqref{:52b}: 
\begin{align*}&
\dlog^{\muairybeta } ( \xs ) = \beta \limi{s} 
\Big( \sum_{\lvert x-y^i\rvert <s}\frac{1}{x-y^i} - \int_{\verty <s} \frac{\rhohat (y)}{-y}dy\Big) 
.\end{align*}
% Let $ b^1 ( \xs ) = (1/2) \dlog^{\muairybeta } ( \xs ) $. 
%
We set $ \mathbf{x} = (x^1,\ldots,x^m) \in \Rm $ and $ \sss =\sum_{j} \delta_{s^j}$. 
For $ m \in \mathbb{N}$, let 
 \begin{align} \notag 
 b^1 ( \xs ) = &\frac{1}{2}\dlog^{\muairybeta } ( \xs ) 
,\quad m=1 
,\\ \label{:60a}
 b ^m (\xxs ) = &
\Big( 
\frac{\beta }{2} \sum_{j \ne i}^m \frac{1}{ x^i-x^j} + 
\frac{1 }{2} \dlog^{\muairybeta } ( x^i , \sss )
\Big)_{i=1}^m 
,\quad m \ge 2 
 .\end{align}
Then the time-inhomogeneous SDEs in \eqref{:22u} are given by 
\begin{align*}&
d\mathbf{Y}_t^{m} = d\mathbf{B}_t^m + 
b^m (\mathbf{Y}_t^m, \sum_{j > m}^{\infty} \delta_{X_t^j}) dt 
.\end{align*}

We introduce an approximation of $ \Rm \ts \SSS $ consisting of compact sets. Let 
\begin{align}
\notag &% \label{:40x}&
\Ssi ^{[m]} = \{ (\mathbf{x},\sss )\in \Rm \ts \SSS \, ;\, \ulab (\mathbf{x}) + \sss \in \Ssi \} 
,\end{align}
where $ \mathbf{x}=(x^1,\ldots,x^m)$ and $ \ulab (\mathbf{x}) = \sum_{i=1}^m \delta_{x^i}$. 
For $\sss =\sum_i \delta_{s^i} \in \Ssi $, we set 
\begin{align}\notag &%\label{:41Y}&
\RRprsm = \bigcap_{j,k,l=1 \atop j \ne k }^m 
 \big\{ \mathbf{x} =(x^n) _{n=1}^m %;\, \lvert x \rvert \le \rr \}^m %\Stwoqm 
\, ;\, \lvert x^j \rvert \le \rr , \, \lvert x^j-x^k \rvert \ge 2^{-\pp } ,\ 
\inf_{i} \lvert x^l-s^i \rvert \ge 2^{-\pp } \big\} 
.\end{align}
%Here and below, we suppress m from the notations. 
We set 
\begin{align}\notag &%\label{:41y} &
\Hpqrm = \big\{ (\mathbf{x},\sss ) \in \Ssi ^{[m]} \, ;\, \ \mathbf{x} \in \RRprsm ,\ \sss \in \Kakk \big\} 
,\\ \notag &%\label{:41z}& 
\Ham = \bigcup_{\rr =1}^{\infty} \Ham _{\rr } 
, \quad 
 \Ham _{\rr } = \bigcup_{\qq =1}^{\infty} \Ham _{\qq ,\rr }, 
\quad 
\Ham _{\qq ,\rr } = \bigcup_{\pp =1}^{\infty} \Hpqrm 
.\end{align}
%Although $ \RRprsm , \Hpqrm , \Ham _{\qq ,\rr } , \Ham _{\qq } $, and $ \Ham $ depend on $ m \in \mathbb{N}$, we suppress $ m $ from the notation. 
We set $\Sr ^m = \{ \lvert x \rvert < \rr \}^m $. 
Let $ \RRprCs $ be the open kernel of $ \RRprsm $: 
\begin{align} \notag % \label{:43n}
\RRprCs = \bigcap_{j,k,l=1 \atop j \ne k }^m \big\{ \mathbf{x} \in % \{ \lvert x \rvert < \rr \}^m 
\Sr ^m 
 \, ;\, \ &
\inf_{j\not=k } \vert x^j-x^k \vert > 2^{-\pp } ,\ 
\inf_{l,i} \vert x^l-s^i\vert > 2^{-\pp } \big\} 
.\end{align}
For $ (\pqr ) $, we set %let $ \HanC $ be such that 
\begin{align} \notag & %\label{:43m}&
\Ha _{\pqr } ^{[m], \circ } := \big\{ (\mathbf{x},\sss ) \in \Ssi ^{[m]} \, ;\, \ 
 \sss \in \mathbf{x} \in \RRprCs ,\, \sss \in \Kakk \big\} 
,\\ \notag % \label{:43M}%\notag 
&\Ha ^{[m], \circ }= \bigcup_{\rr =1}^{\infty} \Ha _{\rr }^{[m], \circ } 
, \quad 
 \Ha _{\rr }^{[m], \circ } = \bigcup_{\qq =1}^{\infty} \Ha _{\qq ,\rr }^{[m], \circ }, 
\quad 
\Ha _{\qq ,\rr }^{[m], \circ } = \bigcup_{\pp =1}^{\infty} \HpqrmC 
.\end{align}

Let $ \NNN = \NNN _1 \cup \NNN _2 \cup \NNNthree $, where 
\begin{align}
\label{:42u}&
\NNN _1 = \mathbb{N}, \ 
\NNN _2 = \{ (\qq , \rr ) \ ;\, \qq , \rr \in \mathbb{N}\} , \ 
\NNNthree = \{ (\pqr ) \ ;\, \pqr \in \mathbb{N}\}
\end{align}
and for $ \nnNNN $ we define $ \nn + 1 \in \NNN $ such that 
\begin{align}
\label{:42v}&
 \nn + 1 = 
\begin{cases}
(\pp + 1, \qq , \rr ) &\text{ for $ \nn = (\pqr ) \in \NNNthree $, }\\
(\qq +1, \rr )&\text{ for $ \nn = (\qq , \rr ) \in \NNN _2$, }\\
 \rr +1 &\text{ for $ \nn = \rr \in \NNN _1$. }
\end{cases}
\end{align}
We shall take a limit in $ \nn $ along with the order $ \nn \mapsto \nn + 1$. 

We write $ \Hanm = \Ham _{\pqr }$ for $ \nn = (\pqr ) \in \NNNthree $. 
We set $ \Hanm $ for each $ \nn \in \NNN _2 \cup \NNN _1 $ similarly. 
Note that $ \Han ^{[m], \circ } \subset \Hanm $ and $ \Han ^{[m], \circ } \cup \partial \Han ^{[m], \circ } 
= \Hanm $. 

We remark that $ \Hanm $ is compact for each $ \nn \in \NNNthree $. 
% This property has critical importance in the proof of \pref{l:42}. 

Let $ \muAm $ be the $ m $-Campbell measure of $ \muA $ defined by 
\begin{align*}&
\muAm (A \times B ) = \int_{A} \rho_{\beta } ^{m} (\mathbf{x}) \muAxx ( B ) 
d\mathbf{x}
,\end{align*}
where $ \rho_{\beta } ^{m} $ is the $ m $-point correlation function of $ \muA $ and 
$ \muAxx $ is the reduced Palm measure of $ \muA $ conditioned at $ \mathbf{x} \in \Rm $. 

Let $ \{ \IIm \}_{\mm \in \mathbb{N} } $ be an increasing sequence of closed sets in $ \SmSS $. 
Let $ \map{\Pitwo }{\SmSS }{\SSS }$ be the projection such that 
$ (\mathbf{x},\sss ) \mapsto \sss $. Then by definition 
\begin{align}\notag % \label{:Hmn}
\HmnPc = &\{ \sss \in \SSS \, ; \, \Hmnc \cap \big(\Rm \ts \{ \sss \} \big) \not=\emptyset \}
.\end{align}
For $ \nnnn =(\pqr ) \in \NNNthree $ and $ \mm \in \mathbb{N} $, let $ \Ct \label{;Y3a} (\nnnn ,\mm ) $ 
be a constant such that % $ 0\le \cref{;Y3a} \le \infty $ and that 
\begin{align}\notag % \label{:97l} 
\cref{;Y3a} = \sup \{
\frac{\lvert \fg (\mathbf{x},\sss ) - \fg (\mathbf{y},\sss ) \rvert }
{\lvert \mathbf{x}-\mathbf{y} \rvert } & ; 
 \mathbf{x}\not=\mathbf{y},\, \sss \in \HmnPc 
, \\ \notag & 
\ (\mathbf{x},\sss ) , (\mathbf{y},\sss ) \in \RRprCs ,\ 
(\mathbf{x},\sss ) \simpr (\mathbf{y},\sss ) 
\} 
.\end{align}
Here $ \tbm $ is a $ \muAm $-version of $ b ^m $ and $ (\mathbf{x},\sss ) \simpr (\mathbf{y},\sss ) $ means 
$ \mathbf{x} $ and $ \mathbf{y} $ are in the same connected component of $ \RRprCs $. 
We set 
\begin{align}\notag &%\label{:97m}&
\HImnC = \bigcup_{\sss \in \Pitwo (\Hmnc )} \RRprCs \ts \{ \sss \} 
.\end{align}
We assume the following. 

\medskip 

\noindent 
\As{B2} For each $ m \in \mathbb{N} $, there exist a $ \muAm $-version $ \tbm $ of $ b ^m $ 
and an increasing sequence of closed sets $ \{ \IIm \}_{\mm \in \mathbb{N} } $ 
such that, for each $ \nnnn \in \NNNthree $ and $ \mm \in \mathbb{N} $, 
\begin{align}\label{:97n}&
\cref{;Y3a} (\nnnn ,\mm ) < \infty 
,\\ \label{:97k}&
\limi{\mm } \mathrm{Cap}^{\muAm } \Big( \HIiC \Big) = 0 
.\end{align}

The following result is a special case of \cite[Theorem 11.1]{o-t.tail}. 
The assumptions \As{A1}--\As{A4} in \cite{o-t.tail} correspond to 
\As{A4}, \As{A1}, \As{A2}, and \As{A3} in the present paper. 
All these assumptions in the present paper are stronger than or special cases of those in \cite{o-t.tail}. 
\begin{lemma}[{\cite[Theorem 11.1]{o-t.tail}}] \label{l:66}
Assume that $\muA $ %and $\lab_{\rm path}(\mathbb{X})$ under $\mathbb{P}_\muA $ 
satisfies 
\As{TT}, \As{A1}--\As{A4}, and \As{B1}--\As{B2}. 
Then, the ISDE 
\begin{align} \tag{\ref{:53a}}&
dX_t^i=dB_t^i + \frac{1}{2} \dmuA (X_t^i,\mathbb{X}^{i\diamondsuit }_t)dt \quad (i\in \mathbb{N}) 
\quad \text{ and }\quad \mathbf{X} _0=\sbm 
\end{align}
has a family of unique strong solutions $\{ \Fs \}$ starting at $\sbm =\lab (\sss) $ for $\muA $-a.s.\,$\sss$ under the constraints of \As{MF}, \As{IFC}, \As{AC} for $\muA $, \As{SIN}, and \As{NBJ}. 
\end{lemma}

\begin{remark}\label{r:66}
\lref{l:66} leads to the existence of a family of strong solutions that satisfies the conditions \As{MF}, \As{IFC}, \As{AC} for $\muA $, \As{SIN}, and \As{NBJ}. 
The uniqueness holds only under the these constraints. 
Thus, \lref{l:66} does not exclude the possibility that a solution not satisfying these constraints exist. 
It is a challenging open problem to construct such a solution. 
Indeed, we used assumption \As{AC} for $\muA $ to define the coefficient $ \frac{1}{2} \dmuA (x , \sss ) $ in ISDE \eqref{:53a}. 
\end{remark}

We shall use \lref{l:66} to prove \tref{l:23}. Hence, we check that $ \muairybeta $ satisfy the assumptions in \lref{l:66}. 
For a constant $\kappa \in ({3}/{2}, 2)$, we set 
\begin{align}
 \label{:61b}&
\Kq=
 \Big\{ \sss \in \SSS ;\, \Big( \sup_{\rR \in \mathbb{N}} 
 \frac{\sss( \SR ) }{\rR ^{\kappa}}\Big) \le q \Big\}
 \quad \mbox{and} \quad
 \K=\bigcup_{\kk =1}^{\infty} \Kq 
.\end{align}
%Throughout this section, we assume $\kappa \in ({3}/{2}, 2)$. 
%
\begin{lemma} \label{l:61} 
Let $ \beta=1,2,4$. Then, $ \muairybeta (\K )=1 $. 
In particular, $ \muairybeta $ satisfies \As{B1} for 
 $ \mathbf{a} = \{\ak \} $ such that $ \ak (r) = q r ^{\kappa }$. 
\end{lemma}

\begin{proof}
From \eqref{:34b}, we have $ \rairybeta ^1 (x) \le \cref{;35} \sqrt{\vertx +1} $. By $ 3/2 < \kappa $ 
\begin{align}\label{:61c}&
E^{\muairybeta }\Big[ \int_{S_1^c} \frac{1}{\vertx ^{\kappa} } \, \sss(dx) \Big]
=\int_{S_1^c} \frac{1}{\vertx ^{\kappa} } \, \rairybeta ^1 (x)dx <\infty 
.\end{align}
Then, from \eqref{:61c} and Fubini's theorem, we see that for $\muairybeta$-a.s.\,$\sss$, 
\begin{align}\label{:61d}&
 \int_{S_1^c} \frac{1}{\vertx ^{\kappa} } \, \sss(dx) 
 <\infty 
.\end{align}
Hence, from \eqref{:61d} and $ \SR \setminus S_1 \subset S_1^c $, 
we deduce for $\muairybeta$-a.s.\,$\sss$ that 
\begin{align} \label{:61e}
 \sup_{\rR \in \mathbb{N}} \Big \{ \frac{\sss( \SR ) }{\rR ^{\kappa}}
\Big\} 
& \le 
\sup_{\rR \in \mathbb{N}}
\Big\{ \sss (S_1) + \int_{ \SR \setminus S_1}\frac{1}{\vertx ^{\kappa} } \sss (dx) \Big\} 
 < \infty 
.\end{align}
Combining \eqref{:61b} and \eqref{:61e} yields 
\begin{align} & \notag %\label{:61f}&
\muairybeta (\K )=\lim_{q\to\infty} \muairybeta (\Kq )
=
 \lim_{q\to\infty} \muairybeta 
 \Big( \Big\{ \sss \in \SSS ;\, \Big\{ 
 \sup_{\rR \in \mathbb{N}} \frac{\sss( \SR ) }{\rR ^{\kappa}} 
 \Big\} \le q \Big\} \Big) =1
.\end{align}
This completes the proof of \lref{l:61}. 
\end{proof}

%\footnote{$ m $Campbellに変えないといけないかもしれない}
Let $ \mathcal{E}^{\muAm } $ be the bilinear form on $ \SmSS $ 
with domain $ \mathcal{D}_{\circ}^{\muAm } $ defined by 
\begin{align} \notag 
&
\mathcal{E}^{\muAm } (f,g)=
\int_{ \SmSS } 
( 
\frac{1}{2} \nabla_x f \cdot \nabla_x g + \mathbb{D} [f,g] 
) 
\, d\muAm 
,\\ \notag %\label{:64x} 
& 
\mathcal{D}_{\circ}^{\muAm } = \{ f \in C_0^{\infty}(\mathbb{R} )\ot \di \, ;\, 
\mathcal{E}^{\muAm } (f,f) < \infty ,\, f \in L^2(\muAm ) \} 
.\end{align}
We assume that $ \muA $ satisfies \As{A1}, \As{A2}, and \As{A3}. 
Then $ ( \mathcal{E}^{\muAm } , \mathcal{D}_{\circ}^{\muAm } )$ 
is closable on $ L^2(\muAm ) $ from \As{A1} and \As{A2}, that is, the quasi-Gibbs property and 
the local boundedness of the $ m $-point correlation functions for all $ m \in \mathbb{N}$. We can prove this in a similar fashion to the case $ m=0 $ (see \cite[pp. 44--46]{o.rm}). 
We denote by $ ( \mathcal{E}^{\muAm } , \mathcal{D}^{\muAm } )$ 
the closure of $ ( \mathcal{E}^{\muAm } , \mathcal{D}_{\circ}^{\muAm } )$ on 
$ L^2(\muAm ) $.

We proceed with the proof of $ \As{B2}$ for $ \muairybeta $. 
For this purpose, using Taylor expansion, we give a sufficient condition for $ \As{B2}$. 
Let 
\begin{align}&\notag % \label{:U2z}&
\mathbf{J}^{[l]} = 
\{ \mathbf{j}= (j_{k})_{k=1}^m ;\ j_{k} \in \{ 0 \} \cup \mathbb{N},\, \sum_{k=1}^m j_{k}=l \} 
.\end{align}
We set 
$ \partial _{\mathbf{j}} = \prod_{k} (\partial /\partial x_{k} )^{j_{k}} $ for 
$ \mathbf{j}=( j_{k})_{k=1}^m \in \mathbf{J}^{[l]} $, where $ (x_{k})_{k=1}^m \in \Rm $ and 
$ ({\partial }/{\partial x_{k}})^{j_{k}}$ denotes the identity if $ j_{k}=0 $. 
For $ \ellell \in \mathbb{N}$, we introduce the following:

\smallskip\noindent 
\As{C1} 
$ \chin \partial _{\mathbf{j}} b ^m \in \dom ^{\muAm } $ 
for each $ \mathbf{j} \in \cup_{l=0}^{\ellell } \mathbf{J}^{[l]} $ and $ \nn \in \NNNthree $. 

\noindent 
\As{C2} There exists a $ \muAm $-version $ \fg $ of $ b ^m $ such that 
\begin{align} \notag &% \label{:U2c} & 
\sup \big\{ \lvert \partial _{\mathbf{j}}\fg (\xxs ) \rvert ; \, (\xxs ) \in \Hanm \big\} < \infty 
\end{align}
for each $ \mathbf{j} \in \mathbf{J}^{[\ellell ]}$ and $ \nn \in \NNNthree $, where 
$ b^m $ was given by \eqref{:60a}. 

The next result is a special case of \cite[Proposition 11.2]{o-t.tail}, which is a general theorem on a sufficient condition for \As{B2}. 

\begin{lemma}[{\cite[Proposition 11.2]{o-t.tail}}] \label{l:U2}
Assume that $ \muA $ satisfies \As{A1}, \As{A2}, and \As{A3}. 
Assume that $ \muA $ satisfies \As{C1} and \As{C2} for some $ \ellell \in \mathbb{N}$. 
Then $ \muA $ satisfies \As{B2}. 
\end{lemma}

Thus, \As{B2} is reduced to \As{C1} and \As{C2}. 
We shall prove \As{C1} and \As{C2} for $ \muairybeta $ by preparing a sequence of lemmas. 
Below, we prove \As{C1} and \As{C2} only for $m=1$; the general case is similar but more involved. 

For $ p,q,r \in \mathbb{N}$, we set the subset $ \mathfrak{H}_{(p,q,r) }^{[1]}$ of $ \mathbb{R} \ts \SSS $ such that 
\begin{align}& \label{:62p}
\mathfrak{H}_{(p,q,r) }^{[1]} =\{ (x,\sss) \in \Sr \ts \Kq \, ;\, 
 2^{-p } \le \lvert x-s^i\rvert \text{ for all } i \} 
.\end{align}
Here we write $ \sss=\sum_i \delta_{ s^i }$. 
Note that $ \muairybeta (\sss (\{ x \} ) \le 1 \ \text{for all } x ) = 1 $ by \eqref{:20e}. 
From this, \eqref{:61b}, \lref{l:61}, and \eqref{:62p}, we have 
\begin{align}&\notag % \label{:62q}
 \muairybeta ^{[1]} 
\Big(\Big( \bigcup _{p,q,r \in \mathbb{N}}\mathfrak{H}_{(p,q,r) }^{[1]} \Big)^c\Big)=0 
.\end{align}
\begin{lemma} \label{l:62} 
Let $\beta=1,2,4$. The derivative in $ x $ of $ \dlog^{\muairybeta } $ is then given by 
\begin{align}\label{:62a}&
\nabla_x \dlog^{\muairybeta }(x,\sss)=
 - \beta \sum_{i}\frac{1}{(x- s^i )^2} 
.\end{align}
 The convergence in \eqref{:62a} takes place 
 in $ L^2(\mathfrak{H}_{(p,q,r) }^{[1]}, \muairybeta ^{[1]} ) $ for all $ p , q , r \in \mathbb{N}$. 
\end{lemma}
\begin{proof}
Note that 
$ \mathfrak{H}_{(p,q,r) }^{[1]} \subset \{ ( x , \sss ) ;\, x\in \Sr ,\, 2^{-p } \le \lvert x -s^i\rvert \text{ for all } i \} $. 
Then, using \eqref{:50b} and a standard calculation of correlation functions, we obtain 
\begin{align}\notag & 
\int_{\mathfrak{H}_{(p,q,r) }^{[1]} }\Big\lvert \sum_{i}\frac{1}{(x- s^i )^2}\Big\rvert ^2 d\muairybeta ^{[1]}
		%\\\notag \le \, &
\le \, 
\int_{\Sr }\rairybeta ^1 (x) dx \Big( 
\int_{\SSS } \Big\lvert \sum_{2^{- p } \le \lvert x- s^i \rvert }\frac{1}{(x-s^i)^2}\Big\rvert ^2 
d\muairybetax \Big)
\\ \notag =\, &
\int_{\Sr }\rairybeta ^1 (x) dx 
\Big( 
\int_{ 2^{- p } \le \lvert x-y \rvert , \, \lvert x-z\rvert } \frac{\rairybetax ^2 (y,z) }{(x-y)^2(x-z)^2} dydz 
+\, 
\int_{2^{- p } \le \lvert x-y \rvert } \frac{ \rairybetax ^1 (y)}{(x-y)^4} dy \Big) 
\\ \notag =\, &
\int_{ \{x \in \Sr ,\, 2^{- p } \le \lvert x-y \rvert , \, \lvert x-z\rvert \}} 
\frac{\rairybeta ^3 (x,y,z) }{(x-y)^2(x-z)^2} dxdydz 
%	\\ \notag &
 + 
\int_{\{x \in \Sr ,\, 2^{- p } \le \lvert x-y \rvert \}} \frac{ \rairybeta ^2 (x,y)}{(x-y)^4} dxdy 
.\end{align}
Using \cite[Theorem 2]{Dy.70} and \cite[(4.24), (4.27)]{ST03}, we easily see that 
\begin{align}& \label{:62c}
 \rairybeta ^2 (x,y) \le \rairybeta ^1 (x) \rairybeta ^1 (y) 
 ,\quad 
 \rairybeta ^3 (x,y,z) \le \rairybeta ^1 (x) \rairybeta ^1 (y) \rairybeta ^1 (z) 
.\end{align}
Applying \eqref{:34b} and \eqref{:62c} to the first display equation in the proof, we see that 
\begin{align*}&
\int_{\mathfrak{H}_{(p,q,r) }^{[1]} }\Big\lvert \sum_{i}\frac{1}{(x- s^i )^2}\Big\rvert ^2 d\muairybeta ^{[1]}
		%\\\notag \le \, &
< \infty 
.\end{align*}
This implies the convergence of the sum in \eqref{:62a} 
in $ L^2(\mathfrak{H}_{(p,q,r) }^{[1]}, \muairybeta ^{[1]} ) $. 
\end{proof}

\begin{lemma} \label{l:63} 
Let $\beta=1,2,4$. We can then take a $ \muairybeta ^{[1]} $-version of $ \nabla_x \dlog^{\muairybeta }(x,\sss) $ such that 
\begin{align} \label{:63z}&
\sup \{ \lvert \nabla_x \dlog^{\muairybeta }(x,\sss) \rvert ; (x,\sss)\in \mathfrak{H}_{(p,q,r) }^{[1]} \}<\infty 
\quad \text{ for each $p,q,r \in \mathbb{N} $}
.\end{align}
Here we take a $ \muairybeta ^{[1]} $-version for supremum because 
$ \nabla_x \dlog^{\muairybeta } \in L^2(\mathfrak{H}_{(p,q,r) }^{[1]}, \muairybeta ^{[1]} ) $ by \lref{l:62}. 
\end{lemma}
\begin{proof}
Write $ \sss=\sum_{i} \delta _{s^i }$ and set $ \Ct=2^{2p} (r+1)^2 \label{;63}$. Then, 
 for $ (x,\sss)\in \mathfrak{H}_{(p,q,r) }^{[1]}$, 
\begin{align} \notag 
\sum_{i} \frac{1}{(x-s^i)^2}& \le \cref{;63}
\Big( \sss (\Sr ) + 
\sum_{\rR=r }^{\infty} \sum_{s^i \in \SRR \backslash \SR } \frac{1}{ \lvert s^i \rvert ^{2}}
\Big) 
\\ \notag &
\le \cref{;63}
\Big( \sss (\Sr ) + 
\sum_{\rR=r}^{\infty} 
\frac{\sss (\SRR ) - \sss ( \SR )}{\rR ^{2 }} 
\Big)
\\ \notag 
&
= \cref{;63}
\Big( \sss (\Sr ) - \frac{\sss (\Sr ) }{ r ^2} + 
\sum_{\rR=r+1}^{\infty} \sss (\SR ) 
\big( \frac{1}{ (\rR -1)^{2 }} - \frac{1}{\rR ^{2 }} \big) \Big) 
\\ \label{:63b} 
 & \le \cref{;63}
\Big( q r^{\kappa } + \sum_{\rR=r+1 } ^{\infty} q \rR ^{\kappa } 
\big( 
\frac{1}{ (\rR -1)^{2 }} - \frac{1}{ \rR ^{2 }} 
\big) \Big) 
< \infty 
.\end{align}
Here, we used $ \sss (\Sr ) \le q r^{\kappa } $ for $ (x,\sss)\in \mathfrak{H}_{(p,q,r) }^{[1]}$. 
Thus, \eqref{:62a} and \eqref{:63b} yield \eqref{:63z}. 
\end{proof}

% \footnote{$ m $Campbellに変えないといけないかもしれない}
Let $ \mathcal{E}^{\muairybeta ^{[1]} } $ be the bilinear form on $ \mathbb{R} \ts \SSS $ 
with domain $ \mathcal{D}_{\circ}^{\muairybeta ^{[1]} } $ defined by 
\begin{align} \notag 
&
\mathcal{E}^{\muairybeta ^{[1]} } (f,g)=
\int_{ \mathbb{R} \ts \SSS } 
( 
\frac{1}{2} \nabla_x f \cdot \nabla_x g + \mathbb{D} [f,g] 
) 
\, d\muairybeta ^{[1]} 
,\\ \label{:64x} 
& 
\mathcal{D}_{\circ}^{\muairybeta ^{[1]} } = \{ f \in C_0^{\infty}(\mathbb{R} )\ot \di \, ;\, 
\mathcal{E}^{\muairybeta ^{[1]} } (f,f) < \infty ,\, f \in L^2(\muairybeta ^{[1]} ) \} 
.\end{align}
We see that $ ( \mathcal{E}^{\muairybeta ^{[1]} } , \mathcal{D}_{\circ}^{\muairybeta ^{[1]} } )$ 
is closable on $ L^2(\muairybeta ^{[1]} ) $ from \As{A1} and \As{A2}. 
% the local boundedness of $ m $-point correlation functions for all $ m \in \mathbb{N}$. 
We can prove this in a similar fashion to the case $ m=0 $ (see \cite[pp. 44--46]{o.rm}). 
We then denote by $ ( \mathcal{E}^{\muairybeta ^{[1]} } , \mathcal{D}^{\muairybeta ^{[1]} } )$ 
the closure of $ ( \mathcal{E}^{\muairybeta ^{[1]} } , \mathcal{D}_{\circ}^{\muairybeta ^{[1]} } )$ on 
$ L^2(\muairybeta ^{[1]} ) $. 
\begin{lemma}[{\cite[Lemma 11.4]{o-t.tail}}] \label{l:64}
For $ p,q,r \in \mathbb{N}$, there exists a function $ \chi_{(p,q,r)}$ with 
constants $\Ct \label{;64} $ and $\Ct \label{;65} $ such that 
$ \chi_{(p,q,r)} \in \mathcal{D}^{\muairybeta ^{[1]} }$ and that 
\begin{align} \label{:64z}
& 0 \le \chi_{(p,q,r)} \le 1, \quad
%\\\label{:64a}
\chi_{(p,q,r)}(x,\sss) 
=\begin{cases}
 1 &\text{for $ (x,\sss) \in \mathfrak{H}_{(p,q,r) }^{[1]} $ }
\\
 0 &\text{for $(x,\sss) \notin \9 $}
\ ,\end{cases}
\\ \label{:64b} & 
\cref{;64}=
 \sup \{ \lvert \nabla_x \chi_{(p,q,r)}(x,\sss)\rvert ^2 \,;\, (x,\sss) \in \mathbb{R} \ts \SSS \} < \infty 
,\\ \label{:64c}
& \cref{;65}= 
 \sup \{ \mathbb{D}[\chi_{(p,q,r)}, \chi_{(p,q,r)}](x,\sss) \,;\, (x,\sss) \in \mathbb{R} \ts \SSS \} < \infty 
.\end{align}
\end{lemma}
\begin{lemma} \label{l:65} 
Let $\beta=1,2,4$. Then, 
$ \chi _{(p,q,r) } \dlog^{\muairybeta }, 
\chi _{(p,q,r) } \nabla_x \dlog^{\muairybeta } \in \mathcal{D}^{\muairybeta ^{[1]}} $ 
 for all $p,q,r \in \mathbb{N} $. 
In particular, $ \muairybeta $ satisfies \As{C1} with $ \ellell = 1 $. 
\end{lemma}
\begin{proof} 
By \tref{l:52}, $ \dlog^{\muairybeta } \in L_{\mathrm{loc}}^2 (\muairybeta ^{[1]} ) $. 
Then, from \eqref{:64z} and \eqref{:64b}, we have 
\begin{align}\label{:65c}
 \chi _{(p,q,r) } \dlog^{\muairybeta }\in L^2( \muairybeta ^{[1]})
.\end{align}

We write $ \mathcal{E}^{\muairybeta ^{[1]} } (f)=\mathcal{E}^{\muairybeta ^{[1]} } (f,f)$ 
and $ \mathbb{D} [f]=\mathbb{D}[f,f]$. Then, from \eqref{:64x} 
\begin{align} \label{:65f}
\mathcal{E}^{\muairybeta ^{[1]} } (\chi _{(p,q,r) } \dlog^{\muairybeta }) 
 &= 
\int_{\9 } 
\frac{1}{2} \lvert \nabla_x ( \chi _{(p,q,r) } \dlog^{\muairybeta } ) \rvert ^2 
+ \mathbb{D} [\chi _{(p,q,r) } \dlog^{\muairybeta }] \, d\muairybeta ^{[1]} 
.\end{align}
Here, we used $ \chi _{(p,q,r)} (x,\sss )=0 $ if $ (x , \sss ) \notin \9 $. 
By \lref{l:63} and \eqref{:64b}, 
\begin{align}\label{:65g}&
\int_{ \9 } 
\frac{1}{2} \lvert \nabla_x ( \chi _{(p,q,r) } \dlog^{\muairybeta } ) \rvert ^2 \, d\muairybeta ^{[1]}
\\ \notag & \quad \quad \le 
\int_{ \9 } 
 \lvert \nabla_x \chi _{(p,q,r) } \rvert ^2 \lvert \dlog^{\muairybeta } \rvert ^2
+ 
 \lvert \chi _{(p,q,r) } \rvert ^2 \lvert \nabla_x \dlog^{\muairybeta } \rvert ^2 
 \, d\muairybeta ^{[1]}
<\infty
.\end{align}
From \eqref{:64z} and \eqref{:64c}, we obtain 
\begin{align} \notag %\label{:65i}
\int_{ \9 } &
\mathbb{D} [\chi _{(p,q,r) } \dlog^{\muairybeta }] \, d\muairybeta ^{[1]} 
%\\ \nonumber & 
\le 2 
\int_{ \9 } 
\mathbb{D} [\chi _{(p,q,r) }] \lvert \dlog^{\muairybeta }\rvert ^2 
 +
\lvert \chi _{(p,q,r) } \rvert ^2
\mathbb{D}[\dlog^{\muairybeta }]
\, d\muairybeta ^{[1]}
\\ \notag 
&\le 2 
\int_{ \9 } 
 \cref{;65} \lvert \dlog^{\muairybeta }\rvert ^2 
+ 
\frac{\beta^2}{2}\sum_{i} \frac{1}{ \lvert x-s^i\rvert ^4}
\, d\muairybeta ^{[1]}
\\ \notag &
= 2 
\int_{ \9 } 
 \cref{;65} \lvert \dlog^{\muairybeta }\rvert ^2 
\, d\muairybeta ^{[1]}
+ 
\int_{ \Srr \ts \mathbb{R} } 
\frac{\beta^2}{2} 
\rairybeta ^1 (x) \frac{1}{\lvert x- y \rvert ^4}\rairybetax ^1 (y) 
\, dxdy
< \infty
.\end{align} 
This together with \eqref{:65f} and \eqref{:65g} yields 
$ \mathcal{E}^{\muairybeta ^{[1]} } (\chi _{(p,q,r) } \dlog^{\muairybeta } ) <\infty $. 
Combining this with \eqref{:65c}, we obtain 
$ \chi _{(p,q,r) } \dlog^{\muairybeta }\in \mathcal{D}^{\muairybeta ^{[1]} }$. 

The second claim can be proved similarly. Hence, we omit the proof. 

The third claim is clear from the first two claims and the definition of \As{C1}. 
\end{proof}

\noindent {\em Proof of \tref{l:23}}: 
We check the assumptions in \lref{l:66} for $ \muairy $, that is, 
 \As{TT}, \As{A1}--\As{A4}, and \As{B1}--\As{B2} for the Airy$_{2} $ RPF $ \muairy $. 

All determinantal RPFs are tail trivial \cite{bqs, ly.18, o-o.tail}. 
Hence, $ \muairy $ satisfies \As{TT}. 
We have already checked \As{A1}--\As{A4} by \tref{l:47}, \lref{l:48}, and \tref{l:52}. 
% 
%We derive \As{IFC} from \cite[Propositions 11.1 and 11.2]{o-t.tail} by checking assumptions \As{B1}, \As{C1}, and \As{C2} in \cite[Section 11]{o-t.tail}. 
We obtain \As{B1} from \lref{l:61}. 
Thus, it only remains for \thetag{1} to prove \As{B2} for $ \muairy $. 
From \lref{l:U2}, we check \As{C1} and \As{C2} for this purpose. 

For $ m = 1$, \As{C1} and \As{C2} follow from \lref{l:65} and \lref{l:63}, respectively. 
The proof of the general case with $ m \in \mathbb{N}$ is the same and thus omitted. 
Hence, we have verified all the assumptions of \lref{l:66}. 
We thus obtain \thetag{1} from \lref{l:66}. 

Let $ \XB $ under $ P $ and $ (\hat{\mathbf{X}},\hat{\mathbf{B}})$ under $ \hat{P} $ 
be weak solutions of \eqref{:12t} with the same initial distribution. 
Suppose that both satisfy \As{SIN}, \As{AC} for $ \muairy $, and \As{IFC}.

Let $ \mathbf{P}_{\mathbf{s}} = P (\cdot \vert \mathbf{X}_0 = \mathbf{s})$ and 
$ \hat{\mathbf{P}}_{\mathbf{s}}= \hat{P} (\cdot \vert \hat{\mathbf{X}}_0 = \mathbf{s})$. 
Let $ \mathbf{m} = P\circ \mathbf{X}_0^{-1}= \hat{P} \circ \hat{\mathbf{X}}_0^{-1} $. 
Then 
$ \XB $ under $ \mathbf{P}_{\mathbf{s}} $ and 
$ ( \hat{\mathbf{X}}, \hat{\mathbf{B}} ) $ under $ \hat{\mathbf{P}}_{\mathbf{s}}$ 
are weak solutions to \eqref{:12t} starting at $ \mathbf{s}$ for $ \mathbf{m}$-a.s.\,$ \mathbf{s}$. 

We denote by $ \Fs $ the strong solution obtained by \thetag{1}. 
Note that both $ \mathbf{B}$ and $ \hat{\mathbf{B}}$ are Brownian motions under 
$ \mathbf{P}_{\mathbf{s}} $ and $ \hat{\mathbf{P}}_{\mathbf{s}}$, respectively. 
Hence from \dref{dfn:43} \thetag{i}, we see that 
$ \mathbf{X} = \Fs (\mathbf{B})$ $ \mathbf{P}_{\mathbf{s}}$-a.s.\,and 
$ \hat{\mathbf{X}} = \Fs (\hat{\mathbf{B}}) $ $ \hat{\mathbf{P}}_{\mathbf{s}} $-a.s.\,hold. 
Hence, 
$ \mathbf{X} $ and $ \hat{\mathbf{X}} $ satisfy \As{MF} because 
\begin{align} &\notag %\label{:69b}&
P( \Fs (\BBB ) \in A) = P ( \mathbf{X} \in A \vert \mathbf{X}_0 = \mathbf{s}), 
P( \Fs (\hat{\BBB }) \in A) = P ( \hat{\mathbf{X}} \in A \vert \hat{\mathbf{X}}_0 = \mathbf{s})
.\end{align}
From \dref{dfn:43} \thetag{ii}, we see that both 
$ \Fs (\mathbf{B}) \text{ under }\mathbf{P}_{\mathbf{s}}$ and 
$ \Fs (\hat{\mathbf{B}}) \text{ under } \hat{\mathbf{P}}_{\mathbf{s}} $ 
% \begin{align} \notag &%\label{:69a}&
% \Fs (\mathbf{B}) \text{ under }\mathbf{P}_{\mathbf{s}} \text{ and }
% \Fs (\hat{\mathbf{B}}) \text{ under } \hat{\mathbf{P}}_{\mathbf{s}} 
% \end{align}
are strong solutions of \eqref{:12t} satisfying \As{SIN}, \As{AC} for $ \muairy $, and \As{IFC}. 
Hence, using the uniqueness in \thetag{1}, we see that $ \mathbf{X}$ and $ \hat{\mathbf{X}} $ have the same distribution, which implies \thetag{2}. 

The third claim \thetag{3} follows from \eqref{:22a} and \thetag{2}. 
\qed 

\smallskip

We proceed with the proof of \tref{l:24}. 
Recall that $ \muairybeta $ is not necessarily tail trivial for $ \beta = 1,4$. 
Hence, we use the tail decomposition of $ \muairybeta $ given by \eqref{:23A}. 
Let $ \muAt (\cdot )=\muairybeta (\cdot \vert \TailS ) (\mathfrak{t} ) $ 
be the regular conditional probability of $ \muairybeta $ with respect to the tail $ \sigma $-field $ \TailS $. 
Then, by \eqref{:23A}, we have 
\begin{align*}&
\muairybeta (\cdot ) = \int_{\SSS } \muAt (\cdot ) \muairybeta (d\mathfrak{t})
\end{align*}
and $ \muAt $ satisfy \eqref{:23c}. 
In particular, $ \muAt $ is tail trivial for $ \muairybeta $-a.s.\,$ \mathfrak{t}$. 
 
%	To prove \tref{l:24}, we prepare the following. 
The next result is a special case of the general theorem \cite[Theorem 12.1]{o-t.tail}. 
\begin{lemma}[{\cite[Theorem 12.1]{o-t.tail}}] \label{l:68}
Assume that $\muA $ %and $\lab_{\rm path}(\mathbb{X})$ under $\mathbb{P}_\mu$ 
satisfies \As{A1}--\As{A4} and \As{B1}--\As{B2}. 
Then, for $ \muA $-a.s.\,$ \mathfrak{t}$, the ISDE 
\begin{align} \tag{\ref{:53a}}&
dX_t^i=dB_t^i + \frac{1}{2} \dmuA (X_t^i,\mathbb{X}^{i\diamondsuit }_t)dt \quad (i\in \mathbb{N}) 
\quad \text{ and }\quad \mathbf{X} _0=\sbm 
\end{align}
has a family of unique strong solutions $\{ \Fs \}$ starting at $\sbm =\lab (\sss) $ for $\muAt $-a.s.\,$\sss $ under the constraints of \As{MF}, \As{IFC}, \As{AC} for $\muAt $, \As{SIN}, and \As{NBJ}. 
\end{lemma}

\noindent {\em Proof of \tref{l:24}}. 
We have seen that $ \muA $ satisfies \As{A1}--\As{A4} and \As{B1}--\As{B2} in the proof of \tref{l:23}. 
Hence, we obtain \tref{l:24} \thetag{1} from \lref{l:68}. 
The proofs of \thetag{2} and \thetag{3} are the same as \tref{l:23}. Hence, we omit them. 
\qed

\subsection{Proof of \tref{l:25}: Girsanov's formula }\label{s:7}
In this section, $ \mathbf{X} ^m $ and $ \Xmstar $ are restrictions of the original $ \mathbf{X} ^m $ and $ \Xmstar $ to time interval $ [0,T]$. 
We set $ \sbm ^m =(s^i)_{i=1}^m $ for $ \sbm =(s^i)_{i\in \mathbb{N}}$ as before. 

Let $ \mathbf{X} =\lpath (\mathbb{X})$ under $ \mathbb{P}_{\sss } $ 
 be the unique strong solution of \eqref{:12t} in \tref{l:23}. 
We recall that condition \As{IFC} is satisfied by the proof of \tref{l:23} in \sref{s:6}. 
%Let $ \mathbf{X} $ be the strong solution of \eqref{:12t} in \tref{l:23}. 
Hence, for $ \mathbb{P}_{\sss } $-a.s.\,$ \Xmstar $, the first $ m $-components 
$ \mathbf{X} ^{ m } $ become the pathwise unique solution of SDE \eqref{:22u} 
 on the interval $ [0,T]$ for $\mu\circ \lab^{-1}$-a.s.\,$\sbm $. 
That is, $ \mathbf{X} ^{ m } $ under $ \mathbb{P}_{\sss } $ is a solution of the SDE:
\begin{align}\notag %
 dY_t^{m,i} =dB_t^i &+ \frac{\beta}{2} \sum_{ j \not= i}^m \frac{1}{Y_t^{m,i} - Y_t^{m,j} } dt 
\\&\label{:71b}%\notag 
+ \frac{\beta}{2} \lim_{r\to\infty} \Big( 
\sum_{ \lvert X_t^j \rvert <r ,\, j > m }\frac{1}{Y_t^{m,i} -X_t^{m,j} } 
-\int_{\vertx <r}\frac{\rhohat (x)}{-x}dx \Big) dt 
,\end{align}
where the third term vanishes for $ m=1$. We emphasize that $ \Xmstar $ is regarded as a part of the coefficients. % and $ \mathbf{X} ^m $ as a solution. 
It is known that the pathwise uniqueness of solutions implies the uniqueness of solution; that is, the distribution of solutions of SDE with the same initial condition coincide \cite[p.166, Corollary]{IW}. 

Let $ \bXmi $ be as defined before \tref{l:25}. 
Using $ \bXmi $, we rewrite \eqref{:71b} as 
\begin{align}\label{:71c}&
dY_t^{m,i} = dB_t^i + \bXmi (\mathbf{Y} _t^m , t ) dt 
\quad \text{ ($ i=1,\ldots,m$)}
.\end{align}
We set $ \mathbf{P}_{\sbm }=\mathbb{P}_{\sss} \circ \lpath ^{-1}$ as before and 
$ \Pms $ as the same as \eqref{:25x}: 
\begin{align}&\notag %\label{:71a}
 \Pms (\ \cdot \ ) =\mathbf{P}_{\sbm } (\mathbf{X} ^{m } \in \cdot \ \vert \Xmstar )
.\end{align}
Then $ \Pms $ is the distribution of the solution of \eqref{:71c}. 

Let $ \tauh ( \mathbf{w} ^m ) =
\inf \Big\{ t \wedge T \, ;\, h \le \int _0^t \lvert 
\bXm ( \mathbf{w} _u^m,u) \rvert ^2 du \Big\} $ be the stopping time given by \eqref{:25a}. 
We easily see that $ \int _0^t \lvert \bXm ( \mathbf{w} _u^m,u) \rvert ^2 du $
 is a continuous process under $ \Pms $. Hence, $ \Pms $-a.s. 
\begin{align}\notag % \label{:71C}&
\tauh ( \mathbf{w} ^m )=
\min \Big\{ t \wedge T \, ;\, h = \int _0^t \lvert \bXm ( \mathbf{w} _u^m,u) \rvert ^2 du \Big\} 
.\end{align} 
In particular, we have 
\begin{align}\notag % \label{:71f}&
 \int _0^{t \wedge \tauh ( \mathbf{w} ^m ) } \lvert \bXm ( \mathbf{w} _u^m,u) \rvert ^2 du \le h 
.\end{align}
 This implies the Novikov condition \cite[Theorem 5.3, p.152]{IW}
 \begin{align}\label{:71g}&
 E \Big[\exp \Big( \frac{1}{2} \int _0^{t \wedge \tauh ( \mathbf{w} ^m ) } \lvert \bXm ( \mathbf{w} _u^m,u) \rvert ^2 du \Big) \Big] < \infty 
 \quad \text{ for all }t
 .\end{align}
Let 
\begin{align}\label{:71h}&
M = 
\exp \Big( 
 \int_0^{\cdot \wedge \tauh } \bXm ( \mathbf{w} _u^m,u)d \mathbf{w} _u^m - 
 \frac{1}{2}\int_0^{\cdot \wedge \tauh } \lvert \bXm ( \mathbf{w} _u^m,u) \rvert ^2 du 
 \Big) 
.\end{align}
Then, from \eqref{:71c}--\eqref{:71g}, and a standard argument of Girsanov's formula (cf. \cite[p.192]{IW}), we see 
$ \mathbf{w}_{t\wedge \tauh ( \mathbf{w} ^m )}^m $ under $ M \Wsh ^{m } $ is a weak solution of 
$ \mathbf{Y}_{0} = \mathbf{s}^m $ and 
\begin{align}%& \notag 
d Y _{t\wedge \tauh ( \bYm )} ^{m,i}
 = dB_{t\wedge \tauh ( \bYm )}^i + 
\bXmi (\mathbf{Y} _{t\wedge \tauh ( \bYm )}^{m} , t\wedge \tauh ( \bYm ) ) dt 
&\label{:71i}
% \mathbf{Y}_{0\wedge \tauh ( \bYm )} = \mathbf{s}^m 
.\end{align}
As we saw at the beginning of this section, the distribution of solutions of the SDE \eqref{:71c} is unique. 
% $ \mathbf{X} ^{ m } $ is a pathwise unique solution of 
% We proved \As{IFC} in the proof of \tref{l:23}. 
Hence, the distribution of solutions of the SDE \eqref{:71i} is also unique. 
Using these, we obtain 
\begin{align}\label{:71j}&
 M \Wsh ^{m } (*) = \Pms ( \ \mathbf{X} _{\cdot \wedge \tauh ( \mathbf{X} ^m ) } ^m \in * )
.\end{align}
From \eqref{:71h} and \eqref{:71j}, we obtain \eqref{:25c}.

We prove the second statement in the case $ m =1$. The general case is derived in the same way. 
Thus we consider the drift coefficient 
\begin{align}&\notag %\label{:71d}&
\bXonei (x,t) = \frac{\beta}{2} \lim_{r\to\infty} \Big( 
\sum_{ \lvert X_t^j \rvert <r ,\, j \ge 2 }\frac{1}{x -X_t^j } 
-\int_{\verty <r}\frac{\rhohat (y)}{-y}dy \Big)
\end{align}
and the process $ \mathbf{X}^1=X^1$ satisfies the SDE 
\begin{align} &\notag %\label{:71k}&
d X_t^1 = dB_t^1 + \bXonei (X_t^1 ,t)dt 
.\end{align}
Let $ \mathbf{X}^{[1]} = (X^1,\sum_{i\ge 2}^{\infty} \delta_{X^i})$ be the one-labeled process. Then the associated Dirichlet form $ ( \mathcal{E}^{\muairybeta ^{[1]} } , \mathcal{D}^{\muairybeta ^{[1]} } )$ on $ L^2(\muairybeta ^{[1]} ) $ was given by \eqref{:64x}. 
Let $ \capaone $ be the $ 1 $-capacity given by the Dirichlet form 
 $ ( \mathcal{E}^{\muairybeta ^{[1]} } , \mathcal{D}^{\muairybeta ^{[1]} } )$ on $ L^2(\muairybeta ^{[1]} ) $.

Let $ \Kq $ and $ \mathfrak{H}_{(p,q,r) }^{[1]}$ be as in \eqref{:61b} and \eqref{:62p}, respectively. 
Then 
\begin{align}&\notag %\label{:}&
\Kq=
 \Big\{ \sss \in \SSS ;\, \Big( \sup_{\rR \in \mathbb{N}} 
 \frac{\sss( \SR ) }{\rR ^{\kappa}}\Big) \le q \Big\}
 %\quad \mbox{and} \quad \K=\bigcup_{\kk =1}^{\infty} \Kq 
,\\\notag &
\mathfrak{H}_{(p,q,r) }^{[1]}=\{ (x,\sss) \in \Sr \ts \Kq \, ;\, 
 2^{-p } \le \lvert x-s^i\rvert \text{ for all } i \} 
.\end{align}
We set $ \mathfrak{H}^{[1]} = \cup_{p,q,r \in \mathbb{N}}\mathfrak{H}_{(p,q,r) }^{[1]}$
From the condition \As{B2} in \cite[Subsection 11.2]{o-t.tail}, which is derived from \As{C1} and \As{C2} by \cite[Proposition 11.2]{o-t.tail}, we see that 
%an $ \mathfrak{H} \in \mathcal{B}(\SSS ) $ such that 
%$ \dlog ^{\muairybeta } (x,\sss)$ takes a finite value for all $ (x,\sss) \in \mathfrak{H}^{[1]} $, and that 
$ \{\delta_x + \sss ; (x,\sss) \in \mathfrak{H}^{[1]} \} \subset \SSSz $ and 
\begin{align} &\label{:71n}
\capaone ((\mathfrak{H} ^{[1]})^c )=0
,\\& \label{:71o}
\mathfrak{H}^{[1]} \subset \bigcup_{k=1}^{\infty} 
\Big\{ (x,\sss) ; 
 \Big\lvert \lim_{r\to\infty} \Big( 
 \sum_{\lvert s^j \rvert <r}\frac{1}{x - s^j } \Big) 
-\int_{\verty <r}\frac{\rhohat (y)}{-y}dy \Big\rvert < k 
\Big\} 
.\end{align}
Using \eqref{:71n} and \eqref{:71o}, we immediately obtain 
\begin{align} \notag &%\label{:71p}&
\capaone \Big( 
\Big\{
\bigcup_{k=1}^{\infty} 
\Big\{ (x,\sss) ; 
 \Big\lvert \lim_{r\to\infty} \Big( 
 \sum_{\lvert s^j \rvert <r}\frac{1}{x - s^j } \Big) 
-\int_{\verty <r}\frac{\rhohat (y)}{-y}dy \Big\rvert < k 
\Big\} \Big\}^c 
 \Big) =0
.\end{align}
Hence for each $ \rR \in \mathbb{N}$
\begin{align}\notag 
\limi{k} \capaone 
\Big( 
\Big\{ (x,\sss) ; &
 \Big\lvert \lim_{r\to\infty} \Big( 
 \sum_{\lvert s^j \rvert <r}\frac{1}{x - s^j } \Big) 
-\int_{\verty <r}\frac{\rhohat (y)}{-y}dy \Big\rvert < k 
\Big\}^c
\\\ \label{:71q}&
\cap (\SR \ts \SSS ) \Big)
= 0 
.\end{align}
Recall that the one-labeled process 
$ (\mathbf{X}^1, \sum_{i=2}^{\infty} \delta_{X^i})$ under $ \mathbb{P}_{\sss } $ is associated with the Dirichlet form 
$ ( \mathcal{E}^{\muairybeta ^{[1]} } , \mathcal{D}^{\muairybeta ^{[1]} } )$ 
 on $ L^2(\muairybeta ^{[1]} ) $ introduced before \lref{l:64} and that 
$ \capaone $ is the $ 1 $-capacity given by the Dirichlet form 
$ ( \mathcal{E}^{\muairybeta ^{[1]} } , \mathcal{D}^{\muairybeta ^{[1]} } )$ 
 on $ L^2(\muairybeta ^{[1]} ) $. Hence from \eqref{:71q}, $ \mathbf{X}^1 $ under $ \mathbb{P}_{\sss } $ satisfies 
for each $ \rR \in \mathbb{N}$ and $ \mathbb{P}_{\sss } $-a.s.
\begin{align}\label{:71r}&
 \int _0^{t \wedge \tauh ( \mathbf{X}^1 ) \wedge \sigma_R (\mathbf{X}^1) } \lvert \bXone ( \mathbf{X}_u^1,u) \rvert ^2 du \le h 
\quad \text{ for some } h \in \mathbb{N}
.\end{align}
Here $ \sigma_R (\mathbf{w}) = \inf\{ w^1 \notin \SR \} $ is the exit time of the first component $ w^1$ of $ \mathbf{w}$ from $ \SR $. 
Recall that $ \mathbf{X}^1$ does not explode. Hence using \eqref{:71r}, we deduce
\begin{align}\notag &%\label{:71s}&
\int _0^{t \wedge \tauh ( \mathbf{X}^1 ) } \lvert \bXone ( \mathbf{X}_u^1,u) \rvert ^2 du \le h 
\quad \text{ for some } h \in \mathbb{N}
.\end{align} 
This implies
\begin{align}\notag &% \label{:71t}&
\limi{h} \tauh ( \mathbf{X} ^m ) = T 
.\end{align}
Because $ \Pms = \mathbb{P}_{\sss } \circ ( \mathbf{X}^1)^{-1}$, 
% coincides the distribution of $ \mathbf{X}^1 $ under $ \mathbb{P}_{\sss } $, 
we obtain \eqref{:25d}. 
\qed

%%%%%%%%%%%%%%%%%%%%%%%
\section{Appendix 1: Estimates of the Airy function} \label{s:8}
%%%%

In this section, we collect estimates of the Airy function $ \mathrm{Ai}(x)$, 
the one-point correlation function $ \rho_{\beta}^{1}$, 
and the Airy kernels $\Kairybeta $ $ (\beta=1,2,4)$. 
% and related quantities. 

\subsection{Airy functions and related quantities}\label{s:81}
We recall that the Airy function satisfies the differential equation 
\begin{align}\label{:81z}&
\Ai''(x) - x \Ai(x)=0 \quad \text{ for all }x \in \mathbb{R} 
.\end{align}
We use the asymptotic expansions of the Airy function and its integrals \cite{Olv54a, Olv54b, VS04}. 
\begin{lemma} \label{l:81}
We see that as $ x \to \infty $, 
\begin{align} \notag 
&\Ai (x)= \frac{e^{-\frac{2}{3}x^{3/2}}}{2\pi^{1/2}x^{1/4}}
\big(1 + {\mathcal{O}}(x^{-3/2})\big) 
, \quad%\\ \notag 
\Ai '(x)= -\frac{x^{1/4}e^{-\frac{2}{3}x^{3/2}}}{2\pi^{1/2}}
\big(1 +{\mathcal{O}}(x^{-3/2})\big)
,\\ \notag 
&\Ai (-x)= \frac{1}{\pi^{1/2}x^{1/4}}
\Big( 
\cos\Big(\frac{2}{3}x^{3/2}-\frac{\pi}{4}\Big)
\left(1+{\mathcal{O}}(x^{-3/2})\right) \Big) ,
\\ \notag 
&{\Ai}'(-x)= \frac{x^{1/4}}{\pi^{1/2}}
\Big( 
\sin \Big(\frac{2}{3}x^{3/2}-\frac{\pi}{4}\Big)
\left(1+{\mathcal{O}}(x^{-3/2})\right) 
 \Big) 
.\end{align}
\end{lemma}

\begin{lemma} \label{l:82}
As $ x \to \infty $, we have 
\begin{align*} 
&\int_0^x \Ai(u)du =\frac{1}{3}- \frac{\exp (-\frac{2}{3}x^{3/2}) } {2\pi ^{1/2} x^{3/4}}
\left(1+ o(1) \right)
,\\ &\int_{-x}^0 \Ai (u)du =\frac{2}{3} +{\mathcal{O}}\left(x^{-3/4}\right) 
, \quad %\\ & 
 1-\int_{-x}^\infty \Ai(u)du={\mathcal{O}}\left(x^{-3/4}\right)
.\end{align*}
\end{lemma}

We apply the above asymptotic behaviors to examine the one-point correlation function $ \rho_{\beta}^{1}$ 
and the Airy kernels $\Kairybeta $ $ (\beta=1,2,4)$. 

From \eqref{:11g}, \eqref{:81z}, and the continuity of $ \Kairy (x,y) $, we see that 
\begin{align} &\label{:83a} 
\rho_{ 2}^{1}(x)=K_{ 2}(x,x) 
= \Ai'(x)^2 - \Ai''(x) \Ai (x)
= {\Ai }'(x) ^2 - x {\Ai }(x)^2
.\end{align}
Using these, we obtain asymptotics of the one-point correlation function of $ \muairy $. 
\begin{lemma} \label{l:83} 
As $ x \to \infty $, we have 
\begin{align*}
&\rho_{ 2}^{1}(x) ={\mathcal{O}}\big(e^{-\frac{4}{3}x^{3/2}} \big), 
 \quad \rho_{ 2}^{1}(-x) =\frac{\sqrt{x}}{\pi} \big( 1 + {\mathcal{O}}(x^{-3/2}) \big) 
.\end{align*}
\end{lemma}
\begin{proof} 
Combining \lref{l:81} and \eqref{:83a} yields \lref{l:83}. 
\end{proof}

Using \lref{l:81} and \lref{l:82}, we obtain the estimates of the determinantal kernel $ \Kairy $. 
We use \lref{l:84} in the proof of \lref{l:37}. 
\begin{lemma} \label{l:84}
For each $ r \in \mathbb{N}$, there exists a constant $ \Ct \label{;84}$ such that 
% 	Let $x\in \mathbb{R} $ be fixed. It then holds that as $ \verty \to\infty $, 
\begin{align}\label{:84a}
 \Supr \lvert \Kairy (x,y) \rvert &\le \cref{;84}\frac{1}{(\yvee )^{3/4}} \quad \text{ for all } \ y \in \mathbb{R}
,\\
 \Supr \lvert \partialtwo \Kairy (x,y) \rvert & \le \cref{;84}\frac{1}{(\yvee )^{1/4}} \quad \text{ for all } y \in \mathbb{R}
\label{:84b}
.\end{align}
\end{lemma}
\begin{proof} 
Note that $ \Kairy (x,y)$ and $ \partialtwo \Kairy (x,y) $ are bounded on each compact set. 
Then \eqref{:84a} follows from \eqref{:11g} and \lref{l:81}. 
From \eqref{:81z}, we deduce that 
\begin{align*}&
\partialtwo \Kairy (x,y)
=\frac{-\Ai'(x)\Ai'(y) +y\Ai(x)\Ai(y)}{x-y} +\frac{\Kairy(x,y)}{x-y}
.\end{align*}
Hence, \eqref{:84b} follows from \lref{l:81}. 
\end{proof}

We need the following estimates for the asymptotics of the one-point correlation function of $ \muairybeta $, $ \beta = 1,4$. We use \lref{l:85} in the proof of \lref{l:37}. 
\begin{lemma}\label{l:85}
There exists $\Ct \label{;T4Na} > 0 $ such that
\begin{align}
&
\sup_{\xX \in \mathbb{R} }\Big\lvert \int_\yY ^\xX \Kairy(u,\yY )du\Big\rvert \le 
 \cref{;T4Na} ( 1 \vee \log \lvert \yY \rvert ), \quad \yY \in\mathbb{R} 
\label{:85a}
.\end{align}
\end{lemma}

Since the proof of \lref{l:85} is long, we will prove \lref{l:85} in the next two subsections. 

%%%%%%%%%%%%%%%%%%%%%%%%%%%%%%%%%%%%%%%%%%%%%%%
\subsection{Subsidiary lemmas}\label{s:16} 

In Subsection \ref{s:16}, we prepare two general lemmas for the rest of the paper. 
\begin{lemma} \label{l:161}
Let $ f $ be non-negative decreasing (or increasing) on $[a,b]\subset \mathbb{R} $. Then, 
\begin{align}\label{:161b}
\Big\lvert \int_a^b f(v) \sin v dv \Big\rvert &\le 2\pi (f(a)\vee f(b)) 
.\end{align}
\end{lemma}
\begin{proof} 
Let $ f $ be decreasing. 
Let $k\in\mathbb{N}\cup\{ 0 \} $ such that $ 2\pi k \le b - a < 2\pi (k+1)$. Then, 
\begin{align}\label{:161c}&
 \Big\lvert \int_{a+2\pi k }^b f(v) \sin v dv\Big\rvert 
 \le ( b- a-2\pi k ) f(a+2\pi k ) \le 2\pi f(a+2\pi k ) 
,\\\notag 
& 
 \Big\lvert \int_a^{a+2\pi k } f(v) \sin v dv \Big\rvert 
\\ \notag 
 = & 
 \Big\lvert 
\int_{0}^{\pi} \sum_{\ell=0}^{k-1} ( f(v+a+2\pi \ell ) - f(v+a+2\pi \ell +\pi )) 
\sin (v+a) dv \Big\rvert 
\\ \label{:161d} 
 \le & \int_{0}^{\pi} ( f(v+a) - f(v+a+2\pi k - \pi ) ) dv \le \pi (f(a)-f(a+2\pi k ))
.\end{align}
Combining \eqref{:161c} and \eqref{:161d} and using the property that $ f $ is decreasing, we obtain 
\begin{align*}&
\Big\lvert \int_{a}^b f(v) \sin v dv\Big\rvert \le 2\pi f(a+2\pi k ) + \pi (f(a)-f(a+2\pi k )) 
\le 2\pi f(a)
.\end{align*}
This implies \eqref{:161b}. The proof for increasing $ f $ is similar and thus omitted. 
\end{proof}

\begin{lemma} \label{l:162}
Let $g \in C^1 ([\yY ,\yY +1])$. Then, for any $ 0 < \delta \le 1 $, 
\begin{align}\label{:162a}
\sup_{\xX \in [\yY , \yY + 1]}
\Big\lvert \int_\yY ^{\xX } \frac{\int_\yY ^u g(v) dv}{u-\yY }du \Big\rvert \le A(1 - \log \delta) + B\delta
,\end{align}
where $A= \sup_{u \in [\yY ,\yY +1] } \lvert \int_\yY ^u g(v) dv \rvert $ and 
$B= (1/2)\sup_{u \in [\yY ,\yY +1] } \lvert g' (u) \rvert $. 
\end{lemma} 
\begin{proof}
We see $\int_\yY ^u g(v) dv = g(u)(u-\yY ) - (1/2)g'(z)(u-\yY )^2$ for some $ z \in [\yY ,u]$. Then, 
\begin{align*}&
\Big\lvert \int_\yY ^{ \min\{\yY + \delta , \xX \} } 
\frac{\int_\yY ^u g(v) dv}{u-\yY }du \Big\rvert \le 
 \Big\lvert \int_\yY ^{ \min\{\yY + \delta , \xX \} } g(u)du \Big\rvert + 
\frac{1}{2}\sup_{z \in [\yY , \yY +1]} \lvert g' (z)\rvert \delta 
 \le A + B\delta
,\\ \notag &
\Big\lvert \int_ { \min\{\yY + \delta , \xX \} } ^{\yY + \xX }
\frac{\int_\yY ^u g(v) dv}{u-\yY }du \Big\rvert \le A 
\Big\lvert \int_ { \min\{\yY + \delta , \xX \} } ^{\yY + \xX }\frac{du}{u-\yY } \Big\rvert \le 
 A( - \log \delta ) 
.\end{align*}
From these two inequalities, we obtain \eqref{:162a}. 
\end{proof}

\subsection{Proof of \lref{l:85}} \label{s:17}

In Subsection \ref{s:17}, we complete the proof of \lref{l:85}. 

We begin by introducing three necessary quantities. 
\begin{align}\notag % \label{:17a}
&\HHa (x,y)= 
\Big(\int_\yY ^\xX \mathbf{1}(\lvert u-\yY \rvert \le 1) \frac{\int_{\yY }^{u} \Ai'(v)dv}{u-\yY }du \Big) 
\Ai'(\yY ),
\\\notag %\label{:17b}
& \HHb (x,y)= -
\Big( \int_{\yY }^{\xX } \mathbf{1}(\lvert u-\yY \rvert \le 1)
 \frac{ \int_\yY ^u \Ai'' (v)dv }{u-\yY } du \Big)
\Ai(\yY ),
\\\notag % \label{:17c}
& \HHc (x,y)=\int_{\yY }^\xX \mathbf{1}(\lvert u-\yY \rvert > 1) K_2(u,\yY )du 
.\end{align}
The next lemma is critical in the proof of \lref{l:85} and \lref{l:196}. 
\begin{lemma} \label{l:17}
There exists a constant $\Ct \label{;I_1} >0$ such that for $ p = 1,2,3 $, 
\begin{align}\label{:17e}& 
\sup_{ \xX \in \mathbb{R} } \lvert \HHp (x,y)\rvert \le \cref{;I_1}( 1 \vee \log \lvert \yY \rvert ), \quad \yY \in (-\infty, -1] 
.\end{align}
\end{lemma}
\begin{proof}
Throughout the proof, we assume that $\yY + 1 \le \xX $. 
We can prove the other cases $\xX \le \yY -1$ and $\yY -1 < \xX < \yY + 1 $ in the same fashion 
and omit the proof for them. Thus, instead of \eqref{:17e}, we shall prove 
\begin{align}\label{:17f}&
\sup_{ \xX \in [ \yY + 1, \infty ) }
 \lvert \HHp (x,y)\rvert \le \cref{;I_1}( 1 \vee \log \lvert \yY \rvert ), \quad \yY \in (-\infty, -1] 
.\end{align}

We begin by proving \eqref{:17f} for $ p = 1 $. From $ \yY + 1 \le \xX $, we have 
\begin{align}\label{:17h}
&\HHa (x,y)= \HHa (y)= 
\Big(\int_\yY ^{\yY + 1 } \frac{\int_{\yY }^{u} \Ai'(v)dv}{u-\yY }du \Big) \Ai'(\yY ) 
.\end{align}
Using \lref{l:81}, we obtain for $\uU \in \yyone $, 
\begin{align}\notag 
\int_{\yY }^{\uU } \Ai'( v ) dv 
&= \int_{\yY }^{\uU } \frac{\vVV {1/4}}{\pi ^{1/2} } \sin \Big(\frac{2}{3}\vVV {3/2} - \frac{\pi}{4}\Big)
\Big(1 + \mathcal{O}(\vVV {-3/2}) \Big) dv
\\ \nonumber
&= \int_{\yY }^{\uU } \frac{\vVV {1/4}}{\pi ^{1/2} } \sin \Big(\frac{2}{3}\vVV {3/2} - \frac{\pi}{4}\Big)dv 
+(\uU - \yY )\mathcal{O}(\lvert \yY \rvert ^{-5/4}), \quad \yY \to-\infty 
\\ \label{:17k}
 &
= - \frac{1}{\pi ^{1/2} }
\int_{\frac{2}{3}\lvert \yY \rvert ^{3/2}}^{\frac{2}{3} 
\lvert \uU \rvert ^{3/2}} f (w) \sin \Big( w-\frac{\pi}{4} \Big)dw 
+(\uU - \yY )\mathcal{O}(\lvert \yY \rvert ^{-5/4}) 
,\end{align}
%%%%%%%%%%%%%%%%%%%%%%%%%%%%%%%%%%%%%%
where $f (w) = (({3}/{2}) w )^{-1/6}$. 
Note that $ f $ is nonnegative and decreasing on $(0,\infty)$. 
Hence, using \lref{l:161}, we deduce that 
\begin{align}\label{:17m}
\sup_{\uU \in \yyone }
\Big\lvert 
\int_{\frac{2}{3}\lvert \yY \rvert ^{3/2}}^{\frac{2}{3}\lvert \uU \rvert ^{3/2}} f (w) 
\sin \Big( w-\frac{\pi}{4} \Big)dw 
\Big\rvert 
= \mathcal{O}(\lvert \yY \rvert ^{-1/4}), 
\quad \yY \to-\infty
.\end{align}
From \eqref{:17k} and \eqref{:17m}, we obtain 
\begin{align}\label{:17n}
\sup_{ \uU \in \yyone } \Big\lvert 
 \int_{\yY }^{\uU } \Ai'(v) dv \Big\rvert 
= \mathcal{O}(\lvert \yY \rvert ^{-1/4}), 
\quad \yY \to-\infty 
.\end{align}
From \eqref{:81z} and \lref{l:81}, we obtain $ \Ai '' (\yY ) = \mathcal{O}(\lvert \yY \rvert ^{3/4}) $ as $ \yY \to -\infty $. 
Using this and applying \lref{l:162} to $ \gdelta = (\Ai' , \lvert \yY \rvert ^{-1}) $, 
we deduce from \eqref{:17n} that 
\begin{align}\label{:17N}
\int_{\yY }^{\yY + 1} \frac{\int_{\yY }^{\uU } \Ai'(v)dv}{\uU -\yY }d\uU = 
\mathcal{O} (\lvert \yY \rvert ^{-1/4} (1 \vee \log \lvert \yY \rvert ) ) 
,\quad \yY \to -\infty 
.\end{align}
From \lref{l:81}, we see that $ \Ai'(\yY ) = \mathcal{O} (\lvert \yY \rvert ^{1/4}) $ as $ \yY \to -\infty$. 
Combining this with \eqref{:17h} and \eqref{:17N}, we obtain \eqref{:17f} for $ p = 1 $.

We next prove \eqref{:17f} for $ p = 2$. Using \eqref{:81z} and \lref{l:81} and applying \lref{l:161} with $f (w) =(({3}/{2})w)^{1/6}$ using the same argument used in obtaining \eqref{:17n}, 
we see that 
\begin{align}\label{:17o}
\sup_{ \uU \in \yyone } \int_{\yY }^{\uU } \Ai''(v) dv = \mathcal{O}( \lvert \yY \rvert ^{1/4}), 
\quad \yY \to-\infty 
.\end{align}
Take $ \gdelta = (\Ai'' , \lvert \yY \rvert ^{-1}) $ in \lref{l:162}. 
Then, from \lref{l:81}, \eqref{:81z}, and \eqref{:17o}, 
\begin{align}\label{:17O}
\int_{\yY }^{\yY + 1} \frac{\int_{\yY }^{\uU } \Ai''(v)dv}{\uU -\yY }d\uU = 
\mathcal{O} (\yY ^{1/4} (1 \vee \log \yY ) ), \quad \yY \to-\infty 
.\end{align}
Using \eqref{:17O} and \lref{l:81}, 
we obtain \eqref{:17f} for $ p = 2$. 

We finally prove \eqref{:17f} for $ p = 3 $. From \lref{l:81}, we deduce that, as $\yY \to -\infty$, 
\begin{align} \notag 
\Big( \sup_{ \yY + 1 < \xX < \infty } 
\Big\lvert \int_{\yY + 1 }^\xX \frac{\Ai (u) }{u-\yY } du\Big\rvert \Big) 
\Ai' (\yY ) 
 = & 
\mathcal{O}( \int_{\yY + 1 }^{\infty}\frac{\lvert u \rvert ^{-1/4}}{ \lvert u - \yY \rvert } d u ) 
\mathcal{O}( \lvert \yY \rvert ^{1/4}) 
\\ \label{:17r} = & 
\mathcal{O}(1 \vee \log \lvert \yY \rvert ) 
.\end{align}
Through integration by parts, we have 
\begin{align}\notag
\int_{\yY + 1 }^\xX \frac{\Ai' (u) }{u-\yY } du 
&=
\Big[ \frac{\Ai(u)}{u-\yY }\Big]_{\yY + 1 }^\xX + \int_{\yY + 1 }^\xX \frac{\Ai (u)}{(u-\yY )^2}du 
.\end{align}
Hence, applying \lref{l:81} to $ \Ai $, we see that 
\begin{align}\label{:17s}
 \sup_{ \yY + 1 < x < \infty }\Big\lvert 
\int_{\yY + 1 }^\xX \frac{\Ai' (u) }{u-\yY } du 
\Big\rvert 
\Ai(\yY ) 
= \mathcal{O}(1), \quad \yY \to -\infty 
.\end{align}
Combining \eqref{:17r} and \eqref{:17s}, we have \eqref{:17f} for $ p = 3 $.
This completes the proof.
\end{proof}

\noindent 
{\em Proof of \lref{l:85}}: 
Let $ I , J = (-\infty, -1] $ or $ (-1, \infty) $. 
We divide \eqref{:85a} into four cases:
\begin{align}&\label{:80a}
\sup_{\xX \in I }\Big\lvert \int_\yY ^\xX \Kairy(u,\yY )du\Big\rvert \le 
 \cref{;T4Na} ( 1 \vee \log \lvert \yY \rvert ), \quad \yY \in J 
.\end{align}
We only present the proof for the case that $ I = J = (-\infty, -1] $. 
Indeed, other cases can be proved similarly and more easily 
because $ \Ai $ and $ \Ai '$ decay exponentially on $ (-1, \infty) $. 
Note that 
\begin{align} & \label{:172z}
\Ai (u) \Ai'(\yY ) - \Ai'(u)\Ai (\yY )= 
\Big( \int_\yY ^u \Ai'(v)dv \Big) \Ai'(\yY ) 
- \Big( \int_\yY ^u \Ai'' (v)dv \Big) \Ai(\yY )
.\end{align}
Let $ \HHp $, $ p = 1,2,3$ be as in \lref{l:17}. Then, using \eqref{:11g} and \eqref{:172z}, we have 
\begin{align}&\notag %\label{:172a}&
\int_\yY ^\xX K_2(u,\yY )du = \HHa + \HHb + \HHc 
.\end{align}
We therefore deduce \eqref{:80a} for $ I = J = (-\infty, -1] $ from \lref{l:17}. 
\qed

%%%%%%%%%%%%%%%%%%%%%
\section{Appendix 2: Estimates of the normalized oscillator functions $\psi _{N } (x)$}\label{s:Ap2}
In \sref{s:Ap2}, we collect estimates of the normalized oscillator functions $\psi _{N } (x)$ given in \eqref{:31a} and related quantities. 
We use the results of Plancherel and Rotach \cite{PR29} and Deift et al. \cite{DKMVZ} to estimate $ \psi _{N }$. The results in these sections are used in \sref{s:Ap3}. 

\subsection{Plancherel--Rotach asymptotics of $\psi _{N } (x)$}\label{s:9}
%%%%%%%%%%%%%%%%%

In Subsection \ref{s:9}, we collect Plancherel--Rotach asymptotics of $\psi _{N } (x)$. 

We introduce the subsidiary functions $ \mathbb{F} $ and $\mathbb{U} $ on $ \mathbb{N}\times \mathbb{R} $ defined as 
\begin{align}\notag &%\label{:91a}& 
\mathbb{F} (N,w) = ( {N+1} )^{-1/3} ( N+1 - ({w}/{2} )^2 )
,\\ \label{:91b}&
 \mathbb{U} (N,w) = N^{1/6} (w-2\sqrt{N})=2N^{1/6} ( ({w}/{2}) - \sqrt{N} ) 
.\end{align}

In addition to the variable $x$ and parameter $N$, we introduce the auxiliary variable $\theta$ to evaluate the function $\psi_N$ in Plancherel-Rotach \cite{PR29} . Then, $ x $ is expressed as a function of $(N, \theta)$. 
The other functions of $ ( N , \theta ) $ such that 
$ \mathsf{f}, \mathsf{g}, \h , \mathsf{f} ^{\rm h} , \mathsf{g} ^{\rm h}$, and $ \mathsf{u} ^{\rm h}$ below are introduced to express the asymptotic behavior of $\psi_N (x) $. 
Let $ \mathsf{f} $, $ \mathsf{g} $, and $ \h $ be the functions on $\mathbb{N}\times [0, \pi/2]$ defined by 
\begin{align}\label{:91f} & 
 \mathsf{f} (N, \theta ) = \mathbb{F} (N, 2\sqrt{N+1} \cos \theta) = (N+1)^{2/3}\sin^2 \theta ,
\\ \label{:91g}
& \mathsf{g} (N, \theta) =2^{-1}{(N +1)} (2\theta-\sin 2\theta)
,\\ \label{:91h}
& \h (N, \theta ) = \mathbb{U} (N, 2\sqrt{N+1} \cos \theta) = 
2N^{1/6}(\sqrt{N+1}\cos \theta -\sqrt{N} )
.\end{align}
\begin{lemma}\label{l:91} %\notag 
Let $ x= \h (N,\theta) $. Then the following hold for all $ N \in \mathbb{N}$%$\in [-2N^{2/3}, 0] $, 
\begin{align}& \notag 
- \frac{1}{2 (\oneN )^{1/3} }x + 
 \frac{1}{N^{1/3} (\oneN )^{1/3}}
 \le\mathsf{f} (N, \theta)
\\& \label{:91p}
\quad \quad \quad \quad \quad \quad \quad 
 \le - x + \frac{1}{N^{1/3} (\oneN )^{1/3}} 
\quad \text{ for } x \in [-2N^{2/3}, 0] 
, \\&\label{:91P}
- \frac{1}{4}x 
 \le\mathsf{f} (N, \theta) \le - 2 x
\quad \text{ for } x \in [-2N^{2/3}, -1] 
.\end{align}
\end{lemma}
\begin{proof}
From $ x= \h (N,\theta) $ and \eqref{:91h}, 
$ (x + 2N^{2/3})^2=4 N^{1/3}(N+1) \cos^2 \theta $. Then, 
\begin{align} \label{:U1c}& 
 \mathsf{f} (N,\theta)=- \frac{1}{(1+N^{-1} )^{1/3}}x - 
\frac{x^2}{4 N^{2/3}(1+N^{-1})^{1/3}} + \frac{1}{N^{1/3}(1+N^{-1})^{1/3}}
.\end{align}
The estimate \eqref{:91p} follows from \eqref{:U1c}. 
For $ N \in \mathbb{N}$ and $ x \in [-2N^{2/3}, -1] $, it holds that 
\begin{align}& \label{:91Q}
\frac{1}{2} \le 
\frac{1}{ (\oneN )^{1/3}} , \quad \frac{1}{N^{1/3} (\oneN )^{1/3}} \le -x 
.\end{align}
Hence, \eqref{:91P} follows from \eqref{:91p} and \eqref{:91Q} immediately. 
\end{proof}

\begin{lemma} \label{l:W4} 
We regard $ \uu $ as a function of $ \theta $ and omit $ N $ from the notation. Then, 
\begin{align}& \label{:W4a}%\notag & 
 \supN \sup_{x \in [-2N ^{2/3}, -1)} 
 \Big\lvert 
\int_x^{-1}
\frac{1}{\mathsf{f} (N,\theta )^{1/4}} \sin \Big(\mathsf{g} (N,\theta ) +
 \frac{\pi}{4} - \frac{\theta }{2}\Big) 
 \Big\rvert _{\theta = \uu ^{-1} (\uU ) }
 d \uU 
 \Big\rvert 
 <\infty 
.\end{align}
\end{lemma}
%%%%%%%%%%%%%%%%%%%%%%%%%%%%%%%%%%%%%%
%We shall prove \lref{l:W4} in \sref{s:18}. 
%\section{Appendix 3: Proof of \lref{l:W4}} \label{s:18}

\begin{proof}
Let $ \mathsf{f} $, $ \mathsf{g} $ and $ \h $ be as in \eqref{:91f}--\eqref{:91h}. Let 
\begin{align} \label{:18e}&
 \vv =\mathsf{g} + \frac{\pi}{4} -\frac{\theta}{2} =\frac{N +1}{2}(2\theta-\sin 2\theta) + \frac{\pi}{4} -\frac{\theta}{2}
.\end{align}
For each $ N $, we regard $ \uu $ and $ \vv $ as functions of $ \theta $. 
We set $ \uu (\theta ) = \uu ( N , \theta ) $ and $ \vv (\theta ) = \vv ( N , \theta ) $. 
We thus omit $ N $ from the notation of $ \uu (\theta ) $ and $ \vv (\theta )$. 
Simple calculation shows that 
\begin{align} \label{:18f} 
&\ODod{ \uu }{ \theta} ( \theta )
= -2N ^{1/6} (N +1)^{1/2}\sin\theta 
,\\ \label{:18g} 
&\ODod{ \vv }{ \theta} ( \theta )
=(N +1)(1-\cos 2\theta)-\frac{1}{2}= 
2(N +1)\sin^2 \theta-\frac{1}{2} 
.\end{align}
We set 
\begin{align} \label{:18d}&
\dN = \{ \theta \in (0, \pi/2]; \,(N+1)^{-1/2 } < 2 \sin \theta \}
.\end{align}
Then, from \eqref{:18f}--\eqref{:18d}, we see that for each $ N \in \mathbb{N}$, 
\begin{align} \label{:18h} &
\ODod{ \uu }{ \theta} ( \theta ) <0 , \quad 
\ODod{ \vv }{ \theta} ( \theta ) >0 
\quad \text {for $ \theta \in \dN $}
.\end{align}

The functions $ \uu $ and $ \vv $ are injective on $\dN $ by \eqref{:18h}. 
Hence, let $\hat{ \uu } = \uu \circ \vv ^{-1}$ and 
$\hat{ \vv } = \vv \circ \uu ^{-1}$ be the functions defined on 
$ \vv (\dN ) $ and $ \uu (\dN ) $, respectively. 
By definition, 
$ \hat{\uu } (\vv (\theta )) = \uu (\theta )$ and $ \hat{\vv } (\uu (\theta )) = \vv (\theta )$ 
for $ \theta \in \dN $. 
From \eqref{:18f} and \eqref{:18g}, we have 
\begin{align}\label{:18m}
 \frac{ d \hat{ \vv }}{ d \uU } (\uU ) = - 
\frac
{\sqrt{N +1} \sin^2 \theta-(4\sqrt{N +1})^{-1}} 
{N ^{1/6}\sin\theta} 
\Big\rvert _{\theta = \uu ^{-1} (\uU ) }
.\end{align}
Let $ \map{\FN }{\dN }{\mathbb{R} }$ such that 
\begin{align} &\label{:18n}
\FN ( \theta )= 
\Big(\frac{N}{N+1}\Big)^{1/6}
\frac{\sqrt{\sin\theta}}
{\sqrt{N +1} \sin^2 \theta-(4\sqrt{N +1})^{-1}}
.\end{align}
From \eqref{:91f}, \eqref{:18m}, and \eqref{:18n}, we obtain 
\begin{align}\label{:18l}&
\frac{1}{\mathsf{f} ( N , \uu ^{-1}(\uU ) )^{1/4}} = 
- \FNuU \frac{ d \hat{ \vv }}{ d \uU }(\uU ) 
.\end{align}
Using \eqref{:18e} and \eqref{:18l}, we deduce that 
\begin{align} \notag 
& 
\int_x^{-1} 
\frac{1}{\mathsf{f} ( N , \theta )^{1/4}} 
\sin \Big(\mathsf{g} \Ntheta + \frac{\pi}{4} - \frac{\theta }{2}\Big) 
 \Big\rvert _{\theta = \uu ^{-1} (\uU ) }
 d \uU 
\\ \notag &=
\int_x^{-1} 
 - \FNuU 
 \frac{ d \hat{ \vv }}{ d \uU }(\uU ) 
\sin ( \hat{ \vv }(\uU )) d \uU 
 \\ \label{:18k} 
 &
= \int_{\hat{ \vv }(-1)}^{\hat{ \vv }(x)} \FNvV 
 (\sin \vV) d \vV 
.\end{align}

From \eqref{:18d}, \eqref{:18h} and \eqref{:18n}, we deduce that 
for $ \vV \in \vv (\dN ) $
\begin{align} \label{:18O}& 
 \FN ( \vv ^{-1}(\vV ) )> 0 
,\\ \label{:18o} &
\frac{ d \FNvV }{ d \vV } = 
\Big(
\frac{ d \FN }{ d \theta}
	\Big) 
(\vv ^{-1}(\vV ) ) 
	\Big (\frac{ d \vv ^{-1} }{ d \vV }\Big) (\vV ) < 0 
.\end{align}
% From \eqref{:18j}--\eqref{:18n}, we have

From \eqref{:91h}, we see that, for each $ N \in \mathbb{N}$, $ \uu = \uu (\theta ) $ is decreasing, 
$\uu (\pi / 2) = -2N^{2/3}$, and 
\begin{align}&\notag % \label{:18j} &
 \sin^2 \theta = (N+1)^{-1}(N^{1/3} + 1 -\QM N^{-1/3}) \quad \text{if} \quad 
\uu (\theta) = -1
.\end{align}
Hence, we deduce that 
\begin{align} \label{:18p}
 [-2N ^{2/3}, -1] \subset \uu ( \dN ) 
, \quad 
\lim_{N\to\infty} \FNvV \Big\rvert _{ \vV= \hat{ \vv }(-1)} = 1
.\end{align}
Using \eqref{:18O}--\eqref{:18p}, we apply \lref{l:161} to $ \FNvV $, which yields 
\begin{align}\label{:18s} 
\sup_{N\in\mathbb{N} }\sup_{x\in[-2N ^{2/3}, -1]} \Big\lvert 
\int_{\hat{ \vv }(-1)}^{\hat{ \vv }(x)} \FNvV \sin \vV \ d \vV 
\Big\rvert 
&< \infty
.\end{align}
Hence, we obtain \eqref{:W4a} from \eqref{:18k} and \eqref{:18s}. 
This completes the proof. 
\end{proof}
%\qed

Let $ \mathsf{f} ^{\rm h} $, $ \mathsf{g} ^{\rm h} $, and $ \h ^{\rm h} $ be the functions on $ \mathbb{N}\times [0, \infty) $ defined by 
 \begin{align} \label{:92x} &
\mathsf{f} ^{\rm h}(N, \theta) = -\mathbb{F} (N, 2\sqrt{N+1} \cosh \theta)
= (N+1)^{2/3} \sinh^2 \theta,
\\ \label{:92y}
&\mathsf{g} ^{\rm h} (N, \theta) = 2^{-1}{(N +1)} (2\theta-\sinh 2\theta)
,\\ \label{:92z}
& \h ^{\rm h}(N, \theta) = \mathbb{U} (N, 2\sqrt{N+1} \cosh \theta)
=2N^{1/6}(\sqrt{N+1}\cosh \theta -\sqrt{N} )
.\end{align}

\begin{lemma} \label{l:92} For $ x= \h ^{\rm h}(N, \theta)$, 
\begin{align}\label{:92a}& 
 \mathsf{f} ^{\rm h}(N,\theta) =
\Big( \frac{N}{N+1} \Big)^{1/3}x + \frac{x^2}{4 N^{1/3}(N+1)^{1/3}} - 
\Big( \frac{1}{N+1} \Big)^{1/3} 
.\end{align}
\end{lemma}
\begin{proof}
We deduce \eqref{:92a} from \eqref{:92x} and \eqref{:92z}. 
\end{proof}

Let $a_{\KM }$ be the constants in \cite[p. 228]{PR29} such that 
\begin{align*}&
\mathsf{P}_k (\xi)= \sum_{m=0}^k a_{km}\xi^m
,\\&
\exp[\xi \sum_{m=3}^{\infty}\frac{(-1)^m}{m} t^{m-2} ]
= 
\sum_{k=0}^{\infty}\mathsf{P}_k (\xi) t^k 
.\end{align*}
Note that $a_{\KM }$ are independent of $N $ and $\theta$. For instance, 
$$
a_{00}=1, \quad a_{10}=0, \quad a_{11}= - 1/3 , \quad a_{20}=0, \quad a_{21}= 1/4 
\ \text{ and } \ a_{22}= 1/18 
.$$
For $k,m \in \zN $ with $ m \le k$, we put 
\begin{align} \label{:93u}
& 
\ckm (\theta)=\frac{\pi}{4}-\frac{\theta}{2}
-\Big(2m+k\Big)\Big(\frac{\pi}{4}+\frac{\theta}{2}\Big)
, \; \theta \in (0, \frac{\pi}{2}]
,\\ \label{:93v} &
\Ckm (N ,\theta)=\frac{1+(-1)^k}{2}
\frac{\Gamma (m+\frac{k+1}{2})}{(N +1)^{k/2}(\sin\theta)^{m+k/2}}a_{\KM },
\; \theta \in (0, \frac{\pi}{2}]
,\\ \notag % \label{:93w}
&\CC ^{\rm h}_{k m}(N ,\theta)= 
\frac{1+(-1)^k}{2}
\frac{
\Gamma (m + \frac{k+1}{2})}{(N +1)^{k/2}}
\Big( \frac{2}{e^{-2\theta} - 1 } \Big)^{m+k/2}a_{k m},
\; \theta \in (0, \infty)
.\end{align}

\begin{lemma} \label{l:9X}
 Let $ k,m \in \zN $, $m\le k$. 
We find $ \Ct \label{;9X}$ independent of $ N $ such that 
\begin{align}&\notag 
\CC ^{00} (N ,\theta) = \pi ^{1/2} ,\quad 
\Ckm (N ,\theta) = 0 \ (k\text{ is odd}), 
\\\label{:U1d}&
 \big\lvert \Ckm (N ,\theta) \big\rvert \le \cref{;9X}\vertx ^{-3k/4} \ (k\text{ is even}) 
\quad \text{for $ x= \h (N,\theta) \in [-2N^{2/3}, -1 ] $}
.\end{align}
\end{lemma}

\begin{proof}
\eqref{:U1d} is clear if $ k = 0 $ or $ k $ is odd. 
Let $k$ be even and $k \ge 2$. From \eqref{:93v}, 
\begin{align}\label{:U1e}
\Ckm (N ,\theta) \le 
\frac{\Gamma (m+\frac{k+1}{2}) a_{\KM }}{(N +1)^{k/2}(\sin\theta)^{3k/2}}
=\frac{\Gamma (m+\frac{k+1}{2}) a_{\KM }}{ \mathsf{f} (N,\theta)^{3k/4}}
.\end{align}
Here, we used $m\le k$ for the inequality and \eqref{:91f} for the equality. 
From \eqref{:91P} and \eqref{:U1e}, we obtain \eqref{:U1d} for even $k \ge 2$. 
\end{proof}

Let $ \DDN $ be such that 
\begin{align} &\notag % \label{:93a} &
\DDN =\Big( \frac{N}{N+1} \Big)^{1/12}\displaystyle{\frac{\sqrt{(N+1)!}}{( 2\pi (N+1) )^{1/4}}
\Big( \frac{e}{N+1} \Big)^{(N+1)/2}}
.\end{align}
Then, by a simple calculation using the Stirling formula, we derive 
\begin{align}\label{:93b}
&\DDN =1+ {\mathcal{O}} ( N ^{-1}) , \text{\quad $N\to\infty$}
.\end{align}

For functions $ f $ and $ g $ defined on 
$ U \subset \mathbb{R} ^n $, $n \in \mathbb{N}$, we write 
\begin{align}&\label{:93c}
 f ( \bmu )=\OL ( g ( \bmu)) \quad \text{ on \ $U$} 
\end{align}
if there exist $C \in (0,\infty)$ and $c \in (0,\infty)$ such that 
$ \lvert f (\bmu)\rvert \le C \lvert g (\bmu)\rvert $ 
for any $ \bmu \in U_{ g }$, where $ U_{ g }$ is such that 
$ U_{ g }=\{ \bmu \in U; \lvert g (\bmu)\rvert \le c \}$.

\begin{remark}\label{r:81}
The asymptotic behavior of \cite{PR29} uses auxiliary variable $\theta$ in addition to the variable $x$ and the parameter $N$. 
%Since those variables are related to each other, it is necessary how to evaluate using constants that do not depend on the variables $(x,N,\theta)$. To this end, we introduce the symbol $\OL (\cdot)$, which is a modified version of Landau's symbol. 
%
Since those variables are related, it is necessary to evaluate using constants that do not depend on the variables $(x,N,\theta)$. To this end, we introduce the symbol $\OL (\cdot)$, a modified version of Landau's symbol. 

For example, if $ g (x,N,\theta) $ is constant in $ (x,\theta) $ such that $ g (x,N ,\theta) = N^{-1} $, then 
\begin{align}& \label{:93D}
f(x,N,\theta) = \OL (N^{-1}) 
\end{align}
implies that there exist $C, c >0$ independent of $(x,N,\theta)$ such that
\begin{align}& \label{:93E}
f(x, N, \theta) \le CN^{-1} \text{ for $N> c^{-1}$}
.\end{align}
This implies that, for each $ (x , \theta )$, $ f(x, N, \theta) = O (N^{-1})$ in the usual Landau notation. 
We note that \eqref{:93D} is a stronger statement than $ f(x, N, \theta) = O (N^{-1})$ because $ C $ in \eqref{:93E} is independent of $ (x,\theta )$. 
\end{remark}

We set 
\begin{align}& \label{:93f}
\Om =\{ (x, N, \theta) \in [-2N^{2/3}, 0]\times \mathbb{N}\times (0, \frac{\pi}{2}];
x=\h (N,\theta) \}
,\\\label{:93g} & 
\Omh =\{ (x, N, \theta) \in [1, \infty) \times \mathbb{N}\times (0, \infty); x=\h ^{\rm h}(N,\theta)\}
.\end{align}
From \eqref{:91f} and \eqref{:91p}, we deduce that 
\begin{align} \label{:93G}& 
\lvert N \sin \theta\rvert ^{-1}=\6 \text{ on \ $\Om$}
.\end{align}
 Let $w= 2\sqrt{N} + N^{-1/6}x$. Then $\mathbb{U} (N,w)=x$. Hence, from \eqref{:31y} and \eqref{:31a}, 
\begin{align}\label{:93h}
\psi_N (x)=N^{1/12} \phi_N (w)=
{N^{1/12}}{(\sqrt{2\pi}N!})^{-1/2} e^{-w^2/4} 
 \widehat{H}_N (w) 
.\end{align}
 The asymptotic behaviors of Hermite polynomials $\widehat{H}_N $ given in \cite{PR29} 
 provide the following estimates from \eqref{:93h} on the domains $ \Om $ and $ \Omh $. 

\begin{lemma}[{\cite[Sect. 2, I, II, pp. 229--230]{PR29}}] \label{l:95}
For each $L\in\mathbb{N} $, we have the following. 
\begin{align} \notag 
 \psi_{N } (x) & ={ \pi^{-1}}{\DDN }
{\mathsf{f} (N, \theta)}^{-1/4}
\sumL 
\Ckm (N ,\theta)\sin \big( \mathsf{g} (N, \theta) + \ckm (\theta)\big)
\\ \label{:93p} 
& \ \qquad\qquad\qquad \qquad\qquad
+ \OL ( {\mathsf{f} (N, \theta )}^{-(3L+1)/4} ) \text{\quad on \ $\Om$} 
,\\ \notag 
\psi_{N }( x ) 
&= 
{ \pi^{-1}}{\DDN }
{( 8e^{2\theta}\mathsf{f} ^{\rm h}(N, \theta) )}^{-1/4} 
\exp ( \mathsf{g} ^{\rm h} (N, \theta) )
\\ \label{:93q} 
& 
\times \Big( \sumL \CC ^{\rm h}_{k m}(N ,\theta) + 
 \OL (( e^{2\theta}\mathsf{f} ^{\rm h}(N,\theta) )^{-3L/4} )
\Big) \text{\quad on \ $\Omh $}
.\end{align}
Here $ \mathsf{f} ^{\rm h}$ and $ \mathsf{g} ^{\rm h}$ are given by \eqref{:92x} and \eqref{:92y}, 
and $ \CC ^{\rm h}_{k m}$ is defined before \lref{l:9X}. 
\end{lemma}

%%%%%%%%%%%%%%%%%%%%%
	\subsection{Estimates of functions in Plancherel--Rotach asymptotics} 
%\subsection{Estimates of functions in \lref{l:95}}
\label{s:U}

In Subsection \ref{s:U}, we present estimates of functions in Plancherel--Rotach asymptotics in \lref{l:95}. 

% Let $ \OL (\cdot ) $ be as in \eqref{:93c}. 
Let $\Om $ be as in \eqref{:93f}. 
For $\varepsilon \in (0, 2/3)$, let
\begin{align} \label{:U1y}
\Ome =\{ (x,N,\theta) \in \Om; x\in [-2N^{2/3}, -N^{\varepsilon}]\}
.\end{align}
For $ \lpm , 0 $ and $(x,N, \theta) \in \Ome$, let $\eta_\ell=\eta_\ell (N, \theta)$ be a number such that 
%$ \theta + \eta_\ell \in (0, \pi / 2] $ and
\begin{align} \label{:U1z}
& \sqrt{N+1}\cos \theta=\sqrt{N+\ell +1} \cos (\theta + \eta_\ell) ,\quad \theta + \eta_\ell \in (0, \pi / 2]
.\end{align}
%%%
\begin{lemma} \label{l:U1} 
Let $ \OL (\cdot ) $ as in \eqref{:93c}. 
It holds that $ \eta _0 = 0 $ and that for $\lpm $ 
\begin{align}\label{:U1a} & 
 \eta_\ell=\frac{\ell \cos \theta}{2N \sin \theta} + \OL \Big(\frac{1}{N^2 \sin^3 \theta}\Big)
\text{\quad on $\Ome$}
.\end{align}
\end{lemma}
%%%%%%%%%%%%%%%%%%%%%%%%%%%%%%%%%%%%%%%%%%%%%%%%%%
\begin{proof}
The claim $ \eta _0 = 0 $ is clear. 
From \eqref{:U1z} and Taylor's theorem, we deduce that 
\begin{align*}
\cos (\theta +\eta_\ell ) - \cos \theta &= 
\frac{\sqrt{N+1}-\sqrt{N+\ell+1}}{\sqrt{N+\ell+1}}\cos \theta 
= 
- \frac{\ell \cos \theta}{2( N + 1 )} + \OL \Big(\frac{1 }{N^2 } \Big)
\text{\quad on $\Ome$}
.\end{align*}
Hence, we see that $\eta_\ell \to 0$, as $N\to \infty$.
Using Taylor's theorem, we readily see that 
\begin{align*}&
\cos (\theta +\eta_\ell ) - \cos \theta =- (\sin \theta ) \eta _{\ell } + (\cos \theta ) \OL (\eta_{\ell }^2)
\text{\quad on $\Ome$}
.\end{align*}
Combining these, we deduce \eqref{:U1a}. 
\end{proof}

%%%%%%%%%%%%%%%%%%%%%%%%%%%%%%%%%%%%%

Let $(x, N,\theta)$ and $\eta_\ell$ be as in \lref{l:U1}. We set for $ \lpm , 0 $, 
\begin{align}\notag &
 \fell (N,\theta)=\mathsf{f} (N+\ell, \theta +\eta_\ell), \quad 
 \gell (N,\theta)= \mathsf{g} (N+\ell, \theta+\eta_\ell) 
,\\ \label{:U2p} 
&
 \ct _\ell^{\KM } (N,\theta)=\ckm (\theta+\eta_\ell) 
, \quad 
 \Clkm (N,\theta)=\Ckm (N+\ell ,\theta+\eta_{\ell}) 
.\end{align}
For $ \ell = 0$, 
we obviously have $ \mathsf{f} _0 = \mathsf{f} $, $ \mathsf{g} _0 = \mathsf{g} $, $ \ct _0^{\KM } = \ct ^{\KM }$, and $ \Ckm _0 = \Ckm $. 
Let 
\begin{align} \notag 
&\delta_1(N,\theta)=\frac{2+\cos^2 \theta}{3N^{1/3}} 
, \qquad
\delta_2^{\KM } (N,\theta)= - \frac{(1+2m+k)\cos \theta}{4N\sin \theta} 
 ,\\ \label{:U2r}
&\delta_3^{\KM } (N,\theta)= - \frac{\Ckm (N,\theta)}{2N}
\Big( 
k + \Big(m+\frac{k}{2}\Big)\frac{\cos^2 \theta}{\sin^2 \theta}
\Big) \mathbf{1}(k\ge 2) 
.\end{align}

\begin{lemma} \label{l:67} \thetag{1}
For $\lpm $, 
\begin{align} \label{:U2a} & 
 \fell (N,\theta) - \fz (N,\theta) 
= \delta_1(N,\theta)\ell 
+ \5 \text{\quad on $\Ome$}, 
\\ \label{:U2b} & 
 \gell (N,\theta) - \gz (N,\theta) 
=\theta\ell + \6 
 \text{\quad on $\Ome$}
.\end{align}
\thetag{2} 
For $ k,m \in \zN $ such that $m\le k$, 
\begin{align} \label{:U2d} &
\ct _\ell^{\KM } (N,\theta) 
-
\ckm (\theta) 
= \delta_2^{\KM } (N,\theta) \ell + \OL ({N^{-1}\vertx ^{-3/2}} ) \text{\quad on $\Ome$},
\\ \label{:U2e} &
 \Clkm (N,\theta) - \Ckm (N,\theta) 	= 
\displaystyle{\delta_3^{\KM } (N,\theta)}\ell
+ \OL ( N^{-2/3} \vertx ^{-(3k/4 +2)} ) \text{\quad on $\Ome$}
.\end{align}
For $k=0$ or $k$ is odd, $ \Clkm (N,\theta) -\Ckm (N,\theta)=0$. 
\end{lemma}

%%%%%%%%%%%%%
\begin{proof} 
We begin by proving \eqref{:U2b}. 
From \eqref{:91g} and \eqref{:U2p}, we have 
\begin{align} \label{:U2i} & 
 \gell - \gz =(I_1 + I_2) /2
,\end{align}
where $I_1= (N+1)( 2\eta_\ell -\sin 2(\theta+\eta_\ell) + \sin 2\theta ) $ 
and $I_2= \ell ( 2(\theta + \eta_\ell) -\sin 2(\theta+\eta_\ell) )$.
Put $\widehat{\eta}={1}/{N\sin\theta}$. 
From Taylor's theorem with \eqref{:U1a}, we obtain on $\Ome$ that 
\begin{align} \notag 
 I_1 
 &=(N+1) ((2 -2\cos 2\theta) \eta_{\ell} +2(\sin 2\theta )\eta_\ell^2 + 
\OL (\widehat{\eta}^3) ) 
\\ \nonumber
&= 4(N+1) ( \sin^2 \theta ) \Big( \frac{\ell \cos \theta}{2N \sin \theta} + 
\OL \Big(\frac{1}{N^2 \sin^3 \theta}\Big) \Big) 
+ \OL (N (\sin 2\theta ) \widehat{\eta}^2 ) + \OL (N\widehat{\eta}^3)
\\ \nonumber
&=\ell\sin 2\theta + \OL \Big(\frac{1}{N}\Big)
+\OL \Big(\frac{1}{N\sin\theta }\Big)
+\OL \Big(\frac{1}{N^2 \sin^3 \theta}\Big)
.\end{align}
From \eqref{:91f} and \eqref{:91P}, we can replace $ N ^{2/3}\sin ^2 \theta $ by $ \vertx $ 
in the term $ \OL ({1}/{N^2 \sin^3 \theta})$. 
From this and applying \eqref{:93G} to the term 
$ \OL ({1}/{N\sin\theta } )$, we obtain on $\Ome$ that 
\begin{align} & \label{:U2j} 
 I_1=\ell\sin 2\theta + \6 
.\end{align}
From Taylor's theorem with \eqref{:U1a} and using \eqref{:91P}, we obtain on $\Ome$ that 
\begin{align} \label{:U2k}
I_2 &= 2 \ell\theta - \ell \sin 2\theta +
 \OL \Big(\frac{1}{N\sin \theta}\Big)
% \\\notag &
=2 \ell\theta - \ell \sin 2\theta +
 \6 
.\end{align}
From \eqref{:U2i}--\eqref{:U2k}, we have \eqref{:U2b}.

The proofs of other claims are similar to the proof of \eqref{:U2b} as follows, and 
we omit the details. 
We have \eqref{:U2a} by applying Taylor's theorem to 
\eqref{:91f} with \eqref{:91P} and \eqref{:U1a}. 
We obtain \eqref{:U2d} by applying Taylor's theorem to 
\eqref{:93u} with \eqref{:93G} and \eqref{:U1a}. 
We derive \eqref{:U2e} by applying Taylor's theorem to \eqref{:93v} with 
 \eqref{:93G}, \eqref{:U1a}, and \eqref{:U1d}. 
 \end{proof}

%%%%%%%%%%%%%%%%%%%%%%%%%%%%%%%%%%%%%%%%%%%%%%%%%%%%%%%%%%%%%%%%%%%%
\subsection{Deift et al.\,strong asymptotics of $ \psi_N (x) $}\label{s:V}

In this subsection, we recall strong asymptotics in \cite{DKMVZ} (\lref{l:V6}) and give an expansion of $ \psi_N $ with respect to the Airy function (\lref{l:V7}). 

Throughout this paper, we take 
\begin{align}\notag &%\label{:V1o}&
\text{$ 0 < \delta < 1 $ and $ 0 < \varepsilon < 1/6 $}
,\\ \notag &%\label{:V1p}&
 \OmAd =\{ (x,N,z) \in [-2N^{2/3}, 0) \times \mathbb{N}\times [1-\delta,1);
x=\mathbb{U} (N, 2\sqrt{N}z) \}
,\\& \label{:V1q}
\OmAe = \{ (x,N,z) \in \OmAdz ; x\in (-N^{\varepsilon}, -1) , \ N\ge 2\}
.\end{align}
Here, $ \mathbb{U} $ is as in \eqref{:91b}. 
Although $ \OmAe $ depends on $ \delta $, we omit $ \delta $ from the notation. 

Let $ x = \mathbb{U} (N , 2\sqrt{N}z)$. Then, $ (x,N,z)$ satisfies 
\begin{align} \label{:V1w} &
2\sqrt{N}z = 2\sqrt{N} + N^{-1/6}x 
.\end{align}
If $ (x,N,z) \in \OmAe $, then $ (x, N,z) $ satisfies $ x = \mathbb{U} (N , 2\sqrt{N}z)$. 
When $ x = \mathbb{U} (N , 2\sqrt{N}z)$ holds, 
we regard $ x $ as a function in $ z $, and also $ z $ as a function in $ x $ for each $ N $. 
We set for $ \lpm , 0 $, 
\begin{align}\label{:V1x}&
x_\ell = (N+\ell )^{1/6} ( 2\sqrt{N}z -2\sqrt{N+\ell}), 
\quad z_\ell = ({1+N^{-1}\ell})^{-1/2} z 
.\end{align}
We note that, if \eqref{:V1w} and \eqref{:V1x} hold, then we find that $ x_0 = x $, $ z_0 = z $, and 
\begin{align} \label{:V1y}
 z_\ell 
=1 + &\frac{x_\ell}{2(N+\ell)^{2/3}}
= \frac{1}{(\oneNell )^{1/2}} + \frac{x}{2 N^{2/3} (\oneNell )^{1/2}} 
 , \ \lpm , 0 
.\end{align}

\begin{lemma} \label{l:V1}
Let $ \OL (\cdot ) $ be as in \eqref{:93c}. 
Let \eqref{:V1w} and \eqref{:V1x} hold. Then, 
\begin{align} \label{:V1a}&
 x_\ell - x= - {N^{-1/3}}\ell +\8 \text{\quad on $\OmAe $},\ \lpm 
,\\
\label{:V1b} &
\sum_{\lpm } (x_\ell - x )= 6^{-1} N^{-4/3} + \OL (N^{-2} x ) 
\text{\quad on $\OmAe $} 
.\end{align}
\end{lemma}
%%%%%%%%%%%%%
\begin{proof} 
 From \eqref{:V1w} and \eqref{:V1x}, we deduce that 
\begin{align}\notag 
x_\ell &= - (N+\ell)^{1/6}( 2\sqrt{N+\ell} - 2\sqrt{N}z )
\\\notag 
&= - (N+\ell)^{1/6}( 2\sqrt{N+\ell} - 2\sqrt{N} - N^{-1/6}x )
\\ \label{:V1B}
&= - (N+\ell)^{1/6}( 2\sqrt{N+\ell} - 2\sqrt{N}) + (1+N^{-1}\ell)^{1/6} x 
.\end{align}
From \eqref{:V1B}, we deduce that 
\begin{align}&\notag 
x_\ell - x = - 2 N^{2/3}(1+N^{-1}\ell)^{1/6}(\sqrt{1+N^{-1}\ell} - 1) + ( (1+N^{-1}\ell)^{1/6} -1)x 
.\end{align}
Using this, we obtain \eqref{:V1a} and \eqref{:V1b}. 
\end{proof}

Let $ \sigma_{\rm semi}$ be as in \eqref{:11a}. 
Then, $ \sigma_{\rm semi}(2z)=({1}/{\pi}) (1-z)^{1/2}(1+z)^{1/2}$ and 
\begin{align*}
 \int_1^z 2 \sigma_{\rm semi}(2y) dy &= 
 \frac{1}{\pi}\big( z(1-z)^{1/2}(1+z)^{1/2} - \arccos z \big) 
 \quad \text{ for $z\in [-1, 1]$}
.\end{align*}
We set for $z\in [-1, 1]$, 
\begin{align}\notag & %\label{:93y} &
\alpha _N(z) =- \big(-3N \int_1^z (1-y)^{1/2}(1+y)^{1/2} dy \big)^{2/3}
.\end{align}
Taking $y=1- (2N^{2/3})^{-1}u$, we obtain from \eqref{:V1w} 
\begin{align} 
 \label{:V2a}
& (-\alpha _N(z))^{3/2}=\frac{3}{2}\int_0^{-x} u^{1/2} (1-\frac{u}{4N^{2/3}} )^{1/2}du 
.\end{align}

\begin{lemma} \label{l:V2} 
% Let $ z \in [-1,1]$ and $ x=2 N^{2/3}(z-1) $. Then, 
It holds that $ \alpha _N $ satisfies 
\begin{align}
\label{:V2b}
& \alpha _N(z)= 
 x +{20}^{-1} N^{-2/3}{ x ^2} + 
 \OL (N^{-4/3} x ^3 ) \text{\quad on $\OmAe $}
.\end{align} 
\end{lemma}
\begin{proof}
From \eqref{:V2a}, we see that 
\begin{align} \notag %\label{:V2f}&
(-\alpha _N(z))^{3/2} = (-x)^{3/2} - \frac{3 (-x)^{5/2}}{40 N^{2/3}} + \OL \Big(\frac{\vertx ^{7/2}}{N^{4/3}}\Big) 
\text{\quad on $\OmAe $}
.\end{align}
We deduce \eqref{:V2b} from this by a simple calculation with Taylor's theorem. 
\end{proof}

\begin{lemma} \label{l:V3}
Make the same assumptions as in \lref{l:V1}. We set for $ \lpm $, 
\begin{align}\label{:V3a}&
\aNl = \alpha _{N+\ell}(z_\ell)
.\end{align}
Then, we find for $ \lpm $ that 
\begin{align}
 \label{:V3b} &
 \aNl = x_{\ell } +{20}^{-1} N^{-2/3 } x_{\ell } ^2+ \OL (N^{-4/3 } x_{\ell }^3 ) \text{\quad on $\OmAe $}
,\\\label{:V3c} &
 \aNl - \aN = - {N^{-1/3}} \ell + \OL (N^{-1+\varepsilon }) \text{\quad on $\OmAe $},
\\ \label{:V3d}&
\sum_{\lpm } ( \aNl - \aN ) = \OL (N^{-1+\varepsilon }) \text{\quad on $\OmAe $} 
.\end{align}
\end{lemma}
\begin{proof}
% From $ z_0 = z $, $ \aN = \aNz $ is clear. 
From \eqref{:V1y}, $ x_\ell=2(N+\ell)^{2/3} (z_\ell - 1)$. 
Then, using \eqref{:V1q} and \eqref{:V2b}, we deduce \eqref{:V3b}. 
 From this and \eqref{:V1a}, we obtain 
 \begin{align} \label{:V3f}&
 \aNl = x - {N^{-1/3}} \ell +{20}^{-1} N^{-2/3}{ x ^2} + \8 \text{\quad on $\OmAe $} 
.\end{align}
We have \eqref{:V3c} from \eqref{:V2b}, \eqref{:V3f}, and $ \varepsilon < 1/6 $. 
We deduce \eqref{:V3d} from \eqref{:V3c}. 
\end{proof}
%%%%%%%%%% Lemma 12.6 %%%%%%%%%%%

We set $ \aNz = \aN $. 
\begin{lemma} \label{l:V4}
Make the same assumptions as in \lref{l:V1}. 
We set 
\begin{align}\label{:V4x} &
% \aN = \aNz ,\quad 
\gaNl =
 \frac{1}{4} 
 \frac{z_\ell +1}{z_\ell -1}
 \frac{ \aNl }{(N+\ell)^{2/3}} \text{ for $ \lpm , 0 $}, \quad 
\gaN = \gaNz 
.\end{align}
Then, for $ \lpm , 0 $, 
\begin{align}
\label{:V4z}&
\gaNl = \frac{z_\ell+1}{2} \frac{\aNl }{x_\ell}, 
\quad 
 \gaN =\frac{(z + 1) \aN }{ 2x} 
,\\\label{:V4b}& 
\gaNl 
= 1 + \frac{3 x }{10N^{2/3}} - \frac{3 \ell}{10 N} 
+ \OL \Big(\frac{ x ^2 }{N^{4/3}}\Big) \text{\quad on $\OmAe $} 
.\end{align}
\end{lemma}
%%%%%%%%%%%%%%%%%%%%%%%%%%%%%%%%%%%%%%%%%%%%%%
\begin{proof}
From \eqref{:V1y}, we see that $z_\ell -1 = x_\ell /2(N+\ell)^{2/3}$. 
From this, we obtain the first equation in \eqref{:V4z}. 
Recall that $ x_0 = x $, $ z_0 = z $, and $ \aNz = \aN $. 
Hence, we obtain the second equation in \eqref{:V4z} from the first. 
From \eqref{:V1y} and Taylor's theorem, we see that 
\begin{align}&\notag %\label{:}&
z_\ell=
\Big( 1 + \frac{\ell }{N} \Big)^{-1/2} \Big( 1 + \frac{x}{2N^{2/3}}\Big) =
 1 + \frac{x}{2N^{2/3}} - \frac{\ell }{2N } + 
 \OL \Big( \frac{ x }{N^{5/3}}\Big) \text{\quad on $\OmAe $}
.\end{align}
Hence, we have 
\begin{align}\label{:V4f}&
 \frac{z_\ell+1}{2}=
 1 + \frac{x}{4N^{2/3}} - \frac{\ell }{4N } + 
 \OL \Big( \frac{ x }{N^{5/3}}\Big) 
 \text{\quad on $\OmAe $}
.\end{align}
From \eqref{:V1a} and \eqref{:V3b}, we deduce on $\OmAe $ that 
\begin{align} \notag 
 \frac{\aNl }{x_\ell} & 
 = 1 + \frac{x_{\ell } }{20N^{2/3}} + \OL \Big( \frac{x_{\ell }^2}{N^{4/3}}\Big) 
\quad \text{ by \eqref{:V3b}}
 \\ \notag & 
 =
 1 + \frac{x}{20N^{2/3}} - \frac{\ell}{20N}
 + 
 \OL \Big( \frac{ x }{N^{5/3}}\Big) + \OL \Big( \frac{ x ^2}{N^{4/3}}\Big) 
\quad \text{ by \eqref{:V1a}}
 \\\label{:V4g} 
 &=
 1 + \frac{x}{20N^{2/3}} - \frac{\ell}{20N} + \OL \Big( \frac{ x ^2}{N^{4/3}}\Big) 
\quad \text{ by } \vert x \vert > 1 
.\end{align}
Multiplying both sides of \eqref{:V4f} and \eqref{:V4g} with each other, we obtain \eqref{:V4b}. 
\end{proof}

\begin{lemma} \label{l:V5}
Make the same assumptions as in \lref{l:V1}. Then, for $ \lpm $, 
\begin{align} \notag 
\mathbf{I}:= &
\Big( \gaN ^{1/4} \Ai (\aN )
- \frac{\Ai' (\aN )}{2N^{1/3}\gaN ^{1/4} } \Big)^2 
- \prod_{\lpm } \Big( \gaNl ^{1/4} \Ai (\aNl )
- \frac{\Ai' (\aNl )}{2N^{1/3}\gaNl ^{1/4} } \Big)
\\ \label{:V5a} 
&
- ( \gaN ^{1/4} \Ai (\aN ) )^2 + \prod_{\lpm } \gaNl ^{1/4} \Ai (\aNl ) 
= 
\8 \text{\quad on $\OmAe $}
.\end{align}
\end{lemma}

\begin{proof}
By straightforward calculations, we obtain 
\begin{align} \notag 
\mathbf{I}&= 
\frac{-1}{2N^{1/3}}\Big[
 2\Ai (\aN )\Ai' (\aN )
- \sum_{\ell=\pm 1}\Big( \frac{\gaNl }{\gaNll } \Big)^{1/4}\Ai (\aNl )\Ai'(\aNll )
\Big]
\\ \notag
&+ \frac{1}{4N^{2/3}}
\Big[ \frac{\Ai'( \aN)^2}{\gaN ^{1/2}} - \prod_{\ell=\pm 1} \frac{ \Ai' (\aNl )}{\gaNl ^{1/4}}
\Big]
\\ \notag
& = 
\frac{-1}{2N^{1/3}}\Big[
 2\Ai (\aN )\Ai' (\aN )
- \sum_{\ell=\pm 1}\Ai (\aNl )\Ai'(\aNll )
\Big]
\\ \notag
&+ \frac{-1}{2N^{1/3}}\Big[
 \sum_{\ell=\pm 1}\Big(1 -\big( \frac{\gaNl }{\gaNll } \big)^{1/4}\Big) 
\Ai (\aNl )\Ai'(\aNll )
\Big]
\\ \label{:V5f} 
&+ \frac{1}{4N^{2/3}}
\Big[ \frac{\Ai'( \aN)^2}{\gaN ^{1/2}} - \prod_{\ell=\pm 1} \frac{ \Ai' (\aNl )}{\gaNl ^{1/4}}
\Big]
.\end{align}

Using \eqref{:81z}, \lref{l:81}, and \eqref{:V3c}, we have 
\begin{align}&\notag 
 \mathrm{Ai}(\alpha_{N,\ell})-\mathrm{Ai}(\alpha_N) 
=
 \mathrm{Ai}'(\alpha_N)(\alpha_{N,\ell}-\alpha_N) + \underline{\mathcal{O}}(\lvert x \rvert ^{3/4}N^{-2/3}),
\\ \label{:V51}
&\mathrm{Ai}'(\alpha_{N,\ell})-\mathrm{Ai}'(\alpha_N)=
\mathrm{Ai}''(\alpha_N)(\alpha_{N,\ell}-\alpha_N)+ \underline{\mathcal{O}}(\lvert x \rvert ^{5/4}N^{-2/3})
.\end{align}
Note that 
\begin{align}\notag 
&2\mathrm{Ai}(\alpha_N) \mathrm{Ai}'(\alpha_N) 
 - \sum_{\ell=\pm1} 
\mathrm{Ai}(\alpha_{N,\ell})\mathrm{Ai}'(\alpha_{N,-\ell})
\\\notag 
= & \sum_{\ell=\pm1} 
\Big( \mathrm{Ai}(\alpha_N) - \mathrm{Ai}(\alpha_{N,\ell}) \Big) \mathrm{Ai}'(\alpha_N) 
+ \sum_{\ell=\pm1}
 \mathrm{Ai}(\alpha_{N,\ell}) 
\Big( \mathrm{Ai}'(\alpha_N) - \mathrm{Ai}'(\alpha_{N,-\ell}) \Big)
\\\notag 
= &\sum_{\ell=\pm1} 
\Big( \mathrm{Ai}(\alpha_N) - \mathrm{Ai}(\alpha_{N,\ell}) \Big) \mathrm{Ai}'(\alpha_N) 
+ 
 \mathrm{Ai}(\alpha_{N}) \sum_{\ell=\pm1} 
\Big( \mathrm{Ai}'(\alpha_N) - \mathrm{Ai}'(\alpha_{N,-\ell}) \Big) 
\\ \label{:V52} 
+ & \sum_{\ell=\pm1} 
\Big( \mathrm{Ai}(\alpha_{N,\ell}) - \mathrm{Ai}(\alpha_N) \Big) 
\Big( \mathrm{Ai}'(\alpha_N) - \mathrm{Ai}'(\alpha_{N,-\ell}) \Big)
.\end{align}
Hence from \lref{l:81}, \eqref{:V3c}, \eqref{:V3d}, \eqref{:V51} and \eqref{:V52}, we have
\begin{align}\label{:V5g}
\frac{1}{2N^{1/3}}\Big[
 2\Ai (\aN )\Ai' (\aN )
- \sum_{\ell=\pm 1}\Ai (\aNl )\Ai'(\aNll )
\Big] = \8 
.\end{align}

Using Taylor's theorem with \eqref{:V4b}, we have
\begin{align} & \label{:V5G} 
1- \frac{\gaNl }{\gaNll } = \frac{1}{\gaNll }( \gaNll - \gaNl ) = 
\4 
.\end{align} 
Hence using Taylor's theorem to \eqref{:V5G}, \lref{l:81}, and \eqref{:V3b}, we see that 
\begin{align}\label{:V5h}
\frac{1}{2N^{1/3}}\Big[
- \sum_{\ell=\pm 1}\Big(1 -\big( \frac{\gaNl }{\gaNll } \big)^{1/4}\Big) \Ai (\aNl )\Ai'(\aNll )
\Big]= 
 \OL (N^{-4/3}) 
.\end{align}
In the same way, we obtain 
\begin{align}\label{:V5i}
\frac{1}{4N^{2/3}}
\Big[ \frac{\Ai'( \aN)^2}{\gaN ^{1/2}} - 
\prod_{\ell=\pm 1} \frac{ \Ai' (\aNl )}{\gaNl ^{1/4}} \Big] = \4 
.\end{align}
Putting \eqref{:V5g}, \eqref{:V5h}, and \eqref{:V5i} into \eqref{:V5f}, we obtain \eqref{:V5a}. 
\end{proof}

\begin{lemma}[{\cite[Theorem 2.2 (iii)]{DKMVZ}}] \label{l:V6}
There exists $ 0 < \deltaz < 1 $ such that, on \ $\OmAdz $, 
\begin{align}\notag &%\label{:93b}
\psi_N(x)=\gaN (z)^{1/4} \mathrm{Ai}( \alpha _N(z) ) ( 1+ \4 ) 
% - (1/2) N^{-1/3}\gaN (z)^{-1/4}{\mathrm{Ai}'( \alpha _N(z))} ( 1+ \4 )
	- \frac{\mathrm{Ai}'( \alpha _N(z))}{2N^{1/3}\gaN (z)^{1/4}} ( 1+ \4 )
.\end{align}
\end{lemma}
\begin{lemma}\label{l:V7}
Make the same assumptions as in \lref{l:V1}. Let $ \lpm , 0 $. 
Recall that $ (x_\ell ,z_\ell )$ is a function in $ (x , N , z )$ from the assumptions. 
Then, for $ (x , N , z ) \in \OmAe $, 
\begin{align}\label{:V7a}
\psi_{N+\ell}(x_\ell)= 
 \gaNl ^{1/4} \Ai (\aNl ) - 
	% (1/2) N^{-1/3}\gaNl ^{-1/4} \Ai' (\aNl ) 
	\frac{1}{2 N^{1/3}\gaNl ^{1/4}} \Ai' (\aNl ) 
+ \OL (N^{-1} ) 
.\end{align}
\end{lemma}
\begin{proof} 
Recall that $ \aNl = \alpha _{N+\ell}(z_\ell)$ by \eqref{:V3a}. 
From \eqref{:V4x}, we have 
\begin{align}\label{:V7f}&
\gaNl ^{1/4}=
 \frac{1}{2^{1/2}(N+\ell)^{1/6}} 
 (\frac{z_\ell +1}{z_\ell -1})^{1/4}\aNl ^{1/4}
.\end{align}
Recall that $x= \mathbb{U} (N, 2\sqrt{N}z)$ and $N\ge 2$ implies 
$ x_\ell = \mathbb{U} (N+\ell, 2\sqrt{N+\ell}z_{\ell})$. 
Suppose that $(x_\ell, N+\ell, z_\ell) \in \OmAe $, $ \lpm $. 
From \lref{l:81}, \eqref{:V3b}, and \eqref{:V4b}, 
\begin{align}\label{:V7g}
\gaNl ^{1/4} \Ai (\aNl ) = \OL (1), \quad 
\frac{\Ai' (\aNl )}{2N^{1/3}\gaNl ^{1/4}} = \OL (N^{-7/24})
\quad \text{ on $\OmAe $}
.\end{align}
Using \lref{l:V6}, \eqref{:V7f}, and \eqref{:V7g}, 
we obtain \eqref{:V7a} for $(x_\ell, N+\ell, z_\ell)\in \OmAe $. 

Let $ \varepsilon < \varepsilon'< 1/6$. 
Suppose that $(x, N, z) \in \OmAeone $. 
Recall that $(x, N, z) \in \OmAeone $ implies 
$x \in [-(2\deltaz N^{2/3}\wedge N^{\varepsilon }), -1)$. 
Then, using \eqref{:V1a}, we find that there exists $N_0 \in \mathbb{N} $ depending on $\varepsilon $ and $\varepsilon'$ such that
$x_\ell \in [-(2\deltaz N^{2/3}\wedge N^{\varepsilon'}), -1)$ for $N\ge N_0$. 
Then, $(x_\ell,N+\ell,z_\ell) \in \OmAetwo $ if $N \ge N_0$. 
Thus, we obtain \eqref{:V7a} for $ (x , N , z ) \in \OmAe $ 
because $\OmAeone \cap \{N \le N_0\}$ is relatively compact and $\psi_N$ is uniformly bounded on each compact set. This completes the proof.
\end{proof}

\begin{lemma} \label{l:V8}
Make the same assumptions as in \lref{l:V1}. 
\begin{align} 
\label{:V8a} &
\Ai (\aNl ) - \Ai (\aN ) = 
 - N^{-1/3}\ell \Ai' (x) + \OL (N^{-\frac{2}{3}+ \frac{3}{4}\varepsilon }) 
 \quad \text{on $\OmAe $}
,\\ \label{:V8b} &
\lvert \Ai ' (\aNp ) - \Ai ' (\aNm ) \rvert = 
 \OL (N^{-\frac{1}{3}+ \frac{3}{4}\varepsilon })
 \quad \text{on $\OmAe $}
,\\
\label{:V8c}&
-\psi_{N+1} ( x_1 ) + \psi_{N-1} ( x_{-1}) = 
2 N^{-1/3}\Ai' (x) + 
 \OL (N^{-\frac{2}{3}+ \frac{3}{4}\varepsilon })
 \quad \text{on $\OmAe $}
.\end{align}
\end{lemma}
\begin{proof} 
From \eqref{:81z} and \lref{l:81}, we see that, on $\OmAe $,
\begin{align}\notag &% \label{:V8ap} & 
\Ai (\aNl ) - \Ai (\aN )
=
\Ai' (\aN ) (\aNl - \aN ) + 
\mathcal{O}( \lvert \aN \rvert ^{3/4} (\aNl - \aN )^2) 
\\ \notag
= & - N^{-1/3}\ell \Ai' (\aN )+ 
\OL ( \vertx ^{3/4}N^{-2/3} ) 
\ \ \ \quad \text{by \eqref{:V2b}, \eqref{:V3c}}
\\ \notag
= & - N^{-1/3}\ell \Ai' (x) + 
\OL ( \vertx ^{3/4}N^{-2/3} ) 
\quad \quad \ \ 
\text{by \eqref{:V2b}}
.\end{align}
Hence, we obtain \eqref{:V8a}. 
Similarly, from \eqref{:81z}, \lref{l:81}, \eqref{:V3c}, and \eqref{:V2b}, we find that 
\begin{align} \notag 
\lvert \Ai '(\aNp ) - \Ai '(\aNm ) \rvert \le &
2 N^{-1/3} \lvert \Ai '' (\aN ) \rvert + \OL (N^{-\frac{2}{3}} \vertx ^{\frac{5}{4}}) 
 \\ \notag = &\OL ( N^{-1/3} \vertx ^{3/4} ) 
 \quad \text{on $\OmAe $}
.\end{align}
Hence, we obtain \eqref{:V8b}. 
Using \eqref{:V4b}, we obtain, on $\OmAe $, 
\begin{align}\label{:V8j}&
\lvert \gaNp ^{\pm 1/4} - \gaNm ^{\pm 1/4} \rvert \le \OL ( N^{-1}) ,\quad 
%\\\label{:V8k}&
\lvert \gaNp ^{\pm 1/4} -1 \rvert \le \OL ( N^{-2/3 } \vertx ) 
.\end{align} 
Using \eqref{:V7a}, we obtain, on $\OmAe $, 
\begin{align} \notag 
& - \psi_{N+1} ( x_1 ) + \psi_{N-1} ( x_{-1}) 
 = 
- \gaNp ^{1/4} \Ai (\aNp ) + \gaNm ^{1/4} \Ai (\aNm ) 
\\ \label{:V8k} 
&
+ (1/2)N^{-1/3}\gaNp ^{-1/4}{\Ai' (\aNp )} - (1/2) N^{-1/3}\gaNm ^{-1/4} {\Ai' (\aNm )} 
%+ \frac{\Ai' (\aNp )}{2N^{1/3}\gaNp ^{1/4}}- \frac{\Ai' (\aNm )}{2N^{1/3}\gaNm ^{1/4}} 
+ \OL (N^{-1}) 
.\end{align}
Putting \eqref{:V8a}, \eqref{:V8b}, and \eqref{:V8j} into \eqref{:V8k}, we obtain \eqref{:V8c}. 
\end{proof}

\subsection{Global and uniform estimates of $ \psi_N (x) $} \label{s:W}
In subsection \ref{s:W}, we shall present global and uniform estimates of $\psi_N$ and their derivatives (
Lemmas \ref{l:W0}--\ref{l:W3}) and integrals (\lref{l:W5}). 

Let $ \fell $, $ \gell $, $ \ct _\ell^{\KM }$, and $ \Clkm $ be as in \eqref{:U2p}. 
Let %We put for $ \lpm , 0 $, 
\begin{align} &\notag %\label{:W1P}
 \Sell ^{\KM }=\sin \big( \gell + \ct _\ell^{\KM } \big) , \quad 
\wSell ^{\KM } = \sin \big( \gz + \ckm + \ell\theta \big)
,\\ & \label{:W1p} 
\Bellkm = \Clkm \Sell ^{\KM } 
, \quad %\\ \nonumber
\Bellkmh =\Ckm \wSell ^{\KM } 
.\end{align}
We note that $ \mathsf{B}_0^{\KM } = \widehat{\mathsf{B}}_0^{\KM }$. 
We introduce the function $\psiNellH $, $(x, N, \theta) \in \Om$, defined by
\begin{align}\label{:W1q}
\psiNellH ( x) &= \frac{ 1}{\pi}
% \Azz 
{\mathsf{f} (N, \theta )}^{-1/4}
\sumL \Bellkmh \Ntheta 
.\end{align}
Although $ \psiNellH $ depends on $L\in\mathbb{N} $, we omit $L$ from the notation. 
We take and fix $L$ appropriately in each situation. 
We set 
\begin{align}\label{:W3f}
\widetilde{\psi}_N (x)&= \frac{1}{\pi}
{\mathsf{f} (N, \theta )}^{-1/4}
 \sumL 
\Ckm (N ,\theta)
\cos \big( \mathsf{g} (N,\theta) + \ckm (\theta) \big)
.\end{align}
From \eqref{:W1p}--\eqref{:W3f}, and 
$ \sin(A+\theta) - \sin (A-\theta)= 2\sin \theta \cos A $, we see that 
\begin{align} \label{:W3h}&
 -\psiNpH (x) + \psiNmH (x) = -2 (\sin \theta ) \widetilde{\psi}_N (x) 
.\end{align}

Let $ \mathbb{U} $ be as in \eqref{:91b}. 
Let $ w $ be such that $\mathbb{U} (N,w)=x$. Then, 
\begin{align*}&
w= 2\sqrt{N} + N^{-1/6}x
.\end{align*}
For $ \lpm $ and $N\ge 2$, we put $x_\ell=\mathbb{U} (N+\ell, w)$. We then have from \eqref{:31a} 
\begin{align}\label{:W1x}
\psi_{N+\ell}^N(x)= (\oneNell ) ^{-1/12} \psi_{N+\ell}(x_\ell)
.\end{align}
Let $\Om$ be as in \eqref{:93f} and $(x,N, \theta) \in \Om$. 
Then, from \eqref{:91h}, \eqref{:U1z}, and \eqref{:W1x}, we deduce that 
\begin{align}\label{:W1z}
\psi_{N+\ell}^{N}(x)=(\oneNell )^{-1/12}
\psi_{N+\ell}(\h (N+\ell, \theta +\eta_\ell))
.\end{align} 
\begin{lemma} \label{l:W1} 
Let $ \OL (\cdot ) $ and $ \Ome $ be as in \eqref{:93c} and \eqref{:U1y}. 
% Let $(x, N, \theta) \in \Ome $ be as in \lref{l:U1}. 
We can then take $L=L(\varepsilon)\in \mathbb{N} $ such that for $ \lpm , 0 $, on $\Ome$, 
\begin{align}\notag 
\psi_{N+\ell}^N(x) -\psiNellH (x) 
=& \frac{1}{\pi}\Big(
\Al \sumL \Bellkm 
- \Azz \sumL \Bellkmh 
 \Big) 
\\\label{:W1a}
+& \5 
,\\\label{:W1c}
\psi_{N+\ell}^N (x) - \psiNellH ( x) = &\OL ({N^{-1/3}\vertx ^{-5/4}}) 
%\text{\quad on $\Ome$}.
,\\ \label{:W1d}
\psi_{N} (x)^2 - \psiNzH ( x) ^2 =& \OL ({N^{-2/3}\vertx ^{-5/4}}) , 
\quad \text{ for } \ell = 0
.\end{align}
\end{lemma}
%%%%%%%%%%%%%
\begin{proof} 
Let $L \ge \frac{8}{9\varepsilon}+1$. Then, $(3L+1)/4 \ge \frac{2}{3\varepsilon} + 1$. 
From \eqref{:91P} and \eqref{:U2a}, 
\begin{align} \label{:W1f}
 \fell ^{ - (3L+1)/4}
= \OL ( {\vertx ^{-\frac{2}{3\varepsilon} - 1} } ) 
=\5 \text{\quad on $\Ome$}
.\end{align}
From \eqref{:W1z} and \eqref{:W1f}, we see from \eqref{:93p} and \eqref{:93b} that
\begin{align} \notag 
\psi_{N+\ell}^N(x) & = 
( \oneNell )^{-1/12}
\frac{\DDN }{\pi}
\Al 
\sumL 
\Bellkm 
+ \5 \text{\quad on $\Ome$}
\\ \label{:W1g} & = 
\frac{1 }{\pi} \fell ^{-1/4} \sumL \Bellkm + \5 \text{\quad on $\Ome$}
.\end{align}
Hence, from \eqref{:W1q} and \eqref{:W1g}, we obtain \eqref{:W1a}. 

\bigskip 

We next prove \eqref{:W1c}. From \eqref{:U2e} and \eqref{:W1p}, 
\begin{align} \notag 
 \Bellkm - \Bellkmh 
&= \Clkm \mathsf{S}_\ell^{km} - \Ckm \hat{\mathsf{S}}_\ell^{km} 
\quad \text{ by }\eqref{:W1p} 
\\ \notag 
& =( \Clkm - \Ckm )\mathsf{S}_\ell^{km} 
+ \Ckm (\mathsf{S}_\ell^{km}-\hat{\mathsf{S}}_\ell^{km}) 
\\ \notag 
& =\Big( \delta_3^{km}(N,\theta)\ell + \OL (N^{-2/3}\vertx ^{-(3k/4+2)})\Big)\mathsf{S}_\ell^{km}
\quad \text{ by }\eqref{:U2e} 
\\ \label{:W1G}&\quad 
+ \Ckm \Big( \sin (\mathrm{g}_\ell + \mathsf{c}_\ell^{km}) -\sin(\mathrm{g}+\mathsf{c}^{km}+\ell\theta) 
 \Big) 
\quad \text{ by }\eqref{:W1p} 
\end{align}
From \lref{l:9X}, we have for $ k \ge 2$
\begin{align}\label{:W1S}&
 \big\lvert \Ckm (N ,\theta) \big\rvert \le \cref{;9X}\vertx ^{-3k/4} 
\quad \text{for $ x= \h (N,\theta) \in [-2N^{2/3}, -1 ] $}
.\end{align}
Using \eqref{:U2r} and \eqref{:W1S}, we see for $ k \ge 2 $ 
\begin{align}\notag 
\delta_3^{\KM } (N,\theta) &= - \frac{\Ckm (N,\theta)}{2N}
\Big( k + \Big(m+\frac{k}{2}\Big)\frac{\cos^2 \theta}{\sin^2 \theta} \Big) 
\quad \text{ by } \eqref{:U2r} 
\\& \notag 
= \OL (\frac{1}{N\vertx ^{3k/4}\sin^2\theta }) 
 \text{\quad on $\Ome$} \quad \text{ by } \eqref{:W1S} 
\\&\label{:W1T}
= \OL (\frac{1}{N^{1/3} \vertx ^{3k/4+1}}) \text{\quad on $\Ome$}
.\end{align}
Here the last line follows from \eqref{:91f} and \eqref{:91P}. Indeed, from \eqref{:91f}, we have 
$ \mathsf{f} \Ntheta = (N+1)^{2/3}\sin^2 \theta $. From \eqref{:91P}, 
$- \frac{1}{4}x \le\mathsf{f} (N, \theta) \le - 2 x $. 
Combining these, we obtain the last line. 
Because $\sin$ is a $1$-Lipschitz continuous function, by using \eqref{:U2b} and \eqref{:U2d}, we see 
\begin{align}\notag 
\lvert \sin &(\mathrm{g}_\ell + \mathsf{c}_\ell^{km}) - \sin(\mathrm{g}+\mathsf{c}^{km}+\ell\theta) \rvert 
 \le \lvert 
\mathrm{g}_\ell + \mathsf{c}_\ell^{km} -(\mathrm{g} + \mathsf{c}^{km}+\ell\theta)
 \rvert 
\\\notag 
&= \lvert 
(\mathrm{g}_\ell -\mathrm{g}-\theta\ell) + \mathsf{c}_\ell^{km}-\mathsf{c}^{km}
 \rvert 
\\\notag &= 
 \OL (\frac{1}{N^{2/3} \vertx ^{1/2}}) + \delta_2^{km} \ell + \OL (\frac{1}{N \vertx ^{3/2}}) 
 \text{\quad on $\Ome$} 
\quad \text{ by } \eqref{:U2b},\, \eqref{:U2d} 
\\\notag &= 
 \OL (\frac{1}{N^{2/3} \vertx ^{1/2}}) - \frac{(1+2m+k)\cos\theta}{4N\sin\theta} \ell 
+ \OL (\frac{1}{N \vertx ^{3/2}}) 
 \text{\quad on $\Ome$}
\quad \text{ by }\eqref{:U2r} 
\\\label{:W1U} &= 
 \OL (\frac{1}{N^{2/3} \vertx ^{1/2}}) 
 \text{\quad on $\Ome$}
.\end{align}
Here the last line follows from the derivation similar to that of the last line of \eqref{:W1T}. 
Putting \eqref{:W1S}--\eqref{:W1U} into the right-hand side of \eqref{:W1G}, we have 
\begin{align} \label{:W1h}
 & \Bellkm - \Bellkmh = \OL ( {N^{-1/3} x ^{-1}} ) \quad \text{ on $\Ome$}
.\end{align}
From \eqref{:W1p}, $ \lvert \Bellkm \rvert = \lvert \Clkm \Sell ^{\KM } \rvert \le \lvert \Clkm \rvert $. 
Hence from \eqref{:U1d} and \eqref{:W1S}, 
\begin{align} \label{:W1H}
\Bellkm = \OL (1) \quad \text{ on $\Ome$}
.\end{align}
 
Using \eqref{:91P}, \eqref{:W1f}, and \eqref{:U2a} and applying Taylor's theorem, we deduce that 
\begin{align} \label{:W1i}
& \Al - \Azz =\OL ( {N^{-1/3} \vertx ^{-5/4}}) \text{\ on $\Ome$}
.\end{align}
Using \eqref{:W1a}, we obtain 
\begin{align} & \notag 
\lvert \psi_{N+\ell}^N(x) -\psiNellH (x) \rvert 
\le \frac{1}{\pi}
\Big\lvert \Al - \Azz \Big\rvert 
\Big\lvert \sumL \Bellkm \Big\rvert 
\\\notag 
&\quad \quad + \frac{1}{\pi}\Azz \Big\lvert \sumL (\Bellkm -\Bellkmh ) \Big\rvert 
+ \5 
\text{ on $\Ome$} 
.\end{align}
Combining this with \eqref{:91P} and \eqref{:W1h}--\eqref{:W1i}, we obtain \eqref{:W1c}. 
Recall that $ \mathsf{B}_0^{\KM } = \widehat{\mathsf{B}}_0^{\KM }$. 
The proof of \eqref{:W1d} is then similar to that of \eqref{:W1a} and thus omitted. 
\end{proof} 
%%%%%%%%%%%%

The following properties of $ \psi _{N }$ $ (N \in\mathbb{N} ) $ are well known \cite[Lemma 3.9.9]{AGZ10}. 
\begin{align} \label{:W3y} &
N ^{1/6}\psi_{N }'(x)
=-\frac{\sqrt{N +1}}{2} \psi_{N +1}^{N }(x)
+\frac{\sqrt{N }}{2} \psi_{N -1}^{N }(x)
,\\\label{:W3x}
&\psi_N ''(x)=\Big(
\frac{x^2}{4N ^{2/3}}+x - \frac{1}{2N ^{1/3}}
\Big) \psi_N (x)
.\end{align}
From \eqref{:W3y}, we have 
\begin{align}\label{:W3z}
\psi'_N(x)= \frac{N^{1/3}}{2} \Big( - (\oneN )^{1/2} \psi_{N+1}^N (x) + \psi_{N-1}^N(x) \Big) 
.\end{align}

\begin{lemma} \label{l:W0}
There exists a positive constant $ \Ct \label{;W3e}$ independent of $ N $ such that 
 \begin{align} 
\label{:W3E}&
 \lvert \psi_N(x) - \Ai (x) \rvert \le \cref{;W3e} N^{-\frac{1}{3}+ \frac{1}{4}\varepsilon},
& x \in [-N^{\varepsilon},0]
,\\\label{:W3e}& 
\lvert \psi_N' (x) - \Ai' (x) \rvert \le \cref{;W3e} N ^{-\frac{1}{3} +\frac{3}{4} \varepsilon }, 
& x \in [-N^{\varepsilon},0] 
.\end{align} 
\end{lemma}
\begin{proof}
Recall that $ x = x_0 $ by \eqref{:V1w} and \eqref{:V1x}, 
 $ \aN = \aNz $ by \lref{l:V3}, and $ \gaN = \gaNz $ by \eqref{:V4z}. 
Let $ \OL (\cdot)$ and $ \OmAe $ be as in \eqref{:93c} and \eqref{:V1q}. 
Then, from \eqref{:V2b}, \eqref{:V4b}, and \eqref{:V7a}, we have 
\begin{align}\notag %\label{:W3D}
\aN = & x +{(1/20)} N^{-2/3}{ x ^2} + \OL (N^{-4/3} x ^3 )
&\text{on $\OmAe $}
,\\ \notag 
\gaN = &
1 + (3/10) N^{-2/3} x %\frac{3 x }{10N^{2/3}} 
+ \OL (N^{-4/3} x ^2 )
% \OL \Big(\frac{ x ^2 }{N^{4/3}}\Big) %)^{1/4} 
&\text{on $\OmAe $}
,\\ \notag 
\psi_{N}(x)= & 
 \gaN ^{1/4} \Ai (\aN ) - (1/2) N^{-1/3}\gaN ^{-1/4} \Ai' (\aN ) 
 + 
 \OL (N ^{-1} ) 
&\text{on $\OmAe $}
.\end{align}
Combining these with the estimate of $ \Ai ' $ in \lref{l:81}, we obtain \eqref{:W3E}. 
Let $x_{\pm1}$ be as in \lref{l:V1}. 
From \eqref{:V8c}, \eqref{:W1x}, \eqref{:W3z}, and \eqref{:W3E}, we obtain on $\OmAe $ 
\begin{align} \notag 
 \psi_N'(x) 
&= \frac{N^{1/3}}{2} \Big( - (\oneN )^{1/2} \psi_{N+1}^N (x) + \psi_{N-1}^N(x) \Big) 
\quad \text{ by \eqref{:W3z} }
\\&\notag 
= \frac{N^{1/3}}{2} \Big(-\psi_{N+1} ( x_1 ) + \psi_{N-1} ( x_{-1}) \Big) 
+\OL ( {N^{-2/3}}) % \text{\quad on $\OmAe $}
\quad \text{by \eqref{:W1x}, \eqref{:W3E}}
,\\ \label{:W3d}
 &
=\frac{N^{1/3}}{2}
\Big( \frac{2}{N^{1/3}} \Ai' (x) + \OL (N^{-\frac{2}{3}+ \frac{3}{4}\varepsilon }) 
\Big) 
+\OL ( {N^{-2/3}}) \quad \text{by \eqref{:V8c}}
.\end{align} 
% Here, we used \eqref{:V8c} for the second line. 
From \eqref{:W3d}, we obtain \eqref{:W3e}. 
\end{proof}

\begin{lemma} \label{l:W2}
There exist positive constants $\Ct \label{;Wg} $ and $\Ct \label{;Wh} $ 
 independent of $ N $ such that
\begin{align}
\label{:W2a} & 
\lvert \psi_{N+ \ell}^N(x)\rvert \le \cref{;Wg} \exp (- \xF /\cref{;Wh} ) 
 \text{ for } x \in [-1, \infty )
%\quad x \ge -1 
,\, \lpm , 0 
.\end{align}
There exist positive constants $\Ct \label{;W3n} $ and $\Ct \label{;W3Y} $ independent of $ N $ such that 
\begin{align}
\label{:W2z}&& 
\lvert \psi_{N} (x) \rvert &\le \cref{;W3n} \exp (- \xF /\cref{;W3Y} ) %\quad x \ge -1 
&& \text{ for } x \in [-1, \infty )
,\\ && \label{:W3n} 
\lvert \psi'_N(x)\rvert &\le \cref{;W3n}\exp ({ - \xF /{ \cref{;W3Y}} }) 
&& \text{ for } x \in [-1, \infty )
,\\ \label{:W3Y} &&
\lvert \psi_N ''(x) \rvert &\le \cref{;W3n} \exp (- \xF / \cref{;W3Y}) 
&& \text{ for } x \in [-1, \infty )
.\end{align}
\end{lemma}
%%%%%%%%%%%%%
\begin{proof} 
Because $ \psi_{N}^N $ is uniformly bounded on $[-1,2]$, 
it suffices to prove \eqref{:W2a} for $x \ge 2$. 
Let $\Omh $ be as in \eqref{:93g} and $(x,N,\theta) \in \Omh $ .
From \eqref{:92z} and \eqref{:93g}, for $ \ x \ge 2 $, %we have
\begin{align}\label{:W2b}
x &\le 2N^{1/6}\sqrt{N+1} (\cosh \theta -1 ) 
+ {N^{-1/3}} 
 \le 4N^{1/6}\sqrt{N+1} (\cosh \theta -1 )
.\end{align}
Let $\mathsf{g} ^{\rm h}$ be as in \eqref{:92y}. From \eqref{:W2b}, 
\begin{align}\label{:W2f}
-\frac{\mathsf{g} ^{\rm h}(N,\theta)}{x^{3/2}} 
&\ge \frac{1}{2^{4}} \frac{\sinh 2\theta - 2\theta}{(\cosh \theta -1)^{3/2}}
.\end{align}
Put $ u(\theta)=\sinh 2\theta - 2\theta$ and $v (\theta)=(\cosh \theta -1)^{3/2}$.
We then have 
\begin{align}\label{:W2h}
u(\theta) \ge v(\theta), \quad \theta \ge 0 
.\end{align}
In fact, \eqref{:W2h} follows from $ u(0)=v(0)=u'(0)=v'(0)=0 $ and 
$ u''(\theta) \ge v''(\theta), \ \theta\ge 0 $. 
Using \eqref{:W2f} and \eqref{:W2h}, we obtain 
\begin{align}\label{:W2j}
-\mathsf{g} ^{\rm h} (N, \theta) \ge 2^{-4} x^{3/2}, \quad x\ge 2 
.\end{align}
From \eqref{:92a}, we easily see that 
\begin{align}\label{:W2k}
\mathsf{f} ^{\rm h}(N,\theta) \ge 2^{-1} x , \quad x\ge 2 
.\end{align}
In the same fashion as \eqref{:U1d}, we can check that 
\begin{align}\label{:W2l}
\sup \{\CC ^{\rm h}_{\KM }(N, \theta) ; {N\in \mathbb{N} , \; x\ge 2} \}< \infty 
\quad \text{ for each $k,m \in \zN $}
.\end{align}
Hence, we obtain \eqref{:W2a} for $ \ell = 0 $ from \eqref{:93q} with \eqref{:W2j}--\eqref{:W2l}. 

Let $\eta^{\rm h}_\ell=\eta^{\rm h}_\ell (N,\theta) \le 0$ such that $\theta + \eta^{\rm h}_\ell \ge 0$ and
\begin{align}\label{:W2p}
\sqrt{N+1}\cosh \theta=\sqrt{N+1+\ell} \cosh (\theta + \eta^{\rm h}_\ell) 
.\end{align}
We put $x_\ell=\h ^{\rm h}(N+\ell, \theta + \eta^{\rm h}_\ell)$. 
Then, from \eqref{:92z} and \eqref{:W2p}, we deduce that 
\begin{align*}
x_\ell &= 2(N+\ell)^{1/6}( \sqrt{N+1}\cosh \theta - \sqrt{N+\ell} ) 
.\end{align*}
We see that%, on \ $\Omh \cap \{ x\ge 2 \}$, 
\begin{align} \notag %\label{:W2q}
\frac{x_\ell}{x} & 
= (\oneNell )^{1/6}
\frac{ (\oneN )^{1/2} \cosh \theta - (\oneNell )^{1/2}} {(\oneN )^{1/2} \cosh \theta - 1 }
\\&\notag 
= 1+ \OL (N^{-1/3}) 
\quad \text{ on \ $\Omh \cap \{ x\ge 2 \}$}
.\end{align}
Combining this with \eqref{:W1x} and \eqref{:W2a} for $ \ell = 0$ yields \eqref{:W2a} for $ \lpm $. 

Recall that $ \psi_{N}^N = \psi_N $ by \eqref{:31a}. Thus, \eqref{:W2z} follows from \eqref{:W2a}.

Using \eqref{:W3x} and \eqref{:W2z}, we have a positive constant $ \cref{;W3Y}$ such that 
\begin{align}\notag %\label{:W3Y-} 
\lvert \psi_N ''(x) \rvert &= 
\Big\lvert\Big( \frac{x^2}{4N ^{2/3}}+x - \frac{1}{2N ^{1/3}} \Big) \Big\rvert 
\lvert \psi_N (x) \rvert 
\\&\notag 
\le \cref{;W3n} \exp (- \xF / \cref{;W3Y}) 
\quad \text{ for } -1 \le x < \infty 
.\end{align}
This implies \eqref{:W3Y}. 

From \eqref{:W2z}, we see $ \limi{x} \psi_N (x) = 0 $. From \eqref{:W3Y}, we have a finite limit 
\begin{align}&\label{:W2x}
 \limi{x} \psi_N ' (x) = \psi_N ' (-1) + \int_{-1}^{\infty} \psi_N '' (y) dy < \infty 
.\end{align}
Hence using $ \limi{x} \psi_N (x) = 0 $, we see $ a:= \limi{x} \psi_N ' (x) = 0 $. 
Indeed, if $ a \ne 0 $, say $ a > 0 $, then there exists an $ R $ such that $ \psi_N ' (x) \ge a/2 $ for all $ x \ge R $. 
Then $ \limi{x}\psi_N (x) - \psi_N (R) = \int_R^{\infty} \psi_N ' (y) dy = \infty $, which yields contradiction. 

Thus, we have from \eqref{:W2x} and $ \limi{x} \psi_N (x) = 0 $ 
\begin{align*}&
\psi_N ' (-1) = - \int_{-1}^{\infty} \psi_N '' (y) dy 
.\end{align*}
Hence, 
\begin{align*}&
\psi_N ' (x) = \psi_N ' (-1) + \int_{-1}^x \psi_N '' (y) dy 
= - \int_{x}^{\infty} \psi_N '' (y) dy
.\end{align*}
Using this and \eqref{:W3Y}, we obtain \eqref{:W3n}. 
\end{proof}

\begin{lemma} \label{l:W3} 
There is a positive constant $ \Ct \label{;W66}$ independent of $ N $ such that
\begin{align}\label{:W3a}
&\lvert \psi_{N }(x)\rvert \le \cref{;W66}(\xvee )^{-1/4} , 
&x \in [-2N ^{2/3},\infty), 
\\ \label{:W3b}
&\lvert \psi_{N }'(x)\rvert \le \cref{;W66}(\xvee )^{1/4} , 
&x \in [-2N ^{2/3},\infty), 
\\ \label{:W3c}
&\lvert \psi_{N }''(x)\rvert \le \cref{;W66} (\xvee )^{3/4} , 
& x \in [-2N ^{2/3},\infty) 
.\end{align}
\end{lemma}
\begin{proof} 
For $ x \in [-2N ^{2/3}, -1]$, \eqref{:W3a} follows from \eqref{:91P}, \eqref{:93p}, and \eqref{:U1d}. 
For $ x \in [-1,1]$, \eqref{:W3a} is clear because of the uniform bound of $ \psi_N $ on $ [-1,1]$. 
For $ x \in [1, \infty )$, \eqref{:W3a} follows from \eqref{:W2a}. 
Applying \lref{l:W2} and \eqref{:W3a} to \eqref{:W3x}, we obtain \eqref{:W3c}.

For \eqref{:W3b}, we divide the case into three parts: 
 $ [-2N ^{2/3},-N ^{\varepsilon}]$, 
 $ (-N ^{\varepsilon},0)$, and 
 $ [0,\infty)$. 
Suppose $ [-2N ^{2/3},-N ^{\varepsilon}]$. 
From \eqref{:W3h}, \eqref{:W1c}, and \eqref{:W3z}, we deduce that 
\begin{align} \notag 
\psi_N'(x) & ={N^{1/3}}{2^{-1}} ( -\psiNpH (x) + \psiNmH (x) ) 
+\OL ( \vertx ^{-5/4}) \text{\quad on $\Ome$}
\\ \label{:W3g} & = 
N^{1/3} (\sin \theta ) \widetilde{\psi}_N (x) +\OL ( \vertx ^{-5/4}) \text{\quad on $\Ome$}
.\end{align}
Note that $ x=\h (N,\theta) $ on $\Ome$ by \eqref{:93f} and \eqref{:U1y}. 
From this and \eqref{:91h}, we see that 
\begin{align}\label{:W3k}&
x= 2N^{2/3}(\cos \theta-1 ) + \3 
= -2N^{2/3}\frac{\sin^2 \theta}{1+\cos \theta} + \3 \text{\quad on $\Ome$}
.\end{align}
Using \eqref{:W3g} and \eqref{:W3k}, we deduce that 
\begin{align}\label{:W3l} 
\psi_N'(x) 
&= - 2^{-1/2} (1+\cos \theta)^{1/2} \vertx ^{1/2} \widetilde{\psi}_N(x)
+\OL ({\vertx ^{-5/4}}) \text{\quad on $\Ome$}
.\end{align}
Using \eqref{:91P}, \eqref{:U1d}, and \eqref{:W3f}, we obtain 
\begin{align} \label{:W3m} &%
\widetilde{\psi}_N(x) =\OL ( \mathsf{f} (N, \theta )^{-1} )^{1/4}
= \OL ({\vertx ^{-1/4}}) \text{\quad on $\Ome$}
.\end{align}
Combining \eqref{:W3l} and \eqref{:W3m}, we obtain \eqref{:W3b} for the case that 
$ [-2N ^{2/3},-N ^{\varepsilon}]$. 

Suppose $ ( -N ^{\varepsilon},0) $. We then obtain \eqref{:W3b} from \eqref{:W3e} and \lref{l:81}. 
Suppose $ [0,\infty)$. Then we deduce \eqref{:W3b} from \eqref{:W3n}, which complete the proof. 
\end{proof}

\begin{lemma} \label{l:W5} 
Let $ \varepsilon \psi_{N } $ be as in \eqref{:31g}. For each $a\in \mathbb{R} $, 
\begin{align}
\label{:W5a} &
 \sup_{N \in\mathbb{N} }\sup_{x\in\mathbb{R} } \Big\lvert \int_a^{x} \psi_{N }(\uU ) d\uU \Big\rvert < \infty 
,\\ \label{:W5c} & 
 \sup_{N \in\mathbb{N} }\sup_{x\in\mathbb{R} }\lvert \varepsilon \psi_{N }(x)\rvert < \infty 
.\end{align}
\end{lemma}
\begin{proof}
We only present the proof for the case that $x\in [-2N^{2/3}, \infty)$ because of \eqref{:31b}. Recall that 
$\CC _{00}=\pi ^{1/2} $, $ \CC _{10}=\CC _{11}=0 $ and $\ct _{00}=\pi/4 - \theta/2$. 
From \eqref{:91p}, \eqref{:93p} with $L=2$, and \eqref{:W4a}, we see that 
\begin{align} &\label{:W5g}
\supN \sup_{ x \in [-2N ^{2/3}, -1) } \Big\lvert \int_x^{ - 1 } \psi_{N }(\uU ) d\uU \Big\rvert <\infty 
.\end{align}
From \eqref{:W2a} with \eqref{:W3a}, we have
\begin{align}& \label{:W5h}
\sup_{N \in\mathbb{N} } \sup_{x\in [- 1, \infty)} \Big\lvert \int_{-1}^x \psi_{N }(\uU ) d\uU \Big\rvert <\infty
.\end{align}
From \eqref{:W5g} and \eqref{:W5h}, we deduce \eqref{:W5a}. 
Using \eqref{:W5a} and \eqref{:31c}, we obtain \eqref{:W5c}. 
\end{proof}

\section{Appendix 3: Remaining proofs of \sref{s:3}} \label{s:Ap3}
In this section, we complete the proof of claims in \sref{s:3}. 
In Subsection \ref{s:X}, we prove \pref{l:35}. 
In Subsection \ref{s:Y1}, we complete the proof of \lref{l:38} except \eqref{:Y1g}. 
In Subsection \ref{s:Z}, we prove \lref{l:3R}. 
In Subsection \ref{s:19}, we prove \eqref{:Y1g}.

\SECT{Proof of \pref{l:35}} \label{s:X}
%\section{Appendix 3: Proof of \pref{l:35}} \label{s:X}

First, we consider the case that $\beta=2$ (\lref{l:X6}). 
For this, we prepare Lemmas \ref{l:X1}--\ref{l:X5}. 
We set 
\begin{align}\label{:X1z}&
\PsiN=\psi_N(x)^2 -\psi_{N+1}^{N}(x)\psi_{N-1}^{N}(x) 
.\end{align}
\begin{lemma}\label{l:X1} 
We find a constant $ \Ct \label{;X1}$ 
independent of $ N \in\mathbb{N}$ such that 
\begin{align} \label{:X1a}
\lvert \rho_2^{N,1}(x) - N^{2/3} \PsiN \rvert \le {\cref{;X1}}{N^{-1/3}({\xvee })^{-1/2}}, 
\quad x\in [-2N^{2/3}, \infty) %,\, N \in \mathbb{N}
.\end{align}
\end{lemma}
\begin{proof}
Note that \eqref{:X1a} is clear for $ N=1$. 
Hence, we suppose $ N \ge 2 $. 

From the Christoffel--Darboux formula (see \cite[p. 511]{Meh04}), we have 
\begin{align} \notag 
\rho_2^{N,1}(x) 
&= {N^{-1/3}} ( N\psi_N(x)^2 -\sqrt{N(N+1)}\psi_{N+1}^N (x) \psi_{N-1}^N(x) )
\\ \label{:X1b}
&= N^{2/3} \PsiN + N^{1/6}( \sqrt{N} - \sqrt{N+1})\psi_{N+1}^N (x)\psi_{N-1}^N(x) 
.\end{align}
Then, \eqref{:X1a} is clear for $x \in [0,\infty)$ from \lref{l:W2}. 

Let $x \in [-2N^{2/3}, 0)$. 
Because $ \lvert \sqrt{N} - \sqrt{N+1}\rvert \le \half N^{ - 1/2}$, we have from \eqref{:X1b} 
\begin{align}& \label{:X1c}
\lvert \rho_2^{N,1}(x) - N^{2/3} \PsiN \rvert \le 
{2^{-1} N^{-1/3}} \lvert \psi_{N+1}^N (x)\psi_{N-1}^N(x) \rvert 
.\end{align}
Let $ x_{\ell }$ be as in \eqref{:V1w} and \eqref{:V1x}. 
Let $ \OL (\cdot ) $ and $ \Om $ be as in \eqref{:93c} and \eqref{:93f}. Then, 
\begin{align}\notag %\label{:X1d}
x_\ell -x &= ((N+\ell)^{1/6} - N^{1/6}) w -2 ( (N+\ell)^{2/3} - N^{2/3} ) 
= \OL ( {N^{-1/3}} ) 
\text{\quad on \ $\Om$}
.\end{align}
Combining this with \eqref{:W1x} and \eqref{:W3a}, we find 
 a constant $ \Ct \label{;X2} > 0$ such that 
\begin{align} \label{:X1e}
\lvert \psi_{N+\ell }^N(x)\rvert &=(\oneNell )^{-1/12} \lvert \psi_{N+\ell }(x_\ell)\rvert 
 \le {\cref{;W66}}{( 1 \vee \lvert x_{\ell }\rvert )^{-1/4}}
 \le {\cref{;X2}}{(\xvee )^{-1/4}} 
\end{align}
for all $ x\in [-2N^{2/3}, 0)$. 
Hence, we obtain \eqref{:X1a} from \eqref{:X1c} and \eqref{:X1e}. 
\end{proof}

Let $ \Bellkm $ and $ \Bellkmh $ be as in \eqref{:W1p}. 
We take $L$ as in \eqref{:W1f}. 
Let $ \Ome $ be as in \eqref{:U1y}. 
\begin{lemma} \label{l:X2}
% Let $\varepsilon \in (0,1/3)$. [[[
For each $ (m,k),\, (n,l) \in \{ (i,j)\, ;\, 0\le i \le j \le L-1 \} $, 
 \begin{align}\label{:X2e}&
 & \sum_{\lpm } ( 
 \Bell ^{\KM } \BB _{-\ell }^{\LN } - 
\Bellkmh \wB _{- \ell }^{\LN } ) 
=\6 \text{\quad on $\Ome$}
.\end{align}
\end{lemma}
\begin{proof}
We divide the proof of \eqref{:X2e} into five steps. 

\noindent {\bf Step I: } 
Let $ \Clkm $, $ \Sell ^{\KM } $, and $ \wSell ^{\KM } $ be as in \eqref{:W1p}. For $ \lpm $, let 
\begin{align}\label{:X2f}
&\CCC _{\ell }= \CCl ^{\KM } \CC _{-\ell }^{\LN },
\quad \mathbb{S}_{\ell }=\Sell ^{\KM } \sSS_{-\ell }^{\LN },
\quad \widehat{\mathbb{S}}_{\ell }= 
\wSell ^{\KM } \wS _{-\ell }^{\LN }
.\end{align}
Then, from \eqref{:W1p} 
\begin{align} \label{:X2g} & 
 \sum_{\lpm } ( \Bell ^{\KM } \BB _{-\ell }^{\LN } - \Bellkmh \wB _{- \ell }^{\LN } ) = 
\Cz \sum_{\lpm } ( \mathbb{S}_{ \ell } -\widehat{\mathbb{S}}_{ \ell } ) + 
\sum_{\lpm } ( \Cl -\Cz ) \mathbb{S}_{ \ell } 
.\end{align}

\noindent {\bf Step II: } 
Let $\fell $, $ \gell $, and $ \ct _\ell^{\KM } $ be as in \eqref{:U2p}. 
From \eqref{:U2b} and \eqref{:U2d}, 
\begin{align}\label{:X2h}
\gell +\ct _\ell^{\KM } - (\gz + \ckm + \ell\theta)
= \6 \text{\quad on $\Ome$}
.\end{align}
Using Taylor's theorem with \eqref{:W1p}, \eqref{:X2f}, and \eqref{:X2h}, 
we see that for $ \lpm $, 
\begin{align}\label{:X2i}
\mathbb{S}_{ \ell } - \widehat{\mathbb{S}}_{ \ell } 
=\6 \text{\quad on $\Ome$}
.\end{align}
Because $\Cz $ is bounded from \eqref{:U1d}, we have from \eqref{:X2i}, 
\begin{align}\label{:X2j}
 \Cz \sum_{\lpm } ( \mathbb{S}_{ \ell } -\widehat{\mathbb{S}}_{ \ell } ) =& 
 \6 \text{\quad on $\Ome$}
.\end{align}
{\bf Step III: }
From Taylor's theorem with \eqref{:U2e} and \eqref{:X2f} and using \eqref{:U1d}, 
we see that 
\begin{align} \notag %abel{:X2k} 
&\sum_{\lpm } ( \Cl -\Cz ) 
=\CC _1^{km}\CC _{-1}^{ln} - 2\Ckm \Cln + \CC _1^{ln}\CC _{-1}^{km} 
\\ \nonumber %
=& (\Ckm + \delta_3^{km})(\Cln - \delta_3^{ln}) - 2\Ckm \Cln 
+ (\Cln + \delta_3^{ln})(\Ckm - \delta_3^{km})
+ \OL ({N^{-2/3}\vertx ^{-7/2}}) 
\\ \notag 
=& - 2 \delta_3^{km} \delta_3^{ln} +
 \OL ({N^{-2/3}\vertx ^{-7/2}})% \OL \Big( \frac{1}{N^{2/3}\vertx ^{7/2}}\Big) 
= \OL ({N^{-2/3}\vertx ^{-7/2}})
 \text{\quad on $\Ome$} 
\quad \text{ by \eqref{:U2r}}
.\end{align}
From this and $\lvert \mathbb{S}_{p}\rvert \le 1, p=\pm1$,
 which follows from \eqref{:W1p} and \eqref{:X2f}, we see on $\Ome$ that 
\begin{align} \label{:X2l} 
& \sum_{\lpm } ( \Cl -\Cz ) \mathbb{S}_{ \ell } = (\Cone -\Cz )(\SSSSp )
 + \OL ({N^{-2/3}\vertx ^{-7/2}})
.\end{align}
 
\noindent {\bf Step IV: } 
Suppose that $ k \ge 2 $ or $\llll \ge 2 $. 
Then from \eqref{:U2e}, we deduce % that for $ k \ge 2 $ or $\llll \ge 2 $, 
\begin{align}\label{:X2n}
& \Cone -\Cz =
\OL \Big( \frac{\Ckm \Cln }{N\sin^2 \theta}\Big) 
\text{\quad on $\Ome$}
.\end{align}
From \eqref{:U1d}, we have $ \Ckm \Cln =\OL ( x ^{-3} ) $ on $\Ome$. 
From \eqref{:93G}, we have
\begin{align*}&
{(N\sin^2 \theta)^{-1}}=N \6 ^2 = \OL (N^{-1/3} x ^{-1})
\text{\quad on $\Ome$}
.\end{align*}
Putting these into \eqref{:X2n}, we have a constant $\Ct \label{;X3} > 0$ such that 
for all $ N\in \mathbb{N}$ 
\begin{align}\label{:X2o}
\lvert \Cone -\Cz \rvert \le \cref{;X3} N^{-1/3} x ^{-4},
\quad x\in \None %, \quad N\in \mathbb{N}
.\end{align}
 
Using $\lvert \mathbb{S}_{\ell }\rvert \le 1$, $ \lpm $, again, we deduce from \eqref{:X2o} that,
for $ N\in \mathbb{N}$, 
\begin{align} \label{:X2p} 
\lvert (\Cone -\Cz )(\SSSSp )\rvert 
\le { \cref{;X3} }{N^{-2/3}\vertx ^{-3}}
, \quad x \in (-2N^{2/3}, -N^{1/3}] %\, N\in \mathbb{N}
.\end{align}
Next, suppose that $ x \in [- N^{1/3}, -N^{\varepsilon}] $. Recall that $x=\h (N,\theta)$. 
Then, from \eqref{:91h}, %we see that $ \lvert 1-\cos \theta \rvert \le N^{-2/3}\vert x \vert $. 
\begin{align*}&
 \lvert 1-\cos \theta \rvert \le N^{-2/3}\vert x \vert 
.\end{align*}
Then, $ \lvert \theta\rvert \le c {N^{-1/3}}\vert x \vert ^{1/2} $ with some $c>0$.
From this, \eqref{:W1p}, and \eqref{:X2f}, we see that 
\begin{align} &\notag % \label{:X2q}&
 \SSSSpt = \OL (\theta )= %\3 
 \OL ({N^{-1/3}}\vert x \vert ^{1/2} ) 
\text{\quad on $\Ome$}
.\end{align}
Hence, combining this with \eqref{:X2i}, we have 
\begin{align} \label{:X2t} 
\SSSSp & = \SSSSpt + \6 
= 
 \OL ({N^{-1/3}}\vert x \vert ^{1/2} ) % \3 
 \text{\quad on $\Ome$}
.\end{align}
We therefore obtain for $ k \ge 2 $ or $\llll \ge 2 $ 
\begin{align} \label{:X2v}&
(\Cone -\Cz )(\SSSSp ) = 
 \OL ({N^{-2/3}\vert x \vert ^{-3}}) \text{\quad on $\Ome$}
\end{align}
for $x\in [-2N^{3/2}, -N^{1/3}]$ from \eqref{:X2p} 
and for $x\in [-N^{1/3}, -N^{\varepsilon}]$ from \eqref{:X2t} with \eqref{:X2o}. 

Suppose $ k , \llll \in \{ 0,1 \} $. 
Then $ \Cone -\Cz = 0 $ by definition. 
We thus obtain \eqref{:X2v} for all $ k, \llll \in \zN $.

%Suppose $ k =1 $ or $ \llll = 1$. Then $ \Cone = \Cz = 0 $. Suppose $ k=\llll = 0 $. Then 

\noindent {\bf Step V: }Combining \eqref{:X2l} and \eqref{:X2v}, we obtain 
\begin{align}\label{:X2w}&
\sum_{\lpm } ( \Cl -\Cz ) \mathbb{S}_{ \ell } 
= \OL ({N^{-2/3} x ^{-3}}) \text{\quad on $\Ome$}
.\end{align}
From \eqref{:X2g}, \eqref{:X2j}, and \eqref{:X2w}, we obtain \eqref{:X2e}. 
 This completes the proof. 
\end{proof}

\begin{lemma} \label{l:X3}
Let $ \Bellkm , \Bellkmh $ 
be as in \eqref{:W1p}. Take $L$ as in \eqref{:W1f}. 
Then, on $\Ome$, 
\begin{align}\label{:X3b}&
\sumL \sumLL
 \big( \Alll \Bellp ^{\KM } \Bellm ^{\LN } - \fz ^{-1/2} \wB _1^{\KM } \wB _{-1}^{\LN } \big) = \5 
.\end{align}
\end{lemma}
\begin{proof}
From \eqref{:91P}, $ \Azz \le {\cref{;Xe}}/{\vertx ^{1/4}} $ for $ x \in [-2N^{2/3}, -1]$. 
Hence from \eqref{:U1d} and \eqref{:W1p}, 
there exist positive constants $\Ct \label{;Xe}$ and $\Ct \label{;Xd}$ 
independent of $N$ and $\theta$ such that 
\begin{align} \label{:X3c} &
 \max_{0\le m \le k \le L}\Bellkm \le \cref{;Xd},
\quad \max_{0\le m \le k \le L}\Bellkmh \le \cref{;Xd} 
\quad \text{on $\Ome$}
\end{align}
for $\lpm $. Using Taylor's theorem, we obtain, on $\Ome$, 
\begin{align}\notag & 
\Alll -\fz ^{-1/2}
\\ \notag =& \prod
_{\lpm }\Big(\fz ^{-1/4} - \frac{1}{4}\fz ^{-5/4}(\fell - \fz ) +\frac{5}{32}\fz ^{-9/4}(\fell -\fz )^2 + \OL\Big(\frac{1}{N \vertx ^{13/4}}\Big)\Big) 
-\fz ^{-1/2} 
\\ \notag =& - 
 \frac{1}{4}\fz ^{-3/2} \sum_{\lpm } (\fell - \fz ) 
 + 
 \frac{5}{32} \fz ^{-5/2} \sum_{\lpm } (\fell -\fz )^2 
 + 
 \frac{1}{16} \fz ^{-5/2} \prod_{\lpm }(\fell -\fz ) + \OL\Big(\frac{1}{N \vertx ^{7/2}}\Big)
% \quad \text{on $\Ome$}
.\end{align}% 
Hence from \eqref{:91P}, \eqref{:U2r}, and \eqref{:U2a}, we deduce 
\begin{align}& \label{:X3D}
\Alll -\fz ^{-1/2} = \OL ( {N^{-2/3}\vertx ^{-5/2}}) \quad \text{ on $\Ome$}
.\end{align}
From \eqref{:X2e}, \eqref{:X3c}, and \eqref{:X3D}, we obtain \eqref{:X3b}. 
\end{proof}

Let $\psiNellH $ be as in \eqref{:W1q}. We set 
\begin{align}\label{:X4x}& 
 \PsiH= \psiNzH (x)^2 - \psiNpH (x)\psiNmH (x) 
.\end{align}
\begin{lemma}\label{l:X4} 
Let $(x, N, \theta) \in \Ome$. Then, 
\begin{align} \label{:X4y} & 
 \lvert \PsiN - \PsiH \rvert =\5 \text{\quad on $\Ome$}
.\end{align}
\end{lemma}
\begin{proof} 
From \eqref{:X1z} and \eqref{:X4x}, we deduce that, on $\Ome $, 
\begin{align}\label{:X4a}&
 \Psi _{N }- \widehat{\Psi }_{N} = 
 \psi_N ^2 - \psiNzH ^2 - ( \psi_{N+1}^{N} \psi_{N-1}^{N} - \psiNpH \psiNmH ) 
.\end{align}
Using \eqref{:W1q}, \eqref{:W1g}, and \eqref{:X3c}, we obtain, on $\Ome$, 
\begin{align}\notag &%\label{:X2d} & 
\psi_{N+1}^{N} \psi_{N-1}^{N} - \psiNpH \psiNmH 
\\ \notag &= \frac{1}{\pi^2}
\prod_{\ell = \pm 1 } \big( \sumL 
\fell ^{-1/4} \Bell ^{\KM } + \5 \big) 
%\\ \notag &
 - \frac{1}{\pi^2}
\prod_{\ell = \pm 1 } \big( \sumL 
\mathsf{f} ^{-1/4} \widehat{\mathsf{B}}_\ell ^{\KM } \big)
\\ \notag &
=\frac{1}{\pi^2} 
\sumL \sumLL
 \big( \Alll \Bellp ^{\KM } \Bellm ^{\LN } - 
 \fz ^{-1/2} \wB _1^{\KM } \wB _{-1}^{\LN } \big) + \OL (N^{-2/3} \vertx ^{-5/4})
.\end{align}
Combining this with \eqref{:W1d}, \eqref{:X4a}, and \eqref{:X3b}, we obtain \eqref{:X4y}. 
\end{proof}

\begin{lemma}\label{l:X5} 
Make the same assumptions as in \lref{l:X4}. Then, 
\begin{align} \label{:X5a}
&\PsiH 
={N^{-2/3}}\rhohat ^{N } (x) + \OL ({N^{-1}\vertx ^{-1/2} }) 
\text{\quad on $\Ome$}
.\end{align}
\end{lemma}
\begin{proof} 
From \eqref{:W1p} and \eqref{:W1q}, we have 
\begin{align}\label{:X5f}
\psiNellH ( x) &= \frac{ 1}{\pi} 
 \Azz \sumL \Ckm 
\sin \big( \gz + \ckm + \ell\theta \big) 
\quad \lpm 
.\end{align}
We put 
$ A=\mathsf{g} (N,\theta) + \ckm (\theta)$ and 
$ A'=\mathsf{g} (N,\theta) + \cln (\theta)$ 
for $0\le m \le k \le L-1$ and $0\le \n \le \llll \le L-1$. 
Then, using the identity
\begin{align}\notag & %\label{:X5b}&
2\sin A \sin A' - \sin (A+\theta) \sin (A'-\theta) - \sin (A-\theta) \sin (A'+\theta) 
= 2\sin^2 \theta \cos (A-A')
,\end{align}
we have, from \eqref{:91f}, \eqref{:U1d}, \eqref{:X4x}, and \eqref{:X5f}, 
\begin{align} \nonumber 
\PsiH 
&= \Big( \frac{1}{\pi}+ \OL \Big(\frac{1}{\vertx ^{3/2}}\Big) \Big) 
\frac{\sin^2 \theta}{\mathsf{f} (N,\theta)^{1/2}}
= \Big( \frac{1}{\pi}+ \OL\Big(\frac{1}{\vertx ^{3/2}}\Big) \Big) \frac{\sin \theta}{N^{1/3}} 
\text{\quad on $\Ome$}
.\end{align}
Because of $x= \h (N,\theta)$, we have from Taylor's theorem and \eqref{:12z} that 
\begin{align}\nonumber 
 \sin \theta &=\sqrt{1-\cos^2 \theta} 
=\sqrt{ 1-\frac{1}{N+1}\Big(\frac{x}{2N^{1/6}} + \sqrt{N} \Big)^2 }
\\ \nonumber
&=\Big(\frac{N}{N+1}\Big)^{1/2}\sqrt{ -\frac{x}{N^{2/3}}\Big(1+ \frac{x}{4N^{2/3}} \Big)+\frac{1}{N} }
\\ \nonumber
&=\sqrt{-\frac{x}{N^{2/3}} \Big(1+\frac{x}{4N^{2/3}}\Big)} +\OL \Big(\frac{1}{N^{2/3}\vertx ^{1/2}}\Big)
 \\ \notag &
= \frac{\pi}{N^{1/3}}\rhohat ^{N } (x) + \6 
\text{\quad on $\Ome$}
.\end{align}
Combining the above two equations, we have \eqref{:X5a}.
\end{proof}
\begin{lemma}\label{l:X6}
Let $\beta=2$. Then, \eqref{:13b} holds. 
\end{lemma}

\begin{proof} 
We divide the proof into three cases: \thetag{I} $x\in [-2N ^{2/3},-N ^{\varepsilon}]$, 
\thetag{II} $x\in (-N ^{\varepsilon},-1)$, and \thetag{III} $x\in [-1,\infty)$. 
For (I), we obtain \eqref{:13b} from \lref{l:X1}, \lref{l:X4}, and \lref{l:X5}. 
Case (III) is readily derived from\lref{l:W2}. 
Hence, it only remains to prove \thetag{II}. Clearly, 
\begin{align}\label{:X6p}&
 \lvert \rho_2^{N,1}(x) - \rhohat^N (x) \rvert \le 
 \lvert \rho_2^{N,1}(x) - \rho_2^1 (x) \rvert + 
 \lvert \rho_2^1(x) - \rhohat (x) \rvert + 
 \lvert \rhohat (x) - \rhohat^N (x) \rvert 
.\end{align}
From \eqref{:34a} and \eqref{:33a}, we deduce that 
\begin{align}\label{:X6q}& 
 \lvert \rho_2^1(x) - \rhohat (x) \rvert + \lvert \rhohat (x) - \rhohat^N (x) \rvert 
= \2 
\text{\quad on $\OmAe $}
.\end{align}
Hence, we obtain \eqref{:13b} for case (II) from 
\eqref{:X6p}, \eqref{:X6q}, and 
\begin{align}\label{:X6r}&
 \lvert \rho_2^{N,1}(x) - \rho_2^1 (x) \rvert = \2 
\text{\quad on $\OmAe $}
.\end{align}
Thus, it only remains to prove \eqref{:X6r}. 
In the rest of the proof, $ \OL (\cdot )$ denotes \lq\lq $ \OL (\cdot )$ on $\OmAe $''.

Using Taylor's theorem with \eqref{:V4x}--\eqref{:V4b}, we obtain 
\begin{align} \label{:X6s}&
\gaN ^{1/2}=1+ \OL \Big(\frac{ x }{N^{2/3}}\Big) , \quad 
\gaN ^{1/2} - \gaga =
 \OL\Big( \frac{ x ^2 }{N^{4/3}} \Big) 
 \text{\quad on $\OmAe $}
.\end{align}

Let $z_\ell$ be as in \eqref{:V1x}. 
Let $ \aN $ and $ \aNl $ be as in \eqref{:V3a}. 
Here, we regard $ \aN $ and $ \aNl $ as functions of $ z $. 
Thus, from \eqref{:V1w} and \eqref{:V1x}, 
we can regard $ \aN $ and $ \aNl $ as functions of $ x $. 
We put $ \delta _{\ell }= \aNl - \aN $. 
We then find $ \deltal = - N^{-1/3}\ell + \OL (N^{-1+\varepsilon })$ from \eqref{:V3c}.
Applying Taylor's theorem to $ \Ai (\aNl ) $ at $ \aN $, we see that 
\begin{align}\notag &%\label{:X6t}&
\Ai (\aNl )=
\Ai (\aN ) + \Ai ' (\aN ) \deltal +
\frac{1}{2} \Ai '' (\aN ) \deltal ^2 + 
 \frac{1}{6}
 \Ai ''' (\aN + \ell \theta_{\ell } \deltal )
 \deltal ^3 
\end{align}
for some constants $ \thetal $ such that $ 0 \le \thetal \le 1 $. Hence, we obtain 
\begin{align} \notag 
 \Ai(\aN )^2 - \0 
% \\ \notag 
= &{N^{-2/3}}( \Ai' (\aN )^2 - \Ai(\aN )\Ai''(\aN ) ) + \5 
\\ \label{:X6u}
= &{N^{-2/3}}\rho_2^1 (x) + \5 
.\end{align}
Here, we used \eqref{:83a}, \lref{l:81}, and \eqref{:V2b} for the last line. 
From \lref{l:81}, \eqref{:V1a}, and \eqref{:V3b}, we obtain 
\begin{align}\label{:X6w} 
\0 = \OL\Big( \frac{1}{\sqrt{\xvee }} \Big) 
.\end{align}

Collecting these, we deduce that 
 \begin{align}& \notag %\label{:X6x}&
\Psi_N(x) =\psi_N(x)^2 - 
\prod_{\lpm } \psi_{N+\ell }(x_{\ell }) + 
\4 
\quad \text{by \eqref{:X1z}, \eqref{:W1x}}
\\ \nonumber
&= \gaN ^{1/2} \Ai(\aN )^2 
- {\prod_{\lpm }\gaNl ^{1/4}\Ai (\aNl ) }
+ \8 
\quad \text{by \eqref{:V7a}, \eqref{:V5a}}
\\ \nonumber
&= \gaN ^{1/2} \Big( \Ai(\aN )^2 - \0 \Big) 
+ \Big( \gaN ^{1/2} - \gaga \Big) \Big( \0 \Big) 
\\ \notag & \quad 
+ \8 
\\ \nonumber
&= 
( 1+ \OL (N^{-2/3} x ) ) ( {N^{-2/3}}\rho_2^1 (x) + \5 )
+ \8 
\\ \notag &
= 
{N^{-2/3}}\rho_2^1 (x) + 
 \5 
 .\end{align}
Here, we used \eqref{:X6s}, \eqref{:X6u}, and \eqref{:X6w} for the fourth line and 
\eqref{:34b} and $\varepsilon \in (0, 1/6)$ for the last line. 
Hence, using the equation above and $\varepsilon \in (0, 1/6)$ again, we deduce
\begin{align}\label{:X6y}&
 \lvert N^{2/3} \Psi_N(x) - \rho_2^1 (x) \rvert \le \2 
.\end{align}
Therefore, we obtain \eqref{:X6r} from \eqref{:X1a} and \eqref{:X6y}. 
\end{proof}
%%%%%%%%%%%%%%% Lemma 14.5 %%%%%%%%%%%%%%%% 
\begin{lemma}\label{l:X7} 
Let $\beta=1,4$. Then, \eqref{:13b} holds.
\end{lemma}
\begin{proof}
If $ \beta=1,4 $, we deduce from \eqref{:31l} and \eqref{:32z} that for $N\in\mathbb{N} $,
\begin{align*} \notag 
\rairyone ^{N ,1} (x) &= 
\rairy ^{N ,1} (x) + 
\frac{1}{2} \psi _{N -1}(x)\varepsilon \psi _{N }(x)+\frac{\psi_{N -1}^{N } (x)}
{\int_{\mathbb{R} }\psi_{N -1}^{N } (t)dt} 
{\bf 1}(\mbox{$N $ is odd})
,\\ \notag 
\rairyfour ^{N ,1}(x) &= 
\frac{1}{2^{1/3}} \Big( \rairy ^{2N ,1} (2^{2/3}x) + \frac{\psi_{2N}(2^{2/3}x)^2}{(2N)^{1/3}}
+ \frac{\sqrt{2N +1}}{2(2N) ^{1/2}} 
\psi_{2N }(2^{2/3}x)\varepsilon \psi_{2N +1}^{2N}(2^{2/3}x) \Big) 
.\end{align*}
Hence, using \lref{l:X6} and \lref{l:W5} with $dx_\ell=(\oneNell )^{1/6}dx$, we obtain \lref{l:X7}. 
Here, $ x $ and $ x_\ell $ are as in \eqref{:V1w} and \eqref{:V1x}. 
\end{proof}

\smallskip 

\noindent 
{\em Proof of \pref{l:35}}: 
Equation \eqref{:13b} is derived from \lref{l:X6} for $\beta=2$,
and \lref{l:X7} for $\beta=1, 4$, respectively.
\qed

\SECT{Subsidiary lemmas and proof of \lref{l:38}}\label{s:Y1}
In Subsection \ref{s:Y1}, we prepare four lemmas. 
Using two of them (\lref{l:Y1} and \lref{l:Y4}), we prove \lref{l:38}. 
The proof of \eqref{:Y1g} in \lref{l:Y1} is very long and thus postponed until Subsection \ref{s:19}.

%%%%%%%%%%%%%%%%%%%%%%%%%%%%%%%%%%%%%%%%%
Let $ \map{\iotaN }{\mathbb{R}}{\mathbb{R}}$ be such that %$ \iotaNy = $
\begin{align}\label{:Y1u}&
\iotaNy =-(4N ^{2/3}+ y ) 
.\end{align}
We will repeatedly use the following simple relations. 
\begin{align}\notag
&\lvert x - y \rvert = \lvert \iotaNx - \iotaNy \rvert \quad \text{for $ x , y \in \mathbb{R} $}
,\\\notag 
&\lvert x - y \rvert \le \lvert x - \iotaNy \rvert \quad \text{for $ x , y \in [-2N ^{2/3},\infty) $}
,\\
&\label{:Y1w}
\lvert x - \iotaNy \rvert \le \lvert x - y \rvert 
\quad \text{for $y \le -2N ^{2/3} \le x $}
.\end{align}
Let $ \psi_{N }^{(0)} = \psi_{N }$ and 
$ \psi_{N }^{(m)}= d^{m}\psi_{N } / dx^{m} $ for $ m \in \mathbb{N} $. 
Then, from \eqref{:31b}, 
\begin{align}
\label{:Z2y} &\quad 
\psi_{N }^{(m)}( y ) = (-1)^{N+m}\psi_{N }^{(m)}( \iotaNy ) 
\quad \text{ for } m \in \zN 
.\end{align}
Let $ \Kairy^{N }$ be as in \eqref{:31f}. 
Then we have 
\begin{align} \notag 
\Kairy^{N }(x,y) =& \Kairy^{N }(\iotaN (x), \iotaNy ) , 
\\\label{:Y1X}
\partial_2 \Kairy^N (x,y) =& - (\partial_2 \Kairy^N )(\iotaN (x), \iotaNy ) 
.\end{align}
For $a,b,c \in\mathbb{R} $, we introduce the functions 
$G_{a,b,c}$ and $\widehat{G}_{a,b,c} $ on $ \mathbb{R}^2$ defined as 
\begin{align}&\notag 
G_{a,b,c}(x,y)={(\xvee )^{-a} (\yvee )^{-b} \lvert x-y \rvert ^{-c}} 
,\\ & \label{:Y1y}
\widehat{G}_{a,b,c}(x,y)=G_{a,b,c}(x,y)+G_{a,b,c}(y,x) 
.\end{align}
Note that $ G_{a,b,c}(x,y) \le G_{a,b,c}(\vertx ,\verty ) $ for $ c\ge 0$ for all $ x , y \in \mathbb{R}$. 
Using $ \widehat{G}_{a,b,c}$, we present a set of estimates of kernel functions. 
We first consider the case $ \beta = 2 $. 
Let $ \mathscr{M}_N $ be the set of maps such that 
\begin{align}\label{:Y1Z}&
\mathscr{M}_N = \{ \mathrm{id}, \iota_N \} 
.\end{align}
In the rest of \sref{s:Y1}, we write $ \xhat \in\mathscr{M}_N $ if 
$ \xhat = u (x)$ for $ u \in \mathscr{M}_N $. 
Thus, $ \xhat \in \mathscr{M} $ means that the function $ \xhat $ is either $ \xhat = x$ or $ \xhat = \iotaNx $. 
Furthermore, 
$ \xhat , \yhat , \zhat \in \mathscr{M}_N $ means $ \xhat = u (x)$, $ \yhat = v (v)$, and $ \zhat = w (z)$ for some $ u , v , w \in \mathscr{M}_N $. 
Here, $ u , v $, and $ w $ are not necessarily the same in general. 
\begin{lemma}\label{l:Y1} 
%For $ \xhat = x $, $ \xhat = \iotaNx $, $ \yhat = y $, $ \yhat = \iotaNy $, 
For $ \xhat , \yhat \in \mathscr{M}_N $, 
the following hold. 
\\ 
{\rm (i)} There exists $\Ct \label{;Y1} >0$ independent of $N\in \mathbb{N} $ such that 
for $x,y\in [-2N ^{2/3},\infty) $ 
\begin{align} \label{:Y1a} 
\lvert \Kairy^{N }(\xhat ,\yhat ) \rvert 
&\le \cref{;Y1} \widehat{G}_{-\frac{1}{4},\frac{1}{4},1}(x,y) 
, \\\label{:Y1b} 
 \lvert 
 \big( \partialtwo \Kairy^{N } \big) (\xhat ,\yhat ) 
\rvert 
& \le \cref{;Y1} \Big(
\widehat{G}_{-\frac{3}{4},\frac{1}{4},1}(x,y)
 +\widehat{G}_{-\frac{1}{4},-\frac{1}{4},1}(x,y) \Big) 
.\end{align}
\noindent {\rm (ii)}
For each $r>0$, there exists $\Ct \label{;Y2} >0$ independent of $N\in \mathbb{N} $
such that for $y\in [-2N ^{2/3},\infty)$ 
\begin{align} 
\max_{x\in [-r,r]}\lvert \Kairy^{N }(x, \yhat )\rvert &\le
 {\cref{;Y2}} \ythreeF 
,\label{:Y1c}
\\ \notag 
\max_{x\in [-r,r]} \lvert \big( \partialtwo \Kairy^{N } \big) (x, \yhat ) \rvert 
&\le {\cref{;Y2}} \yoneF 
,\\ \label{:Y1d} 
\max_{x\in [-r,r]} \lvert \big( \partialtwo \Kairy^{N } \big) (\yhat ,x) \rvert 
&\le \cref{;Y2} \yoneF 
.\end{align}

\noindent {\rm (iii)}
There exists $\Ct \label{;Y3} >0$ independent of $N\in \mathbb{N} $ such that 
for $ \xX \in [-2N^{2/3}, \infty) $
\begin{align} \label{:Y1e}%&&
\max_{\yhaT \in [\xhat -1,\xhat +1]}
\lvert \Kairy^{N }(\xhat ,\yhaT )\rvert &\le \cref{;Y3} \sqrt{\xvee } 
%&& x\in [-2N ^{2/3}, \infty)
,\\ \label{:Y1f} %&&
\max_{\yhaT \in [\xhat -1,\xhat +1]}
\lvert 
 \big( \partialtwo \Kairy^{N } \big) (\xhat ,\yhaT ) 
\rvert 
&\le \cref{;Y3}(\xvee ) 
%&& x\in [-2N ^{2/3}, \infty)
,\\ \label{:Y1F} %&&
\max_{\yhaT \in [\xhat -1,\xhat +1]}
\lvert 
 \big( \partialtwo \Kairy^{N } \big) (\yhaT , \xhat ) 
\rvert 
&\le \cref{;Y3}(\xvee ) 
,\\\label{:Y1g} %&&
\sup_{ \yY \in \mathbb{R} }
\lvert \int_{\xhat }^{\yY } K_2^N (u, \xhat )du\rvert &\le \cref{;Y3} ( 1 \vee \log\lvert \xX \rvert ) 
.\end{align}
\end{lemma}
\begin{proof} 
From \eqref{:W2z}--\eqref{:W3Y} in \lref{l:W2} and \eqref{:Z2y}, 
we see for $\xhat \in \mathscr{M}_N $
\begin{align}&\notag 
\lvert \psi_{N} (\xhat ) \rvert \le \cref{;W3n} \exp (- \xF /\cref{;W3Y} ) 
, \\ \notag &
\lvert \psi'_N(\xhat )\rvert \le \cref{;W3n} \exp (- \xF /\cref{;W3Y} ) 
,\\ \label{:W1v}&
\lvert \psi_N ''(\xhat ) \rvert \le \cref{;W3n} \exp (- \xF / \cref{;W3Y}) 
\quad \text{ for $x \in [-1,\infty )$}
.\end{align}
Using \lref{l:W3}, \eqref{:Z2y}, and Taylor's theorem at $x = -2N^{2/3}$ to the right-hand side of \eqref{:W3a}--\eqref{:W3c}, 
we have a constant $ \Ct \label{;W6}$ such that for $\xhat \in \mathscr{M}_N $ 
\begin{align} \notag & 
\lvert \psi_{N }(\xhat )\rvert \le \cref{;W6}(\xvee )^{-1/4}
,\\ \notag &% \quad 
\lvert \psi_{N }'(\xhat )\rvert \le \cref{;W6}(\xvee )^{1/4}, 
\\\label{:W1u}&
\lvert \psi_{N }''(\xhat )\rvert \le \cref{;W6} (\xvee )^{3/4} 
\quad \text{ for $ x \in [-1 -2N ^{2/3}, \infty ) $}
.\end{align}
%
%	From \eqref{:Y1k}, we have for all $ x , y \in \mathbb{R}$
Applying \eqref{:Y1k} to $ \xhat , \yhat \in \mathscr{M}_N $, we have for all $ x , y \in \mathbb{R}$ 
\begin{align}\label{:Y1K} 
\Kairy^{N }( \xhat , \yhat ) & = 
 \frac
{\psi_{N }( \xhat )\psi'_{N }( \yhat )-\psi'_{N }( \xhat )\psi_{N }( \yhat )}
{ \xhat - \yhat } 
-
\frac{1}{2N ^{1/3}}\psi_{N }( \xhat )\psi_{N }( \yhat ) 
.\end{align}

We examine the first term on the right-hand side of \eqref{:Y1K}. 
From \eqref{:Y1w}, \eqref{:Z2y}, \eqref{:Y1y}, and \eqref{:W1u}, we have for $ x, y \in [-2N ^{2/3}, \infty ) $
\begin{align} \notag 
\Big\lvert 
& \frac{\psi_{N }(\xhat )\psi'_{N }(\yhat ) 
-\psi'_{N }(\xhat )\psi_{N }(\yhat )}{\xhat -\yhat } 
\Big\rvert 
\le 
\frac{ \lvert \psi_{N }(\xhat )\psi'_{N }(\yhat ) \rvert +
 \lvert \psi'_{N }(\xhat )\psi_{N }(\yhat ) \rvert }{\lvert {\xhat -\yhat } \rvert } 
\\&\notag 
= \frac{ \lvert \psi_{N }(x) \psi'_{N }(y ) \rvert +
 \lvert \psi'_{N }(x)\psi_{N }(y ) \rvert }{\lvert {\xhat -\yhat } \rvert } 
\quad \text{ by \eqref{:Z2y}}
\\\notag & \le 
 \cref{;W6}^2
\Big( 
\frac{(\xvee )^{-1/4} (\yvee )^{1/4}}{\lvert {\xhat -\yhat } \rvert } + 
\frac{ (\xvee )^{1/4} (\yvee )^{-1/4}}{\lvert {\xhat -\yhat } \rvert }
\Big) 
\quad \text{ by \eqref{:W1u}}
\\ \label{:Y1p} 
& \le \cref{;W6}^2 
 \widehat{G}_{-\frac{1}{4},\frac{1}{4},1}(x,y) 
\quad \text{ by \eqref{:Y1w}, \eqref{:Y1y}}
%\quad \text{ for $ 2N^{2/3} \ge r+1 $}
.\end{align}
%$ \lvert {\xhat -\yhat } \rvert \ge \lvert x - y \rvert $ for $ -2N^{2/3} \le - (r + 1) \le y $. 

We next consider the second term on the right-hand side of equation \eqref{:Y1K}. 
Let $( x,y ) \in [-2N ^{2/3}, 0]^2 $. 
From \eqref{:Z2y}, \eqref{:Y1y}, and \eqref{:W1u}, we have 
\begin{align} \notag %
\Big\lvert \frac{1}{2N ^{1/3}} & \psi_{N }(\xhat )\psi_{N }(\yhat ) \Big\rvert = 
\Big\lvert \frac{1}{2N ^{1/3}}\psi_{N }(x)\psi_{N }( y ) \Big\rvert 
\quad \text{ by \eqref{:Z2y}}
\\&\notag \le \cref{;W6}^2 
\frac{1}{2N ^{1/3}} ( \xvee )^{-1/4} ( \yvee )^{-1/4} 
\quad \text{ by \eqref{:W1u}}
\\ \notag &
\le \cref{;W6}^2 \frac{1}{\lvert x - y \rvert ^{1/2}} ( \xvee )^{-1/4} ( \yvee )^{-1/4}
\quad \text{ by }(x,y ) \in [-2N ^{2/3}, 0]^2 
\\\label{:Y1q}&
\le \cref{;W6}^2 
\widehat{G}_{-\frac{1}{4},\frac{1}{4},1}(x,y)
\quad \text{ by 
$ (\xvee ) \vee (\yvee ) \ge \lvert x - y \rvert $ and \eqref{:Y1y}}
.\end{align}
Let $ (x , y ) \in [0,\infty)^2$. From \eqref{:Z2y}, \eqref{:Y1y}, and \eqref{:W1v}, 
we have $ \Ct \label{;W1s} > 0 $ such that 
\begin{align}\notag 
\Big\lvert \frac{1}{2N ^{1/3}} & \psi_{N }(\xhat )\psi_{N }(\yhat ) \Big\rvert 
=
 \Big\lvert \frac{1}{2N ^{1/3}}\psi_{N }(x)\psi_{N }( y ) \Big\rvert 
\quad \text{ by \eqref{:Z2y} }
\\\notag & \le 
\frac{1}{2N ^{1/3}} \cref{;W3n}^2 \exp (- \xF /\cref{;W3Y} ) \exp (- \yF /\cref{;W3Y} ) 
\quad \text{ by \eqref{:W1v}}
\\\label{:W1s}&
\le \cref{;W1s}
\widehat{G}_{-\frac{1}{4},\frac{1}{4},1} (x,y) 
\quad \text{ by \eqref{:Y1y}}
.\end{align}
The proof for the case $ (x,y) , (y,x)\in [-2N ^{2/3}, 0] \ts [0,\infty ) $ is similar. Hence, we omit it. Putting \eqref{:Y1p}--\eqref{:W1s} into the right-hand side of \eqref{:Y1K} yields \eqref{:Y1a}. 

Differentiating \eqref{:Y1K}, we see 
\begin{align} \notag 
(\partial_2\Kairy^{N } ) (\xhat , \yhat ) &= 
 \frac{\psi_{N }(\xhat )\psi''_{N }(\yhat )-\psi'_{N }(\xhat )\psi_{N }' (\yhat )}{\xhat -\yhat } + 
 \frac{\psi_{N }(\xhat )\psi_{N }'(\yhat )-\psi'_{N }(\xhat )\psi_{N }(\yhat )}{(\xhat -\yhat )^2} 
\\\notag &
-\frac{1}{2N ^{1/3}}\psi_{N }(\xhat )\psi_{N }'(\yhat ) 
\\\label{:W1t}&
=: I_1 (\xhat , \yhat ) + I_2 (\xhat , \yhat ) + I_3 (\xhat , \yhat ) 
.\end{align}
Similarly as \eqref{:Y1p}--\eqref{:W1s}, using \eqref{:W1v} and \eqref{:W1u}, we have for $ x,y\in [-2N ^{2/3}, \infty ) $ 
\begin{align}& \label{:W1(}
\lvert I_p (\xhat , \yhat ) \rvert \le 
\cref{;W6}^2 \Big( \widehat{G}_{-\frac{3}{4},\frac{1}{4},1}(x,y) +
 \widehat{G}_{-\frac{1}{4}, -\frac{1}{4},1}(x,y) 
\Big), \ p=1,3
.\end{align}
For $ p=2$ and $x,y\in [-2N ^{2/3}, \infty ) $ such that $ \lvert \xhat - \yhat \rvert \ge 1 $
\begin{align} \notag 
\lvert I_2 (\xhat , \yhat ) \rvert 
&=
 \Big\lvert 
\frac{\psi_{N }(\xhat )\psi_{N }'(\yhat )-\psi'_{N }(\xhat )\psi_{N }(\yhat )}
{( \xhat - \yhat )^2} 
\Big\rvert 
\\\notag &
\le \frac{ \lvert \psi_N(x)\psi''_N(y) \rvert + \lvert \psi_N'(x)\psi_N'(y) \rvert }
{ \lvert \xhat - \yhat \rvert }
\quad \text{ by \eqref{:Z2y}, $ \lvert \xhat - \yhat \rvert \ge 1 $}
\\ \label{:W1)} &\le 
\cref{;W6}^2 
\widehat{G}_{-\frac{1}{4}, \frac{1}{4},1}(x,y) 
\quad \text{ by \eqref{:Y1w}, \eqref{:Y1y}, \eqref{:W1u}}
.\end{align}
We next consider the case $ p=2$ and $x,y\in [-2N ^{2/3}, \infty ) $ such that $ \lvert \xhat - \yhat \rvert \le 1 $. 
We estimate the numerator of $ I_2 (\xhat , \yhat ) $ using Taylor's theorem. 
Then, there exists $\yhat _1 $ such that 
$ \yhat \le \yhat_1\le \xhat $ or $ \xhat \le \yhat_1 \le \yhat $, where 
$ \yhat _1 = v (y _1 )$ for $ v \in \mathscr{M}_N $ such that $ \yhat = v (y)$, and that 
\begin{align} &\psi_{N }(\xhat )\psi'_{N }(\yhat )-\psi'_{N }(\xhat )\psi_{N }(\yhat )
= (\yhat -\xhat ) ( \psi_N(\xhat )\psi''_N(\yhat _1) - \psi_N'(\xhat )\psi_N'(\yhat _1) ) 
\label{:W1w}
.\end{align}
Because $ \yhat = v (y)$ and $ \yhat_1 = v(y_1)$ for the same $ v \in \mathscr{M}_N $, 
we have 
% Here for the fifth line, we used $ \rvert y - y_1 \rvert \le 1 $, which follows from 
\begin{align}& \label{:W1Z}
 \rvert y - y_1 \rvert = \rvert v (y) - v (y_1) \rvert 
= \lvert \yhat - \yhat _1 \rvert \le \lvert \xhat - \yhat \rvert \le 1
.\end{align}
Hence, for $ x , y \in [-2N ^{2/3}, \infty ) $ such that $ \lvert \xhat - \yhat \rvert \le 1 $, we have a constant $ \Ct \label{;W1!}$ such that 
\begin{align}\notag &
\lvert I_2 (\xhat , \yhat ) \rvert =
 \Big\lvert 
\frac{\psi_{N }(\xhat )\psi_{N }'(\yhat )-\psi'_{N }(\xhat )\psi_{N }(\yhat )}
{( \xhat - \yhat )^2} 
\Big\rvert 
\\&\notag 
= \Big\lvert \frac{\psi_N(\xhat )\psi''_N(\yhat _1) - \psi_N'(\xhat )\psi_N'(\yhat _1) }{ \xhat - \yhat } \Big\rvert 
\quad \text{ by \eqref{:W1w}}
\\&\notag 
\le \frac{ \lvert \psi_N(x)\psi''_N(y_1) \rvert + \lvert \psi_N'(x)\psi_N'(y_1) \rvert }{ \lvert x- y \rvert }
\quad \text{ by \eqref{:Y1u}--\eqref{:Z2y}}
\\&\notag 
\le \cref{;W6}^2\Big(
\frac{(\xvee )^{-1/4} ( 1 \vee \lvert y_1 \rvert )^{3/4}}{\lvert { x - y } \rvert } + 
\frac{ (\xvee )^{1/4} ( 1 \vee \lvert y_1 \rvert )^{1/4}}{\lvert { x - y } \rvert }
\Big)
\quad \text{ by \eqref{:W1u}}
\\&\notag 
\le \cref{;W1!} \Big(
\frac{(\xvee )^{-1/4} (\yvee )^{3/4}}{\lvert { x - y } \rvert } + 
\frac{ (\xvee )^{1/4} (\yvee )^{1/4}}{\lvert { x - y } \rvert }
\Big) 
\quad \text{ by \eqref{:W1Z}}
\\&\label{:W1!} 
= \cref{;W1!} \Big(\widehat{G}_{-\frac{3}{4}, \frac{1}{4},1}(x,y) +
\widehat{G}_{-\frac{1}{4}, -\frac{1}{4},1}(x,y) \Big) 
\quad \text{ by \eqref{:Y1y}}
.\end{align}
Putting \eqref{:W1(}, \eqref{:W1)}, and \eqref{:W1!} into \eqref{:W1t}, we have \eqref{:Y1b}. 
Thus, we obtain \thetag{i}.

From \eqref{:Y1a}, we have for $y\in [-2N ^{2/3},\infty)$ 
\begin{align}\label{:W1'} 
\max_{x\in [-r,r] \cap [-2N^{2/3}, r]} \lvert \Kairy^{N }(x, \yhat )\rvert 
&\le 
\cref{;Y1} \max_{x\in [-r,r] } \widehat{G}_{-\frac{1}{4},\frac{1}{4},1}(x,y) 
.\end{align}
Let $ 2N^{2/3} < r $. Then $ \iotaN ([-r , -2N^{2/3}) ) = (-2N^{2/3}, r - 4N^{2/3} ] $. 
Hence using \eqref{:Y1X} and \eqref{:W1'}, we have for $y\in [-2N ^{2/3},\infty)$ 
\begin{align}\notag 
\max_{x\in [-r,r] \backslash [-2N^{2/3}, r ]} 
\lvert \Kairy^{N }(x, \yhat )\rvert 
& = 
\max_{x\in [-r,r] \backslash [-2N^{2/3}, r ]} 
\lvert \Kairy^{N }( \iotaNx , \iotaN (\yhat ) )\rvert 
\quad \text{ by \eqref{:Y1X} }
\\\notag & = 
\max_{ x \in [-2N^{2/3}, r - 4N^{2/3} ] } \lvert \Kairy^{N }( x , \iotaN (\yhat ) )\rvert 
\\&\label{:W1D}
 \le
\cref{;Y1} \max_{x\in [-r,r] } 
\widehat{G}_{-\frac{1}{4},\frac{1}{4},1}(x,y) 
\quad \text{ by \eqref{:W1'}}
.\end{align}
From \eqref{:W1'} and \eqref{:W1D}, we have for $y\in [-2N ^{2/3},\infty)$ 
\begin{align}\label{:W1E}&
\max_{x\in [-r,r] } \lvert \Kairy^{N }(x, \yhat )\rvert \le 
\cref{;Y1} \max_{x\in [-r,r] } \widehat{G}_{-\frac{1}{4},\frac{1}{4},1}(x,y) 
.\end{align} 

We divide the case into two parts: \thetag{I} $ \lvert y \rvert > r+1 $ and 
\thetag{II} $ \lvert y \rvert \le r+1 $. % and $ \lvert \yhat \rvert > r+1 $. 

Suppose \thetag{I}. Then using $ \lvert x \vert \le r $ and $ r + 1 < \lvert y \rvert $, 
we have $ 1 \le \xy $. Hence, 
\begin{align}\label{:W1"}&
 \lvert y \rvert \le \xy + \lvert x \rvert \le \xy + r \le (1+r) \xy 
.\end{align}
From \eqref{:Y1y}, \eqref{:W1E}, \eqref{:W1"}, and $ 1 < \lvert y \rvert $, 
we obtain for $ y\in [-2N ^{2/3},\infty) $ 
\begin{align} \notag 
\max_{x\in [-r,r] \cap \{ \xy \ge 1 \}} \lvert \Kairy^{N }(x, \yhat )\rvert 
&\le 
 \cref{;Y1} (1+r)(1+ r^{1/4}) \lvert y \rvert ^{-3/4}
% \quad \text{ by \eqref{:Y1y}, \eqref{:W1"}}
\\& \label{:W1A} 
= 
 \cref{;Y1} (1+r)(1+ r^{1/4}) \ythreeF 
.\end{align}
Here we used $ \lvert y \rvert > r+1 $ for the second line. 
Because $ \Kairy^{N }(x, \yhat )$ is locally bounded uniformly in $ N $ from \eqref{:W1v}--\eqref{:Y1K}, we have a constant $ \Ct \label{;W1B}$ independent of $ N $ such that for $ y\in [-2N ^{2/3},\infty) $ 
\begin{align}& \label{:W1B} 
\max_{x\in [-r,r] \cap \{ \xy \le 1 \}} \lvert \Kairy^{N }(x, \yhat )\rvert 
\le \cref{;W1B}
.\end{align}
Using \eqref{:W1A} and \eqref{:W1B}, we obtain \eqref{:Y1c} for \thetag{I}. 
Next, suppose \thetag{II}. Then, similarly as \eqref{:W1B}, we have a constant $ \Ct \label{;W1C}$ independent of $ N $ such that 
\begin{align} \notag & %\label{:W1C}&
\max_{x\in [-r,r] \cap \{\lvert y \rvert \le r+1 \}} 
\lvert \Kairy^{N }(x, \yhat )\rvert 
\le \cref{;W1C}
.\end{align}
Thus we have \eqref{:Y1c} for \thetag{II}. This yields \eqref{:Y1c}. 
We can prove \eqref{:Y1d} in the same fashion as \eqref{:Y1c}. Hence we omit it. 
% We have \eqref{:Y1d} from \eqref{:Y1y} and \eqref{:Y1b}. 
Thus, we obtain \thetag{ii}.

We proceed with the proof of \thetag{iii}. 
Using \eqref{:Y1k} and \eqref{:W1w}, we see 
\begin{align*}&
\max_{\yhaT \in [\xhat -1, \xhat +1]} \lvert \Kairy^{N }(\xhat ,\yhaT ) \rvert
\\
&\le \max_{\yhaT \in [\xhat -1, \xhat +1]} \Big(
 \Big\lvert 
\frac{\psi_{N }(\xhat )\psi_{N }'(\yhaT )-\psi'_{N }(\xhat )\psi_{N }(\yhaT )}
{\xhat -\yhaT } \Big\rvert 
+
\Big\lvert \frac{1}{2N ^{1/3}}\psi_{N }(\xhat )\psi_{N }(\yhaT ) \Big\rvert \Big)
\\ & = 
\max_{\yhaT \in [\xhat -1, \xhat +1]} \Big(
\Big\lvert \psi_N(\xhat )\psi''_N(y_1) - \psi_N'(\xhat )\psi_N'(y_1) \Big\rvert + 
\Big\lvert \frac{1}{2N ^{1/3}}\psi_{N }(\xhat )\psi_{N }(\yhaT ) \Big\rvert \Big)
\\ & \le 
\max_{\yhaT \in [\xhat -1, \xhat +1]} \Big(
\Big\lvert \psi_N(\xhat )\psi''_N(y_1) \Big\rvert + \Big\lvert \psi_N'(\xhat )\psi_N'(y_1) \Big\rvert + 
\Big\lvert \frac{1}{2N ^{1/3}}\psi_{N }(\xhat )\psi_{N }(\yhaT ) \Big\rvert \Big)
.\end{align*}
Hence we obtain \eqref{:Y1e} from this and \eqref{:W1u}.

We next prove \eqref{:Y1f}. Using the notation of \eqref{:W1t}, we have 
\begin{align}\label{:W1k}&
(\partial_2\Kairy^{N } ) (\xhat , \yhat ) = 
I_1 (\xhat , \yhat ) + I_2 (\xhat , \yhat ) + I_3 (\xhat , \yhat ) 
.\end{align}
Applying Taylor's theorem to the numerator in \eqref{:W1t}, we see that there exist 
$ \yhat _1, \yhat _2 \in \mathscr{M}_N $ such that 
$ \lvert \xhat - \yhat _i \rvert 
 \le \lvert \xhat - \yhat \rvert $, 
$ i=1,2$, and that 
\begin{align} &\notag 
\psi_{N }(\xhat )\psi''_{N }( \yhat )-\psi'_{N }(\xhat )\psi'_{N }( \yhat )
\\\notag 
= &\psi_{N }(\xhat )\psi''_{N }(\xhat ) - \psi_N'(\xhat )^2 + 
( \yhat -\xhat ) ( \psi_{N }(\xhat )\psi'''_{N }( \yhat _1) - \psi_N'(\xhat )\psi''_{N }( \yhat _1) )
\end{align}
and 
\begin{align}&\notag 
\psi_{N }(\xhat )\psi'_{N }( \yhat )-\psi'_{N }(\xhat )\psi_{N }( \yhat )
%	\\\notag = &( \yhat -\xhat ) ( \psi_N(\xhat )\psi''_N( \yhat _1) - \psi_N'(\xhat )\psi_N'( \yhat _1) ) 
\\ \notag 
=& ( \yhat -\xhat ) ( \psi_N(\xhat )\psi''_N(\xhat ) - \psi_N'(\xhat )^2 ) 
+ \frac{( \yhat -\xhat )^2}{2} ( \psi_N(\xhat )\psi_N'''( \yhat _2) -\psi'_N(\xhat )\psi_N''( \yhat _2) ) 
.\end{align}
Hence, we have a cancellation such that 
\begin{align}&\notag 
I_1 (\xhat , \yhat ) + I_2 (\xhat , \yhat )
\\ \notag
= &
 \frac{\psi_{N }(\xhat )\psi''_{N }( \yhat )-\psi'_{N }(\xhat )\psi_{N }' ( \yhat )}{\xhat - \yhat } + 
 \frac{\psi_{N }(\xhat )\psi_{N }'( \yhat )-\psi'_{N }(\xhat )\psi_{N }( \yhat )}{(\xhat - \yhat )^2} 
%		\\\notag &
\\\label{:W1l} 
=& - 
( \psi_{N }(\xhat )\psi'''_{N }( \yhat _1) - \psi_N'(\xhat )\psi''_{N }( \yhat _1) ) 
+ 
\frac{1}{2} 
%		\frac{\xhat - \yhat _1 }{\xhat - \yhat } 
( \psi_N(\xhat )\psi_N'''( \yhat _2) -\psi'_N(\xhat )\psi_N''( \yhat _2) )
.\end{align}
From \eqref{:W3x}, we have 
\begin{align}\label{:W1m}&
\psi'''_{N }(z) = \Big( \frac{z}{2N ^{2/3}}+1 \Big) \psi_N (z) + 
 \Big( \frac{z^2}{4N ^{2/3}}+z - \frac{1}{2N ^{1/3}}\Big) \psi_N' (z) 
.\end{align}
We have $ \lvert x - \yy _i \rvert \le 1 $ from 
$ \lvert x - \yy _i \rvert \le \lvert \xhat - \yhat _i \rvert \le \lvert \xhat - \yhat \rvert \le 1 $ for $ i =1,2$. 
Using this, \eqref{:Z2y}, \eqref{:W1u}, and \eqref{:W1m}, we have constants $ \Ct \label{;W1n}$ and $ \Ct \label{;W1N}$ independent of $ N $ such that 
\begin{align}\notag &
\sum_{i=1,2}
\Big\lvert \psi_{N }(\xhat )\psi'''_{N }( \yhat _i ) \Big\rvert + 
\Big\lvert \psi_{N }'(\xhat )\psi''_{N }( \yhat _i ) \Big\rvert 
\\\notag = & 
\sum_{i=1,2}
\Big\lvert \psi_{N }(x )\psi'''_{N }(\yy _i ) \Big\rvert + 
\Big\lvert \psi_{N }'(x )\psi''_{N }(\yy _i ) \Big\rvert 
\quad \text{ by \eqref{:Z2y}}
\\\le & \notag 
\cref{;W1n} \sum_{i=1,2}
\Big( ( \xvee )^{-1/4} ( 1 \vee \lvert \yy _i \rvert )^{5/4} + 
 ( \xvee )^{1/4} ( 1 \vee \lvert \yy _i \rvert )^{3/4} \Big) 
 \quad \text{ by \eqref{:W1u}, \eqref{:W1m}}
\\\label{:W1n} \le &
 \cref{;W1N} (\xvee ) 
\quad \text{ by $ \lvert x - \yy _i \rvert \le 1 $, $ i = 1 , 2$}
.\end{align}
% Here we used \eqref{:W1u} and \eqref{:W1m} for the third line. 
Hence from \eqref{:W1l} and \eqref{:W1n}, 
\begin{align}\label{:W1N}& 
I_1 (\xhat , \yhat ) + I_2 (\xhat , \yhat )
\le 
 \cref{;W1N} (\xvee ) 
.\end{align}
We have $ \lvert x- \yy \rvert \le 1 $ from $ \lvert \xhat - \yhat \rvert \le 1 $. 
From \eqref{:W1u} and $ \lvert x- \yy \rvert \le 1 $, 
we have 
\begin{align}\label{:W1o}& 
 I_3 (\xhat , \yhat ) = 
\Big\lvert \frac{1}{2N ^{1/3}}\psi_{N }(\xhat )\psi_{N }'( \yhat ) \Big\rvert 
 \le 
 \frac{1}{2}\cref{;W6}^2 
 ( \xvee )^{-1/4} (\yvee )^{1/4} \le \cref{;W6}^2
\end{align}
for $ - 2N ^{2/3} \le x < \infty $. 
Putting \eqref{:W1N} and \eqref{:W1o} into \eqref{:W1k}, we obtain \eqref{:Y1f}. 

The proof of \eqref{:Y1F} is similar to that of \eqref{:Y1f}. Hence, we omit it. 

We leave the proof of \eqref{:Y1g} until Subsection \ref{s:19}. 
\end{proof}

We next consider the case $ \beta = 1,4$. As we see in \eqref{:31l}, 
\begin{align}\notag %\label{:31k}&
\Kairyone^N (x,y)&=
\begin{bmatrix} 
J_1^{N }(x,y) & -\partialtwo J_1^{N }(x,y)
\\
 \int_y^x J_1^{N }(\zz ,y) d\zz - \frac{1}{2}{\rm sign}(x-y) & J_1^{N }(y,x)
\end{bmatrix}
,\\ \notag 
\frac{1}{2^{\frac{2}{3}}}
\Kairyfour^{N }\Big( \frac{x}{2^{\frac{2}{3}}}, \frac{y}{2^{\frac{2}{3}}}\Big)
&= 
\frac{1}{2}
\begin{bmatrix} 
J_4^{N }(x,y) & -\partialtwo J_4^{N }(x,y)
\\
\int_y^x J_4^{N }(\zz ,y) d\zz & J_4^{N }(y,x)
\end{bmatrix}
.\end{align}
Here $ J_{\beta }^{N}$ is as in \eqref{:31j}. 
From \eqref{:31j}, we have for all $ x , \yhatyhat \in \mathbb{R} $
\begin{align}\notag %\label{:31i}
J_1^{N }(x,\yhatyhat ) -\Kairy^{N }(x,\yhatyhat )&=
\frac{1}{2}\psi_{N -1}^{N } (x)\varepsilon \psi_{N }(\yhatyhat )
+\frac{\psi_{N -1}^{N } (x)}
{\int_{\mathbb{R} }\psi_{N -1}^{N } (\zz )d\zz } 
{\bf 1}(\mbox{$N $ is odd})
,\\\label{:Y3e}
J_4^{N }(x,\yhatyhat ) - K_2^{2N} (x,\yhatyhat ) 
& = 
\frac{\psi_{2N}(x) \psi_{2N}(\yhatyhat )}{(2N)^{1/3} } + 
\frac{\sqrt{2N +1}}{2(2N )^{1/2}}
 \psi_{2N }(x)\varepsilon \psi_{2N +1}^{2N }(\yhatyhat ) 
.\end{align}
% By \eqref{:31g}, 
% $ \varepsilon \psi_{N }(\yhat )=
% \frac{1}{2}\int_\mathbb{R} \psi_{N } (z)dz -\int_{\yhat }^\infty \psi_{N } (z)dz $. 
Recall that 
 $ (\varepsilon f)(x) =\frac{1}{2}\int_\mathbb{R} f(y)dy -\int_x^\infty f(y)dy $ by \eqref{:31g}. 
Hence $ (\varepsilon f )' = f $. 
% \begin{align}\label{:Y2X}&
% (\varepsilon \psi_{N })' = \psi_{N} ,\quad 
% (\varepsilon \psi_{2N +1}^{2N })' = \psi_{2N +1}^{2N }
% .\end{align}
Differentiating \eqref{:Y3e} in $ y $ and using $ (\varepsilon f )' = f $, 
we have for all $ x , \yhatyhat \in \mathbb{R} $ 
\begin{align} \notag %\label{:Y3j}
\partial_2(J_1^{N } -\Kairy^{N }) (x,\yhatyhat ) &=
\frac{1}{2}\psi_{N -1}^{N } (x) \psi_{N }(\yhatyhat )
,\\\label{:Y3k}
\partial_2(J_4^{N } - K_2^{2N}) (x,\yhatyhat ) & 
= \frac{\psi_{2N}(x) \psi_{2N}' (\yhatyhat )}{(2N)^{1/3} } + 
\frac{\sqrt{2N +1}}{2(2N )^{1/2}} 
\psi_{2N }(x) \psi_{2N +1}^{2N }(\yhatyhat ) 
.\end{align}
In \lref{l:Y1}, we already have the necessary estimate for $ \Kairy ^N $. 
In the following lemma, we establish the required estimate for $ J_{1}^N $ and $ J_4^N $. 
\begin{lemma} \label{l:Y3} 
Let $\beta=1,4$. 
Let $ \mathscr{M}_N $ be as in \eqref{:Y1Z}. 
Let $ \xhat , \yhat \in \mathscr{M}_N $. 

\noindent 
\thetag{i} There exists $\Ct \label{;Y5}>0 $ independent of $N \in \mathbb{N} $, 
such that for $x,y\in [-2N ^{2/3},\infty)$ 
\begin{align} \label{:Y2a}
 \Big\lvert J_{\beta }^{N }(\xhat ,\yhat ) \Big\rvert 
\le & \cref{;Y5} \Big(
\widehat{G}_{-\frac{1}{4},\frac{1}{4},1}(x,y) + ( \xvee )^{-\frac{1}{4}} 
\Big)
,\\ \label{:Y2b}
 \Big\lvert \big( \partialtwo J^{N }_{\beta } \big) (\xhat ,\yhat ) 
 \Big\rvert 
\le & \cref{;Y5} 
\Big(
 \widehat{G}_{-\frac{3}{4},\frac{1}{4},1}(x,y)
 +\widehat{G}_{-\frac{1}{4},-\frac{1}{4},1}(x,y)
 +\widehat{G}_{\frac{1}{4},\frac{1}{4},0}(x,y) 
\Big) 
.\end{align}
\thetag{ii} 
For $r\in \mathbb{N}$, we have $\Ct \label{;Y]} >0$ independent of $N\in \mathbb{N} $
such that for $y\in [-2N ^{2/3},\infty)$ 
\begin{align} \notag 
\max_{x\in [-r,r]} 
\lvert J_{\beta }^{N } (x,\yhat ) \rvert\le &
{\cref{;Y]}} 
,\\\label{:Y2c}
\max_{x\in [-r,r]} \lvert J_{\beta }^{N } (\yhat ,x) \rvert \le &
\cref{;Y]}
\yoneF 
,\\ \notag 
\max_{x\in [-r,r]} \lvert \big( \partialtwo J_{\beta }^{N } \big) (x,\yhat ) \rvert \le &
{\cref{;Y]}} \yoneF , 
\\ \label{:Y2d}
\max_{x\in [-r,r]} \lvert \big( \partialtwo J_{\beta }^{N } \big) (\yhat ,x) \rvert 
\le &
 \cref{;Y]} \yoneF %\verty ^{-1/4}
.\end{align}

\noindent {\rm (iii)}
\; There exists $\Ct \label{;Y[} >0$ independent of $N\in \mathbb{N} $ such that 
for $ x\in [-2N ^{2/3}, \infty)$
\begin{align} \label{:Y2e}
\max_{\yhaT \in [\xhat -1,\xhat +1]}
\lvert J_{\beta }^{N }(\xhat ,\yhaT )\rvert &\le \cref{;Y[} \sqrt{\xvee } 
,\\\label{:Y2E}
\max_{\yhaT \in [\xhat -1,\xhat +1]}
\lvert J_{\beta }^{N }(\yhaT ,\xhat )\rvert & \le \cref{;Y[} \sqrt{\xvee } 
,\\ \label{:Y2f}
\max_{\yhaT \in [\xhat -1,\xhat +1]}
\lvert \big( \partialtwo J_{\beta }^{N } \big) (\xhat ,\yhaT ) \rvert 
& \le \cref{;Y[}(\xvee ) 
,\\ \label{:Y2F}
\max_{\yhaT \in [\xhat -1,\xhat +1]}
\lvert \big( \partialtwo J_{\beta }^{N } \big) (\yhaT ,\xhat ) \rvert 
& \le \cref{;Y[}(\xvee ) 
,\\ \label{:Y2g}
\sup_{\yY \in\mathbb{R}} \lvert \int_{\xhat }^{\yY } J_{\beta }^{N } (u, \xhat )du\rvert 
 &\le \cref{;Y[} ( 1 \vee \log\lvert x \rvert ) 
.\end{align}
\end{lemma}

\begin{proof} 
From \eqref{:W1x} and the second line of \eqref{:V1B}, we have, for $N\ge 2 $ and $ \lpm $, 
\begin{align}\label{:Y3f}&
\psi_{N+\ell}^{N}(x) = \Big(\frac{N+\ell}{N}\Big)^{-1/12}\psi_{N+\ell}\Big(
\Big(\frac{N+\ell}{N}\Big)^{1/6}x - 
\frac{2\ell (N +\ell)^{1/6}}{\sqrt{N} + \sqrt{N+\ell}}
\Big) 
.\end{align}
Recall that $ \iotaNx = -(4N ^{2/3}+ x )$ by \eqref{:Y1u}. 
Hence from \eqref{:Z2y} and \eqref{:Y3f}, we have 
\begin{align} & \notag 
\psi_{N+\ell}^{N}(\iotaNx )
%\\&\notag 
= 
\Big(\frac{N+\ell}{N}\Big)^{-1/12}\psi_{N+\ell}\Big(
\Big(\frac{N+\ell}{N}\Big)^{1/6} \iotaNx - 
\frac{2\ell (N +\ell)^{1/6}}{\sqrt{N} + \sqrt{N+\ell}}
\Big) 
\\&\notag = (-1)^{N+\ell}
\Big(\frac{N+\ell}{N}\Big)^{-1/12}
\\&\notag \quad \ts \psi_{N+\ell}
\Big(
 - \Big( \frac{N+\ell}{N}\Big)^{1/6} \iotaNx + 
\frac{2\ell (N +\ell)^{1/6}}{\sqrt{N} + \sqrt{N+\ell}} 
- 4(N +\ell)^{2/3}
\Big) 
\quad \text{ by \eqref{:Z2y} }
\\ \label{:Y3F} & = 
(-1)^{N+\ell}
\Big(\frac{N+\ell}{N}\Big)^{-1/12} \psi_{N+\ell}
\Big( \Big(\frac{N+\ell}{N}\Big)^{1/6} x + \cref{;Y3F} \Big) 
\quad \text{ by \eqref{:Y1u}}
.\end{align}
Here $ \Ct \label{;Y3F} = \cref{;Y3F} (N)$ is the constant defined by 
\begin{align} \notag 
\cref{;Y3F} = &
4 N ^{2/3} 
\Big( \frac{N+\ell}{N}\Big)^{1/6} 
 + 
\frac{2\ell (N +\ell)^{1/6}}{\sqrt{N} + \sqrt{N+\ell}} 
- 4(N +\ell)^{2/3}
\\ \label{:Y3ff}
= & \mathcal{O}(N^{-1/3}) \quad \text{ as } N \to \infty 
.\end{align}
From \lref{l:W3} and \eqref{:Y3f}--\eqref{:Y3ff} and using the fact that $ \psi_{N+\ell}^{N} $ is locally bounded uniformly in $ N $, 
we have a constant $ \Ct \label{;Y3f}$ independent of $N\ge 2 $ and $ \lpm $ such that 
\begin{align}\label{:Y3G}&
\Big\lvert \psi_{N+\ell}^{N}(\xhat ) \Big\rvert \le \cref{;Y3f} (\xvee )^{-1/4} 
,\ x \in [-2N ^{2/3},\infty) 
.\end{align}
Using \lref{l:W5} and \eqref{:31c}, we have 
\begin{align} \notag 
 \sup_{N \in\mathbb{N} }\sup_{\yhatyhat \in\mathbb{R} } 
\lvert \varepsilon\psi_{2N +1}^{2N } (\yhatyhat ) \rvert &< \infty 
\quad \text{ by \lref{l:W5}}
,\\ \label{:Y3g}
\SupN 
\frac{\mathbf{1}(\mbox{$N $ is odd})}{\int_{\mathbb{R} }\psi_{N -1}^{N } (\zz )d\zz } 
&< \infty 
\quad \text{ by \eqref{:31c}}
.\end{align}
From \eqref{:W1u}, \eqref{:Y3G}, and \eqref{:Y3g}, we have $ \Ct \label{;Y3i} > 0 $ such that for $ x \in [-2N ^{2/3},\infty) $
%independent of $N $ such that 
\begin{align}\notag 
\SupN \sup_{\yhatyhat \in \mathbb{R}}
\Big\lvert \frac{1}{2}\psi_{N -1}^{N } (\xhat )\varepsilon \psi_{N }(\yhatyhat )
+\frac{\psi_{N -1}^{N } (\xhat )}
{\int_{\mathbb{R} }\psi_{N -1}^{N } (\zz )d\zz } 
{\bf 1}(\mbox{$N $ is odd}) \Big\rvert
 &\le \cref{;Y3i} (\xvee )^{-1/4}, 
\\ \notag 
\supN 
\sup_{\yhatyhat \in \mathbb{R}}
\Big\lvert \frac{\psi_{2N}(\xhat ) \psi_{2N}(\yhatyhat )}{(2N)^{1/3} } + 
\frac{\sqrt{2N +1}}{2(2N )^{1/2}}
 \psi_{2N }(\xhat )\varepsilon \psi_{2N +1}^{2N }(\yhatyhat ) \Big\rvert &\le 
 \cref{;Y3i} (\xvee )^{-1/4}
% \\&\label{:Y3i}
% \quad \quad \quad \quad \quad \quad \quad \quad \quad \quad \quad \quad \quad \quad \quad \quad \quad \quad \quad \quad 
% \quad \text{for $ x \in [-2N ^{2/3},\infty) $} 
.\end{align}
From this, \lref{l:Y1}, and \eqref{:Y3e}, we obtain \eqref{:Y2a}, \eqref{:Y2c}, \eqref{:Y2e}, and \eqref{:Y2E}. 

Similarly, % Using \eqref{:Y3G} and \eqref{:Y3g}, 
we find $ \Ct \label{;Y3l} > 0 $ such that, for $ x \in [-2N ^{2/3},\infty) $, 
\begin{align} \notag 
\SupN 
\Big\lvert \frac{1}{2}\psi_{N -1}^{N } (\xhat ) \psi_{N }(\yhat ) \Big \rvert 
&\le \cref{;Y3l} (\xvee )^{-1/4} (\yvee )^{-1/4} 
,\\ \notag 
\supN 
\Big\lvert 
\frac{\psi_{2N}(\xhat ) \psi_{2N}' (\yhat )}{(2N)^{1/3} } + 
\frac{\sqrt{2N +1}}{2(2N )^{1/2}} \psi_{2N }(\xhat ) \psi_{2N +1}^{2N }(\yhat ) 
\Big \rvert 
& \le \cref{;Y3l} (\xvee )^{-1/4} (\yvee )^{-1/4}
.\end{align}
 Combining this with \lref{l:Y1} and \eqref{:Y3k}, we have 
\eqref{:Y2b}, \eqref{:Y2d}, \eqref{:Y2f}, and \eqref{:Y2F}. 

From \eqref{:Y1g}, \eqref{:Y3e}, \eqref{:Y3f}, and \lref{l:W5}, we obtain \eqref{:Y2g}. 
\end{proof} 

Let $ L_{\beta }^{N }$, $ \beta = 1,2,4$, be as in \eqref{:31n}: 
\begin{align}\label{:Y4f}&
L_{\beta }^{N } ( x, y )=
\Kairybeta ^{N }( y , x )\Kairybeta ^{N }( x, y )
.\end{align}
Let $ G_{a,b,c}$ and $ \widehat{G}_{a, b, c}$ be as in \eqref{:Y1y}. 
We easily see that 
\begin{align} \notag &% \label{:Y1z}&
G_{a,b,c}(x,y)G_{a',b',c'}(x,y)=G_{a+a',b+b',c+c'}(x,y) 
,\\&\label{:Y4u}
\widehat{G}_{a, b, c} ( x,y ) \le \widehat{G}_{a',b',c} ( x,y )\text{\ if } (a, b, c) \ge (a', b', c)
.\end{align}
\begin{lemma} \label{l:Y4} 
Let $ \beta=1, 2, 4$. 
Let $ \mathscr{M}_N $ be as in \eqref{:Y1Z}. 
Let $ \xhat , \yhat \in \mathscr{M}_N $. 

\noindent {\rm (i)} There exists $\Ct \label{;Y8}>0$ independent of $N\in \mathbb{N}$ such that 
for $x, y\in [-2N ^{2/3},\infty)$ 
\begin{align} \notag & 
\lvert L_{\beta }^{N }(\xhat ,\yhat )\rvert 
\le \cref{;Y8}\Big(
\widehat{G}_{-\frac{1}{2},\frac{1}{2},2}(x, y)
+\widehat{G}_{0,0,2}(x,y)
\\ & \label{:Y4a} 
 \quad \quad 
+\Big(\widehat{G}_{-\frac{3}{4},\frac{1}{4},1}(x, y)
+ \widehat{G}_{-\frac{1}{4},-\frac{1}{4},1}(x, y)
+\widehat{G}_{\frac{1}{4},\frac{1}{4},0}(x,y)\Big)
( 1 \vee \log \vertx )
\Big)
.\end{align}

\noindent {\rm (ii)}
For each $r>0$, there exists $\Ct \label{;Y10} >0$ independent of $N \in \mathbb{N} $ such that for $ y \in [-2N ^{2/3}, \infty) $ 
\begin{align} & \label{:Y4b}
\max_{x\in [-r,r]} \lvert L_\beta^{N }(x, \yhat )\rvert \le 
{\cref{;Y10}} \yoneF 
.\end{align}

\noindent {\rm (iii)} 
There exists $\Ct \label{;Y9}>0 $ independent of $N\in \mathbb{N} $ such that for $ x\in [-2N ^{2/3}, \infty )$ 
\begin{align} & \label{:Y4c}
\max_{\yhaT \in [\xhat -1, \xhat +1]} \lvert L_\beta^{N }(\xhat ,\yhaT )\rvert 
\le \cref{;Y9}(\xvee ) ( 1 \vee \log \vertx )
.\end{align}
\end{lemma}

\begin{proof} 
Let $\beta=2$. Then we obtain \lref{l:Y4} from \lref{l:Y1}, \eqref{:Y4f}, and \eqref{:Y4u}.

Let $\beta=1, 4$. Then $ \Kairybeta ^{N }$, $ \beta =1,4$, are given by \eqref{:31l}: 
\begin{align}\notag 
\Kairyone^N (\xhat ,\yhat )&=
\begin{bmatrix} 
J_1^{N }(\xhat ,\yhat ) & -\partialtwo J_1^{N }(\xhat ,\yhat )
\\
 \int_{\yhat }^{\xhat } J_1^{N }(\zz ,\yhat ) d\zz - \frac{1}{2}{\rm sign}(\xhat -\yhat ) & J_1^{N }(\yhat ,\xhat )
\end{bmatrix}
,\\ \label{:Y4g} 
\frac{1}{2^{\frac{2}{3}}}
\Kairyfour^{N }\Big( \frac{\xhat }{2^{\frac{2}{3}}}, \frac{\yhat }{2^{\frac{2}{3}}}\Big)
&= 
\frac{1}{2}
\begin{bmatrix} 
J_4^{N }(\xhat ,\yhat ) & -\partialtwo J_4^{N }(\xhat ,\yhat )
\\
\int_{\yhat }^{\xhat } J_4^{N }(\zz ,\yhat ) d\zz & J_4^{N }(\yhat ,\xhat )
\end{bmatrix}
.\end{align}
By \lref{l:31}, $ L_{\beta }^{N } (\xhat ,\yhat ) $ can be regarded as a real number, denoted by the same symbol. Then we have from \eqref{:Y4f} and \eqref{:Y4g} 
\begin{align} \notag 
 L_{1}^{N }(\xhat ,\yhat ) &= 
J_1^{N }(\xhat ,\yhat ) J_1^{N }(\yhat ,\xhat ) - 
\partialtwo J_1^{N }(\xhat ,\yhat ) 
\Big( \int_{\xhat }^{\yhat } J_1^{N }(\zz ,\xhat ) d\zz - \frac{1}{2}{\rm sign}(\yhat -\xhat ) \Big)
\\\label{:Y4j}
\frac{1}{2^{ \frac{4}{3}}}
 L_{4}^{N } \Big( \frac{\xhat }{2^{\frac{2}{3}}}, \frac{\yhat }{2^{\frac{2}{3}}}\Big)
&= 
\frac{1}{4} \Big(J_4^{N }(\xhat ,\yhat ) J_4^{N }(\yhat ,\xhat ) - 
\partialtwo J_4^{N }(\xhat ,\yhat ) \int_{\xhat }^{\yhat } J_4^{N }(\zz ,\xhat ) d\zz \Big)
.\end{align}
Hence, we obtain \thetag{i} from \lref{l:Y3}, \eqref{:Y4u}, and \eqref{:Y4j}. 

We have (ii) and (iii) from \lref{l:Y3} and \eqref{:Y4j}. 
\end{proof} 

\noindent 
{\em\it Proof of \lref{l:38}.} 
From \eqref{:31n}, \eqref{:37b}, and \eqref{:32z}, we deduce 
\begin{align*}&
 \lvert \rairybetax^{N ,1}(y)-\rairybeta^{N ,1}(y)\rvert =
{ \lvert L_{\beta }^{N } (x,y) \rvert }/ {\rairybeta^{N ,1}(x)}
.\end{align*}
Using this and \eqref{:3Si}, we have \eqref{:38a} from \eqref{:Y1c} and 
\lref{l:Y4} \thetag{ii}. 
\qed

Let $ \MbN $ be as in \eqref{:31o}: 
\begin{align*}&
\MbN (x,y,z)=
\Kairybeta ^{N }(y,z)
\Kairybeta ^{N }(z,x)
\Kairybeta ^{N }(x,y)
+
\Kairybeta ^{N }(y,x)
\Kairybeta ^{N }(x,z)
\Kairybeta ^{N }(z,y)
.\end{align*}
\begin{lemma} \label{l:Y4_M} 
Let $ \beta=1, 2, 4$. 
Let $ \mathscr{M}_N $ be as in \eqref{:Y1Z}. 
Let $ \yhat , \zhat \in \mathscr{M}_N $. 

\noindent {\rm (i)} 
\; There exists $\Ct \label{;Z10}>0$ independent of $N\in \mathbb{N}$ such that 
for $y, z\in [-2N ^{2/3},\infty)$ 
\begin{align} \notag 
& \Supr \lvert {M}_{\beta}^{N }(x,\yhat ,\zhat )\rvert \le 
 \cref{;Z10}\Big( 
\widehat{G}_{-\frac{3}{4},\frac{1}{4},1}(y,z) +
 \widehat{G}_{-\frac{1}{4},-\frac{1}{4},1}(y,z) 
\\& \label{:Z2d} 
\quad \quad \quad \quad \quad \quad \quad \quad \quad 
+ 
 \widehat{G}_{\frac{1}{4},\frac{1}{4},0}(y,z) 
( 1 \vee \log \verty + 1 \vee \log \lvert z\rvert )
\Big)
.\end{align}
\noindent {\rm (ii)}
\; For each $r>0$, there exists $\Ct \label{;Z10_2} >0$ independent of $N \in \mathbb{N} $ such that 
\begin{align}\label{:Z2c}
&\sup_{ z \in [\yhat -1, \yhat +1]} 
\Supr \lvert {M}_{\beta}^{N }(x,\yhat , z ) \rvert 
\le \cref{;Z10_2} ( \yvee ) 
\end{align}
\end{lemma}
\begin{proof}

We only consider the case $\beta=1$ since the proof for the other cases is similar. Recall that $ {M}_{1}^{N }$ is a scalar quaternion treated as a real number from \lref{l:31}.
From \eqref{:31l} and \eqref{:31o}, we deduce that 
\begin{align}\notag%\label{:Z2B}
&{M}_{1}^{N }(u_1,u_2,u_3) = 
\\\notag &
{J}_{1}^{N }(u_1,u_2) 
{J}_{1}^{N }(u_2,u_3) 
{J}_{1}^{N }(u_3,u_1) 
+ 
{J}_{1}^{N }(u_2,u_1) 
{J}_{1}^{N }(u_3,u_2) 
{J}_{1}^{N }(u_1,u_3) 
 \\\notag &
 - \sum_{\sigma }
\big( \partialtwo {J}_1^{N } \big) (u_{\sigma(1)},u_{\sigma(2)}) 
\Big( 
 \int_{u_{\sigma(3)}}^{u_{\sigma(2)}}
 {J}_1^{N }(t,u_{\sigma(3)})dt - 
\frac{{\rm sign}(u_{\sigma(2)}-u_{\sigma(3)})}{2}
\Big) 
\\&\notag \quad \quad 
\times 
{J}_1^{N }(u_{\sigma(3)},u_{\sigma(1)})
 \\\notag &
 - \sum_{\sigma }
\big( \partialtwo {J}_1^{N } \big) (u_{\sigma(1)},u_{\sigma(2)}) 
\Big( 
 \int_{u_{\sigma(1)}}^{u_{\sigma(3)}}
 {J}_1^{N }(t,u_{\sigma(1)})dt - 
\frac{{\rm sign}(u_{\sigma(3)}-u_{\sigma(1)})}{2}
\Big) 
\\&\label{:Z1g} \quad \quad 
 \times
{J}_1^{N }(u_{\sigma(3)},u_{\sigma(2)})
,\end{align}
where $ \sigma $ are taken over the permutations such that $ \sigma = \mathrm{id}, (1,2,3), (2,3)$. 

Applying \lref{l:Y3} to each term on the right-hand side of \eqref{:Z1g} and using \eqref{:Y1y} and \eqref{:Y4u}, we obtain (i) and (ii).
For example, 
\begin{align}& \notag 
\sup_{\vert x \vert \le r} \vert J_1^N(x , \yhat )J_1^N(\yhat , \zhat )J_1^N(\zhat , x )\vert 
\\ \notag &\le 
\Big(\sup_{\vert x\vert \le r} \vert J_1^N(x , \yhat ) \vert \Big) 
 \vert J_1^N(\yhat , \zhat )\vert 
\Big( \sup_{\vert x\vert \le r} \vert J_1^N(\zhat , x ) \vert \Big) 
\\\notag 
&\le \cref{;Y5} \cref{;Y]}^2 
\Big( \widehat{G}_{-\frac{1}{4},\frac{1}{4},1}(\yy , \zz ) + 
\yoneF \Big)
\zoneF 
\quad \text{ by \lref{l:Y3}}
\\\notag 
& \le 
 \cref{;Y5} \cref{;Y]}^2
 (G_{0,\frac{1}{4},1}(\yy , \zz )+G_{-\frac{1}{4},\frac{1}{2},1}(\yy , \zz )+ G_{\frac{1}{4},\frac{1}{4},0}(\yy , \zz )) 
\quad \text{ by \eqref{:Y1y}}
\\&\le \cref{;Y5} \cref{;Y]}^2 
\Big( 
 \widehat{G}_{-\frac{1}{4},-\frac{1}{4},1}(\yy , \zz ) + 
 \widehat{G}_{\frac{1}{4},\frac{1}{4},0}(\yy , \zz ) \Big) 
\quad \text{ by \eqref{:Y4u}}
\notag %\label{:Z2h}
.\end{align}
Thus, we have the estimate for the first term on the right-hand side of \eqref{:Z1g}. 
We can similarly prove the estimates for other terms on the right-hand side of \eqref{:Z1g}, and we omit it. 
\end{proof}

\SECT{Proof of \lref{l:3R}} \label{s:Z} 

In Subsection \ref{s:Z}, we prove \lref{l:3R}. We begin with a subsidiary lemma. 
For $a,b,c \in\mathbb{R} $, let $\widehat{G}_{a,b,c}$ be as in \eqref{:Y1y}. 

\begin{lemma} \label{l:Z1}
Let 
$a_{i} ,b_{i}, c_{i}, \kappa_{i}, \lambda_{i} , \nu_{i} \in\mathbb{R} $ be such that 
 $0 < a_{i}+b_{i}+c_{i} $, $ 0 < \nu_{i} $, and 
\begin{align}\label{:Z1c}&
 0 < 1+a_{i}+b_{i} -\kappa_{i} - ( \kappa_{i} + c_{i} -1) \lambda_{i} 
 ,\
 \kappa_i + c_i >1
,\
\max\{0 , \nu_{i} -1\} < \lambda_{i}
,\end{align}
where $ 1 \le i \le m $. 
Suppose that $ g $ is a non-negative function on $\mathbb{R} ^2$ satisfying 
\begin{align} \label{:Z1a}&
 g(u ,v) \le \cref{;Z8} \sum_{i=1}^m 
 \widehat{G}_{a_{i},b_{i},c_{i}} 
(\lvert u \rvert ,\lvert v \rvert ) \SV (\vert u \vert )\SV (\vert v \vert )
\quad \text{for all } u,v \in \mathbb{R} 
,\\ \label{:Z1b} &
\sup_{ v;\, \lvert u-v\rvert \le 1}g(u ,v) \le \cref{;Z8} 
\sum_{i=1}^m (1+\lvert u \rvert ^{\nu_{i}}) \SV (\vert u \vert )
\quad 
\text{for all } u \in \mathbb{R} 
\end{align}
for some positive constant $\Ct \label{;Z8} $ depending on 
$ (a_{i} ,b_{i}, c_{i}, \kappa_{i}, \lambda_{i} , \nu_{i})_{i=1}^m \in\mathbb{R}^{6m} $ and a positive slowly varying function $\SV $. 
Then,
\begin{align}\label{:Z1d}
\lim_{s\to\infty}
\int_{s \le \lvert u \rvert ,\, s \le \lvert v \rvert } \ \frac{g(u, v)}{\lvert uv\rvert } dudv =0 
.\end{align}
Here we call a positive measurable function $SV$ on $[0,\infty)$ slowly varying if 
\begin{align*}&
\lim_{x\to\infty}\frac{SV(\lambda x)}{SV(x)}=1 \text{ for all $\lambda>0$}
.\end{align*}
\end{lemma}
\begin{proof}
Without loss of generality, we can assume $ m=1$. 
Hence, we write $ a=a_1$, $ b=b_1$, $ c=c_1$, and so on. 
In \eqref{:Z1c}--\eqref{:Z1b}, retaking $ \cref{;Z8} $ and $\nu $ 
larger and $ a $ and $ b $ smaller for the same $ c, \kappa$, and $\lambda$, we can assume $\SV \le 1$. 

The assumptions made in \eqref{:Z1a} and \eqref{:Z1b} 
are essentially symmetric in $ u $ and $ v $, and depend only on 
the absolute values $ \lvert u \rvert $ and $ \lvert v \rvert $, and it is thus sufficient for \eqref{:Z1d} to prove that 
\begin{align}\label{:Z1m}&
\lim_{s\to\infty}
\int_{s}^{\infty} du \int_{u}^{\infty} dv 
 \ \frac{g(u, v)}{uv} =0
.\end{align}
Dividing $ [u,\infty)$ into $ [u,u + u^{-\lambda })$ and 
$ [u + u^{-\lambda }, \infty )$ and taking $ w=v-u$, we obtain that 
\begin{align}\label{:Z1n}&
\int_{s}^{\infty} du \int_{u}^{\infty} dv 
 \ \frac{g(u, v)}{uv}=
\int_{s}^{\infty} du 
\Big( \int_{0}^{u^{-\lambda }} dw + \int_{u^{-\lambda }}^{\infty} dw \Big) 
 \ \frac{g(u,u+w)}{u(u+w)} 
.\end{align}
From \eqref{:Z1c} and \eqref{:Z1a}, we deduce that 
\begin{align}\notag &
\int_{s}^{\infty} du
\int_{u^{-\lambda }}^{\infty} 
\frac{g(u,u+w)}{u(u+w)} 
dw 
\\ \notag 
& \le \cref{;Z8}
\int_{s}^{\infty} du \int_{u^{-\lambda }}^{\infty} 
\Big( 
\frac{1}{u(u+w)} 
\frac{1}{u^a (u+w)^b} 
\frac{1}{w^c} 
+ 
\frac{1}{u(u+w)} 
\frac{1}{u^b (u+w)^a} 
\frac{1}{w^c} 
\Big) 
dw 
\\ \notag &
 \le 2\cref{;Z8}
\int_{s}^{\infty} du \int_{u^{-\lambda }}^{\infty} 
\frac{1}{u^{2+a+b-\kappa }}\frac{1}{w^{\kappa + c}} 
dw 
\\ \label{:Z1p} &
= O (s^{-( 1+a+b -\kappa - ( \kappa + c -1) \lambda ) }), 
\quad s\to \infty
.\end{align}
From \eqref{:Z1c} and \eqref{:Z1b}, we deduce that 
\begin{align}\notag 
\int_{s}^{\infty} du 
\int_{0}^{u^{-\lambda }} \frac{g(u,u+w)}{u(u+w)} dw 
& \le 2 \cref{;Z8}
\int_{s}^{\infty} du 
\int_{0}^{u^{-\lambda }} \frac{1+u^{\nu }}{u(u+w)} dw 
\\ \notag &
 \le 2 \cref{;Z8}
\int_{s}^{\infty} du \ \frac{1+u^{\nu }}{u^2 }
\int_{0}^{u^{-\lambda }} dw 
\\ \label{:Z1o} & 
= O (s^{-1 + \nu - \lambda }), \quad s\to\infty
.\end{align}
Hence, \eqref{:Z1m} follows from 
\eqref{:Z1n}--\eqref{:Z1o} combined with \eqref{:Z1c}. 
\end{proof}

We are now ready to prove \lref{l:3R}. 

\smallskip 
\noindent 
{\em Proof of \lref{l:3R}}: 
We assume $s\ge 2r +1$, which implies 
\begin{align}\label{:Z1r}&
\text{$\vert y \vert \ge r+1$ for $ \lvert x \rvert \le r $ and $ \lvert x - y \rvert \ge s $}
.\end{align}
Let $ I^{N }_{\beta,k}$ be as in \eqref{:3Rz}. 
We begin by proving \eqref{:3Rb} for $k=1$. 
Using \lref{l:Y4} \thetag{ii}, \eqref{:Y1u}, \eqref{:Y1w}, and \eqref{:Z1r}, 
we have for $ \vertx \le r \le N^{2/3} $
\begin{align} \notag &
 I^{N }_{\beta,1}(x,s) =
\int_{\lvert x-\yy \rvert >s \atop y < -2 N^{2/3} } 
\frac{ \big\lvert L_{\beta }^{N } (x,\yy ) \big\lvert }
{\lvert x-\yy \rvert ^2}d\yy 
+ 
\int_{\lvert x-\yy \rvert >s \atop -2 N^{2/3} \le y } 
\frac{ \big\lvert L_{\beta }^{N } (x,\yy ) \big\lvert }{\lvert x-\yy \rvert ^2}d\yy 
\\ &\label{:Z3a} 
 \le 
\cref{;Y10} \Big(
\int_{\lvert x-\yy \rvert >s \atop y < -2 N^{2/3} } 
\frac{1} {\lvert x- \yy \rvert ^2 \lvert \iotaN (\yy ) \rvert ^{1/4}}d\yy 
+ 
 \int_{\lvert x-\yy \rvert >s \atop -2 N^{2/3} \le y } 
\frac{1} {\lvert x-\yy \rvert ^2\verty ^{1/4}}d\yy 
\Big)
.\end{align}
Here $ \iotaN (y) =-(4N ^{2/3}+ y ) $ is as in \eqref{:Y1u}. We easily see 
\begin{align}\notag &
\lim_{s\to\infty}\SupN \Supr 
\int_{\lvert x-\yy \rvert >s \atop y < -2 N^{2/3} } 
\frac{1} {\lvert x- \yy \rvert ^2 \lvert \iotaN (\yy ) \rvert ^{1/4}}d\yy 
=0
,\\& \label{:Z3b}
\lim_{s\to\infty}\Supr 
\int_{\lvert x-\yy \rvert >s} \frac{1} {\lvert x-\yy \rvert ^2\verty ^{1/4}}d\yy =0 
.\end{align}
From \eqref{:Z3a} and \eqref{:Z3b}, we obtain \eqref{:3Rb} for $ k=1$. 

We next prove \eqref{:3Rb} for $k=2$ using \lref{l:Z1}. 
So we verify the conditions in \lref{l:Z1}. 
Taking \eqref{:Y4a} into account, we set $ m = 5$, $\nu_i=1$, $1 \le i \le 5 $, and 
\begin{align*}&
((a_i,b_i,c_i))_{i=1}^5 = (
(-\frac{1}{2}, \frac{1}{2}, 2 ), (0, 0,2 ), 
 (-\frac{3}{4}, \frac{1}{4},1 ), (-\frac{1}{4}, -\frac{1}{4}, 1 ), 
 (\frac{1}{4}, \frac{1}{4}, 0 ) 
)
.\end{align*}
We then find constants $\lambda_i$, $\kappa_i$, $1\le i \le 5$, such that \eqref{:Z1c} holds. 

Let $ H_N = [-2N^{2/3},\infty )$. 
Then $ \iotaN (H_N) = H_N^c \cup \{ -2N^{2/3} \} $. 
We set 
\begin{align*}&
H_N^{1} = H_N \ts H_N , \ 
H_N^{2} = H_N^c \ts H_N , \ 
H_N^{3} = H_N \ts H_N^c , \ 
H_N^{4} = H_N^c \ts H_N^c 
.\end{align*}
Let $ \Ixxs = \{ (y,z) \in \mathbb{R}^2\, ;\,\lvert x-\yy \rvert >s ,\, \lvert x' -\zz \rvert >s \} $. 
Then we have 
\begin{align}\notag &
 I^{N }_{\beta,2}(x,s) = \sum_{i=1}^4
\int_{ \Ixxxs \cap H_N^i} 
\frac{\big\lvert L_{\beta }^{N } (y,z)\big\rvert }
{\lvert x-\yy \rvert \lvert x-\zz \rvert } d\yy d\zz 
\\&\notag = 
 \int_{ \Ixxxs \cap H_N^1} 
\frac{\big\lvert L_{\beta }^{N } (y,z)\big\rvert }
{\lvert x-\yy \rvert \lvert x-\zz \rvert } 
d\yy d\zz 
+ 
 \int_{ I ( \iotaN(x) , x , s ) \cap H_N^1 }
\frac{\big\lvert L_{\beta }^{N } ( \iotaNy , z )\big\rvert }
{\lvert x- \iotaNy \rvert \lvert x - z \rvert } d\yy d\zz 
\\\notag &
+ 
 \int_{ I ( x , \iotaN(x) , s ) \cap H_N^1 }
\frac{\big\lvert L_{\beta }^{N } ( y , \iotaN(z ) )\big\rvert }
{\lvert x-\yy \rvert \lvert x - \iotaN(z ) \rvert } d\yy d\zz 
\\& \label{:Z3k} 
+ 
 \int_{ I ( \iotaN (x) , \iotaN (x) , s ) \cap H_N^1 }
\frac{\big\lvert L_{\beta }^{N } ( \iotaNy , \iotaN (z) )\big\rvert }
{\lvert x - \iotaN (\yy ) \rvert \lvert x - \iotaN (\zz) \rvert } d\yy d\zz 
.\end{align}
We take 
\begin{align}\label{:Z3d}
g (y,z)=\supN 1_{H_N^1} (y,z) \max \{\big\{ & 
\lvert L_{\beta }^{N } ( u (y) , v(z) ) \big\rvert ; u , v \in \mathscr{M}_N 
\}
.\end{align}
From \eqref{:Z3k}, \eqref{:Z3d}, and $ \lvert x - \iotaNy \rvert = \lvert \iotaN (x) - \yy \rvert $, we have 
\begin{align}&\notag 
 I^{N }_{\beta,2}(x,s) \le 
 \int_{ \Ixxxs } 
\frac{ g (y,z) }
{\lvert x-\yy \rvert \lvert x-\zz \rvert } d\yy d\zz 
\\&\notag 
+ 
 \int_{ I ( \iotaN(x) , x , s ) }
\frac{ g ( y , z ) }
{\lvert \iotaN (x) - \yy \rvert \lvert x - z \rvert } d\yy d\zz 
\\&\notag 
+ 
 \int_{ I ( x , \iotaN(x) , s ) }
\frac{ g ( y , z ) }
{\lvert x-\yy \rvert \lvert \iotaN (x) - z \rvert } d\yy d\zz 
\\&\notag 
+ 
 \int_{ I ( \iotaN (x) , \iotaN (x) , s ) } 
\frac{ g (y,z) }
{\lvert \iotaN (x) - \yy \rvert \lvert \iotaN (x) - \zz \rvert } d\yy d\zz 
\\\label{:Z3C} &
=: I_{(1)} ( x , s ) + I_{(2)}^N ( x , s ) + I_{(3)}^N ( x , s ) + I_{(4)}^N ( x , s ) 
.\end{align}
%	Here we used \eqref{:Y1w}, \eqref{:Z3c}, and \eqref{:Z3d} for the third line. 
Recall that $ s\ge 2r +1$. Then we have 
\begin{align}& \label{:Z1f}
 \sup \Big\{
\frac{\lvert \yy \rvert \lvert \zz \rvert } { \lvert x-\yy \rvert \lvert x-\zz \rvert }; 
\vertx \le r , (y,z) \in \Ixxxs 
\Big\}
\le \Big(\frac{s}{s - r}\Big)^2
.\end{align}

From \eqref{:Y4a} and \eqref{:Y4c}, $ g $ in \eqref{:Z3d} 
satisfies \eqref{:Z1a} and \eqref{:Z1b}. 
Hence, we can apply \lref{l:Z1} to $ g $. 
Using \eqref{:Z3d}--\eqref{:Z1f} and \lref{l:Z1}, we deduce 
\begin{align}\notag &
\limsupi{s} \Supr I_{(1)} ( x , s ) 
% \\\notag & 
\le 
\limsupi{s}
 \Supr 
 \int_{ \Ixxxs }
\Big(\frac{s}{s-r} \Big)^2 \frac{g (y,z)}{\lvert y \rvert \lvert z \rvert } 
 d\yy d\zz 
\\ \label{:Z1h}&
 \le 
\limsupi{s}
 \int_{ I (0,0,s-r) }
\Big(\frac{s}{s-r} \Big)^2 \frac{g (y,z)}{\lvert y \rvert \lvert z \rvert } 
 d\yy d\zz = 0 
.\end{align}
Here we used $ I ( x , x , s ) \subset I ( 0 , 0 , s - r ) $ for $ \lvert x \rvert \le r $ for the last line. 

From \lref{l:Y4}, for each $ \epsilon > 0 $, there exists an $ N_0 $ such that 
\begin{align}\label{:Z1H}&
 I_{(2)}^N ( x , 2r + 1 ) < \epsilon 
\quad \text{ for all } N \ge N_0 
.\end{align}
It is noteworthy that $ I_{(2)}^N ( x , s ) $ is decreasing in $ s $. 
Hence from \eqref{:Z1H}, we have 
\begin{align}\label{:Z1x}&
\sup_{N \ge N_0 }
\sup_{ s \ge 2r + 1 } I_{(2)}^N ( x , s ) < \epsilon 
.\end{align}
From \lref{l:Y4}, for each $ \epsilon > 0 $, we have an $ s_0 \ge 2r + 1 $ such that 
\begin{align}\label{:Z1y}&
\max_{N \le N_0 }
\sup_{ s \ge s_0 } I_{(2)}^N ( x , s ) < \epsilon 
.\end{align}
From \eqref{:Z1H}--\eqref{:Z1y} and 
$ I ( x , \iotaN (x) , s ) \subset I ( 0 , \iotaN (0) , s - r ) $ for $ \lvert x \rvert \le r $, 
\begin{align}\label{:Z1z}&
\limsupi{s} \SupN \Supr I_{(2)}^N ( x , s ) < 2\epsilon 
.\end{align} 
Because \eqref{:Z1z} holds for arbitrary $ \epsilon > 0 $, we have 
\begin{align}\label{:Z1Z}&
\limsupi{s} \SupN \Supr I_{(2)}^N ( x , s ) = 0 
.\end{align} 
In a similar fashion, we have \eqref{:Z1Z} for $ I_{(3)}^N $ and $ I_{(4)}^N $. 
Hence from \eqref{:Z3C}, \eqref{:Z1h}, and \eqref{:Z1Z}, we obtain \eqref{:3Rb} for $k=2$. 

We finally prove \eqref{:3Rb} for $k=3 $ using \lref{l:Z1}. 
Taking \eqref{:Z2d} and \eqref{:Z2c} in \lref{l:Y4_M} into account, we set $ m=3$ and $\nu_i =1$ for $1 \le i \le 3$. 
We set 
$ 
((a_i,b_i,c_i))_{i=1}^3 = 
((-\frac{3}{4}, \frac{1}{4},1 ), (-\frac{1}{4}, -\frac{1}{4}, 1), (\frac{1}{4}, \frac{1}{4}, 0) ) 
$. 
We then easily see that there exist constants $\lambda_i$, $\kappa_i$, $1\le i \le 3 $ 
such that the condition \eqref{:Z1c} holds.

Similarly as \eqref{:Z3k}, we have 
\begin{align}\notag &
 I^{N }_{\beta,3}(x,s) = \sum_{i=1}^4
\int_{ \Ixxxs \cap H_N^i} 
\frac{\big\lvert M_{1 }^{N } (x,y,z)\big\rvert }
{\lvert x-\yy \rvert \lvert x-\zz \rvert } d\yy d\zz 
\\&\notag = 
 \int_{ \Ixxxs \cap H_N^1} 
\frac{\big\lvert M_{1 }^{N } (x,y,z)\big\rvert }
{\lvert x-\yy \rvert \lvert x-\zz \rvert } 
d\yy d\zz 
+ 
 \int_{ I ( \iotaN(x) , x , s ) \cap H_N^1 }
\frac{\big\lvert M_{1 }^{N } ( x, \iotaNy , z )\big\rvert }
{\lvert x - \iotaNy \rvert \lvert x- z \rvert } d\yy d\zz 
\\\notag & + 
 \int_{ I ( x , \iotaN(x) , s ) \cap H_N^1 }
\frac{\big\lvert M_{1 }^{N } ( x,y , \iotaNz )\big\rvert }
{\lvert x-\yy \rvert \lvert x - \iotaN(z ) \rvert } d\yy d\zz 
\\&\label{:Z3d_2} 
+ 
 \int_{ I ( \iotaN (x) , \iotaN (x) , s ) \cap H_N^1 }
\frac{\big\lvert M_{1 }^{N } (x, \iotaNy , \iotaN (z) )\big\rvert }
{\lvert x - \iotaN (\yy ) \rvert \lvert x - \iotaN (\zz) \rvert } d\yy d\zz 
.\end{align}
Let $ \mathscr{M}_N = \{ \mathrm{id}, \iota_N \} $ be the set of maps as in \eqref{:Y1Z}. We take
\begin{align} \label{:Z3d_3} 
\tilde{g}_x (y,z) =& \SupN 1_{H_N^1}(y,z) 
 \max\{ \vert M_1^N(x, u (y) , v (z) )\vert ; u , v \in \mathscr{M}_N 
\} 
.\end{align}
Using \eqref{:Z3d_2} and \eqref{:Z3d_3}, we have 
\begin{align}
 I^{N }_{\beta,3}(x,s) 
&\le \notag
 \int_{ \Ixxxs } 
\frac{ \tilde{g}_x (y,z) }{\lvert x-\yy \rvert \lvert x-\zz \rvert } d\yy d\zz 
+ 
 \int_{ I ( \iotaNx , x , s ) } 
\frac{ \tilde{g}_x (y,z ) }
{\lvert \iotaNx - \yy \rvert \lvert x - z \rvert } d\yy d\zz 
\\&\notag + 
 \int_{ I ( x , \iotaN(x) , s ) }
\frac{ \tilde{g}_x (y,z ) }
{\lvert x-\yy \rvert \lvert \iotaN (x) - z \rvert } d\yy d\zz 
\\&\notag 
+ 
 \int_{ I ( \iotaN (x) , \iotaN (x) , s ) } 
\frac{ \tilde{g}_x (y,z) }
{\lvert \iotaN (x) - \yy \rvert \lvert \iotaN (x) - \zz \rvert } d\yy d\zz 
\\ \notag%\label{:Z3C} 
&=: \tilde{I}_{(1)} ( x , s ) + \tilde{I}_{(2)}^N ( x , s ) + \tilde{I}_{(3)}^N ( x , s ) 
+ \tilde{I}_{(4)}^N ( x , s ) 
.\end{align}
From \eqref{:Z2d} and \eqref{:Z2c}, $ \tilde{g}_x $ satisfies \eqref{:Z1a} and \eqref{:Z1b}. 
Thus, we can apply \lref{l:Z1} to $ \tilde{g}_x $. 
The quantities $ \tilde{I}_{(1)}, \tilde{I}_{(2)}^N, \tilde{I}_{(3)}^N $, and $ \tilde{I}_{(4)}^N $ correspond to $ I_{(1)} , I_{(2)}^N , I_{(3)}^N $, and $ I_{(4)}^N $ in \eqref{:Z3C}, respectively. 
Hence, using the same argument as for $k=2$, we obtain (3.45) for $k=3$. 
\qed

\subsection{Proof of \1} \label{s:19} 
Let $ \mathsf{f} $, $ \mathsf{g} $, $ \uu $, and $\vv $ be as in \eqref{:91f}--\eqref{:91h} and \eqref{:18e}. 
We often omit $ N $ from the notation of these symbols if 
no confusion occurs. We thus write $ \mathsf{f} (\theta ) = \mathsf{f} \Ntheta $, $ \mathsf{g} (\theta ) = \mathsf{g} \Ntheta $, 
$ \uu (\theta ) = \uu \Ntheta $, and $ \vv (\theta ) = \vv \Ntheta $. 
Let 
\begin{align*}&
 \dNz = \big\{ \theta \in (0,\pi /2] ; \uu (\theta) \in \None \big\}
, \\&
\dNone = \{ \theta \in \dNz ;\, N^{-1} < \cos \theta \}
.\end{align*}
Clearly, $ \dNz \subset \dN $ by \eqref{:18d}. 
From \eqref{:91h}, we see that 
\begin{align} \label{:191r}&
 \uu (\dNz ) = \None 
.\end{align}

\begin{lemma}\label{l:191} %%%%%%%%%%%%%%%
Let $\FN $ be as in \eqref{:18n}. Set 
\begin{align} \label{:191s}&
\FN \zzz ( \theta ) = \FN ( \theta ) (1+\cos \theta)^{1/2} \lvert \uu (\theta ) \rvert ^{1/2}
.\end{align} 
There then exist $N_0\in \mathbb{N} $ and $\Ct \label{;21a} > 0 $ such that, for all $ N\ge N_0 $, 
\begin{align} \label{:191u}&
0 < \FN \zzz ( \uu ^{-1} (\uU ) ) \le \cref{;21a}\lvert \uU \rvert ^{-1/4} \quad \text{ for }
\uU \in \None 
,\\\label{:191t}& 
\ODod{ }{ \vV } \FNvVz <0 
\quad \text{ for } \vV \in \vv (\dNone ) 
.\end{align}
\end{lemma}
\begin{proof} 
From \eqref{:18O} and \eqref{:191s}, 
\begin{align}&\label{:191b}
 \FN \zzz ( \theta ) > 0 \quad \text{for } \theta \in \dN 
.\end{align}
Let $ \mathsf{f} (N, \theta ) =(N+1)^{2/3}\sin^2 \theta$ be as in \eqref{:91f}. 
From this, \eqref{:18n}, and \eqref{:191s}, we see that 
\begin{align}\label{:191z}
\FN \zzz ( \theta ) 
&=\Big(\frac{N}{N+1}\Big)^{1/6} \frac{(1+\cos \theta)^{1/2}\lvert \uu (\theta ) \rvert ^{1/2}}
{ \mathsf{f} ^{3/4}-(4\sqrt{(N +1) \sin \theta })^{-1}}
% \Big\rvert _{\theta = \uu ^{-1}(\uU )}
.\end{align}
Using \eqref{:91P} and taking $ \theta = \uu ^{-1} (\uU )$, we obtain \eqref{:191u} 
from \eqref{:191b} and \eqref{:191z}. 

From \eqref{:91h} and \eqref{:191z}, we deduce that 
\begin{align} \notag 
 \FN \zzz ( \theta ) ^2 &= 
 \Big( \frac{N}{N+1}\Big)^{1/3}
 \frac{4\sin \theta (1+\cos \theta)\lvert \uu (\theta ) \rvert }
 {4(N+1)\sin^4 \theta - 2 \sin^2\theta + (4(N+1) )^{-1}}
 \\ \notag
 &= \frac{2N^{1/2}}{(N+1)^{1/3}} 
\frac{(1+\cos \theta)(\sqrt{N}-\sqrt{N+1}\cos \theta)}
{(N+1)\sin^3 \theta - \frac{\sin \theta}{2} + ( 16(N+1)\sin\theta )^{-1}}
 \\ \label{:191v} 
 &=
\frac{2N^{1/2}}{(N+1)^{5/6}} \frac{\pppW }{\qqqW } 
,\end{align}
where we set $ \pppW = \mathsf{p} (N, \theta )$ and $ \qqqW = \qqq (N, \theta )$ such that 
\begin{align} \notag 
&\pppW = 
\sin^2 \theta - \Big(1-\sqrt{\frac{N}{N+1}}\Big)(1+\cos\theta), 
\\\label{:191w} %\label{:191x}
&\qqqW = 
\sin^3\theta - \frac{\sin \theta}{2(N+1)} + 
\frac{1}{16(N+1)^2}\frac{1}{\sin \theta}
.\end{align}
Let $ \theta \in \dNone $. Then, $ N^{-1} < \cos \theta $ by definition. 
Furthermore, we have $ \uu ( \theta ) \le - 1 $. 
Hence, we find $ \Ct > 0 \label{;192}$ independent of $ N $ such that 
$ \cref{;192} N^{-1/3} \le \sin \theta $. We thus see for $ \theta \in \dNone $ that 
\begin{align}\label{:191x}& 
 N^{-1} < \cos \theta , \quad \cref{;192} N^{-1/3} \le \sin \theta 
.\end{align}
From \eqref{:191w} and \eqref{:191x}, we find $ \Ct \label{;191} > 0 $ and 
$N_0\in \mathbb{N} $ such that, for all $ N\ge N_0 $, 
\begin{align} \notag %\label{:191y}
& \pppWW \qqqW- \pppW \qqqWW 
< 
-\sin^4 \theta \cos \theta + \Big ( 1- \sqrt{\frac{N}{N+1}} \Big ) \sin^4 \theta 
	+ \frac{\cref{;191}}{N} \sin^2 \theta \cos \theta 
\\\notag & =
\Big( - \frac{\sin^2 \theta }{3} + \frac{\cref{;191}}{N}\Big) \sin^2 \theta \cos \theta 
 - 
\Big( \frac{2\cos \theta }{3} -
\Big ( 1 - \sqrt{\frac{N}{N+1}} \Big ) \Big) \sin^4 \theta 
 \\\notag & < 
\Big( - \frac{ \cref{;192} ^2 }{ 3 N ^{2/3}} + \frac{\cref{;191}}{N}\Big) \sin^2 \theta \cos \theta - 
\Big( \frac{2}{3 N} - \Big ( 1 - \sqrt{\frac{N}{N+1}} \Big ) \Big) \sin^4 \theta < 0 
,\quad \theta \in \dNone 
.\end{align}
Combining this with \eqref{:191v}, we obtain for all $ N\ge N_0 $, 
\begin{align}&\label{:191a}
\ODod{}{\theta } \FN \zzz ( \theta ) ^2 = 
\frac{2N^{1/2}}{(N+1)^{5/6}} \frac{\pppWW \qqqW - \pppW \qqqWW }{\qqqW ^2}< 0 
 , \quad \theta \in \dNone 
.\end{align} 
From \eqref{:191b} and \eqref{:191a}, we find that $ d \FN \zzz / d \theta < 0 $ on $ \dNone $. 
From \eqref{:18h}, we deduce that $ { d \vv }/{ d \theta} >0$ on $\dN $. 
Hence, we obtain $ {d \vv ^{-1} } /{ d \vV }>0$ on $ \vv ( \dN ) $. 
Clearly, $ \vv (\dNone ) \subset \vv ( \dN ) $ by \eqref{:18d}. 
Combining these, we deduce that 
\begin{align*}&
\ODod{ \FNvVz }{ \vV } = 
\Big (\ODod{ \FN \zzz }{ \theta}\Big ) ( \vv ^{-1}(\vV ) )
\Big (\ODod{ \vv ^{-1} }{ \vV }\Big ) (\vV ) < 0 
,\quad \vV \in \vv (\dNone )
.\end{align*}
Hence, we obtain \eqref{:191t}. 
\end{proof}
\begin{lemma}\label{l:192} 
Let $ \uu $, $ \vv $, and $\FN $ be as in \eqref{:91h}, \eqref{:18e}, and \eqref{:18n}. 
For $ \iii = 1,2 $, let 
\begin{align}\label{:192a}&
\FN ^{(\iii )}( \theta ) = \FN ( \theta )\lvert \uu (\theta ) \rvert ^{\iii }
.\end{align}
There then exists $N_0\in \mathbb{N} $ and a positive constant $\Ct \label{;22c}$ such that for all $ N\ge N_0 $, 
\begin{align} \label{:192d}& 
 0 < \FN ^{(\iii )} ( \uu ^{-1} (\uU ) ) 
\le \cref{;22c}\lvert \uU \rvert ^{ \iii - ( 3/4 )} 
\quad \text{ for }\uU \in \None 
,\\ \label{:192c}& 
0 < 
\ODod{ \FNvVi }{ \vV }
\quad \text{ for }
\vV \in \vv (\dNz ) 
.\end{align}
\end{lemma}

%	\end{document}

\begin{proof} 
From \eqref{:18O}, we have $ \FN > 0 $ on $ \dN $. 
From this and \eqref{:192a}, we obtain the first inequality in \eqref{:192d}. 
Recall that $\sqrt{(N+1)\sin \theta} = (N+1)^{1/3}\mathsf{f} \Ntheta ^{1/4}$ by \eqref{:91f}. 
From \eqref{:18n} and \eqref{:192a}, we see that $ \FN \Fone = \FN \Fone ( \theta )$ satisfies 
%	\end{document}
\begin{align}\label{:192f}&
\FN \Fone 
= \Big( \frac{N}{N+1} \Big)^{1/6} 
 \frac{ -\mathsf{f}^{1/4} \uu }{\mathsf{f} - \QM (N+1)^{-1/3}}
.\end{align}
Hence, we obtain the second inequality in \eqref{:192d} for $ \iii = 1$ from \eqref{:91P} and \eqref{:192f}. 
%
%	\end{document}
From \eqref{:91f} and \eqref{:192f}, we have for each $ N \in \mathbb{N}$, 
\begin{align} \label{:192g}
\ODod{ \FN \Fone }{ \theta} 
&=\Big( 
\frac{N}{N+1} \Big)^{1/6}
\frac{\mathsf{f}^{1/4}\rrr }{4 (\mathsf{f} - \QM (N+1)^{-1/3} )^2 }
,\end{align}
where $ \rrr (\theta ) = \rrr \Ntheta $ is such that 
\begin{align} \label{:192h} 
\rrr &= \uu \big( 3\mathsf{f}' + \QM (N+1)^{-1/3}(\log\mathsf{f} )' \big)+
 \uu ' ( -4 \mathsf{f} +(N+1)^{-1/3} )
.\end{align}
We shall prove $ \rrr (\theta ) > 0 $ on $ \dNz $ for sufficiently large $ N $. 

%	\end{document}
A straightforward calculation with \eqref{:91f} and \eqref{:91h} yields 
\begin{align}\label{:192i}
\frac{3\uu \mathsf{f} '-4 \uu ' \mathsf{f}}{4N^{1/6}(N+1)^{5/6}\mathsf{f} ^{1/2}}=
 \cos^2 \theta -3\sqrt{1 - (N+1)^{-1}}\cos \theta +2 
.\end{align}
From \eqref{:91f}, \eqref{:91h}, and \eqref{:18f}, we see that 
\begin{align}\notag & 
(N+1)^{-1/3}(\log\mathsf{f} )' = 2 (\cos \theta )\mathsf{f}^{-1/2}
,\quad %\\ \notag &
\uu ' = -2N^{1/6}(N+1)^{1/6} \mathsf{f}^{1/2}
.\end{align}
Hence, using \lref{l:91}, we obtain 
\begin{align}\label{:192j}
\sup_{\theta \in \dNz }\Big\lvert 
\frac{\uu \QM (N+1)^{-1/3} (\log\mathsf{f} )' + \uu '(N+1)^{-1/3}}
{4N^{1/6}(N+1)^{5/6}\mathsf{f}^{1/2}} 
\Big\rvert 
= 
\mathcal{O}(N ^{-1})
.\end{align}
From \eqref{:192h}--\eqref{:192j}, we deduce that 
\begin{align} \label{:192l} &
\sup_{\theta \in \dNz }\Big\lvert 
\frac{\rrr}{4N^{1/6}(N+1)^{5/6}\mathsf{f}^{1/2}}
- (\cos^2\theta -3\cos\theta +2) \Big\rvert 
= \mathcal{O} (N^{-1})
.\end{align}

Note that $ \cos^2\theta -3\cos\theta +2 \ge \cref{;192l} N^{-2/3}$ on $ \dNz $ for some 
$ \Ct \label{;192l} > 0 $ independent of $ N $. 
From this and \eqref{:192l}, we find $N_0\in \mathbb{N} $ such that $ \rrr (\theta ) > 0 $
on $\dNz $ for all $ N\ge N_0$. 
Hence, from \eqref{:192g}, we obtain $ d \FN \Fone/ d \theta > 0$ on $\dNz $ for all $ N\ge N_0$. 
From \eqref{:18h}, we see that $ d \vv / d \theta >0$ on $\dN $. 
Then, $ d \vv ^{-1}/ d \vV >0$ on $ \vv (\dN ) $. 
Hence, for all $ N\ge N_0$, 
\begin{align*}&
\ODod{ \FNvVone }{ \vV } =
\Big (\ODod{ \FN \Fone }{ \theta}\Big ) ( \vv ^{-1}(\vV ) ) 
\Big (\ODod{ \vv ^{-1}}{ \vV }\Big ) (\vV )
 > 0 
\quad \text{ for } v \in \vv (\dNz ) 
,\end{align*}
which yields \eqref{:192c} for $ \iii = 1 $. 
From \eqref{:192a}, we have $ \FN \Ftwo = \FN \Fone \lvert \uu \rvert $. 
Hence, we deduce \eqref{:192d} and \eqref{:192c} for $ \iii =2 $ from the equations for $ \iii = 1 $. 
\end{proof}
%

%	\end{document}

Let $ \mathsf{f} $, $ \mathsf{g} $, $ \ct ^{km}$, and $ \Ckm $ be as in 
\eqref{:91f}, \eqref{:91g}, \eqref{:93u}, and \eqref{:93v}. 
Set 
\begin{align}\label{:194x}
\UNkm ( \uU )= & 
\mathsf{f} \Ntheta ^{-1/4} (1+\cos \theta)^{1/2} \lvert \uU \rvert ^{1/2}
 \Ckm \Ntheta \cos ( \mathsf{g} \Ntheta + \ct ^{km}\Ntheta ) 
,\\ \label{:194y}
\UN ( \uU ) = & \mathsf{f} \Ntheta ^{-1/4} 
\sin \Big( \mathsf{g} \Ntheta + \frac{\pi}{4}-\frac{\theta}{2}\Big), 
\quad \text{ where $ \theta = \uu ^{-1} (\uU )$}
.\end{align}

\begin{lemma} \label{l:194}%%%%%%%%%%%%%%%%%%%%%%%%%%%%%%%%%%%
Let $k,m \in \zN $ with $m\le k$. Then there exist positive constants $ \Ct \label{;19X}$ and $\Ct \label{;Up_11} $ independent of $ N $ such that for $\uU \in \None $,
\begin{align}& \label{:194z}
\lvert \UNkm ( \uU ) \rvert \le \cref{;19X} \lvert \uU \rvert ^{1/4} \ ( k \ge 1), %text{ on $ [-2 N ^{2/3}, -1]}
\lvert \UNz ( \uU ) \rvert \le \cref{;19X} N^{1/6}, 
\lvert \UN ( \uU ) \rvert \le \cref{;19X} \lvert \uU \rvert ^{-1/4} 
,\\&\label{:194a}
\Big\lvert \ODod{}{\uU } \UNkm ( \uU )\Big\rvert \le \cref{;Up_11}\lvert \uU \rvert ^{3/4}, \quad 
\Big\lvert \ODod{}{\uU } \UN ( \uU ) \Big\rvert \le \cref{;Up_11}\lvert \uU \rvert ^{1/4} 
.\end{align}
\end{lemma}%%%%%%%%%%%%%%%%%%%%%
\begin{proof}
From \eqref{:91P}, \eqref{:U1d}, \eqref{:194x}, and \eqref{:194y}, we obtain \eqref{:194z}. 

From \eqref{:194x}, we see that 
\begin{align*}
\ODod{\UNkm ( \uU )}{ \uU } = & \UNkm ( \uU ) 
\Big( 
 \ODod{\log \mathsf{f} ^{- 1/4}}{\uU } +
\ODod{\log (1+\cos \theta )^{1/2}}{\uU } + 
 \ODod{ \log \lvert \uU \rvert ^{1/2}}{\uU } 
 + 
 \ODod{ \log \Ckm }{\uU } \Big) 
 \\ \notag + & \Big( 
 \mathsf{f} ^{-1/4}(1+\cos \theta )^{1/2} \lvert u \rvert ^{1/2} \Ckm \Big) 
 \ODod{\cos (\mathsf{g} + \ckm )}{\uU }=:\mathbb{I}_1+\mathbb{I}_2 
 .\end{align*}
Using \eqref{:91f}, we see $ { d \mathsf{f} }/{ d \theta } = (N+1)^{2/3} 2\sin \theta \cos \theta $. 
From this, $ \theta = \uu ^{-1} (\uU ) $, and \eqref{:18f},% we see 
\begin{align}& \label{:194e}
\ODod{\mathsf{f} ( \uu ^{-1} (\uU ))}{\uU } =
 \frac{ (N+1)^{2/3} 2 \sin \theta \cos \theta }
{-2 N^{1/6}(N+1)^{1/2} \sin \theta }
= - (1 + N^{-1})^{1/6} \cos \theta
.\end{align}
From \eqref{:91P}, \eqref{:194z}, and \eqref{:194e}, 
we find $\Ct \label{;194e} >0$ independent of $ N $ such that 
\begin{align*}&
 \Big\lvert \UNkm ( \uU ) \ODod{ \log \mathsf{f} ^{- 1/4}}{ \uU }\Big\rvert \le \cref{;194e} \lvert \uU \rvert ^{-3/4} 
\quad \text{ for $\uU \in \None $}
.\end{align*}
By direct calculation, we see that other three terms in the curly brackets of $ \mathbb{I}_1$ are bounded. 
Hence, from this and \eqref{:194z}, we deduce that $ \lvert \mathbb{I}_1 \rvert \le \Ct \lvert \uU \rvert ^{1/4} \label{;194E}$ on $ [-2 N ^{2/3}, -1]$, 
where $ \cref{;194E}$ is a positive constant independent of $ N $. 

From \eqref{:91g}, we see that 
$ { d \mathsf{g} }/{ d \theta } = (N+1)(1-\cos 2\theta ) $. 
Using this, \eqref{:18f}, and \eqref{:91f}, we deduce that 
\begin{align}\label{:194f} & 
\ODod{\mathsf{g} (\uu ^{-1} (\uU )) }{\uU } = 
\ODod{\mathsf{g} }{\theta } \ODod{\theta }{\uU } = 
\frac{ (N+1)(1-\cos 2\theta )} {-2 N^{1/6}(N+1)^{1/2} \sin \theta } =
 - (1 + N^{-1})^{1/6}\mathsf{f}^{1/2}
.\end{align}
Combining \eqref{:194f} with \eqref{:91P}, \eqref{:93u}, and \eqref{:93v}, we find 
a constant $ \Ct \label{;19g} > 0 $ independent of $ N $ such that, for $\uU \in \None $, 
\begin{align}&\notag %\label{:194h}&
 \lvert \mathbb{I}_2 \rvert = \Big\lvert \Big( \mathsf{f} ^{-1/4}(1+\cos \theta )^{1/2} \lvert u \rvert ^{1/2} \Ckm \Big) 
 \ODod{\cos (\mathsf{g} + \ckm )}{\uU }\Big\rvert 
 \le \cref{;19g} \lvert u \rvert ^{3/4} 
.\end{align}
We thus prove the first inequality in \eqref{:194a}. 
The proof of the second inequality is similar. 
\end{proof}

\begin{lemma}\label{l:195}%%%%%%%%%%%%%%%%%%%%%%%%%%
Let $ k,m \in \zN $ with $m\le k$ and $ i = 1 , 2 $. 
Then there exists a constant $\Ct \label{;J_11} >0$ independent of $ N $ 
such that for $ \xXx \in \None $, 
\begin{align}\label{:195x}
& 
\sup_{\yYy \in \mathbb{R} }
\Big\lvert \int_\xXx ^\yYy \mathbf{1}( \lvert u-\xXx \rvert \le 1)\frac{\int_{\xXx }^{u} \UNkm ( \vV ) d\vV }{u - \xXx } du \Big\rvert \le 
\cref{;J_11} \lvert \xXx \rvert ^{-1/4}( 1 \vee \log \lvert \xXx \rvert ),
\\ \label{:195y}
& 
\sup_{\yYy \in \mathbb{R} }
\Big\lvert \int_\xXx ^{\yYy } \mathbf{1}( \lvert u-\xXx \rvert \le 1) 
\frac{\int_\xXx ^{u} \lvert \vV \rvert ^{ i }\UN ( \vV ) d\vV }{u - \xXx } du \Big\rvert \le 
 \cref{;J_11} \lvert \xXx \rvert ^{i -(3/4)}( 1 \vee \log \lvert \xXx \rvert )
.\end{align}
\end{lemma}

\begin{proof}
Throughout the proof of \lref{l:195}, we only consider the case that $ \xXx < \yYy $ 
because we can prove the other case in the same fashion. 
We also note that we can replace $ \sup_{\yYy \in\mathbb{R} }$ with 
$ \sup_{ x \le \yYy \le x + 1}$ because of 
$ \mathbf{1}( \lvert u-\xXx \rvert \le 1) $ in the integrand. 
We shall prove \eqref{:195x} and \eqref{:195y} by replacing 
$ \sup_{\yYy \in\mathbb{R} }$ with $ \sup_{\xXx \le \yYy \le \xXx + 1}$. Thus, $ \xXx \le \yYy \le \xXx + 1$ in the following. 

{\em Proof of \eqref{:195x}}: 
If $ k \ge 1$, then we easily obtain \eqref{:195x} from \eqref{:194z}. 

We present the proof for the case that $ k=0 $. 
Let $ \hat{ \vv } = \vv \circ \uu ^{-1}$ be as in the proof of \lref{l:W4}. 
Recall that $\hat{ \vv }(\xXx ) > \hat{ \vv }(u)$ for $ \xXx < u $. 
Similarly as for \eqref{:18k}, we obtain from 
\eqref{:194x}, \eqref{:18l}, and \eqref{:191s} 
\begin{align} \label{:195c}
\int_{\xXx }^{\uU } \UNz (\uU ) d \uU = 
\pi ^{1/2} \int _{\hat{ \vv }(u)} ^{\hat{ \vv }(\xXx )} \FNvVz \cos \vV d \vV 
.\end{align}

Recall that $ \None = \uu (\dNz )$ by \eqref{:191r}. 
We then divide the interval $ \None $ into two parts: $ \uu (\dNone ) $ and $ \uu (\dNz \setminus \dNone ) $. 

We first consider $ \xXx \in \uu (\dNone ) $. 
From \eqref{:191u} and \eqref{:191t}, we see that $\FNvVz $ 
is a positive and decreasing function in $ \vV $ on $ \hat{ \vv }( \uu (\dNone ) ) $. 
Then applying \lref{l:161} to the right-hand side of \eqref{:195c} and using \eqref{:191u}, 
we find a constant $\Ct \label{;IU00} >0$ such that
\begin{align}\label{:195d}
\sup_{u \in [\xXx , \xXx +1]} \Big\lvert \int_{\xXx }^{u} \UNz ( \vV ) d\vV \Big\rvert \le 
\cref{;IU00} \lvert \xXx \rvert ^{-1/4} \quad \text{for $\xXx \in \uu ( \dNone )$}
.\end{align}
From \lref{l:162} with $ \gdelta = (\UNz , \lvert \xXx \rvert ^{-1}) $, 
 \eqref{:194a}, and \eqref{:195d}, we see that \eqref{:195x} holds for $\xXx \in \uu ( \dNone )$. 

We next consider $ \uu ( \dNz \setminus \dNone ) $. 
Let $ \xXx _0 $ be such that $ [-2N^{2/3}, \xXx _0) = \uu ( \dNz \setminus \dNone ) $. Then, 
$$
\xXx _0 = \uu (\theta) \Big{\vert }_{\cos \theta = 1/N} = 2N^{-5/6}\sqrt{N+1} -2N^{2/3}
.$$
Let $\xXx _1= \xXx _0 + N^{-1/3}$. We then find a constant $\Ct \label{;W2z}> 0 $ independent of $ N $ such that 
\begin{align} \label{:195f} 
 \lvert -2N^{2/3} - u\rvert & \le \cref{;W2z} N^{-1/3}
 \quad \text{ for } \uU \in [-2N^{2/3}, \xXx _1]
.\end{align}

Because $\xXx \in \uu ( \dNz \setminus \dNone )$, we have $ \xXx < \xXx _0 < \xXx _1 $ and 
\begin{align} \label{:195F}
&\int_\xXx ^\yYy \frac{\int_{\xXx }^{u} \UNz ( \vV ) d\vV }{u - \xXx } du
= \II _1 + \II _2 + \II _3 
,\end{align}
where 
\begin{align*}
& 
\II _1 = \int_\xXx ^{\xXx _1} \frac{\int_{\xXx }^{u} \UNz ( \vV ) d\vV }{u - \xXx } du, \;
\II _2 = \int_{\xXx _1} ^\yYy \frac{\int_{\xXx }^{\xXx _0} \UNz ( \vV ) d\vV }{u - \xXx } du, \;
\II _3 = \int_{\xXx _1} ^\yYy \frac{\int_{\xXx _0 }^{u} \UNz ( \vV ) d\vV }{u - \xXx } du 
.\end{align*}
From \eqref{:194z} and \eqref{:195f}, we find $ \Ct >0 \label{;195}$ independent of $ N $ such that 
for $\xXx \in \uu ( \dNz \setminus \dNone )$ 
\begin{align} \notag 
 \lvert \II _1 (\xXx ) \rvert &\le \cref{;19X}\cref{;W2z}
N^{-1/3} N^{1/6} \le
\cref{;195} \verty ^{-1/4}
,\\ \label{:195o} 
\sup_{\yYy \in (\xXx _1, \xXx +1)} 
 \lvert \II _2 ( \yYy , \xXx ) \rvert &\le \cref{;19X}\cref{;W2z} 
N^{-1/3} N^{1/6} \int_{\xXx _1} ^\yYy \frac{1}{u - \xXx } du 
\le 
 \cref{;195} \verty ^{-1/4}\log \verty 
.\end{align}
Using \eqref{:194z}, \eqref{:195d}, and \eqref{:195f}, 
we find $ \Ct \label{;196}$ independent of $ N $ such that 
\begin{align} \notag 
&
\sup_{\yYy \in (\xXx _1, \xXx +1)} 
 \lvert \II _3 ( \yYy , \xXx ) \rvert \le 
\cref{;196} \verty ^{-1/4} \log \verty 
\quad \text{ for $\xXx \in \uu ( \dNz \setminus \dNone )$}
,\\ \label{:195p} &
\sup_{\yYy \in (\xXx , \xXx _1]} 
\Big\lvert \int_\xXx ^\yYy \frac{\int_{\xXx }^{u} \UNz ( \vV ) d\vV }{u - \xXx } du \Big\rvert \le 
\cref{;196} \verty ^{-1/4} \log \verty 
\quad \text{ for $\xXx \in \uu ( \dNz \setminus \dNone )$}
.\end{align}
From \eqref{:195F}--\eqref{:195p}, we see that \eqref{:195x} holds for $\xXx \in \uu ( \dNz \setminus \dNone ) $. 

%%%%%%%%%%%%%%%%%%%%%%%%%%%%%%%%
{\em Proof of \eqref{:195y}}: 
Similarly as for \eqref{:18k}, we obtain from 
 \eqref{:18l}, \eqref{:192a}, and \eqref{:194y} 
\begin{align}\notag &
\int_{\xXx }^{\uU } \lvert \vV \rvert ^i \UN ( \vV ) d\vV = 
\int _{\hat{ \vv }(u)} ^{\hat{ \vv }(\xXx )} \FNvVi \sin \vV d \vV 
.\end{align}
From \eqref{:18n} and \eqref{:192a}, $\FNvVi >0$ on $ \vv (\dN ) $. 
Hence, using \lref{l:161} and \lref{l:192}, we find $\Ct \label{;IU0} > 0 $ 
independent of $ N $ such that
\begin{align}\label{:195i}
\sup_{u \in [\xXx , \xXx +1]} 
\Big\lvert \int_{\xXx }^{u}
 \lvert \vV \rvert ^i \UN ( \vV ) d\vV \Big\rvert \le \cref{;IU0} \lvert \xXx \rvert ^{ i - (3/4)} \quad
\text{for $\xXx \in \None $}
.\end{align}
Taking $ \gdelta = ( \lvert \uU \rvert ^i \UN , \lvert \xXx \rvert ^{-1}) $ in \lref{l:162} and using 
\eqref{:194a} and \eqref{:195i}, we obtain \eqref{:195y}. 
\end{proof}

We set 
\begin{align} \label{:196x} 
& \hat{K}_2^N (x,y)= \frac{\psi_{N }(x)\psi'_{N }(y)-\psi'_{N }(x)\psi_{N }(y) }{x-y}
.\end{align}
Then from \eqref{:Y1k}, we have 
\begin{align}\label{:196w}&
K_{2}^{N }(x,y) =\hat{K}_2^N (x,y) - \frac{1}{2N ^{1/3}}\psi_{N }(x)\psi_{N }(y)
.\end{align}
\begin{lemma} \label{l:196}
There exists a constant $\Ct \label{;198} $ independent of $ N $ such that 
\begin{align} \label{:196b} &
\sup_{\yYY \in \mathbb{R}}
 \Big\lvert \int_\xXX ^\yYY 
 \mathbf{1}( \lvert u-\xXX \rvert \le 1) \hat{K}_2^N (u,\xXX )du \Big\rvert 
\le \cref{;198}( 1 \vee \log \lvert \xXX \rvert )
 ,\ \xXX \in \Ninfty 
.\end{align}
%for $ \xXX \in \Ninfty $.
\end{lemma}

\begin{proof} 
We first prove the case $ \xXX \in [-1, \infty ) $. 
Using \eqref{:W1w} and \lref{l:W2}, we have 
\begin{align} &\notag %\label{:1197e}&
\SupN 
\sup_{\xXX \in [-1,\infty)}
\sup_{ \lvert u - \xXX \rvert \le 1 }
 \vert \hat{K}_2^N (u,\xXX ) \vert 
< \infty 
.\end{align}
Hence, \eqref{:196b} follows from the following: 
\begin{align}\notag & %\label{:197g}&
\SupN 
\sup_{\xXX \in [-1,\infty)}
 \int_{\xXX -1 }^{\xXX +1} \lvert \hat{K}_2^N (u,\xXX ) \rvert du 
< \infty 
.\end{align}

We proceed with the case $ \xXX \in \Noneopen $. 
Let 
\begin{align}\notag 
& \JaN (\xX , \yYY ) =
 \Big( \int_\xXX ^\yYY \mathbf{1}( \lvert u-\xXX \rvert \le 1) \frac{\int_{\xXX }^{u} \psi_N'(v)dv}{u-\xXX }du \Big) 
\psi_N'(\xXX ),
\\ \label{:196z}%\\\notag 
&\JbN (\xX , \yYY ) = - 
\Big( \int_{\xXX }^{\yYY } \mathbf{1}( \lvert u-\xXX \rvert \le 1) \frac{ \int_\xXX ^u \psi_N'' (v)dv }{u-\xXX } du 
\Big) \psi_N(\xXX ) 
.\end{align}
Then from \eqref{:196x} and \eqref{:196z} with a simple calculation, we see 
\begin{align} \label{:196i}&
 \int_\xXX ^\yYY \mathbf{1}( \lvert u-\xXX \rvert \le 1) \hat{K}_2^N (u,\xXX )du =
\JaN (\xX , \yYY ) + \JbN (\xX , \yYY )
.\end{align}
Hence, we shall prove that there exists a constant $\Ct \label{;tI_11} >0$ 
independent of $ N $ such that %for $ p = 1,2,3 $
\begin{align}\label{:196a} &\quad \quad 
\sup_{\yYY \in \mathbb{R}}
 \lvert \JpN (x,y)\rvert \le \cref{;tI_11}( 1 \vee \log \lvert \xXX \rvert ), \quad \xXX \in \Noneopen ,\ p = 1,2 
.\end{align}
Using \eqref{:196i} and \eqref{:196a}, we deduce \eqref{:196b} for $ \xXX \in \Noneopen $.

Let $ \Omei $ be as in \eqref{:U1y} and \eqref{:V1q} and $ \Vi =\{ v ; ( v , N , \theta ) \in \Omei \}$. Then, 
\begin{align}\label{:196f}&
\Noneopen \subset \Vone \cup \Vtwo 
.\end{align}
Using \eqref{:W3f}, \eqref{:W3l}, and \eqref{:194x}, we deduce that 
\begin{align}\label{:196l}
\psi_N'(\vV )= - \frac{1}{\sqrt{2}\pi}\sumL \UNkm (\vV ) + 
 \OL ( \vV ^{-1}) \quad \text{on $\Ome $}
.\end{align}
From \eqref{:196l}, \eqref{:195x}, and \eqref{:W3b}, we obtain \eqref{:196a} for $ p = 1 $ on $ \Vone $. 
From \eqref{:W3e} and \eqref{:17e} for $ p = 1$, we have \eqref{:196a} for $ p = 1 $ on $ \Vtwo $. 
Thus, \eqref{:196a} holds for $ p = 1$ by \eqref{:196f}.

From \eqref{:W3x}, we obtain 
\begin{align} \notag 
\psi_N''(\vV ) & = \Big(\frac{\vV ^2}{4N^{2/3}} + \vV - \frac{1}{2N^{1/3}}\Big) \psi_N(\vV )
\\ \label{:196n} &
= \Big(\frac{\vV ^2}{4N^{2/3}} + \vV - \frac{1}{2N^{1/3}}\Big) \Big(
 {\pi}^{-1}\UN (\vV ) + \OL ({ \lvert \vV \rvert ^{-7/4}}) \Big) 
\text{ on $\Ome $}
.\end{align}
Here, we used \eqref{:93p} with $L=2$, \eqref{:91P}, and \eqref{:93b} for the second line. 
Using \eqref{:196n}, \eqref{:195y}, and \eqref{:W3a}, we deduce \eqref{:196a} for $ p = 2$ on $ \Vone $. 
From \eqref{:81z} and \eqref{:W3x}, %we have
\begin{align}\notag
\vert \psi_N '' (x) - \mathrm{Ai}''(x) \vert 
&\le \vert x\vert \vert \psi_N(x)- \mathrm{Ai}(x)\vert
+ \vert \frac{x^2}{N^{2/3}} + \frac{1}{2N^{1/3}}\vert \vert \psi_N(x) \vert
\\ \label{:196q} 
&= \underline{\mathcal{O}}( N^{-1/3+5\varepsilon/4})
\quad \text{on $\Omega_2^\varepsilon$},
\end{align}
where we used \eqref{:W3E} and \eqref{:W3a} in the last line.
Then from \eqref{:17e} and \eqref{:196q} we have \eqref{:196a} for $p=2$ on $V_{N,2}$. 
We thus have \eqref{:196a} for $ p = 2$ by \eqref{:196f}. 
%We have thus obtained \eqref{:196a}, from which we conclude \eqref{:196b}. 
\end{proof}

\noindent {\em Proof of \eqref{:Y1g}}: 
By integration by parts, we have for $ x , x ^{\dagger }, y \in \mathbb{R}$
\begin{align}\notag & %\label{:197B}&
 \int_{x ^{\dagger }}^\yYY \hat{K}_2^N (u,\xXX )du 
= 
\int_{x ^{\dagger }} ^\yYY \frac{\psi_{N }(u)\psi'_{N }(\xXX ) -\psi'_{N }(u)\psi_{N }(\xXX ) }{u - \xXX }du 
\\&\notag 
= 
\Big(\int_{x ^{\dagger }} ^\yYY \frac{\psi_{N }(u)}{u - \xXX } du \Big)\psi'_{N }(\xXX ) 
 - 
\Big(\int_{x ^{\dagger }} ^\yYY \frac{\psi'_{N }(u) }{u - \xXX }du \Big) \psi_{N }(\xXX ) 
\\&\label{:197n} 
= 
\Big(\int_{x ^{\dagger }} ^\yYY \frac{\psi_{N }(u)}{u - \xXX } du \Big)\psi'_{N }(\xXX ) 
-\Big( 
\Big[ \frac{\psi _{N }(u) }{u - \xXX }\Big]_{x ^{\dagger }} ^\yYY 
 + \int_{x ^{\dagger }} ^\yYY \frac{\psi_{N }(u)}{(u - \xXX )^2} du \Big) \psi_{N }(\xXX ) 
.\end{align}
Let $ \iotaN (u) =-(4N ^{2/3}+ u )$ be as in \eqref{:Y1u}. 
We set 
\begin{align}\label{:q2q}&
Q_N(x) = 
\int_{-2N^{2/3} }^{\infty}
 \mathbf{1}( \lvert \iotaN (u) -\xXX \rvert \ge 1) 
 \mathbf{1}( \lvert u-\xXX \rvert < 1) 
\frac{\lvert \psi_{N }(u) \rvert}{ \lvert \iotaN (u) - x \rvert }du
.\end{align}
Then from \eqref{:W3a}, we have a constant $ \Ct \label{;q2q}$ independent of $ N $ such that 
\begin{align}\label{:q2r}&
Q_N(x) \le \cref{;q2q}( \xvee ) ^{-1/4} \quad \text{ for } x \in \Ninfty 
.\end{align}
Note that $ \psi_{N } (u) = (-1)^{N} \psi_{N }(\iotaN (u)) $ for all $ u \in \mathbb{R}$. 
Recall that for $ u , \xXX \in \Ninfty $, we have $ \lvert u -x \rvert \le \lvert \iota_N(u)-x \rvert $. 
Then for $ \xXX \in \Ninfty $
\begin{align}&\notag 
 \int_{\mathbb{R}} \mathbf{1}( \lvert u-\xXX \rvert \ge 1) 
\frac{\lvert \psi_{N }(u) \rvert}{ \lvert u - \xXX \rvert }du 
\\&\notag = 
\int_{-\infty }^{- 2N^{2/3}}
 \mathbf{1}( \lvert u-\xXX \rvert \ge 1) 
\frac{\lvert \psi_{N }(u) \rvert}{ \lvert u - \xXX \rvert }du
+ 
\int_{-2N^{2/3}}^{\infty}
 \mathbf{1}( \lvert u-\xXX \rvert \ge 1) 
\frac{\lvert \psi_{N }(u) \rvert}{ \lvert u - \xXX \rvert }du
\\&\notag = 
\int_{-2N^{2/3}}^{\infty}
 \mathbf{1}( \lvert \iotaN (u) -\xXX \rvert \ge 1) 
\frac{\lvert \psi_{N }(\iotaN (u)) \rvert}{ \lvert \iotaN (u) - \xXX \rvert }du
+ 
\int_{-2N^{2/3}}^{\infty}
 \mathbf{1}( \lvert u-\xXX \rvert \ge 1) 
\frac{\lvert \psi_{N }(u) \rvert}{ \lvert u - \xXX \rvert }du
\\&\notag =
\int_{-2N^{2/3} }^{\infty}
 \mathbf{1}( \lvert \iotaN (u) -\xXX \rvert \ge 1) 
\frac{\lvert \psi_{N }(u) \rvert}{ \lvert \iotaN (u) - x \rvert }du
+ 
\int_{-2N^{2/3}}^{\infty}
 \mathbf{1}( \lvert u-\xXX \rvert \ge 1) 
\frac{\lvert \psi_{N }(u) \rvert}{ \lvert u - \xXX \rvert }du
\\&\notag =
\int_{-2N^{2/3} }^{\infty}
 \mathbf{1}( \lvert \iotaN (u) -\xXX \rvert \ge 1) 
 \mathbf{1}( \lvert u-\xXX \rvert < 1) 
\frac{\lvert \psi_{N }(u) \rvert}{ \lvert \iotaN (u) - x \rvert }du
\\&\notag 
+
\int_{-2N^{2/3} }^{\infty}
 \mathbf{1}( \lvert \iotaN (u) -\xXX \rvert \ge 1) 
 \mathbf{1}( \lvert u-\xXX \rvert \ge 1) 
\frac{\lvert \psi_{N }(u) \rvert}{ \lvert \iotaN (u) - x \rvert }du
\\&\notag 
+ 
\int_{-2N^{2/3}}^{\infty}
 \mathbf{1}( \lvert u-\xXX \rvert \ge 1) 
\frac{\lvert \psi_{N }(u) \rvert}{ \lvert u - \xXX \rvert }du
\\&\label{:197p} \le 
% \int_{-2N^{2/3} }^{\infty} \mathbf{1}( \lvert \iotaN (u) -\xXX \rvert \ge 1) 
% \frac{\lvert \psi_{N }(u) \rvert}{ \lvert u - x \rvert }du 
%Q_N(\xXX )
\cref{;q2q}( \xvee ) ^{-1/4} + 
2 \int_{-2N^{2/3}}^{\infty}
 \mathbf{1}( \lvert u-\xXX \rvert \ge 1) 
\frac{\lvert \psi_{N }(u) \rvert}{ \lvert u - \xXX \rvert }du
.\end{align}
Here we used \eqref{:q2q}, \eqref{:q2r}, and 
$ \lvert u -x \rvert \le \lvert \iota_N(u)-x \rvert $ for the last line. 
Hence using \lref{l:W3} and \eqref{:197p}, 
we have a constant $ \Ct \label{;197h}$ independent of $ N \in \mathbb{N}$ such that 
for $ \xXX \in \Ninfty $ 
\begin{align}&\notag 
\Big( \int_{\mathbb{R}} \mathbf{1}( \lvert u-\xXX \rvert \ge 1) 
\frac{\lvert \psi_{N }(u) \rvert}{ \lvert u - \xXX \rvert }du \Big)
\big\lvert \psi'_{N }(\xXX ) \big \rvert 
\\\notag 
 \le & \Big( \cref{;q2q}( \xvee ) ^{-1/4} + 
2\int_{-2N^{2/3}}^{\infty}
 \mathbf{1}( \lvert u-\xXX \rvert \ge 1) 
\frac{\lvert \psi_{N }(u) \rvert}{ \lvert u - \xXX \rvert }du
\Big) 
\big\lvert \psi'_{N }(\xXX ) \big \rvert 
\\ \notag \le &
 \cref{;q2q} \cref{;W66} + 
2 \cref{;W66}^2 
 \Big(
\int_{-2N^{2/3}}^{\infty}
 \mathbf{1}( \lvert u-\xXX \rvert \ge 1) 
\frac{ ( 1 \vee \lvert u \vert )^{-1/4}}{ \lvert u - \xXX \rvert }du \Big) 
( \xvee ) ^{1/4}
\\ \notag \le &
 \cref{;q2q} \cref{;W66}+ 
2 \cref{;W66}^2 \Big(
\int_{-\infty}^{\infty}
 \mathbf{1}( \lvert u-\xXX \rvert \ge 1) 
\frac{ ( 1 \vee \lvert u \vert )^{-1/4}}{ \lvert u - \xXX \rvert }du \Big) 
( \xvee ) ^{1/4}
\\ \label{:197o} 
 \le &\cref{;197h} 
( 1 \vee \log \lvert \xXX \rvert )
.\end{align}

From \lref{l:W3} and $ \psi_{N } (u) = (-1)^{N} \psi_{N }(-4N^{2/3} - u) $ for all $ u \in \mathbb{R}$, we deduce that 
$ \supN \sup_{\xXX \in \mathbb{R}} \lvert \psi_{N }( \xXX ) \rvert < \infty $. 
Hence, we have 
\begin{align} \notag 
\supN 
\sup_{\xXX \in \mathbb{R}}
\Big( &
\sup_{ \yYY \ge \xXX + 1}
\Big\lvert 
 \frac{ \psi_{N }( \yYY )}{\yYY - \xXX } - 
 \psi_{N }( \xXX + 1 ) 
\Big \rvert 
+ 
\sup_{ \yYY \le \xXX - 1}
\Big\lvert 
 \frac{ \psi_{N }( \yYY )}{\yYY - \xXX } - 
 \psi_{N }( \xXX - 1 ) 
\Big \rvert 
\\ & \label{:197j}
 + \int_{ \lvert u - \xXX \rvert \ge 1 } 
 \frac{\lvert \psi_{N }(u) \rvert }{(u - \xXX )^2} du \Big)
\big\lvert \psi_{N }(\xXX ) \big \rvert 
< \infty 
.\end{align}

From \eqref{:197n}, \eqref{:197o}, and \eqref{:197j}, we have a constant $ \Ct \label{;197f}$ 
independent of $ N $ such that for $ \xXX \in \Ninfty $
\begin{align}\label{:197f}&
\sup_{ \yYY \in \mathbb{R}}
 \Big\lvert \int_{\xXX }^\yYY 
 \mathbf{1}( \lvert u-\xXX \rvert \ge 1) 
 \hat{K}_2^N (u,\xXX )du \Big\rvert 
\le \cref{;197f} ( 1 \vee \log \lvert \xXX \rvert ) 
%\quad \text{ for } \xXX \in \Ninfty 
.\end{align} 

Using \eqref{:W3a} and \eqref{:W5a}, we deduce 
\begin{align}\label{:197d}
& \SupN 
\sup_{\xXX \in \Ninfty }
\sup_{\yYY \in \mathbb{R}}
\Big\lvert \int_\xXX ^\yYY \frac{1}{2N ^{1/3}}\psi_{N }(u )\psi_{N }(\xXX ) du \Big\rvert < \infty
.\end{align}

From \eqref{:196w}, \eqref{:196b}, \eqref{:197f}, and \eqref{:197d}, we obtain \eqref{:Y1g}. 
\qed

\section*{Acknowledgments}

This work was supported by JSPS KAKENHI 
 Grant Numbers JP16K13764, JP16H06338, JP18H03672, JP19H01793, JP20K20885, JP21H04432. 
%%%%%%%%%%%%%%%%%%%%%%%%%%%%%%%%%%%%%%%%%%%%%%%

%%% Please do NOT use ``\bysame'' command in your bibliography list.

 \end{document}